\documentclass[smallextended]{svjour3}
\pdfoutput=1 
\usepackage{algorithm}
\usepackage{algpseudocode}
\usepackage{amssymb}
\usepackage[title,titletoc]{appendix}
\usepackage{bbm}
\usepackage{braket}
\usepackage{bookmark}
\usepackage{booktabs}
\usepackage{caption}
\usepackage{enumitem}
\usepackage[T1]{fontenc} 
\usepackage{float}
\usepackage{indentfirst}
\usepackage[utf8]{inputenc}
\usepackage[autostyle,italian=guillemets]{csquotes}
\usepackage{graphicx}
\graphicspath{{Immagini/}}
\DeclareGraphicsExtensions{{.png},{.pdf},{.jpg}}
\graphicspath{{./}}
\usepackage{lmodern}
\usepackage{lipsum} 
\usepackage{listings}
\usepackage{yhmath}
\usepackage{mathtools}
\usepackage{mathrsfs}
\usepackage{microtype} 
\usepackage{multirow}
\usepackage{setspace}
\usepackage{subcaption}
\usepackage{tabularx}
\usepackage[table]{xcolor}
\usepackage{wasysym}
\usepackage{url} 
\usepackage[overload]{empheq}

\graphicspath{{Images/}}
\begin{document}

\title{A reduced order model for domain decompositions with non--conforming interfaces}

\author{Elena Zappon \and Andrea Manzoni \and Paola Gervasio \and Alfio Quarteroni}
\institute{E. Zappon, A. Manzoni, A. Quarteroni \at MOX -- Dipartimento di Matematica,\\
	 Politecnico di Milano,\\
	  P.zza Leonardo da Vinci 32, \\
	  I-20133 Milano, Italy\\
	  \email{andrea1.manzoni@polimi.it}, \email{alfio.quarteroni@polimi.it}
	  \and
	  E. Zappon \at Gottfried Schatz Research Center, Biophysics, \\
	  Medical University of Graz,\\
	  Neue Stiftingtalstrasse 6/D04,\\
	  8010 Graz, Austria\\
	  \email{elena.zappon@medunigraz.at}	  
	  \and
	  P. Gervasio \at DICATAM,\\
	  Universit\`a degli Studi di Brescia,\\
	  Via Branze, 38, \\
	  25123 Brescia, Italy\\
	  \email{paola.gervasio@unibs.it}
	  \and
	  A. Quarteroni \at 
	  Institute of Mathematics,\\
	   Ecole Polytechnique Federale de Lausanne,\\
	    Station 8, \\
	     CH-1015 Lausanne, Switzerland \\(Professor Emeritus)}
	
\maketitle
	
\begin{abstract}
In this paper, we propose a reduced--order modeling strategy for two--way Dirichlet--Neumann parametric coupled problems solved with domain--decomposition (DD) sub--structuring methods. We split the original coupled differential problem into two sub--problems with Dirichlet and Neumann interface conditions, respectively. After discretization by, \emph{e.g.}, the finite element method, the full--order model (FOM) is solved by Dirichlet--Neumann iterations between the two sub--problems until interface convergence is reached. We then apply the reduced basis (RB) method to obtain a low--dimensional representation of the solution of each sub--problem. Furthermore, we apply the discrete empirical interpolation method (DEIM) at the interface level to achieve a fully reduced--order representation of the DD techniques implemented. To deal with non--conforming FE interface discretizations, we employ the INTERNODES method combined with the interface DEIM reduction. The reduced--order model (ROM) is then solved by sub--iterating between the two reduced--order sub--problems until the convergence of the approximated high--fidelity interface solutions. The ROM scheme is numerically verified on both steady and unsteady coupled problems, in the case of non--conforming FE interfaces.
\end{abstract}

\section{Introduction}

Reduced order modeling (ROM) techniques encompass a wide class of numerical methods able to solve differential problems several orders of magnitude faster than conventional, high--fidelity full order models (FOMs), achieving real--time decision--making and operational modeling in several contexts, ranging from fluid dynamics \cite{Amsallem2010,Ballarin2017,Ballarin2016,FortiRozza2014,Lorenzi2016,noack2016} to biomedical engineering \cite{BallarinFaggiano2016,Fresca2021,Pagani2018}. 

In several cases, when dealing with complex -- possibility multi--physics -- problems, domain decomposition (DD) \cite{QuarteroniValli1999,Schwarz1869} techniques are necessary to {\em (i)} split the domain into two or more regions in which either the same or different methods are used to approximate the solution \cite{Confalonieri2013,Glowinski1083,Paz2010,Weilun1998}, to {\em (ii)}  solve systems that stem from the assembly of independently generated meshes \cite{Hartmann2009,McGee2005,Oliver2009} or to {\em (iii)} frame coupled problems where the physical nature of the involved sub--models is very different \cite{deboer2007}. In the last two cases, interface non--conformity issues can arise, and \emph{ad--hoc} techniques, such as \emph{e.g.} the MORTAR method \cite{Chan1996,Bernardi2005,Hesch2014} or the INTERNODES method \cite{Deparis2016,GervasioQuarteroni2018,Gervasio2019} have to be implemented in order to ensure the correct exchange of information at the interface.

Domain decomposition schemes coupled with ROM techniques, especially the reduced basis (RB) method \cite{Benner2017,Hesthaven2016,Quarteroni2016}, have been first used in \cite{lucia2003}, where ROMs were applied only on small regions of the domains, for instance for PDEs with discontinuous solutions in those regions where the discontinuity occurs. Similar strategies have been implemented to tackle the numerical simulation of problems in fluid dynamics \cite{Baiges2013,IAPICHINO201263,martini2015reduced,PEGOLOTTI2021113762,Washabaugh2012,Wicke2009,Xiao2019}, aerospace engineering \cite{Legresley2006} and structural mechanics  \cite{Kerfriden2013,Corigliano2015}, as well as for the optimization of complex systems \cite{Antil2012,Antil2011}. In order to reduce the global parametrized problem, other methods such as the reduced basis element methods (RBEM) \cite{Maday2004,Lvgren2006ARB}, the reduced basis hybrid (RBHM) method \cite{IAPICHINO201263}, and the static condensation methods (SCRBEM) \cite{Eftang2014,Huynh2013} have been developed.
A common feature shared by several of these works is the application of a reduced order model in small parts of the domains, through an \emph{ad--hoc} definition of a small set of basis functions in each subdomain able to preserve the solution continuity along the internal subdomain interfaces. An effort to include different interface conditions comes with the RDF method \cite{IAPICHINO2016408}, where local basis functions are employed on each selected sub--problem and interface conditions are treated in high--fidelity form through the application of different techniques, \emph{e.g.} Lagrangian multipliers or Fourier basis functions, for sub--iterating schemes and/or solving the coupled model in monolithic form, naturally imposing interface constraints. In \cite{MAIER201431} also Dirichlet and Neumann data on conforming interfaces are considered in reduced form, while the Dirichlet--Neumann DD iterations structure is preserved by the reduced algorithm, whereas an algebraic--splitting approach between internal and interface nodes is employed in the domain--decomposition least--squares Petrov--Galerkin (DD--LSPG) approach \cite{HOANG2021113997}. However, to the best of our knowledge, interface non--conformity has never been considered, given its intrinsic complexity.

In this work, we present a Dirichlet--Neumann DD--ROM relying on the RB method able to transfer interface data across non--conforming interfaces. In particular, we consider parametric second--order elliptic and/or parabolic problems solved with Dirichlet--Neumann sub--structuring DD algorithms. Therefore, we split the problem domain into two non--overlapping subdomains with a common interface and we define two parametrized sub--problems with Dirichlet or Neumann interface conditions. RB methods are then applied at the sub--problems level to approximate both the sub--problems solution together with the Dirichlet or Neumann interface conditions.

To solve the coupled problem and compute the snapshots (for different parameter instances) required
to train our ROM, we rely on the finite element (FE) method as a high--fidelity FOM. In particular, the
FOM solution is sought by sub--iterating between the two sub--models solutions until convergence, being
this latter is reached when the difference between the solution at both sides of the interface falls under a
prescribed tolerance. In the presence of a non--conforming interface, FOM solutions can be sought through FE--based methods such as the MORTAR method or the INTERNODES. In our code, however, we surrogate the non--conforming coupling by solving for each selected parameter instance the
coupled problems twice: indeed, we define two possible FE discretizations in the complete domain (each one being obtained by extending to the other subdomain the spatial discretization set in the first one) and compute the solution of each sub--problem using both discretizations, one for each simulation, featuring interface conformity. 
Then, we extract the solution of each of the sub--problems (in the discretization originally set for the corresponding subdomain), therefore obtaining two sets of solution snapshots as well as the corresponding Dirichlet and Neumann traces at the common interface that can be seen as the original model solution when interface nonconformity is considered (see Remark \ref{rem:complete_ROM} and Fig. \ref{Fig:domain_representation}). 

The RB method is then applied to define a low--dimensional representation of the solution in each subdomain, while the discrete empirical interpolation method (DEIM) is applied to both Dirichlet and Neumann interface data in order to achieve a fully reduced--order representation of the DD method adopted. In the online phase, a Galerkin projection is used to reduce the sub--problems' dimension obtaining two reduced--order sub--problems. The solution of the coupled problem, for any new parameter instance, is then found by iterating between the solutions of the two reduced sub--problems until (a suitable norm of) the difference between the two solutions at the interface, once traced back at the high--fidelity level, falls below a prescribed tolerance. In this phase, we transfer Dirichlet and Neumann interface data by applying the INTERNODES method and by using the DEIM as well as a low--order piecewise constant interpolation. This approach, which extends the work presented in \cite{Zappon2022} is still modular, and allows us to achieve a complete reduction of the model at hand, which can be seen as a two--way coupled model, including interface non--conforming grid cases. 

The organization of the paper is the following: in Section \ref{sec:Two_way_models} we present the formulation and the high--fidelity discretization of the parametrized problem, considering both steady and unsteady cases. Section \ref{sec:reduced_order_formulation} is devoted to the reduced--order formulation of the two sub--problems, while Section \ref{sec:interface_reduction} is dedicated to the interface Dirichlet and Neumann reduced formulation. In Section \ref{sec:numerical_results}, the algorithm is numerically verified by means of two test cases dealing with second--order linear PDEs, considering both an elliptic and a parabolic problem. Section \ref{sec:conclusion} then reports some final remarks and possible perspectives of this work.

\section{Two--way coupled problem and its FE discretization}
\label{sec:Two_way_models}
%
We introduce the parameter--dependent two--way coupled problem. Such a problem can arise from the application of splitting domain decomposition methods to single--physics or multi--physics models. For simplicity, let us consider the single--physics case -- the same procedure can be considered for the multi--physics case -- defined on an open bounded domain $\Omega \subset \mathbb{R}^n$, $(n=2,3)$ with Lipschitz boundary $\partial \Omega$. $\partial \Omega_D$ and $\partial \Omega_N$ denote the Dirichlet and the Neumann boundary, respectively, such that $\overline{\partial\Omega_D} \cup \overline{\partial\Omega_N} = \partial\Omega$ and $\partial\Omega_D \cap \partial\Omega_N =\emptyset$. Given the set of parameters $\boldsymbol{\mu} \in \mathscr{P}^d \subset \mathbb{R}^d$, $d\geq 1$, we search for $u(\boldsymbol{\mu})$ on $\Omega$ such that \vspace{-0.05cm}
\begin{equation}
\label{eq:FOM_global}
\begin{cases}
\mathcal{L}(\boldsymbol{\mu})u(\boldsymbol{\mu}) = f(\boldsymbol{\mu}) &\text{in }\Omega\\
u(\boldsymbol{\mu}) = g_D(\boldsymbol{\mu}) &\text{on }\partial \Omega_{D}\\
\partial_{\mathcal{L}(\boldsymbol{\mu})}u(\boldsymbol{\mu}) = g_N(\boldsymbol{\mu}) &\text{on } \partial \Omega_{N},
\end{cases}
\end{equation}
where $\mathcal{L}(\boldsymbol{\mu})$ is a second--order elliptic operator, $f(\boldsymbol{\mu})$, $g_D(\boldsymbol{\mu})$ and $g_N(\boldsymbol{\mu})$ are functions defined on $\Omega$, $\partial \Omega_D$ and $\partial \Omega_{N}$, respectively, and $\partial_{\mathcal{L}(\boldsymbol{\mu})}u(\boldsymbol{\mu})$ is the conormal derivative associated with the operator $\mathcal{L}(\boldsymbol{\mu})$ on $\partial\Omega$. 

Now, we split the computational domain $\Omega$ into two non--overlapping subdomains $\Omega_1$ and $\Omega_2$ with Lipschitz boundary $\partial \Omega_i$, $i = 1,2$ and a common interface $\Gamma := \partial \Omega_1 \cap \partial\Omega_2$. For each $i=1,2$, we denote $\partial\Omega_{i,D} = \partial\Omega_i \cap \partial \Omega_D$ and $\partial\Omega_{i,N} = \partial \Omega_i \cap \partial \Omega_N$. 

The two--domains formulation of model \eqref{eq:FOM_global} reads \cite{QuarteroniValli1999}:
\begin{equation}
\label{Eq:generale_DD_model}
\begin{cases}
\mathcal{L}_i(\boldsymbol{\mu}) u_i(\boldsymbol{\mu}) = f(\boldsymbol{\mu}) &\text{in }\Omega_i, ~i=1,2\\
u_i(\boldsymbol{\mu}) = g_D &\text{on }\partial\Omega_{i,D}\\
\partial_{\mathcal{L}_i(\boldsymbol{\mu})}u_i(\boldsymbol{\mu}) = g_N(\boldsymbol{\mu}) &\text{on }\partial\Omega_{i,N}\\
u_1(\boldsymbol{\mu}) = u_2(\boldsymbol{\mu}) &\text{on }\Gamma\\
\partial_{\mathcal{L}_1(\boldsymbol{\mu})}u_1(\boldsymbol{\mu}) +\partial_{\mathcal{L}_2(\boldsymbol{\mu})}u_2(\boldsymbol{\mu}) = 0  &\text{on }\Gamma,
\end{cases}
\end{equation}
where $\mathcal{L}_i$ is the differential operator acting on functions defined on the sub--domain $\Omega_i$. 

The \emph{Dirichlet-Neumann iterative scheme} \cite{Bjorstad2018,QuarteroniValli1999} is applied to solve problem \eqref{Eq:generale_DD_model}. Therefore, starting from the initial guesses $u_1^0(\boldsymbol{\mu})$ and $u_2^0(\boldsymbol{\mu})$, for each $k\geq 0$, we search for  $u_1^{k+1}(\boldsymbol{\mu})$ on $\Omega_1$ and $u_2^{k+1}(\boldsymbol{\mu})$ on $\Omega_2$ such that  
\begin{subequations}
\label{eq:master_FOM}
\begin{align}[left = { \empheqlbrace}]
&\mathcal{L}_1(\boldsymbol{\mu})u_1^{k+1}(\boldsymbol{\mu}) = f(\boldsymbol{\mu}) &\text{in }\Omega_1\\
&u_1^{k+1}(\boldsymbol{\mu}) =  u_2^{k}(\boldsymbol{\mu}) &\text{on } \Gamma \label{eq:dirichlet_master}\\
&u_1^{k+1}(\boldsymbol{\mu}) = g_{D}(\boldsymbol{\mu}) &\text{on }\partial\Omega_{1,D} \label{eq:general_Dirichlet_master}\\ &\partial_{\mathcal{L}_1(\boldsymbol{\mu})}u_1^{k+1}(\boldsymbol{\mu}) = g_{N}(\boldsymbol{\mu}) &\text{on }\partial\Omega_{1,N}  
\end{align}
\end{subequations}
and  
\begin{subequations}
\label{eq:slave_FOM}
\begin{align}[left = { \empheqlbrace}]
&\mathcal{L}_2(\boldsymbol{\mu})u_2^{k+1}(\boldsymbol{\mu}) = f(\boldsymbol{\mu}) &\text{in }\Omega_2\\
&\partial_{\mathcal{L}_2(\boldsymbol{\mu})}u_2^{k+1}(\boldsymbol{\mu})+\partial_{\mathcal{L}_1(\boldsymbol{\mu})}u_1^{k+1}(\boldsymbol{\mu}) = 0 &\text{on }\Gamma \label{eq:interface_neumann_slave}\\
&u_2^{k+1}(\boldsymbol{\mu}) = g_{D}(\boldsymbol{\mu}) &\text{on }\partial\Omega_{2,D}\label{eq:general_Dirichlet_slave}\\ &\partial_{\mathcal{L}_2(\boldsymbol{\mu})}u_2^{k+1}(\boldsymbol{\mu}) = g_{N}(\boldsymbol{\mu}) &\text{on }\partial\Omega_{2,N},
\end{align}
\end{subequations}
where $\partial_{\mathcal{L}_i(\boldsymbol{\mu})}\mathbf{u}(\boldsymbol{\mu})$ is the conormal derivative associated with the operator $\mathcal{L}_i(\boldsymbol{\mu})$ on $\partial\Omega_i$. Moreover, a relaxation technique \cite{QuarteroniValli1999} is usually applied to ensure and accelerate the scheme convergence.

Then, we denote the problem and the corresponding solution on $\Omega_1$ as \emph{slave model} and \emph{slave solution}, and the ones in $\Omega_2$ as \emph{master model} and \emph{master solution}, respectively.

Defining $V = H_0^1(\Omega)$, for each $i=1,2$, we first define the local spaces  
\begin{equation}
\label{eq:v_i_spaces}
V_i = \{ v \in H^1(\Omega_i)\ | \ v = 0 \text{ on } \partial \Omega_{i,D}\} \quad \text{and} \quad V_i^0 = \{ v \in V_i~|~ v = 0 \text{ on }\Gamma \}.  \end{equation}

We consider two a--priori independent discretizations $\mathcal{T}_{h_1}$ and $\mathcal{T}_{h_2}$ on the domains $\Omega_1$ and $\Omega_2$ that can imply a mesh non--conformity at the interface. For instance, $\mathcal{T}_{h_i}$ can be made of simplices (triangles or tetrahedra) or quads (quadrilaterals or hexahedra), depending on the mesh size, a positive parameter $h_i > 0$. Moreover, different mesh sizes $h_1$ and $h_2$, or different polynomial degrees $p_1$ or $p_2$, can be selected. Then, we denote the internal interfaces of $\Omega_1$ and $\Omega_2$ induced by $\mathcal{T}_{h_1}$ and $\mathcal{T}_{h_2}$ as $\Gamma_1$ and $\Gamma_2$, respectively: 
we talk of geometrical conformity if $\Gamma_1 = \Gamma_2$ and of geometrical non--conformity if $\Gamma_1 \not = \Gamma_2$.
Finally, we assume that for any $T_{i,m} \in \mathcal{T}_{h_i}$, $\partial T_{i,m} \cap \partial \Omega$ fully belongs to $\partial \Omega_{i,D}$ or $\partial \Omega_{i,N}$. According to the test cases in Section \ref{sec:numerical_results}, hereon we only consider quads elements. 
 
For each partition $\mathcal{T}_{h_i}$ we define the finite element approximation spaces as
$$ X_{h_i}^{p_i} = \{ v \in C^0(\overline{\Omega_i}) : v_{|T_{i,m}}\circ F_{i,m} \in \mathbb{Q}_{p_i}, \forall T_{i,m} \in \mathcal{T}_{h_i}\},$$
in which $F_{i,m}$ is a smooth bijection that maps the reference quad $(-1,1)^n$ into the quad $T_{i,m}$, and $p_i$ are chosen integers. 
The finite--dimensional spaces to define the discrete formulation of the exploited problems will be
\begin{equation}
\label{eq:V_h_i}
V_{h_i} = \{ v \in X_{h_i}^{p_i} : v_{|\partial \Omega_{i,D}} = 0 \}, \quad V_{h_i}^0 = \{ v \in V_{h_i}, v_{|\Gamma} = 0\}, \quad i =1,2,\end{equation}
while the spaces of traces on $\Gamma$ are
\begin{equation}
\label{eq:Y_h_i} Y_{h_i} = \{ \lambda = v_{|\Gamma}, v \in X_{h_i}^{p_i} \} \quad \text{ and } \quad \Lambda_{h_i} = \{ \lambda = v_{|\Gamma}, v \in V_{h_i} \}. \end{equation}
When the Dirichlet boundaries and the interface share common DoFs, the spaces $Y_{h_1}$ and $\Lambda_{h_1}$ (as well as $Y_{h_2}$ and $\Lambda_{h_2}$) are different and should be addressed separately, whereas when $\Gamma_i \cap \Omega_{i,D} = \emptyset$, the two spaces $Y_{h_1}$ and $\Lambda_{h_1}$ (as well as $Y_{h_2}$ and $\Lambda_{h_2}$) coincide. Since the problem we are interested in falls in the first situation, from now on we assume that $\Lambda_{h_i}=Y_{h_i}$.

We  also introduce two independent transfer operators able to exchange information between the independent grids on the interface $\Gamma$, namely
$$ \Pi_{12} : Y_{h_2} \rightarrow Y_{h_1} \quad \text{and} \quad \Pi_{21} : Y_{h_1}\rightarrow Y_{h_2}.$$
In the non--conforming case, if $\Gamma_1$ and $\Gamma_2$ coincide, such operators could be the classical Lagrange interpolation operators, while $\Pi_{jk}$ are the identity operators when the meshes are conforming. Instead, if the meshes are non--conforming and $\Gamma_1 \not = \Gamma_2$, $\Pi_{12}$ and $\Pi_{21}$ could be, \emph{e.g.}, Rescaled Localized Radial Basis Function operators, as for the INTERNODES \cite{Deparis2016,Gervasio2018,Gervasio2019}. 

To exploit a FE--Galerkin approximation to set the high--fidelity FOM and get the algebraic formulation of problems \eqref{eq:master_FOM} and \eqref{eq:slave_FOM}, it is useful to consider local vectors and matrices. In particular, we define the following set of indices associated with the nodes $\mathbf{x}_j$ of the mesh in $\Omega_i$:
\begin{equation}
\begin{split}
&\mathcal{I}_{\overline{\Omega}_i} = \{1,\dots,N_{i}\}, \qquad \mathcal{I}_{\Gamma_i} = \{ j \in \mathcal{I}_{\overline{\Omega}_{i}}: \mathbf{x}_j \in \overline{\Gamma_i}\}, \qquad   \mathcal{I}_{D_i} = \{ j \in \mathcal{I}_{\overline{\Omega}_{i}}: \mathbf{x}_j \in \overline{\partial\Omega_{D,i}} \}\\ 
&\mathcal{I}_{1} = \{j \in \mathcal{I}_{\overline{\Omega}_{1}} : \mathbf{x}_j \in \overline{\Omega}_{1} \backslash (\overline{\partial \Omega_{D,1}}\cup \mathring{\Gamma_1})\} \qquad \mathcal{I}_{2} = \{j \in \mathcal{I}_{\overline{\Omega}_{2}} : \mathbf{x}_j \in \overline{\Omega}_{2} \backslash \overline{\partial \Omega_{D,2}}\},\\
\end{split}
\end{equation}
being $N_i$ the cardinality of $\mathcal{I}_{\overline{\Omega}_i}$. Moreover, we denote by $\tilde{N}_i$ the cardinality of $\mathcal{I}_i$. Note that the definitions of $\mathcal{I}_1$ and $\mathcal{I}_2$ are different: indeed, $\mathcal{I}_1$ represents all the nodes in $\Omega_1$ minus the nodes on both the interface and the Dirichlet portion of the boundary of $\Omega_1$, whereas $\mathcal{I}_2$ contains all the nodes in $\Omega_2$ minus only the nodes on the Dirichlet portion of the boundary of $\Omega_2$. 

Then, for each $i = 1,2$, we set the local stiffness matrices $\mathbb{A}_i(\boldsymbol{\mu})$
so that  
$$\mathbb{A}_{i,i}(\boldsymbol{\mu})= \mathbb{A}_i(\mathcal{I}_i,\mathcal{I}_i;\boldsymbol{\mu})$$
is the submatrix of $\mathbb{A}_i(\boldsymbol{\mu})$ of the rows and columns of $\mathbb{A}_i(\boldsymbol{\mu})$ whose indices belong to $\mathcal{I}_i$, \emph{i.e.} the internal nodes of $\Omega_i$ or those on $\partial\Omega_{i,N}$. Similarly, we can define $\mathbb{A}_{\Gamma_i,\Gamma_i}(\boldsymbol{\mu}) = \mathbb{A}_i(\mathcal{I}_{\Gamma_i},\mathcal{I}_{\Gamma_i};\boldsymbol{\mu})$, $\mathbb{A}_{i,\Gamma_i}(\boldsymbol{\mu}) = \mathbb{A}_i(\mathcal{I}_{i},\mathcal{I}_{\Gamma_i};\boldsymbol{\mu})$, $\mathbb{A}_{\Gamma_i,i}(\boldsymbol{\mu}) = \mathbb{A}_i(\mathcal{I}_{\Gamma_i},\mathcal{I}_{i};\boldsymbol{\mu})$, and $\mathbb{A}_{i,D}(\boldsymbol{\mu}) = \mathbb{A}_i(\mathcal{I}_{i},\mathcal{I}_{D_i};\boldsymbol{\mu})$.

Moreover, if $\mathbf{f}_{N,i}(\boldsymbol{\mu})$ and $\mathbf{u}_{N,i}(\boldsymbol{\mu})$ are the right--hand side vector (including the Neumann data $g_N(\boldsymbol{\mu}$)) and the vector of degrees of freedom of the approximated solution in $\overline{\Omega_i}$, respectively, we set
$$\mathbf{f}_{i}(\boldsymbol{\mu}) = \mathbf{f}_{N,i}(\mathcal{I}_{i};\boldsymbol{\mu}), \quad \mathbf{f}_{\Gamma_{i}}(\boldsymbol{\mu}) = \mathbf{f}_{N,i}(\mathcal{I}_{\Gamma};\boldsymbol{\mu}),$$
$$\mathbf{u}_i(\boldsymbol{\mu}) = \mathbf{u}_{N,i}(\mathcal{I}_{i};\boldsymbol{\mu}), \quad \mathbf{u}_{\Gamma_i}(\boldsymbol{\mu}) = \mathbf{u}_{N,i}(\mathcal{I}_{\Gamma_i};\boldsymbol{\mu}).$$
Then we can explicit the algebraic form of \eqref{eq:master_FOM} as: for each $k\geq 0$, find $\mathbf{u}_1^{k+1}(\boldsymbol{\mu})$ solution of
\begin{equation}
\label{eq:master_steady_algebraic}
\begin{cases}
\mathbb{A}_{1,1}(\boldsymbol{\mu}) \mathbf{u}_{1}^{k+1}(\boldsymbol{\mu}) = \mathbf{f}_1(\boldsymbol{\mu}) - \mathbb{A}_{1,D}(\boldsymbol{\mu})\mathbf{g}_{D,1}(\boldsymbol{\mu}) - \mathbb{A}_{1,\Gamma_1}(\boldsymbol{\mu})\mathbf{u}_{\Gamma_1}^{k+1}(\boldsymbol{\mu}) \\
\mathbf{u}_{\Gamma_1}^{k+1}(\boldsymbol{\mu}) = \mathbb{R}_{12} \mathbf{u}_{\Gamma_2}^{k}(\boldsymbol{\mu}).
\end{cases}
\end{equation}
where $\mathbb{R}_{12}$ is the rectangular matrix associated with $\Pi_{12}$, and $\mathbf{g}_{D,1}(\boldsymbol{\mu})$ the vector whose elements are the evaluation of $g_D(\boldsymbol{\mu})$ on the Dirichlet boundary nodes of $\Omega_1$. 

The algebraic formulation of problem \eqref{eq:slave_FOM} reads as: for each $k\geq 0$, find $\mathbf{u}_2^{k+1}(\boldsymbol{\mu})$ such that 
\begin{equation}
\label{eq:slave_steady_algebraic}
\begin{cases}
\mathbb{A}_{2,2}(\boldsymbol{\mu}) \mathbf{u}_{2}^{k+1}(\boldsymbol{\mu}) = \mathbf{f}_{2}(\boldsymbol{\mu})+{\mathbb E}_2\mathbf{r}_{\Gamma_2}^{k+1}(\boldsymbol{\mu})-\mathbb{A}_{2,D}(\boldsymbol{\mu})\mathbf{g}_{D,2}(\boldsymbol{\mu})\\
 \mathbf{r}_{\Gamma_2}^{k+1}(\boldsymbol{\mu}) = - \mathbb{M}_{\Gamma_2}\mathbb{R}_{21}\mathbb{M}_{\Gamma_1}^{-1} \mathbf{r}_{\Gamma_1}^{k+1}(\boldsymbol{\mu}),\end{cases}
\end{equation}
where, for $i=1,2$, $\mathbb{M}_{\Gamma_i}$ is the \emph{interface mass matrix} on $\Gamma_i$, $\mathbb{R}_{21}$ is the matrix associated with $\Pi_{21}$, $\mathbf{g}_{D,2}(\boldsymbol{\mu})$ is the vector whose elements are the evaluation of $g_D(\boldsymbol{\mu})$ at the Dirichlet boundary nodes of $\Omega_2$, \begin{equation} 
\label{eq:interface_residual_def_steady} \mathbf{r}_{\Gamma_i}^{k+1}(\boldsymbol{\mu}) = \left(\mathbb{A}_i(\boldsymbol{\mu})\mathbf{u}_{N,i}^{k+1}(\boldsymbol{\mu}) - \mathbf{f}_{N,i}(\boldsymbol{\mu})\right)_{|\Gamma_i},
\end{equation} 
is the residual at the interface $\Gamma_i$, and 
\begin{eqnarray*}
\mathbf{u}_{N,1}^{k+1}(\boldsymbol{\mu}) = \left [ \begin{matrix}
\mathbf{u}_1^{k+1}(\boldsymbol{\mu})\\
\mathbf{u}_{\Gamma_1}^{k+1}(\boldsymbol{\mu})\\
\mathbf{g}_{D,1}(\boldsymbol{\mu})
\end{matrix} \right ], \qquad
\mathbf{u}_{N,2}^{k+1}(\boldsymbol{\mu}) = \left [ \begin{matrix}
\mathbf{u}_2^{k+1}(\boldsymbol{\mu})\\
\mathbf{g}_{D,2}(\boldsymbol{\mu})
\end{matrix} \right ]
\end{eqnarray*}
are the complete solutions of the sub--problems \eqref{eq:master_FOM} and \eqref{eq:slave_FOM}, and, finally, ${\mathbb E}_2\in{\mathbb R}^{\tilde N_2\times N_{2,\Lambda}}$ is the extension matrix from the interface to all the domain DoFs in ${\mathcal I}_2$ such that $({\mathbb E}_2)_{ij}=1$ only if $i\in  \mathcal{I}_2 $, $j \in {\mathcal I}_{\Gamma_2}$, otherwise it is null.

The meaning of \eqref{eq:slave_steady_algebraic} is, therefore, that of: 
\begin{enumerate}[noitemsep]
\item moving the residual vector ${\bf r}_{\Gamma_1}^{k+1}(\boldsymbol{\mu})$ from the dual space to the primal one by the product ${\mathbb M}^{-1}_{\Gamma_1}{\bf r}_{\Gamma_1}^{k+1}(\boldsymbol{\mu})$; 
\item interpolating on the primal space by applying ${\mathbb R}_{12}$; 
\item coming back to the dual space by applying the mass matrix ${\mathbb M}_{\Gamma_2}$. 
\end{enumerate}
This procedure is the one employed by the INTERNODES method \cite{Gervasio2019}. Indeed, INTERNODES is a general--purpose method to deal with non--conforming discretizations
of partial differential equations on 2D and 3D regions partitioned into two (or several) disjoint subdomains. Differently from Mortar methods, which are based on projection,
the idea of INTERNODES consists in exchanging the information between the two subdomains by resorting to both two independent interpolation operators ${\mathbb R}_{12}$ and ${\mathbb R}_{21}$, and two local mass matrices at the interfaces ${\mathbb M}_{\Gamma_1}$ and ${\mathbb M}_{\Gamma_2}$. To transfer the trace ${\bf u}_{\Gamma_2}(\boldsymbol{\mu})$ of the solution from $\Gamma_2$ to $\Gamma_1$, only the interpolation operator ${\mathbb R}_{12}$ is needed. On the contrary, for what concerns the balance of the residuals (which are the algebraic counterparts of the fluxes at the interfaces) the rule ${\bf r}_{\Gamma_2}^{k+1}(\boldsymbol{\mu})=-{\mathbb R}_{12} {\bf r}_{\Gamma_1}^{k+1}(\boldsymbol{\mu})$ would lead to a suboptimal method (see, e.g., \cite{bbdmkmp,Bernardi1994}). To recover the optimal convergence rate, INTERNODES employs \eqref{eq:slave_steady_algebraic} instead of ${\bf r}_{\Gamma_2}^{k+1}(\boldsymbol{\mu})=-{\mathbb R}_{12} {\bf r}_{\Gamma_1}^{k+1}(\boldsymbol{\mu})$. 
Further details on the derivation of these systems,  which are indeed quite standard in the DD literature, can be found, e.g.,  in  \cite{Gervasio2019,QuarteroniValli1999}. 

The residual vector is the algebraic counterpart of an element of the dual space 
$Y'_{h_i}$ of $Y_{h_i}$ -- see, e.g., \cite[Chapter 3]{QuarteroniManzoniNegri2016} -- then we define by $\mathbf{z}_{\Gamma_i}(\boldsymbol{\mu})$ the element obtained from $\mathbf{r}_{\Gamma_i}(\boldsymbol{\mu})$ by solving  
\begin{equation}
\label{eq:def_dual_residual}
\mathbb{M}_{\Gamma_i} \mathbf{z}_{\Gamma_i}(\boldsymbol{\mu}) = \mathbf{r}_{\Gamma_i}(\boldsymbol{\mu}).
\end{equation}
Here $\mathbf{z}_{\Gamma_i}(\boldsymbol{\mu})$ is the algebraic counterpart of the Riesz' element associated with the residual $\mathbf{r}_{\Gamma_i}(\boldsymbol{\mu})$. In other words, 
 the interface mass matrix becomes the transfer matrix from the Lagrange basis to the dual one and vice--versa \cite{Brauchli1971,Gervasio2019} and the array $\mathbf{z}_{\Gamma_i}(\boldsymbol{\mu})$ represents the residual in primal form. 

Therefore, the matrix $\mathbb{R}_{21}$ transfers the function of $Y_{h_1}$ whose nodal values are stored in $\mathbf{z}_{\Gamma_1}$ 
(corresponding to the interface residual vectors $\mathbf{r}_{\Gamma_1}$) at the nodes on $\Gamma_2$. 
%
%
Note that the conforming interface case can be recovered by taking $\mathbb{R}_{12}$ and $\mathbb{M}_{\Gamma_2}\mathbb{R}_{21}\mathbb{M}_{\Gamma_1}^{-1}$ both equal to the identity matrix. \smallskip

\begin{remark}
	When $\partial\Gamma_i \cap \partial \Omega_{i,D} \not = \emptyset$, the residual $\mathbf{r}_{\Gamma_i}$ should be corrected to take into account the interpolation process on all the degrees of freedom of $\Gamma_i$, including those on $\partial \Gamma_i$ (see \emph{e.g.} \cite{Gervasio2019}). Even if the reduced technique presented in this paper will work in both cases, hereon we will consider only $\partial\Gamma_i \cap \partial \Omega_{i,D} = \emptyset$.\end{remark}

\smallskip 

\begin{remark}
	The above formulation can be easily extended to time--dependent second--order parabolic PDE problems. In such cases, suitable numerical schemes have to be implemented to handle the time discretization, and Dirichlet--Neumann subdomains iterations must be applied for each time step of the approximated solution \cite{QuarteroniValli1999}. The application of our method to a time--dependent test case will be addressed in Section \ref{sec:numerical_results}. 
\end{remark}

\begin{figure}[!h]
\centering
\includegraphics[width=0.8\textwidth]{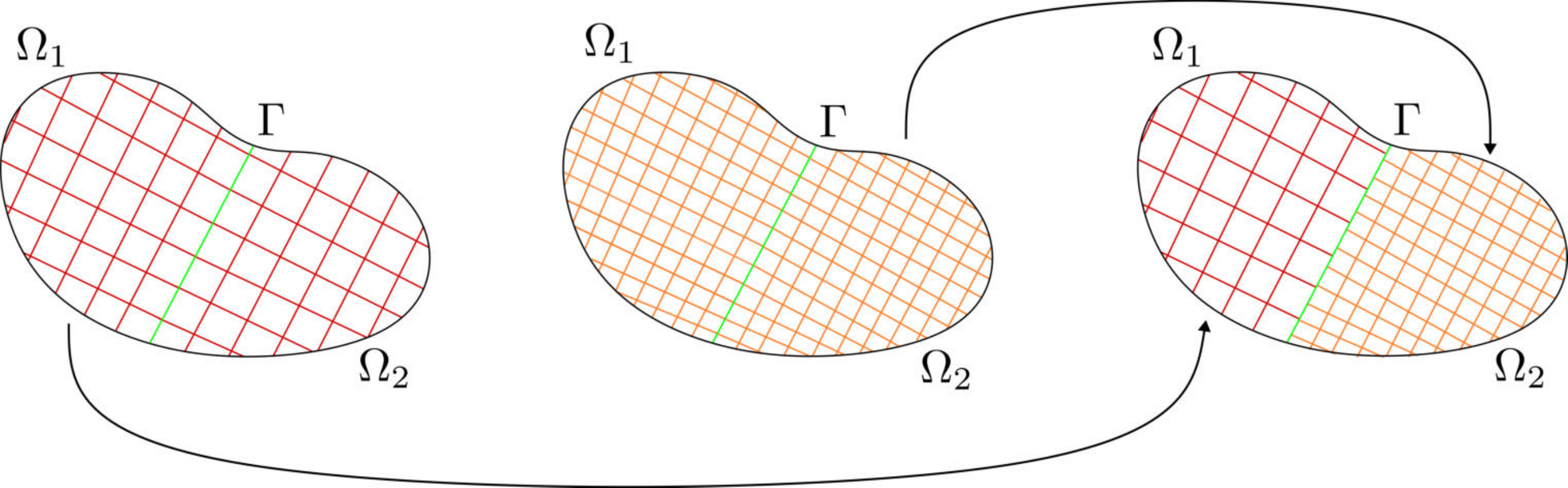}
\caption{Schematic representation of the two discretizations of the domain $\Omega$ used to compute the FOM snapshots (left and center) and the discretization of the domain $\Omega$ used to compute the ROM snapshots (right).}
\label{Fig:domain_representation}
\end{figure}

\begin{remark}
	\label{rem:complete_ROM}
	The master and the slave solution snapshots can be directly collected from the FOM computations by solving \eqref{eq:master_steady_algebraic}--\eqref{eq:interface_residual_def_steady} when \emph{(i)} conforming discretizations are considered in the two subdomains or \emph{(ii)} when interpolation/projection methods are implemented to handle non--conforming grids, \emph{e.g.} MORTAR methods or INTERNODES. However, since their implementation is not trivial, both MORTAR and INTERNODES methods are not easy to find in standard scientific libraries. Moreover, RB methods rely on snapshot data, whose generation can be performed with any desired numerical methods. For the sake of generality, here we assume to not have either INTERNODES or MORTAR method implemented and therefore compute the snapshots in an alternative -- perhaps naive -- way by solving for each parameter instance the FOM problem twice, \emph{i.e.} with two different conforming discretizations. In particular, we define two possible FE grids in the global domain by setting on each first subdomain a chosen spatial discretization and extending it (conformingly) to the corresponding second one. Now, both coupled problems feature interface conformity and can be solved with Dirichlet--Neumann iterations. Then, we collect the snapshots of each sub--problem and the relative interface data in the discretization set in the non--conforming case for the corresponding subdomain, obtaining two sets of solution and their Dirichlet and Neumann data. Indeed, these snapshots can be seen as the model solution when interface non--conformity is considered (see Fig. \ref{Fig:domain_representation} for a schematic sketch of the used procedure). 
	In this paper, such a technique is used to collect the snapshots for the test cases reported in Section \ref{sec:numerical_results}.  
\end{remark}

\section{Master and slave reduced order problems}
\label{sec:reduced_order_formulation}
The strategy we propose here aims both to reduce the two sub--problems separately and to employ reduced techniques to also address the Dirichlet and Neumann interface conditions arising from the Dirichlet--Neumann subdomains iterations. In particular, this ROM  technique is an extension of the one proposed in \cite{Zappon2022}, and combines different RB methods, one set for each sub--problem, and the DEIM to treat both Dirichlet and Neumann interface conditions, thus defining independent reduced order representations of the involved quantities.

\begin{figure}[!h]
	\centering
	\includegraphics[width=1\textwidth]{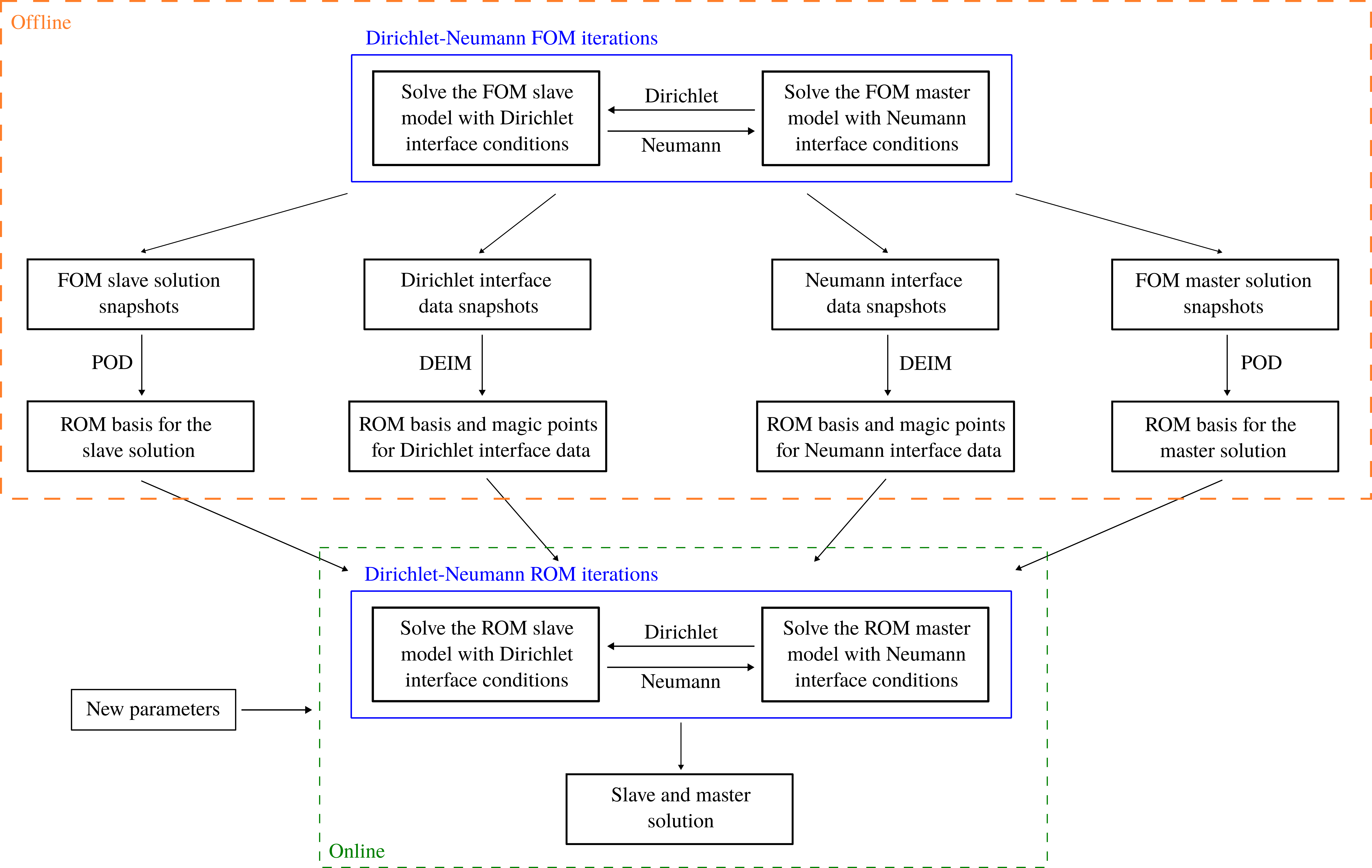}
	\caption{Schematic representation of the reduced order Dirichlet--Neumann domain decomposition algorithm.}
	\label{Fig:schemes_algorithm}
\end{figure}

We first approximate the FOM solution of the master and slave models by means of a POD--derived small number of basis functions defined in the corresponding subdomain $\Omega_i$. Moreover, employing the DEIM, we identify a suitable set of basis functions for the master and slave interface snapshots, and we use them to transfer Dirichlet and Neumann data across conforming or non--conforming interface grids. Lastly, considering the same Dirichlet--Neumann iteration scheme of the high--fidelity FOM, we iterate between the reduced solutions of the two sub--problems by imposing the continuity of both the interface solutions and fluxes at each iteration (see Fig. \ref{Fig:schemes_algorithm}). 

The reduced forms of the master and the slave problems are described in this section, while we derive the procedure to reduce parameter--dependent Dirichlet and Neumann interface conditions in Section \ref{sec:interface_reduction}.
	
We define the reduced--order version of problems $\eqref{eq:master_steady_algebraic}$ and $\eqref{eq:slave_steady_algebraic}$ relying on a POD--Galerkin approach \cite{QuarteroniManzoniNegri2016}. Therefore, in the offline stage, we collect the set of snapshots solving the sub--FOMs for a suitable set of parameter values. In particular, we choose as snapshots the FOM slave and master solutions \emph{at convergence} of the sub--iterations,  \emph{i.e.} $\mathbf{S}_1 = \{\mathbf{u}_1(\boldsymbol{\mu}_{\ell}),\boldsymbol{\mu}_{\ell} \in \mathscr{P}^d_{train}\}$ and $\mathbf{S}_2 = \{\mathbf{u}_2(\boldsymbol{\mu}_{\ell}),\boldsymbol{\mu}_{\ell} \in \mathscr{P}^d_{train} \}$, respectively. A sampling of the parameter space is here done considering a Latin hypercube sampling (LHS) method \cite{Iman2006,Mckay1979}. 

\begin{remark}
 From now on, the index $k$ is omitted when quantities at the convergence of the Dirichlet--Neumann FOM iterations are considered.
\end{remark}

The POD technique is applied to each set of snapshots $\mathbf{S}_1$ and $\mathbf{S}_2$ and a corresponding set of reduced basis functions is computed and stored for the approximation of the solution on each subdomain. We denote by $n_i$ the cardinality of the set of reduced basis functions. Defining 
$\mathbb{V}_i \in \mathbb{R}^{\tilde{N}_i \times n_i}$, $n_i \ll  \tilde{N}_i$, for $i=1,2$,
the matrices whose columns yield the obtained basis functions, the ROM seeks an approximation of the FOM solutions under the form
\begin{equation}
\label{eq:FOM_rec_1}
\mathbf{u}_1(\boldsymbol{\mu}) \approx \mathbb{V}_1 \mathbf{u}_{n,1}(\boldsymbol{\mu})
\end{equation}
and 
\begin{equation} 
\label{eq:FOM_rec_2}\mathbf{u}_2(\boldsymbol{\mu}) \approx \mathbb{V}_2 \mathbf{u}_{n,2}(\boldsymbol{\mu}).\end{equation}

Projecting problems $\eqref{eq:master_steady_algebraic}$ and $\eqref{eq:slave_steady_algebraic}$ onto the reduced spaces defined by $\mathbb{V}_i$, starting from an initial guess $\mathbf{u}_{n,i}^0(\boldsymbol{\mu})$ and $\mathbf{u}_{n,2}^0(\boldsymbol{\mu})$, in the online phase, for each $k \geq 0$, we search for the reduced solutions $\mathbf{u}_{n,1}^{k+1}(\boldsymbol{\mu}) \in \mathbb{R}^{n_1}$ and $\mathbf{u}_{n,2}^{k+1}(\boldsymbol{\mu})\in \mathbb{R}^{n_2}$ such that 

\begin{equation}
\label{eq:master_reduced_problem}
\begin{cases}
\mathbb{A}_{n,1}(\boldsymbol{\mu}) \mathbf{u}_{n,1}^{k+1}(\boldsymbol{\mu}) = \mathbf{f}_{n,1}(\boldsymbol{\mu}) - \mathbb{V}_1^T\mathbb{A}_{1,D}(\boldsymbol{\mu})\mathbf{g}_{D,1}(\boldsymbol{\mu}) - \mathbb{V}_1^T\mathbb{A}_{1,\Gamma_1}(\boldsymbol{\mu})\mathbf{u}_{\Gamma_1}^{k+1}(\boldsymbol{\mu}) \\
\mathbf{u}_{\Gamma_1}^{k+1}(\boldsymbol{\mu}) = \mathbb{R}_{12}\mathbf{u}_{\Gamma_2}^{k}(\boldsymbol{\mu})
\end{cases}
\end{equation}
and 
\begin{equation}
\label{eq:slave_reduced_problem}
\begin{cases}
\mathbb{A}_{n,2}(\boldsymbol{\mu})  \mathbf{u}_{n,2}^{k+1}(\boldsymbol{\mu}) = \mathbf{f}_{n,2}(\boldsymbol{\mu}) +\mathbb{V}_2^T{\mathbb E}_2\mathbf{r}_{\Gamma_2}^{k+1}(\boldsymbol{\mu})-\mathbb{V}_2^T\mathbb{A}_{2,D}(\boldsymbol{\mu})\mathbf{g}_{D,2}(\boldsymbol{\mu})\\
\mathbf{r}_{\Gamma_2}^{k+1}(\boldsymbol{\mu}) = -\mathbb{M}_{\Gamma_2}\mathbb{R}_{21}\mathbb{M}_{\Gamma_1}^{-1} \mathbf{r}_{\Gamma_1}^{k+1}(\boldsymbol{\mu}),
\end{cases}
\end{equation}
where 
\[
\mathbb{A}_{n,i}(\boldsymbol{\mu})  = \mathbb{V}_i^T\mathbb{A}_{i,i}(\boldsymbol{\mu})\mathbb{V}_i, \qquad \mathbf{f}_{n,i}(\boldsymbol{\mu}) = \mathbb{V}_i^T\mathbf{f}_i(\boldsymbol{\mu}), \qquad  i=1,2.
\]
Notice that the second equations in problems \eqref{eq:master_reduced_problem} and \eqref{eq:slave_reduced_problem} (the interface equations) are defined in the FOM space, whereas the first ones are in the ROM space. The reduced version of the interface equations will be derived in the following section.

\begin{remark}
For simplicity, in this paper, we only consider the case of linear PDE problems. In the case of non--linear problems, the presence of non--linear terms in the master and slave formulations can be handled through suitable hyper--reduction techniques like, \emph{e.g}, the DEIM \cite{Barrault2004,Chaturantabut2009,Chaturantabut2010}. The ROM approach would be still modular, requiring an \emph{ad--hoc} reduction of each sub--problem, whereas the processing of interface data would not be affected by the operator non--linearity, and would be treated with the approach shown in the next section.
\end{remark}

\begin{remark}
	Time--dependent problems can be reduced with RB methods considering the time variable as an additional parameter of the model. Indeed, in such a case the FOM solutions at each time step of the simulation are collected in the snapshots set, and the reduced basis allows to approximate the time--dependent solution of each sub--problem involving vectors of (reduced) degrees of freedom  $\mathbf{u}_{n,1}$, $\mathbf{u}_{n,2}$ that are also time--dependent (see Section \ref{sec:numerical_results}).
\end{remark}

\section{Parametric interface data reduction}
\label{sec:interface_reduction}

Dealing with interface conditions, especially when using non--conforming grids, requires special care. Since the sub--problems \eqref{eq:master_reduced_problem} and \eqref{eq:slave_reduced_problem} are parameter--dependent, the interface data naturally inherit the parameter dependency and the DEIM \cite{Barrault2004,Chaturantabut2009,Chaturantabut2010,farhat20205,Grepl2007,Maday2008,Negri2015} can be applied to both reduce the dimension of such data and to transfer the information across the interface grids. 

Indeed, using DEIM requires: \emph{(i)} to compute a set of basis functions for the quantity of interest employing POD, \emph{(ii)} to use a greedy algorithm to identify a small number of DoFs to compute the weights for the corresponding basis functions (instead of the weights used in a {\em simple}  POD). The nodes corresponding to such DoFs define the so--called reduced mesh, effectively determining a relation between the FE grids and the reduced space. Therefore, when considering interface data reduction, a small number of interface nodes can be selected through the DEIM to describe the complete vector of parametrized interface data. In the conforming case, DEIM can be used directly on the quantity of interest, \emph{i.e.} the interface solution in the case of Dirichlet data, or the interface residual in the case of Neumann data. For non--conforming interface grids, the primal form of the residuals must be employed to properly treat the Neumann terms \cite{Oliver2009}. 

Even if using the DEIM at the interface between the two subdomains could seem a rather involved approach compared to more classic interpolation algorithms -- like, e.g., piecewise--constant interpolation -- DEIM offers two main advantages. First, in this context, DEIM also aims at approximating interface data that could in principle depend on parameters unrelated to the problem solution. Therefore, DEIM adjusts to variations in the interface data given to changes in the interface parametrization, independently of the model solution. Second, DEIM can be applied to larger domains with a substantially higher number of DoFs at the interface, compared to the one employed in our test cases, only involving a small subset of the interface DoFs, thus significantly accelerating the overall interpolation process.

In Subsection \ref{subsec:Dirichlet_reduction} and \ref{subsec:Neumann_reduction}, we define the reduction of the Dirichlet and Neumann interface conditions, respectively, when non--conforming interface grids are considered. These new ROM interface conditions are then used to replace the interface equations of problems \eqref{eq:master_reduced_problem} and \eqref{eq:slave_reduced_problem}.

\subsection{Parameter--dependent Dirichlet data}
\label{subsec:Dirichlet_reduction}

The parametric interpolation method of the Dirichlet data used in this work is similar to the one introduced in \cite{Zappon2022}. Such a technique relies on the DEIM and can be applied in the case of both conforming and non--conforming interface grids. 

First, in the offline phase, we collect from the slave domain $\Omega_1$ the interface snapshots, \emph{i.e.}, we extract the interface (Dirichlet) degrees of freedom obtained for different instances of the parameter vector from the FOM computation. Notice that, as for the solution reduction, we  only select the interface DoFs at the convergence of the FOM Dirichlet--Neumann iterations, namely
$$ \mathbf{S}_D = \{\mathbf{u}_{\Gamma_1}(\boldsymbol{\mu}_{\ell}), ~\boldsymbol{\mu}_{\ell} \in \mathscr{P}^d_{train} \}.$$

Let us denote by $N_{1,\Lambda}$ the number of FOM DoFs on $\Gamma_1$. A low--dimensional representation of the interface DoFs can then be computed by determining a set of $M_1 \ll N_{1,\Lambda}$ POD basis functions from $\mathbf{S}_D$ that we store in the matrix $\boldsymbol{\Phi}_D\in {\mathbb R}^{N_{1,\Lambda}\times M_1}$, with the purpose of getting
$$ \mathbf{u}_{\Gamma_1}(\boldsymbol{\mu}) \approx \boldsymbol{\Phi}_D \mathbf{u}_{1,M_1}(\boldsymbol{\mu}),$$
where $\mathbf{u}_{1,M_1}(\boldsymbol{\mu})$ is a vector of $M_1$ coefficients. Furthermore, with a greedy algorithm \cite{Maday2008}, we select iteratively $M_1$ indices in $\{1,\dots,N_{1,\Lambda}\}$, by minimizing the interpolation error over the interface snapshots set $\mathbf{S}_D$, according to the maximum norm. The set of such indices is denoted by 
\begin{equation}
\label{eq:set_indices_Dirichlet}
\mathcal{I}_{1,D} \subset \{1,\dots,N_{1,\Lambda}\}, \quad \text{with cardinality } |\mathcal{I}_{1,D}| = M_1.
\end{equation}
The points corresponding to the indices of $\mathcal{I}_{1,D}$ are usually referred to as \emph{magic points} on $\Gamma_1$, and are used to impose the Dirichlet interface conditions on $\Gamma_1$ for the reduced online problem. Let us denote by $\mathbf{u}_{\Gamma_{1,|\mathcal{I}_{1,D}}}(\boldsymbol{\mu})$ the vector of the FOM DoFs at the magic points. 

In the online phase, at each Dirichlet--Neumann iteration $k$, we ask that the reduced interface vector  $\mathbf{u}_{1,M_1}^{k+1}(\boldsymbol{\mu})$ satisfies the relation
$$ \boldsymbol{\Phi}_{D_{|\mathcal{I}_{1,D}}} \mathbf{u}_{1,M_1}^{k+1}(\boldsymbol{\mu}) = \mathbf{u}_{\Gamma_{1_{|\mathcal{I}_{1,D}}}}^{k+1}(\boldsymbol{\mu}),$$ 
where $\boldsymbol{\Phi}_{D_{|\mathcal{I}_{1,D}}} \in \mathbb{R}^{M_1 \times M_1}$ is the sub--matrix of $\boldsymbol{\Phi}_{D}$ containing the $\mathcal{I}_{1,D}$ rows (see \emph{e.g.} \cite{Chaturantabut2010} for the well--posedness of the above procedure). The FOM interface DoFs on $\Gamma_1$ can be then approximated as
\begin{equation}
\label{eq:approx_dirichlet_data_1}
\mathbf{u}_{\Gamma_1}^{k+1}(\boldsymbol{\mu}) \approx \boldsymbol{\Phi}_D\boldsymbol{\Phi}^{-1}_{D_{|\mathcal{I}_{1,D}}}\mathbf{u}_{\Gamma_{1_{|\mathcal{I}_{1,D}}}}^{k+1}(\boldsymbol{\mu}).\end{equation}

Now we replace $\mathbf{u}_{\Gamma_{1_{|\mathcal{I}_{1,D}}}}^{k+1}$ with the values of the master solution $\mathbf{u}_{N,2}^k$ at the points on $\Gamma_2$ corresponding to the magic points on $\Gamma_1$. Thus, given the position $\mathbf{p}_1$ of the magic point corresponding to the index $i_{1,D} \in \mathcal{I}_{1,D}$ in Cartesian coordinates, we search for the corresponding node in the master interface $\Gamma_2$, \emph{i.e.} for the point $\mathbf{p}_2 \in \Gamma_2$ such that
$$\mathbf{p}_2 = \underset{\mathbf{p}_2^j \in \text{nodes}_{\Gamma_2}}{\arg\min}(\text{dist}(\mathbf{p}_1 - \mathbf{p}_2^j)),$$
where $\text{dist}$ represents the Euclidean distance.

\begin{remark}
When the solution of the above minimization problem is not unique, we choose as $\mathbf{p}_2$ the last node found by the algorithm that satisfies such a relation. 
\end{remark}

Then, we can define the set $\mathcal{I}_{2,D}$ of the indices on the master grid corresponding to the indices in $\mathcal{I}_{1,D}$, \emph{i.e.} 
$$\mathcal{I}_{2,D} = \{i_{2,D}(i_{1,D})\}_{i_{1,D} \in \mathcal{I}_{1,D}}.$$
Notice that $\mathcal{I}_{2,D}$ is computed in the offline phase.

Finally, in the online phase, we replace the FOM interface values at the magic points on $\Gamma_1$ with the values of the master solution at the corresponding points on $\Gamma_2$, \emph{i.e.},
\begin{equation}
\label{eq:approx_dirichlet_data_2}
\left( \mathbf{u}_{\Gamma_1}^{k+1}(\boldsymbol{\mu}) \right)_{i_{1,D}} = \left( \mathbf{u}_{\Gamma_2}^{k}(\boldsymbol{\mu}) \right)_{i_{2,D}(i_{1,D})}, \quad i_{1,D} \in \mathcal{I}_{1,D}.
\end{equation}
More briefly, we write
$\mathbf{u}_{\Gamma_{1_{|\mathcal{I}_{1,D}}}}^{k+1}(\boldsymbol{\mu}) = \mathbf{u}_{\Gamma_{2_{|\mathcal{I}_{2,D}}}}^{k}(\boldsymbol{\mu})$, so that \eqref{eq:approx_dirichlet_data_1} becomes
$$\mathbf{u}_{\Gamma_1}^{k+1}(\boldsymbol{\mu}) \approx \boldsymbol{\Phi}_D\boldsymbol{\Phi}^{-1}_{D_{|\mathcal{I}_{1,D}}}\mathbf{u}_{\Gamma_{2_{|\mathcal{I}_{2,D}}}}^{k}(\boldsymbol{\mu}).$$

\begin{figure}[h!]
	\centering
	\includegraphics[width=0.6\textwidth]{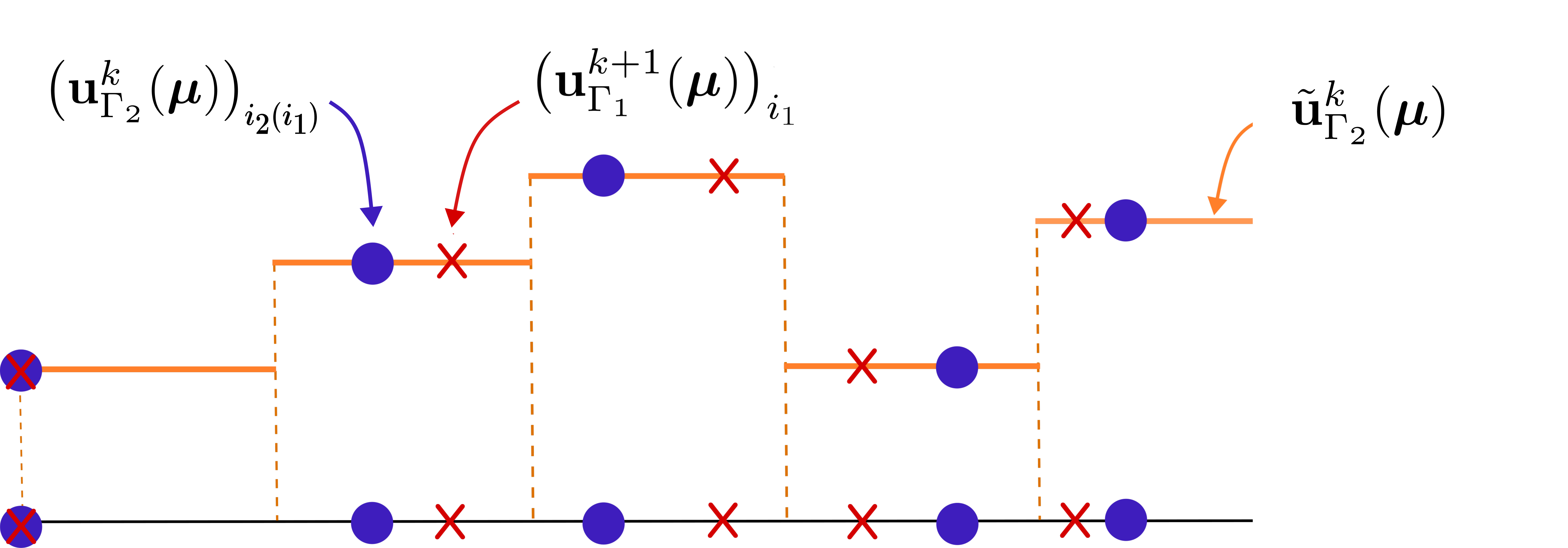}
	\caption{Geometrical interpretation of the interface reduction with the DEIM. The blue points represent the magic points on $\Gamma_2$ and the values of $\mathbf{u}_{\Gamma_2}^k$ at these magic points, the orange lines represent the piecewise constant interpolating function $\tilde{\mathbf{u}}_{\Gamma_2}^k$, while the red crosses are the points on $\Gamma_1$ corresponding to the magic points of $\Gamma_2$ and the values of $\mathbf{u}_{\Gamma_1}^{k+1}$ at these points.}
	\label{Fig:interpretation_interpolation}
\end{figure}

\begin{remark}
	\label{rem:interporetation_DEIM_interface}
	The substitution \eqref{eq:approx_dirichlet_data_2} can be interpreted as a low--order interpolation process: first, we build the piecewise constant function $\tilde{\mathbf{u}}_{\Gamma_2}^k(\boldsymbol{\mu})$ (the orange one in Fig. \ref{Fig:interpretation_interpolation}) that interpolates the values of $\mathbf{u}_{\Gamma_2}^k(\boldsymbol{\mu})$ at the magic points on $\Gamma_2$ (the blue dots); then the values $\mathbf{u}_{\Gamma_1}^{k+1}(\boldsymbol{\mu})$ at the corresponding points of $\Gamma_1$ (the red symbols) are obtained by evaluating the function $\tilde{\mathbf{u}}_{\Gamma_2}^k(\boldsymbol{\mu})$ at such points.
\end{remark}

Note that $\mathbf{u}_{\Gamma_{2_{|\mathcal{I}_{2,D}}}}^{k}(\boldsymbol{\mu})$ refers to the approximation of the FOM solution of the master problem that must, therefore, be computed from the ROM solution $\mathbf{u}_{n,2}^{k}(\boldsymbol{\mu})$ during the online phase. However, only a part of the approximated FOM master solution is needed, i.e., the one at the magic points on $\Gamma_2$. Therefore, one can store the FOM solution at the interface only at the magic points by pre--multiplying the ROM master solution by those rows of $\mathbb{V}_2$ corresponding to the magic points. To this end, we introduce a restriction matrix $\mathbb{U}_{12} \in \mathbb{R}^{M_1 \times \tilde N_2}$ whose only entries different from zero (and equal to one) are those whose column index belongs to $\mathcal{I}_{2,D}$, so that
$$ \mathbf{u}_{\Gamma_{2_{|\mathcal{I}_{2,D}}}}^k(\boldsymbol{\mu}) = \mathbb{U}_{12} \mathbb{V}_2 \mathbf{u}_{n,2}^k(\boldsymbol{\mu}).$$ 

Therefore, the action of the operator $\mathbb{R}_{12}$ in \eqref{eq:master_reduced_problem} (in \emph{both} conforming and non--conforming cases), can be summarized by
$$ \mathbb{R}_{12}\mathbf{u}_{\Gamma_{2}}^{k}(\boldsymbol{\mu}) = \boldsymbol{\Phi}_D\boldsymbol{\Phi}^{-1}_{D_{|\mathcal{I}_{1,D}}}\mathbb{U}_{12} \mathbb{V}_2\mathbf{u}_{n,2}^k(\boldsymbol{\mu}),$$
and the interface condition in  \eqref{eq:master_reduced_problem} becomes 
$$ \mathbf{u}_{\Gamma_1}^{k+1}(\boldsymbol{\mu}) =  \boldsymbol{\Phi}_D\boldsymbol{\Phi}^{-1}_{D_{|\mathcal{I}_{1,D}}}\mathbb{U}_{12} \mathbb{V}_2\mathbf{u}_{n,2}^k(\boldsymbol{\mu}).$$
In this way, the last term of  \eqref{eq:master_reduced_problem} can be approximated as
\begin{equation}
\label{eq:lifting_term_DEIM}
\mathbb{V}_1^T\mathbb{A}_{1,\Gamma_1}(\boldsymbol{\mu})\mathbf{u}_{\Gamma_1}^{k+1}(\boldsymbol{\mu}) = \mathbb{V}_1^T\mathbb{A}_{1,\Gamma_1}(\boldsymbol{\mu})\mathbb{R}_{12}\mathbf{u}_{\Gamma_2}^{k}(\boldsymbol{\mu}) = \mathbb{V}_1^T\mathbb{A}_{1,\Gamma_1}(\boldsymbol{\mu}) \boldsymbol{\Phi}_D\boldsymbol{\Phi}^{-1}_{D_{|\mathcal{I}_{1,D}}}\mathbb{U}_{12} \mathbb{V}_2\mathbf{u}_{n,2}^{k}(\boldsymbol{\mu}),
\end{equation}
 where the matrix $\boldsymbol{\Phi}_D\boldsymbol{\Phi}^{-1}_{D_{|\mathcal{I}_{1,D}}}\mathbb{U}_{12} \mathbb{V}_2$ does not depend on the solution and can be pre--computed and stored in the offline phase. Note that if $\mathbb{A}$ is parameter independent or depends affinely on $\boldsymbol{\mu}$, in the online phase we can precompute and store the entire product $\mathbb{V}_1^T\mathbb{A}_{1,\Gamma_1}\boldsymbol{\Phi}_D\boldsymbol{\Phi}^{-1}_{D_{|\mathcal{I}_{1,D}}}\mathbb{U}_{12} \mathbb{V}_2$.

Note that even if the interface snapshots stored in $\boldsymbol{\Phi}_D$ contain the DoFs on $\Gamma_2$ at the convergence of the Dirichlet--Neumann iterations, the reduced coupled model is solved by iterating between the reduced master and slave models. For this reason, the index $k$ appears on the quantities computed online. 

\begin{remark}
As for the FOM computation, an initial guess of the Dirichlet interface conditions must be considered, but only at those points on $\Gamma_2$ corresponding to the magic points on $\Gamma_1$. Therefore, for $k = 0$, the approximated FOM solution $\mathbb{U}_{12} \mathbb{V}_2\mathbf{u}_{n,2}^{0}$ 
can be replaced with the FOM initial guess $\mathbf{u}^{0}_{\Gamma_{2|\mathcal{I}_{2,D}}}$.
\end{remark}

\begin{remark}
Note that if the coupled problem is unsteady, to take into account the time variations of the solution, the interface data at the convergence of the Dirichlet--Neumann iterations for each time instant $n = 1,\dots, N_t$ must be collected in the set of snapshots. Once the POD basis has been computed on this set of snapshots, the interpolation of the Dirichlet data can be performed as in the steady case. 
\end{remark}

\subsection{Parameter--dependent Neumann data}
\label{subsec:Neumann_reduction}
The DEIM used to interpolate the parametric Dirichlet interface conditions can be also applied to the parametric Neumann interface conditions; as before, let us detail the case of a steady problem. 
If the interface grids are conforming, for each $k \geq 0$,
$$\mathbf{r}_{\Gamma_2}^{k+1}(\boldsymbol{\mu}) = - \mathbf{r}_{\Gamma_1}^{k+1}(\boldsymbol{\mu})$$
so that the DEIM can be used on the interface residual. Instead, in the non--conforming case, the interface mass matrices are involved (see Section \ref{sec:Two_way_models}). However, recalling definition $\eqref{eq:def_dual_residual}$ of $\mathbf{z}_{\Gamma_i}(\boldsymbol{\mu})$, \emph{i.e.} the algebraic counterpart of the Riesz' element associated with the residual, the second equation of \eqref{eq:slave_reduced_problem} can be replaced by 
\begin{equation} \label{eq:dual_interpolation}\mathbf{z}_{\Gamma_2}^{k+1}(\boldsymbol{\mu}) =- \mathbb{R}_{21}\mathbf{z}_{\Gamma_1}^{k+1}(\boldsymbol{\mu}).\end{equation}
Therefore, the vector $\mathbf{z}_{\Gamma_2}^{k+1}$ is the quantity to be reduced and reconstructed using the DEIM. 

For the sake of generality, in this subsection we derive the Neumann data approximation for the non--conforming case, however, the same procedure also holds for the conforming case. 
Starting from \eqref{eq:interface_residual_def_steady}--\eqref{eq:def_dual_residual}, for each $\boldsymbol{\mu}$ we compute the interface residual snapshots, 
$$\mathbf{S}_{N} = \{ \mathbf{z}_{\Gamma_2}(\boldsymbol{\mu}_{\ell}), \boldsymbol{\mu}_{\ell} \in \mathscr{P}^d_{train}\};$$
we remind that these snapshots are saved at the convergence of the FOM Dirichlet--Neumann iterations.

Applying POD, a set of $M_2$ basis functions is selected and stored in $\boldsymbol{\Phi}_N\in{\mathbb R}^{N_{2,\Lambda} \times M_2}$, being $M_2 \ll N_{2,\Lambda}$ (recall that $ N_{2,\Lambda}$ is the number of mesh points on $\Gamma_2$). Given a generic $\boldsymbol{\mu}$, the vector $\mathbf{z}_{\Gamma_2}$ is, therefore, approximated as
$$\mathbf{z}_{\Gamma_2}(\boldsymbol{\mu}) = \boldsymbol{\Phi}_N \mathbf{z}_{2,M_2}(\boldsymbol{\mu}),$$
where $\mathbf{z}_{2,M_2}$ is a vector of $M_2$ coefficients.
Then, $M_2$ magic points on $\Gamma_2$ are selected through a greedy algorithm and their indices in the master grid numbering are collected in the set
$$ \mathcal{I}_{2,N} \subset\{1,\dots, N_{2,\Lambda}\}, \quad \text{with cardinality } |\mathcal{I}_{2,N}| = M_2.$$ In the online phase, at each Dirichlet--Neumann iteration $k$, we need to find $\mathbf{z}_{2,M_2}^{k+1}(\boldsymbol{\mu})$ such that
$$ \boldsymbol{\Phi}_{N_{|\mathcal{I}_{2,N}}} \mathbf{z}_{2,M_2}^{k+1}(\boldsymbol{\mu}) = \mathbf{z}_{\Gamma_{2|\mathcal{I}_{2,N}}}^{k+1}(\boldsymbol{\mu}),$$
where $\boldsymbol{\Phi}_{N_{|\mathcal{I}_{2,N}}}\in{\mathbb R}^{M_2\times M_2}$ is the restriction of $\boldsymbol{\Phi}_{N}$ to the indices associated with the magic points identified by $\mathcal{I}_{2,N}$, while $\mathbf{z}_{\Gamma_{2|\mathcal{I}_{2,N}}}^{k+1}(\boldsymbol{\mu})$ is the restriction of $\mathbf{z}_{\Gamma_2}^{k+1}(\boldsymbol{\mu})$ at such magic points. Then, we can approximate
$$\mathbf{z}_{\Gamma_2}^{k+1}(\boldsymbol{\mu}) \approx \boldsymbol{\Phi}_N \boldsymbol{\Phi}_{N_{|\mathcal{I}_{2,N}}}^{-1} \mathbf{z}_{\Gamma_{2|\mathcal{I}_{2,N}}}^{k+1}(\boldsymbol{\mu}).$$

Now, $\mathbf{z}_{\Gamma_{2|\mathcal{I}_{2,N}}}^{k+1}(\boldsymbol{\mu})$ is replaced by the values of $\mathbf{z}_{\Gamma_1}^{k+1}(\boldsymbol{\mu})$ extracted at the points on $\Gamma_1$ corresponding to the magic points on $\Gamma_2$. As done in the previous subsection, to find such points on $\Gamma_1$, we first need to select the set of indices $\mathcal{I}_{1,N}$ on the slave interface $\Gamma_1$ corresponding to the magic points identified by $\mathcal{I}_{2,N}$. Therefore, denoting by $\mathbf{p}_2$ the node in cartesian coordinates corresponding to the index $i_{2,N}$ in $\mathcal{I}_{2,N}$, we search for 
$$ \mathbf{p}_1 = \underset{\mathbf{p}_1^j \in \text{nodes}_{\Gamma_1}}{\arg\min} (\text{dist}(\mathbf{p}_2-\mathbf{p}_1^j)),$$
where $\text{dist} (\cdot)$ denotes the Euclidean distance.
Then, denoting by $i_{1,N}(i_{2,N})$ the index of $\mathbf{p}_1$ in the numbering of the slave grid, we define the set of indices in the slave grid corresponding to the nodes defined by the indices in $\mathcal{I}_{2,N}$, \emph{i.e.}
$$ \mathcal{I}_{1,N} = \{i_{1,N}(i_{2,N})\}_{i_{2,N}\in \mathcal{I}_{2,N}}.$$

Thus, in the online phase, we impose that 
$$ \left( \mathbf{z}_{\Gamma_2}^{k+1}(\boldsymbol{\mu})\right)_{i_{2,N}} = \left( \mathbf{z}_{\Gamma_1}^{k+1}(\boldsymbol{\mu})\right)_{i_{1,N}(i_{2,N})}, \quad \text{with }i_{2,N}\in \mathcal{I}_{2,N}$$
or more briefly $\mathbf{z}_{\Gamma_{2|\mathcal{I}_{2,N}}}^{k+1}(\boldsymbol{\mu}) = \mathbf{z}_{\Gamma_{1|\mathcal{I}_{1,N}}}^{k+1}(\boldsymbol{\mu})$, so that
$$\mathbf{z}_{\Gamma_2}^{k+1}(\boldsymbol{\mu}) \approx \mathbb{R}_{21}\mathbf{z}_{\Gamma_1}^{k+1}(\boldsymbol{\mu}) \approx  \boldsymbol{\Phi}_N \boldsymbol{\Phi}_{N_{|\mathcal{I}_{2,N}}}^{-1} \mathbf{z}_{\Gamma_{1|\mathcal{I}_{1,N}}}^{k+1}(\boldsymbol{\mu})$$
and, similarly to the Dirichlet interpolation, the operator $\mathbb{R}_{21}$ here is represented by
$$\mathbb{R}_{21}\mathbf{z}_{\Gamma_1}^{k+1}(\boldsymbol{\mu}) = \boldsymbol{\Phi}_N \boldsymbol{\Phi}_{N_{|\mathcal{I}_{2,N}}}^{-1}\mathbf{z}_{\Gamma_{1|\mathcal{I}_{1,N}}}^{k+1}(\boldsymbol{\mu}).$$

Finally, recalling formula \eqref{eq:def_dual_residual}, we can recover the (Neumann) interface condition in the second equation of \eqref{eq:slave_reduced_problem} as  
\begin{equation}\label{eq:Neumm_ROM_term}  \mathbf{r}_{\Gamma_2}^{k+1}(\boldsymbol{\mu}) \approx  - \mathbb{M}_{\Gamma_2} \boldsymbol{\Phi}_N \boldsymbol{\Phi}_{N_{|\mathcal{I}_{2,N}}}^{-1}\mathbb{U}_{21}\mathbb{M}_{\Gamma_1}^{-1}\mathbb{U}_{21}^T\mathbf{r}_{\Gamma_{1|\mathcal{I}_{1,N}}}^{k+1}(\boldsymbol{\mu}),\end{equation}
with $\mathbb{U}_{21} \in \mathbb{R}^{M_2 \times N_{1,\Lambda}}$ the restriction matrix whose entries different from zero and equal to 1 are only those with column index in $\mathcal{I}_{1,D}$. 
Note that the matrix  $\mathbb{V}_2^T {\mathbb E}_2\mathbb{M}_{\Gamma_2} \boldsymbol{\Phi}_N \boldsymbol{\Phi}_{N_{|\mathcal{I}_{2,N}}}^{-1} \mathbb{U}_{21}\mathbb{M}_{\Gamma_{1}}^{-1}\mathbb{U}_{21}^T$ (that appears in \eqref{eq:slave_reduced_problem} to evaluate $\mathbb{V}_2^T{\mathbb E}_2\mathbf{r}_{\Gamma_2}^{k+1}(\boldsymbol{\mu})$) does not depend on $\boldsymbol{\mu}$ and can be therefore computed and stored offline.

Moreover, recalling the definition \eqref{eq:interface_residual_def_steady} of the interface residual, we can recover the dependency of such vectors on the reduced slave solution, computing $\mathbf{r}_{\Gamma_{1|\mathcal{I}_{1,N}}}^{k+1}$ as the vector of dimension $M_1$ whose only non--zero entries are those at the magic points, i.e., 
\begin{equation}\label{eq:master_residual_ROM} \mathbf{r}_{\Gamma_{1|\mathcal{I}_{1,N}}}^{k+1}(\boldsymbol{\mu}) \approx \left(\mathbb{A}_1(\boldsymbol{\mu})\begin{bmatrix}\mathbb{V}_1\mathbf{u}_{n,1}^{k+1}(\boldsymbol{\mu})\\
\mathbf{u}^{k+1}_{\Gamma_1}(\boldsymbol{\mu}) \\
\mathbf{g}_{D,1}(\boldsymbol{\mu})\\
\end{bmatrix} 
- \mathbf{f}_{N,1}
(\boldsymbol{\mu})\right)_{|\mathcal{I}_{1,N}}.
\end{equation}

We notice that the magic points on $\Gamma_1$ found during the Dirichlet phase do not necessarily coincide with the corresponding points on $\Gamma_2$ found  during the Neumann phase, i.e, ${\cal I}_{1,D}$ and ${\cal I}_{1,N}$ can differ, as well as 
 ${\cal I}_{2,D}$ and ${\cal I}_{2,N}$ can do.

We summarize the coupled problem training procedure in Algorithms \ref{alg:DEIMInterface}--\ref{alg:DEIMInterface_ROM_arrays}, including the interface DEIM reduction of both Dirichlet and Neumann processing, while the complete reduction of the two--way coupled model can be found in Algorithms \ref{alg:ROM_complete1}--\ref{alg:ROM_complete2} concerning the offline training and the online query of the ROM.

We remark that when the interface grids are conforming, a perfect matching between the corresponding nodes on the master and slave interface is found. However, this does not happen in the non--conforming case, where interpolation errors arise for both the Dirichlet and Neumann interface approximations, especially when the grids differ substantially. Considering the numerical tests of Section \ref{sec:numerical_results}, to minimize such errors a possible remedy is to consider a finer discretization on the master domain than in the slave one since the Dirichlet approximation seems to suffer more from the interface difference than the Neumann one. Another remedy consists in employing more accurate interpolation operators, like Lagrange interpolation when the interfaces are geometrically conforming, or Radial Basis functions interpolation in the presence of geometrical non--conforming interfaces; this will be the subject of future work.

Moreover, given the smaller number of DoFs in the slave interface than in the master one -- considering a coarser discretization in the slave domains, as done above -- we check  the continuity of the interface solution using the $\ell_2$ norm of the difference between the approximated solutions (expressed in a high--fidelity format) on the DoFs of the slave interface, \emph{i.e.} 
\begin{equation}
\label{eq:interface_continuity_check}
\|\mathbf{u}_{\Gamma_1}^{k+1}(\boldsymbol{\mu}) -{\mathbb R}_{12}\mathbf{u}_{\Gamma_2}^{k+1}(\boldsymbol{\mu})\|_{\ell_2} < \epsilon,
\end{equation}
where ${\mathbb R}_{12}\mathbf{u}_{\Gamma_2}^{k+1}(\boldsymbol{\mu})=\boldsymbol{\Phi}_D\boldsymbol{\Phi}^{-1}_{D_{|\mathcal{I}_{1,D}}}\mathbb{U}_{12} \mathbb{V}_2\mathbf{u}_{n,2}^{k+1}(\boldsymbol{\mu})$ effectively represents the interpolation of the master solution on the slave interface grids by means of the DEIM (see Subsection \ref{subsec:Dirichlet_reduction}).

\begin{remark} If the coupled problem is unsteady, similarly to the case of parameter--dependent Dirichlet data, to take into account the time variations of the solution, the interface data at the convergence of the sub--iterations for each time instant $n = 1,\dots,N_t$ must be collected in the set of snapshots. Once the POD basis has been computed on this set of snapshots, the interpolation of the Neumann data can be performed as in the steady case. 
\end{remark}

\begin{remark}
\label{Rem:RBF_formula}
The nearest neighbor interpolation adopted in exchanging the information across the interfaces in formulas \eqref{eq:lifting_term_DEIM} and \eqref{eq:Neumm_ROM_term} could be replaced by more sophisticated interpolation techniques, such as, e.g., that based on Rescaled Localized Radial Basis Function (RBF). In such a case, if we denote by \(\mathbb{P}_{12}\) and \(\mathbb{P}_{21}\) the RBF interpolation matrices (see \cite{Gervasio2019,Deparis2016a}) of size \((M_1 \times N_{\Gamma_2})\) and \((M_2 \times N_{\Gamma_1})\), respectively, and define a matrix $\mathbb{U}_{\Gamma_2,2}$ of size
\((N_{\Gamma_2} \times N_2)\), used to extract the solution of the master model on \(\Gamma_2\), then equation \eqref{eq:lifting_term_DEIM} for the interpolation of the Dirichlet data on the interfaces can be reformulated as:
\[ 
\mathbb{V}_1^T \mathbb{A}_{1,\Gamma_{1}}(\boldsymbol{\mu})
{\bf u}_{\Gamma_1}^{k+1}(\boldsymbol\mu)=
\mathbb{V}_1^T \mathbb{A}_{1,\Gamma_{1}}(\boldsymbol{\mu})\mathbf{\Phi}_{D}\mathbf{\Phi}_{D_{|\mathcal{I}_{1,D}}}^{-1} \mathbb{P}_{12}\mathbb{U}_{\Gamma_2,2}\mathbb{V}_2\mathbf{u}_{n,2}^k(\boldsymbol{\mu}). 
\]
Similarly, for the interpolation of the residuals, equation \eqref{eq:Neumm_ROM_term} can be reformulated as:
$$
{\bf r}_{\Gamma_2}^{k+1}(\boldsymbol\mu)\approx - \mathbb{M}_{\Gamma_2} \boldsymbol{\Phi}_{N}\boldsymbol{\Phi}_{D_{|\mathcal{I}_{2,N}}}^{-1}\mathbb{P}_{21}\mathbb{M}_{\Gamma_1}^{-1}\mathbb{U}_{21}^T\mathbf{r}_{\Gamma_{1|\mathcal{I}_{1,N}}}^{k+1}(\boldsymbol{\mu}).
$$
Notice that it is sufficient to compute the matrix products appearing in the right--hand sides of the above formulas only once during the offline phase.
At least for the coupled problems exploited in Section \ref{sec:numerical_results}, we experienced that, keeping fixed all the other approximations (such as the FOM finite element spaces and the POD/DEIM techniques), the interpolation approach across the interfaces affects very mildly the convergence rate of the ROM Dirichlet--Neumann method, as well as the error of the ROM solution with respect to the FOM one. In Appendix \ref{App_B} we report some numerical results obtained with RBF interpolation that corroborate this remark.
\end{remark}

\begin{algorithm}[h!]
	\caption{ROM training procedure -- Snapshots computation} \label{alg:DEIMInterface}
	\begin{algorithmic}[1]
		\Procedure{[SNAPSHOTS arrays] = Snapshots}{FOM arrays, $\mathscr{P}^d_{train}$, tol}
		\State \emph{Dirichlet and Neumann data snapshots}
		\State Given the set of $\mathscr{P}^d_{train} \subset \mathscr{P}^d$ parameters:
		\For{$\boldsymbol{\mu} \in \mathscr{P}^d_{train}$}
		\State $\mathbf{u}_1,~\mathbf{u}_2 \gets$ solve the coupled problem \eqref{eq:master_steady_algebraic}--\eqref{eq:interface_residual_def_steady} with Dirichlet--Neumann iterations, the convergence is achieved when $\|\mathbf{u}_{\Gamma_1}^{k+1}-{\mathbb R}_{12}\mathbf{u}_{\Gamma_2}^{k+1}\|_2 <$tol 
		\State $\mathbf{u}_{\Gamma_1} \gets$ extract the slave interface solution;
		\State $\mathbf{z}_{\Gamma_2} \gets$ extract the master primal residual.
        \State $\mathbf{S}_1 = [\mathbf{S}_1, \mathbf{u}_1]$; 
        \State $\mathbf{S}_2 = [\mathbf{S}_2, \mathbf{u}_2]$;
		\State $\mathbf{S}_{D} = [\mathbf{S}_D, \mathbf{u}_{\Gamma_1}]$;
		\State $\mathbf{S}_{N} = [\mathbf{S}_N, \mathbf{z}_{\Gamma_2}]$;
		\EndFor
		\EndProcedure
			\end{algorithmic}
	\end{algorithm}	

 \begin{algorithm}[h!]
	\caption{Interface DEIM training procedure -- ROM arrays} \label{alg:DEIMInterface_ROM_arrays}
	\begin{algorithmic}[1]
		\Procedure{[ROM arrays] = ROM arrays}{$\mathbf{S}_D$, $\mathbf{S}_N$, $\epsilon_{tol_{D}}$, $\epsilon_{tol_{N}}$}
		\State \emph{DEIM reduced--order arrays:}
		\State $\mathbf{\Phi}_D \gets$ POD($\mathbf{S}_{D},\epsilon_{tol_{D}}$);$\quad \mathcal{I}_{1,D} \gets$ DEIM--indices($\mathbf{\Phi}_D$);
		\State $\mathbf{\Phi}_N \gets$ POD($\mathbf{S}_{N},\epsilon_{tol_{N}}$);$\quad \mathcal{I}_{2,N} \gets$ DEIM--indices($\mathbf{\Phi}_N$);
		\State \emph{Dirichlet magic points:}
		\For{$i_{1,D} \in \mathcal{I}_{1,D}$}
		\State $\mathbf{p}_1 \gets$ get Cartesian coordinates of $i_{1,D}$ node; 
		\State $\mathbf{p}_2 = \arg\min_{\mathbf{p}_2^j \in nodes_{\Gamma_2}}(\text{dist}(\mathbf{p}_1 - \mathbf{p}_2^j)) \gets$ search in $\Gamma_2$ the nearest node to $\mathbf{p}_1 \in nodes_{\Gamma_1}$;
		\State $i_{2,D} \gets$ get the Dirichlet index for $\mathbf{p}_2$;
		\State $\mathcal{I}_{2,D} = [\mathcal{I}_{2,D}, i_{2,D}]$;
		\EndFor
		\State \emph{Neumann magic points:}
		\For{$i_{2,N} \in \mathcal{I}_{2,N}$}
		\State $\mathbf{p}_2 \gets$ get Cartesian coordinates of $i_{2,N}$ node; 
		\State $\mathbf{p}_1 = \arg\min_{\mathbf{p}_1^j \in nodes_{\Gamma_1}}(\text{dist}(\mathbf{p}_2 - \mathbf{p}_1^j)) \gets$ search in $\Gamma_1$ the nearest node to $\mathbf{p}_2 \in nodes_{\Gamma_2}$;
		\State $i_{1,N} \gets$ get the Neumann index for $\mathbf{p}_1$;
		\State $\mathcal{I}_{1,N} = [\mathcal{I}_{1,N}, i_{1,N}]$;
		\EndFor
		\EndProcedure
			\end{algorithmic}
	\end{algorithm}	
\begin{algorithm}[h!]
	\caption{Complete ROM training procedure} \label{alg:ROM_complete1}
	\begin{algorithmic}[1]
		\Procedure{[ROM arrays] = Offline training}{FOM arrays, $\mathscr{P}^d_{train}$,$\epsilon_{tol_1}$,$\epsilon_{tol_2}$, $\epsilon_{tol_{D}}$, $\epsilon_{tol_{N}}$, tol}
		\State \emph{Solution, Dirichlet, and Neumann data snapshots}
        \State $[\mathbf{S}_1, \mathbf{S}_2, \mathbf{S}_D, \mathbf{S}_N]$ = SNAPSHOTS(FOM arrays, $\mathscr{P}^d_{train}$, tol); 
		\State \emph{POD reduced--order arrays:}
		\State $\mathbb{V}_1\gets$ POD($\mathbf{S}_1, \epsilon_{tol_1}$);
		\State $\mathbb{V}_2\gets$ POD($\mathbf{S}_2, \epsilon_{tol_2}$);
		\State $\{\mathbb{A}_{n,1},\mathbf{f}_{n,1}\}\gets$ Galerkin projection of the FOM slave arrays onto $\mathbb{V}_1$;
		\State $\{\mathbb{A}_{n,2},\mathbf{f}_{n,2}\}\gets$ Galerkin projection of the FOM master arrays onto $\mathbb{V}_2$;
		\State \emph{DEIM reduced--order arrays:}
		\State $[\boldsymbol{\Phi_D}, \boldsymbol{\Phi_N},\mathcal{I}_{1,D}, \mathcal{I}_{2,N}]$ = ROM ARRAYS($\mathbf{S}_D$, $\mathbf{S}_N$, $\epsilon_{tol_{D}}$, $\epsilon_{tol_{N}}$);
		\State $\boldsymbol{\Phi}_D\boldsymbol{\Phi}^{-1}_{D_{|\mathcal{I}_{1,D}}}{\mathbb U}_{12}\mathbb{V}_{2} \gets $ save matrix product for the slave trace term;
		\State $\mathbb{V}_2^T {\mathbb E}_2\mathbb{M}_{\Gamma_2} \boldsymbol{\Phi}_N \boldsymbol{\Phi}_{N_{|\mathcal{I}_{2,N}}}^{-1}{\mathbb U}_{21}\mathbb{M}_{\Gamma_1}^{-1}{\mathbb U}_{21}^T \gets$ save matrix product for the master residual term;
		\EndProcedure
		\end{algorithmic}
\end{algorithm}	
\begin{algorithm}[h!]
	\caption{ROM query} \label{alg:ROM_complete2}
	\begin{algorithmic}[1]
		\Procedure{[$\mathbf{u}_{N,1}$,$\mathbf{u}_{N,2}$] = Online Query}{ROM arrays, FOM arrays, $\boldsymbol{\mu}$, tol}
		\State given the index $k \geq 0$ and the initial guess $\mathbf{u}_2^0$
		\While{($\|\mathbf{u}_{\Gamma_1}^{k+1} - 
        {\mathbb R}_{12}\mathbf{u}_{\Gamma_2}^k(\boldsymbol{\mu})
        \|_{\ell_2}>$tol)}
		\State $\mathbf{u}_{\Gamma_{2|\mathcal{I}_{2,D}}}^{k+1}(\boldsymbol{\mu}) \gets $ extract Dirichlet magic points values on $\Gamma_2$;
		\State$\mathbb{V}_1^T\mathbb{A}_{1,\Gamma_1}(\boldsymbol{\mu}) \boldsymbol{\Phi}_D\boldsymbol{\Phi}^{-1}_{D_{|\mathcal{I}_{1,D}}}\mathbf{u}_{\Gamma_{2|\mathcal{I}_{2,D}}}^{k+1}(\boldsymbol{\mu}) \gets$ assemble the trace term;
		\State $\mathbf{u}_{n,1}^{k+1}(\boldsymbol{\mu}) \gets$ solve the slave problem \eqref{eq:master_reduced_problem};
		\State $\mathbf{r}_{\Gamma_{1|\mathcal{I}_{1,N}}}^{k+1}(\boldsymbol{\mu}) \gets$ extract the slave interface residual on the magic points $\mathcal{I}_{1,N}$;
		\State $-\mathbb{V}_2^T \mathbb{E}_2\mathbb{M}_{\Gamma_2} \boldsymbol{\Phi}_N \boldsymbol{\Phi}_{N_{|\mathcal{I}_{2,N}}}^{-1} \mathbb{U}_{21}\mathbb{M}_{\Gamma_1}^{-1}\mathbb{U}_{21}^T \mathbf{r}_{\Gamma_{1|\mathcal{I}_{1,N}}}^{k+1}(\boldsymbol{\mu}) \gets$ assemble the interface residual term;
		\State $\mathbf{u}_{n,2}^{k+1}(\boldsymbol{\mu}) \gets$ solve the master problem \eqref{eq:slave_reduced_problem};
		\State relax the master trace on $\Gamma_2$ (see e.g. \eqref{eq:relaxation_step}).
		\State $k \gets k + 1$
		\EndWhile
        \State Recover $\mathbf{u}_{N,1}$ and $\mathbf{u}_{N,2}$, using \eqref{eq:FOM_rec_1} and \eqref{eq:FOM_rec_2}, respectively.
		\EndProcedure
	\end{algorithmic}
\end{algorithm}	

\section{Numerical results}
\label{sec:numerical_results}
In this section, we present numerical results obtained solving {\em (i)} a steady problem, namely a Dirichlet boundary value problem for a linear diffusion--reaction equation, and {\em (ii)} a time--dependent problem, namely an initial--boundary value problem for the heat equation with Neumann boundary conditions. In particular, we aim to investigate the performances of the proposed algorithm by comparing FOM and ROM results in terms of both efficiency and accuracy.  

The mathematical models and numerical methods presented in this section have been implemented in C++ and Python languages and based on life$^\text{x}$ (\url{https://lifex.gitlab.io}) \cite{africa2022lifex}, an in--house high--performance C++ FE library mainly focused on cardiac applications based on deal.II FE core \cite{dealII93} (\url{https://www.dealii.org}). Both online and offline stages of the simulations have been performed in serial on a notebook with an Intel Core i7--10710U processor and 16GB of RAM. 

In what follows, spatial domains are discretized employing Q1-FEM. We would like to point out that the geometric error induced by the approximation of curve surfaces by planes is of order $h^2$ and is totally consistent with the approximation properties of Q1-FEM. The gaps and overlaps on the interface are again of order $h^2$ and they do not downgrade the FOM approximation error. Moreover, the ROM procedure is totally independent of the FE element used, since the reduced basis is not related to spatial properties
of the mesh. 

\subsection{Test\#1. Steady case: diffusion--reaction equation}
\label{Subsect:steady_case}
\begin{figure}[h!]
	\centering
    \includegraphics[width = 1\textwidth]{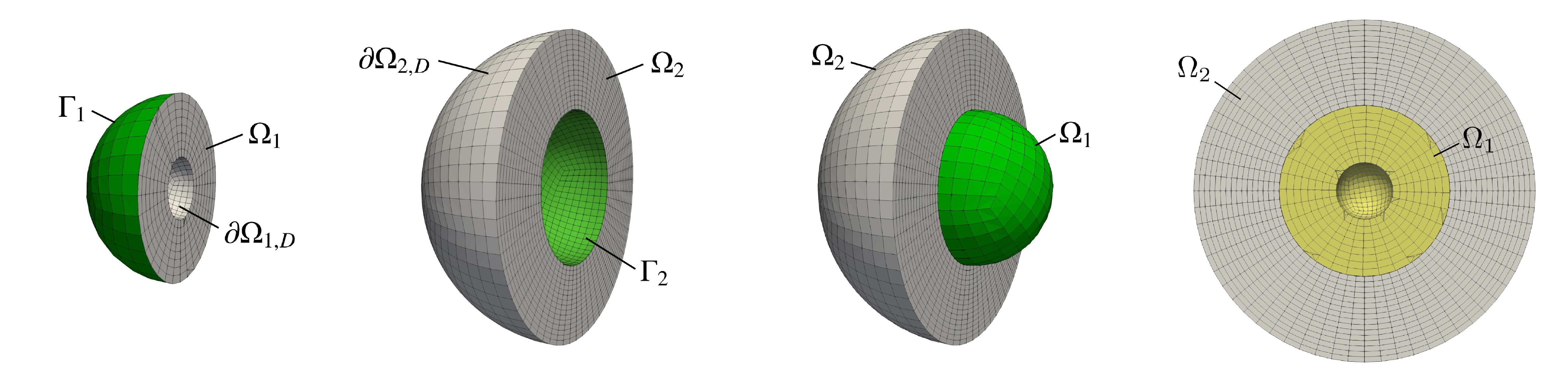}
	\caption{\emph{Test\#1.} From left to right: Half of the slave domain, half of the master domain, lateral view, and cross--section of the two half subdomains. In green the interface $\Gamma$.}
	\label{Fig:Laplace_domains}
\end{figure}

In this first test case, we solve the following boundary value problem for a diffusion--reaction equation: find $u \in \Omega$ such that
\begin{equation}
\label{Eq:test_case_1}
\begin{cases}
-\nabla \cdot (\alpha \nabla u) + \beta u = f &\text{in }\Omega \\
u = g_D &\text{on }\partial\Omega_D 
\end{cases}
\end{equation}
where the domain $\Omega$ is a hollow spheroid (see Fig. \ref{Fig:Laplace_domains}) with inner and outer radius equal to $0.5 m$ and $3.0 m$, respectively, and parameters vector $\boldsymbol{\mu} = (\alpha,\beta)$. The interface $\Gamma$ is here represented by a sphere of radius $1.5m$. Moreover, $f(x,y,z) = \frac{\pi}{4} y x^2 \sin\left(\
\frac{\pi}{2} y \right) \exp(z - 1)$, while 
\[g_D = \begin{cases}
0.01 &\text{on the internal sphere}\\
0 &\text{on the external sphere}.
\end{cases}\]
We split the domain into two hollow spheroids with a common interface (see Fig. \ref{Fig:Laplace_domains}), where $\Omega_1$ is the inner and $\Omega_2$ is the outer spheroid, and apply the Dirichlet--Neumann iterative scheme so that, given $\lambda_2^0$, for each $k\geq 0$, we solve the two following sub--problems
\[
\begin{cases}
-\nabla \cdot (\alpha \nabla u_1^{k+1}) + \beta u_1^{k+1} = f &\text{in }\Omega_1 \\
u_1^{k+1} = 0.01 &\text{on }\partial\Omega_{1,D}\\
u_1^{k+1} = \lambda_2^{k}   &\text{on } \Gamma_1
\end{cases}
\]
and
\[
\begin{cases}
-\nabla \cdot (\alpha \nabla u_2^{k+1}) + \beta u_2^{k+1} = f &\text{in }\Omega_2 \\
u_2^{k+1} = 0 &\text{on }\partial\Omega_{2,D}\\
\alpha \nabla u_2^{k+1} \cdot \mathbf{n}_2 = \alpha \nabla u_1^{k+1} \cdot \mathbf{n}_2   &\text{on } \Gamma_2,
\end{cases}
\]
while
\begin{equation}
\label{eq:relaxation_step}
\lambda_2^{k+1} = \omega u_{2_{|\Gamma_2}}^k + (1-\omega)\lambda_2^k.\end{equation}

The acceleration parameter $\omega = 0.25$ is fixed among the iterations, for both FOM and ROM problems, whereas we choose a tolerance $\epsilon = 10^{-10}$ to check at each iteration $k$ the continuity of the solutions at the subdomains interface, according to \eqref{eq:interface_continuity_check}.
 
As FOM dimension, we consider $N = 3474$ in the slave domain $\Omega_1$ and $N = 26146$ in the master domain $\Omega_2$. We choose to vary the parameters $\alpha$ and $\beta$ in $[1,10]$ using an LHS to compute the solution snapshots.  We select $N_{\text{train}} = 150$ to get a sufficiently rich snapshot set to build the reduced bases, while additional $N_{\text{test}} = 20$ values of the parameter vectors are chosen to test the method.

\begin{figure}[h!]
	\centering
	 \hspace{-40pt} $\alpha = 2.35, ~\beta = 9.55$ \qquad  \qquad  \qquad	$\alpha = 6.85, ~\beta = 3.25$ \qquad \qquad \qquad	$\alpha = 4.15, ~\beta = 5.05$ \\
	\vspace{7pt}
	\includegraphics[width=0.3\textwidth]{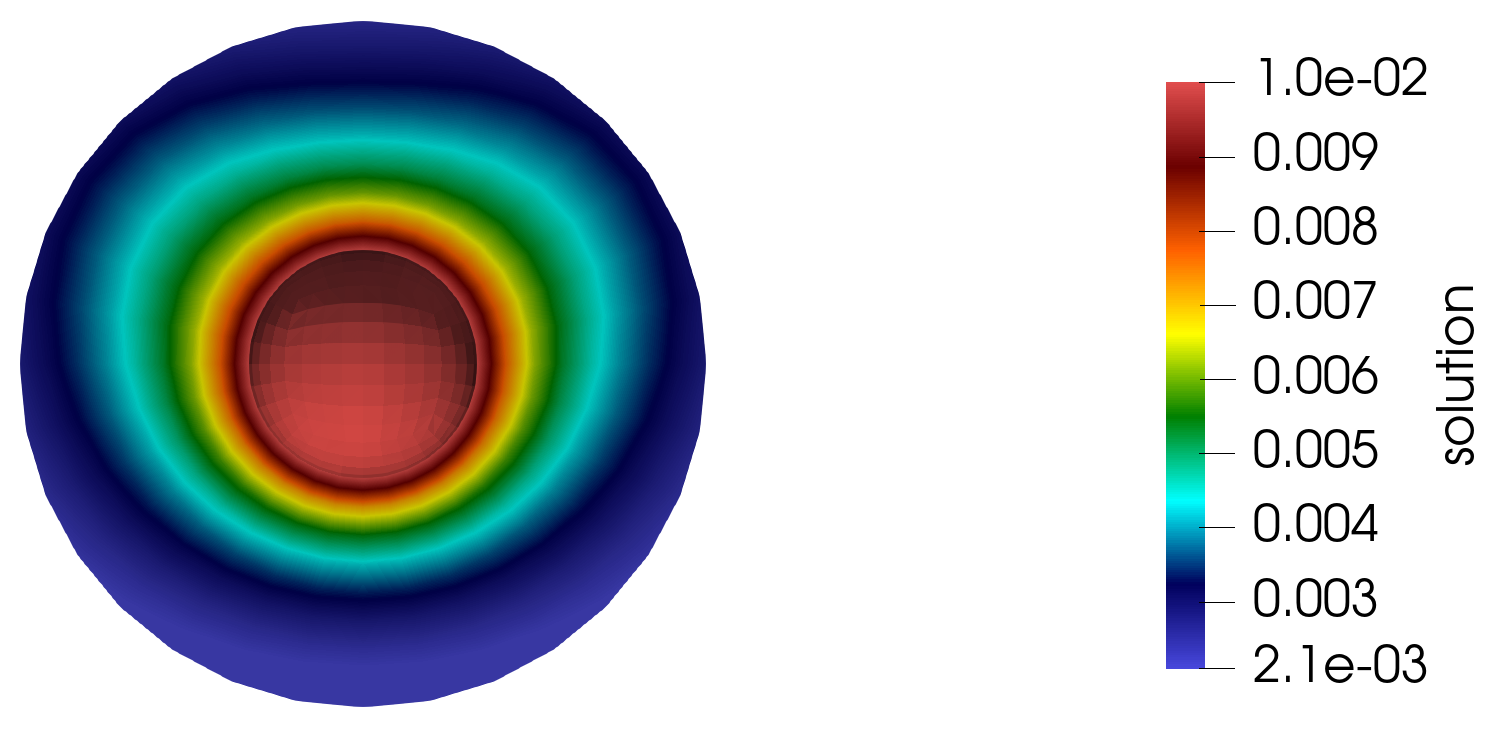}	
	$\quad$
	\vspace{7pt}
	\includegraphics[width=0.3\textwidth]{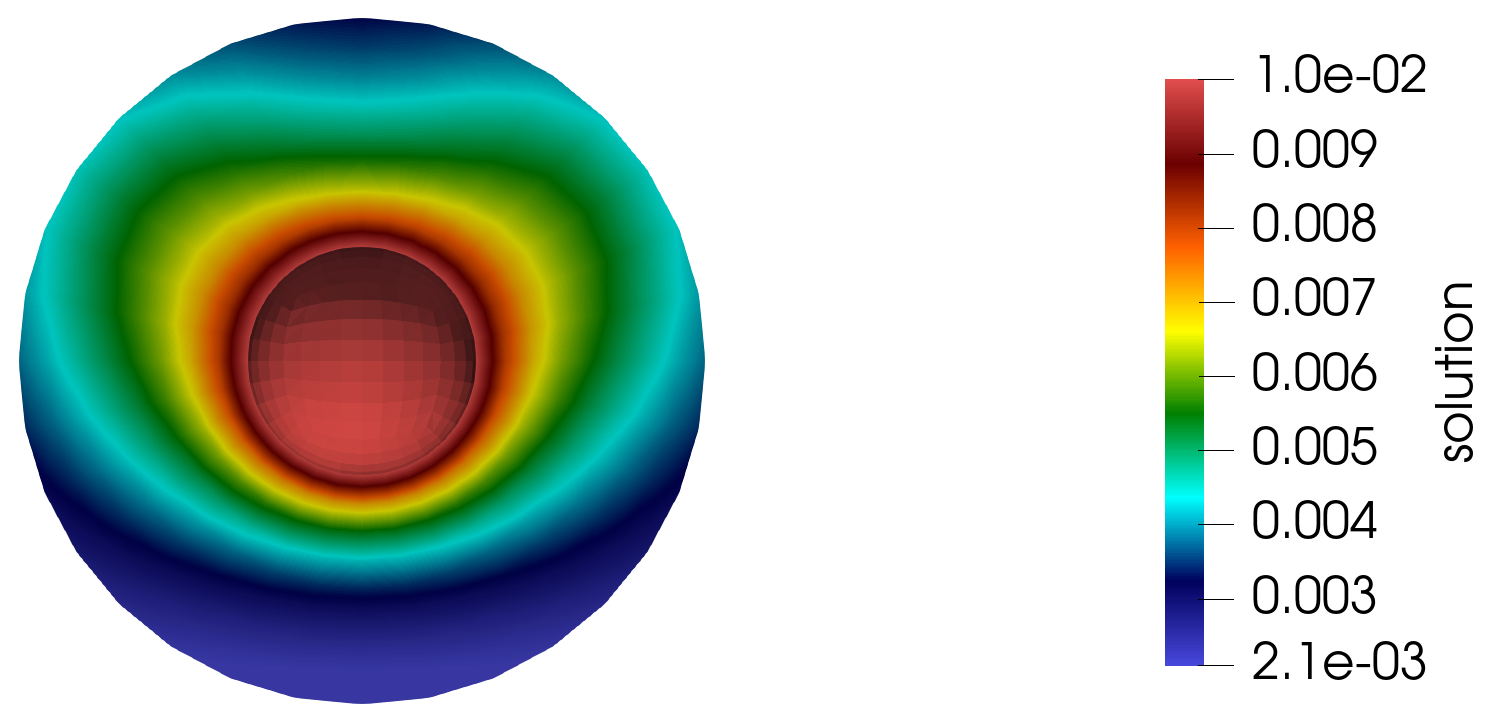}	
	$\quad$
	\vspace{7pt}
	\includegraphics[width=0.3\textwidth]{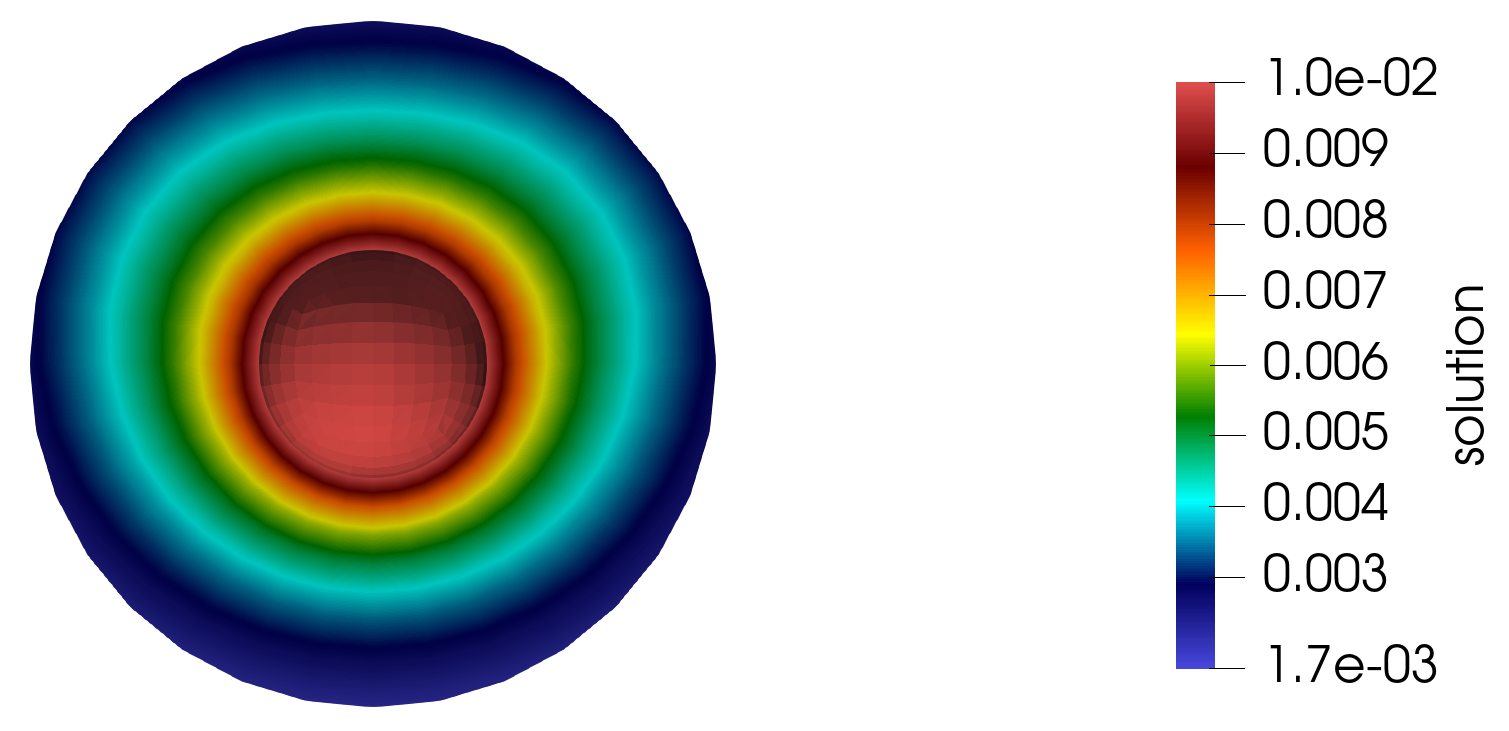}
	 \\
	\includegraphics[width=0.3\textwidth]{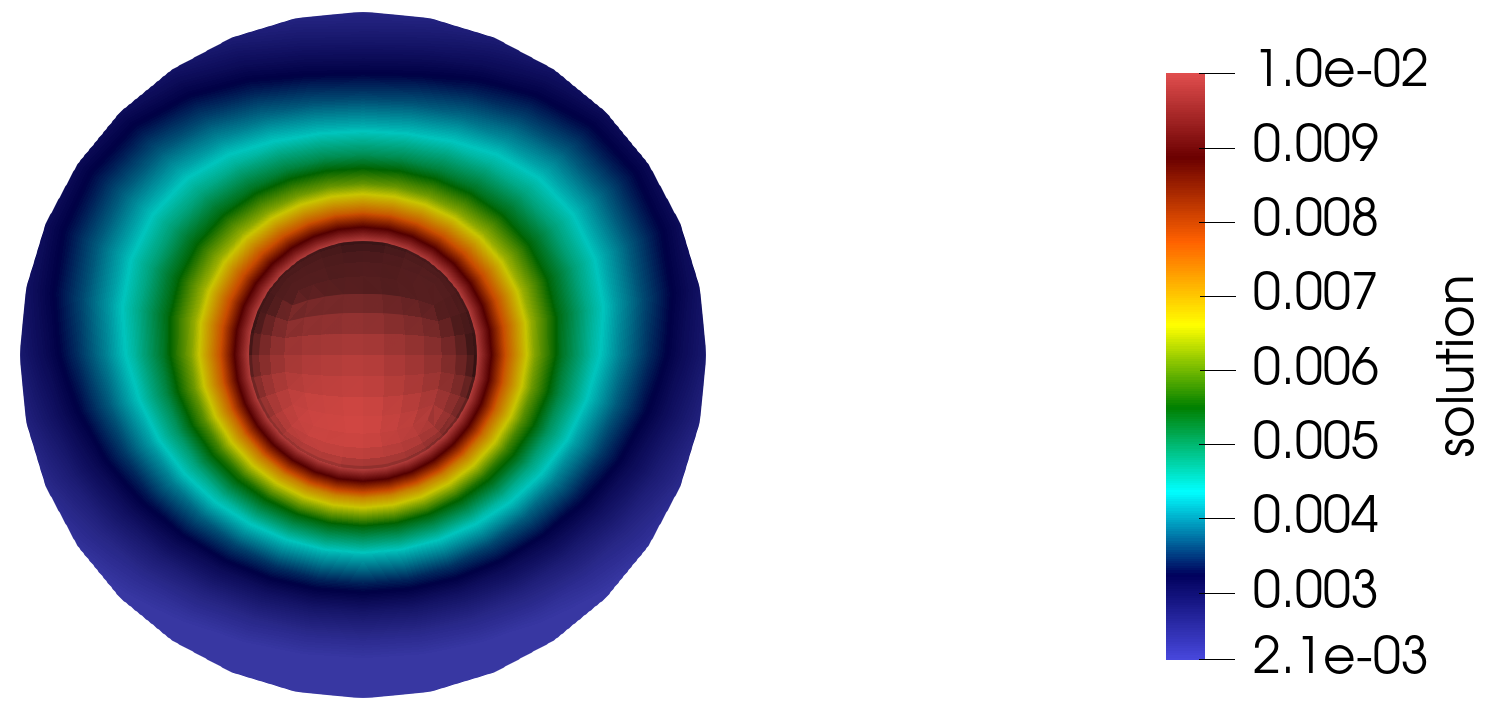}
	$\quad$
	\includegraphics[width=0.3\textwidth]{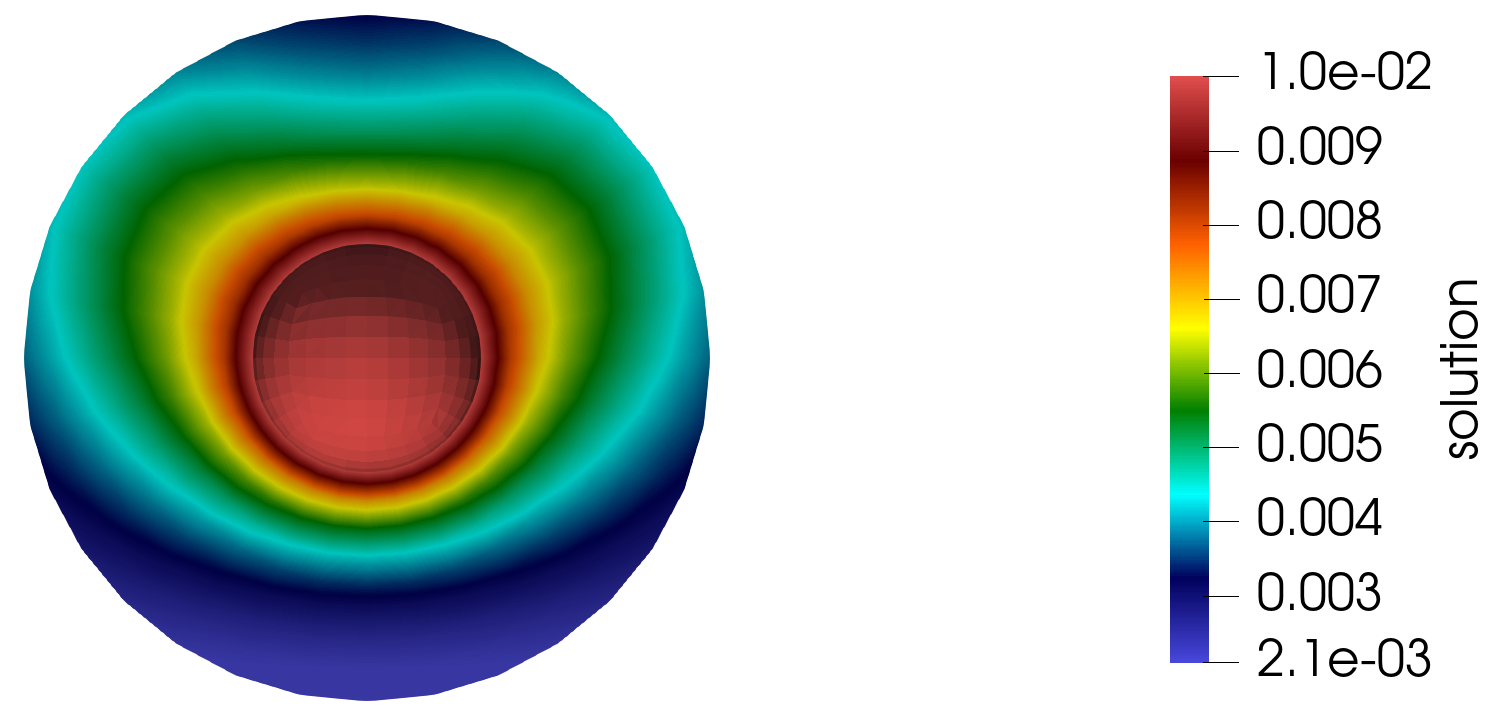}
	$\quad$
	\includegraphics[width=0.3\textwidth]{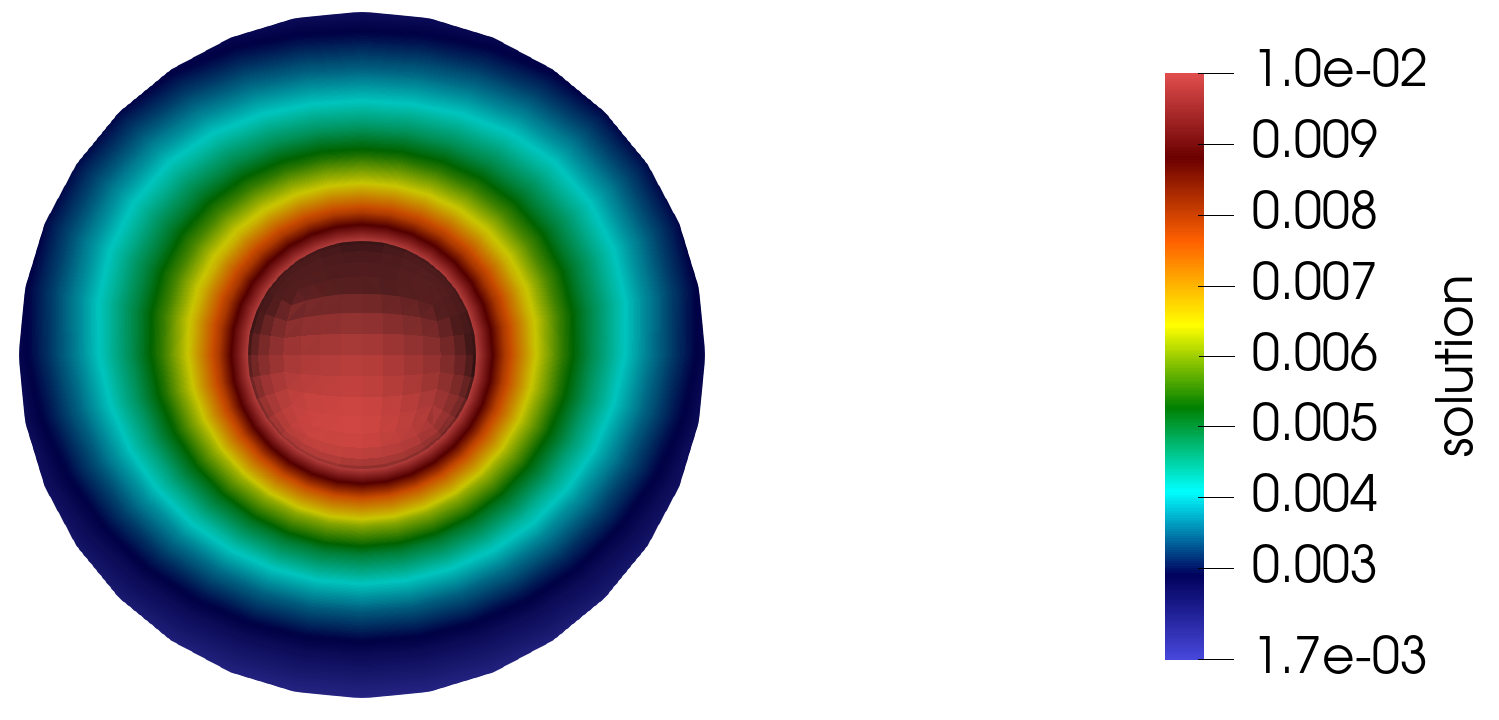}
	\\
	\vspace{12pt}
	\includegraphics[width=0.3\textwidth]{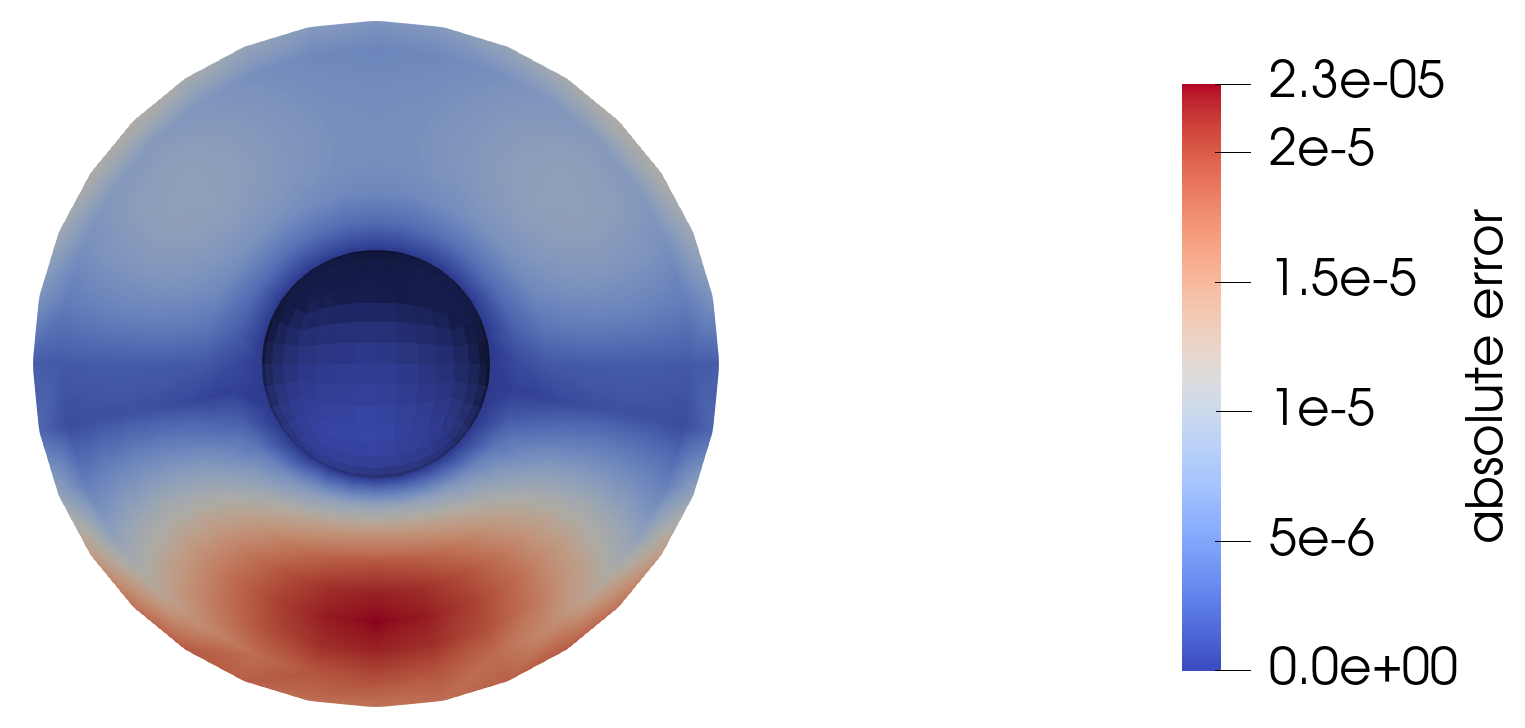}	
	$\quad$
	\includegraphics[width=0.3\textwidth]{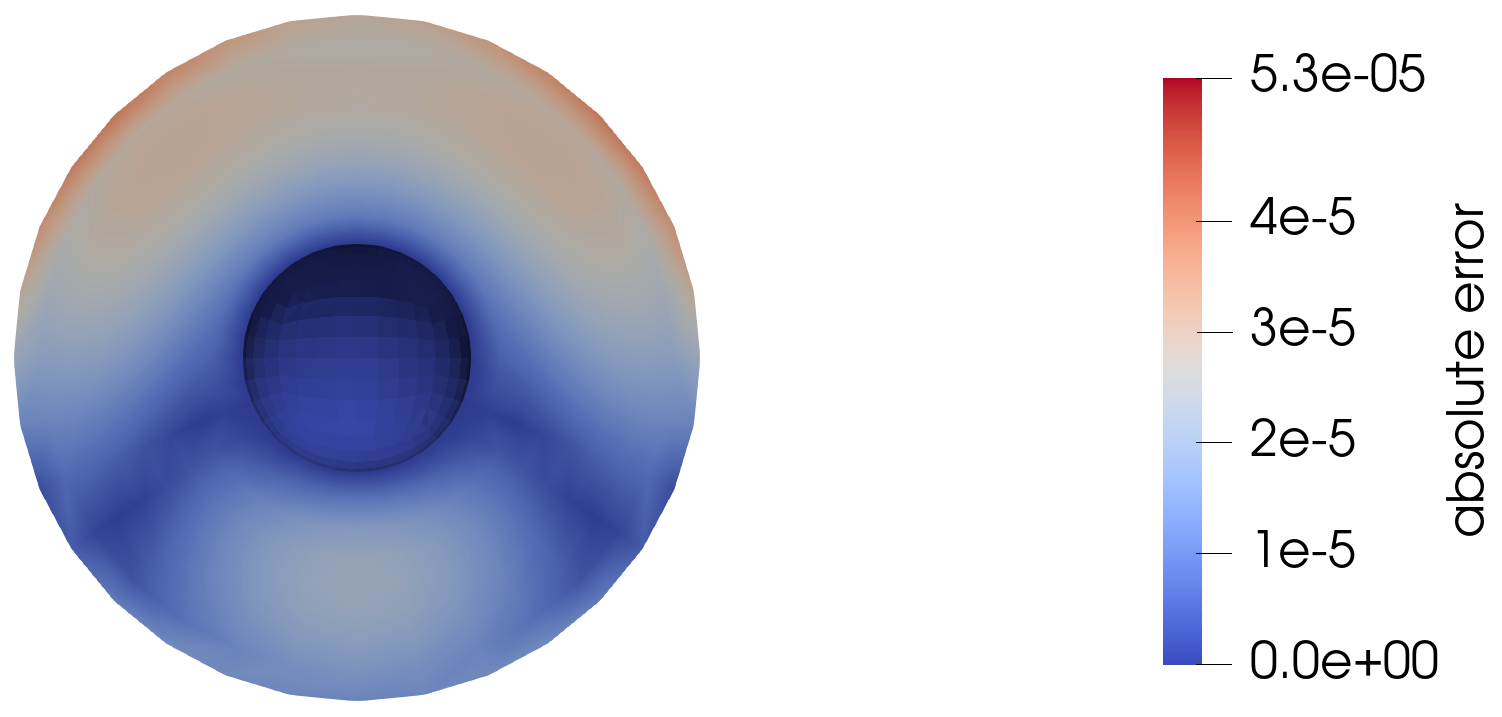}
	$\quad$
	\includegraphics[width=0.3\textwidth]{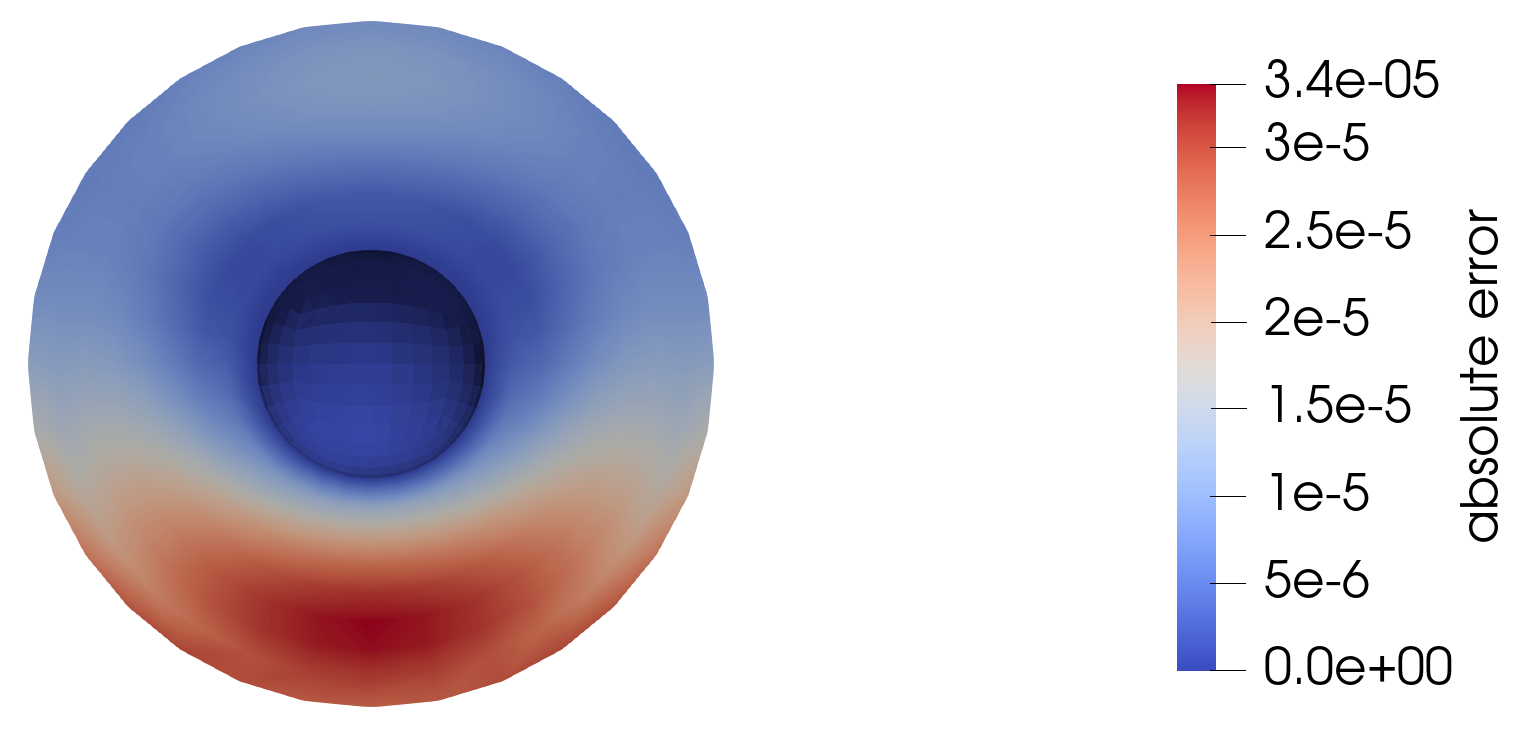}
	
	\caption{\emph{Test\#1.} Slave solution FOM (top), ROM (center) solutions, and absolute error (bottom) for three different vectors of testing parameters.}
	\label{fig:snapshots_laplace_master}
\end{figure}

\begin{figure}[h!]
	\centering
	 \hspace{-40pt} $\alpha = 2.35, ~\beta = 9.55$ \qquad \qquad \qquad $\alpha = 6.85, ~\beta = 3.25$ \qquad \qquad \qquad $\alpha = 4.15, ~\beta = 5.05$ \\
	\vspace{7pt}
	\includegraphics[width=0.3\textwidth]{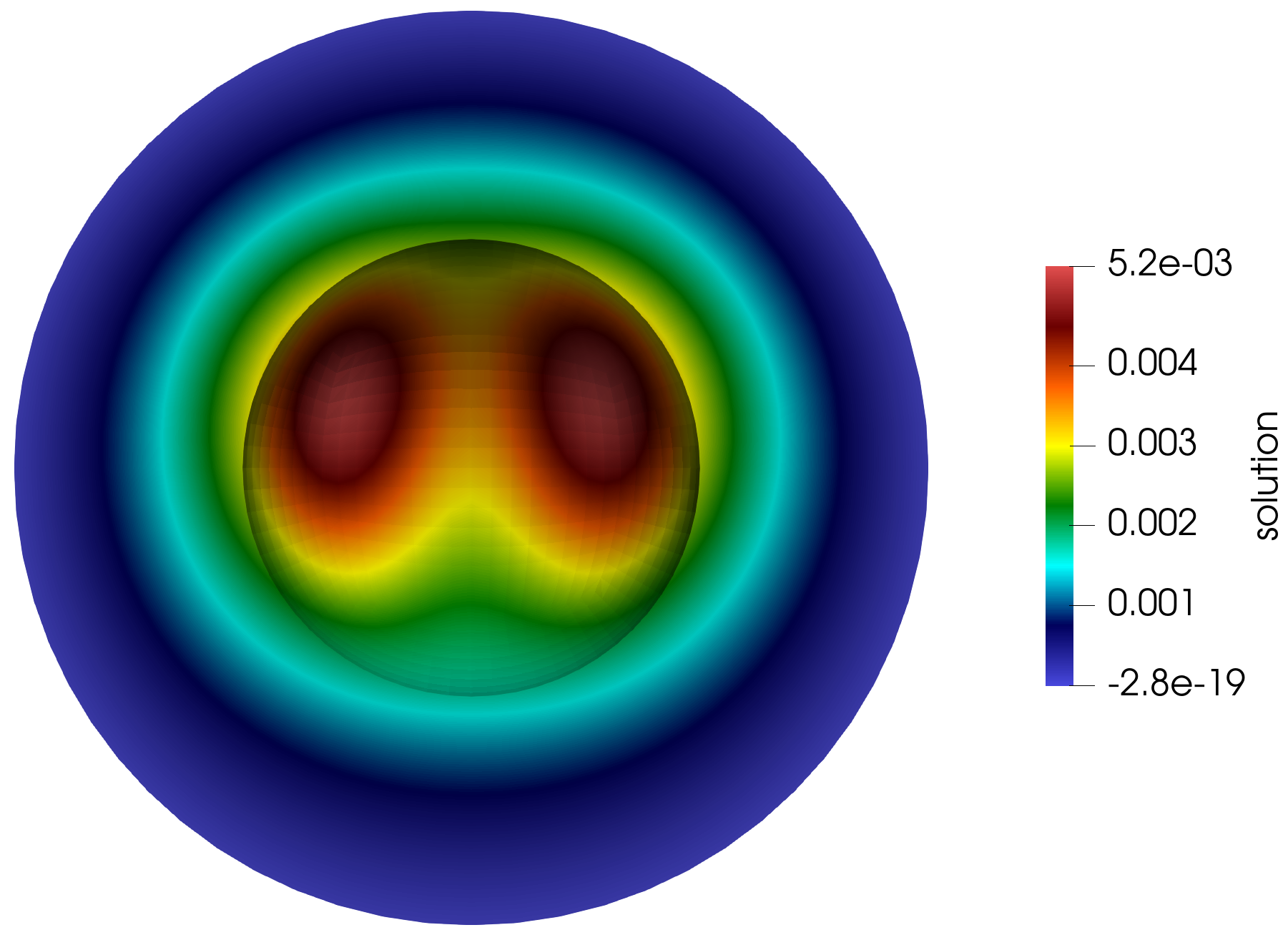}
	$\quad$
	\vspace{7pt}
	\includegraphics[width=0.3\textwidth]{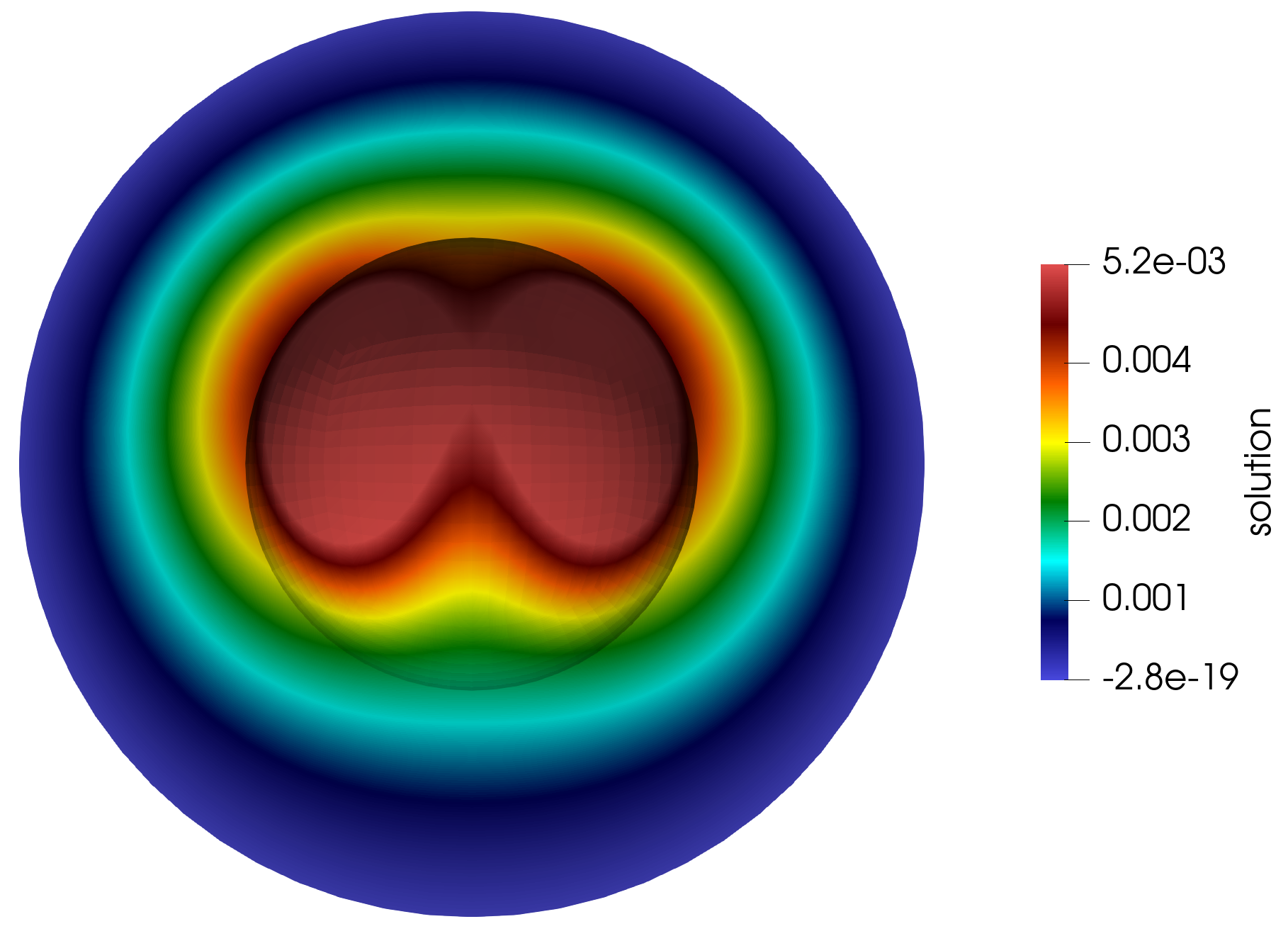}%
	$\quad$
	\vspace{7pt}
	\includegraphics[width=0.3\textwidth]{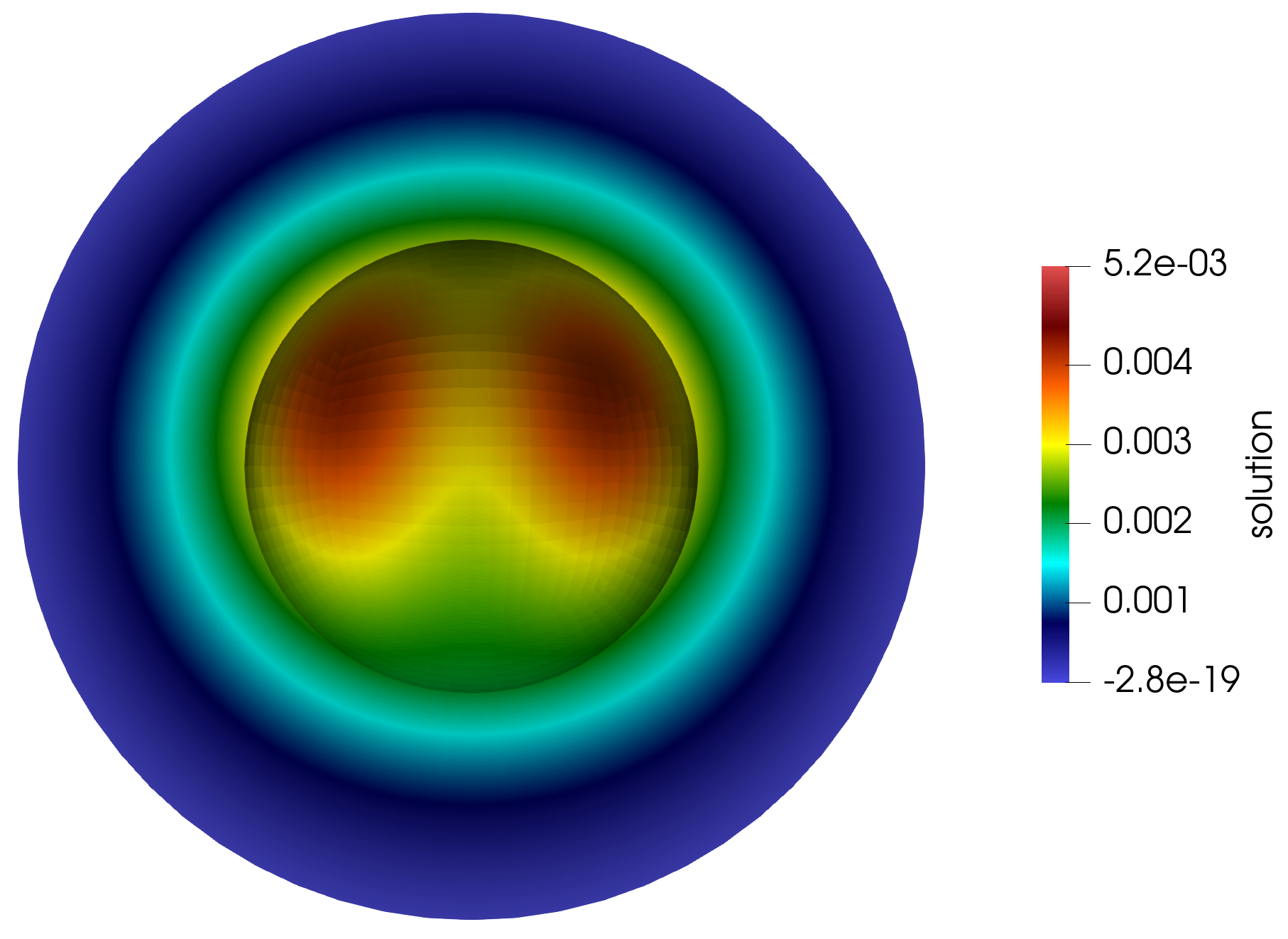}
	\\
	\includegraphics[width=0.3\textwidth]{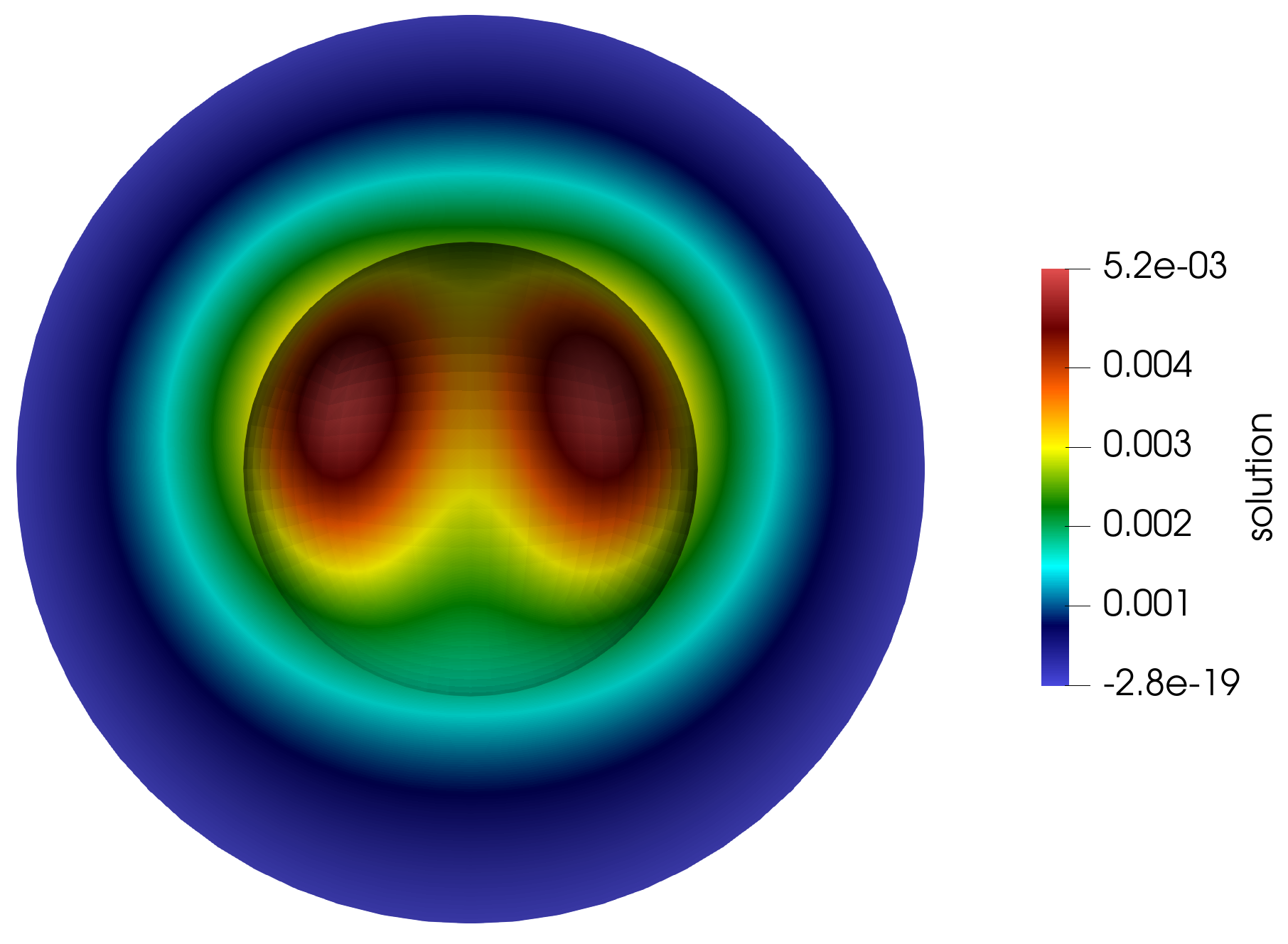}
	$\quad$
	\includegraphics[width=0.3\textwidth]{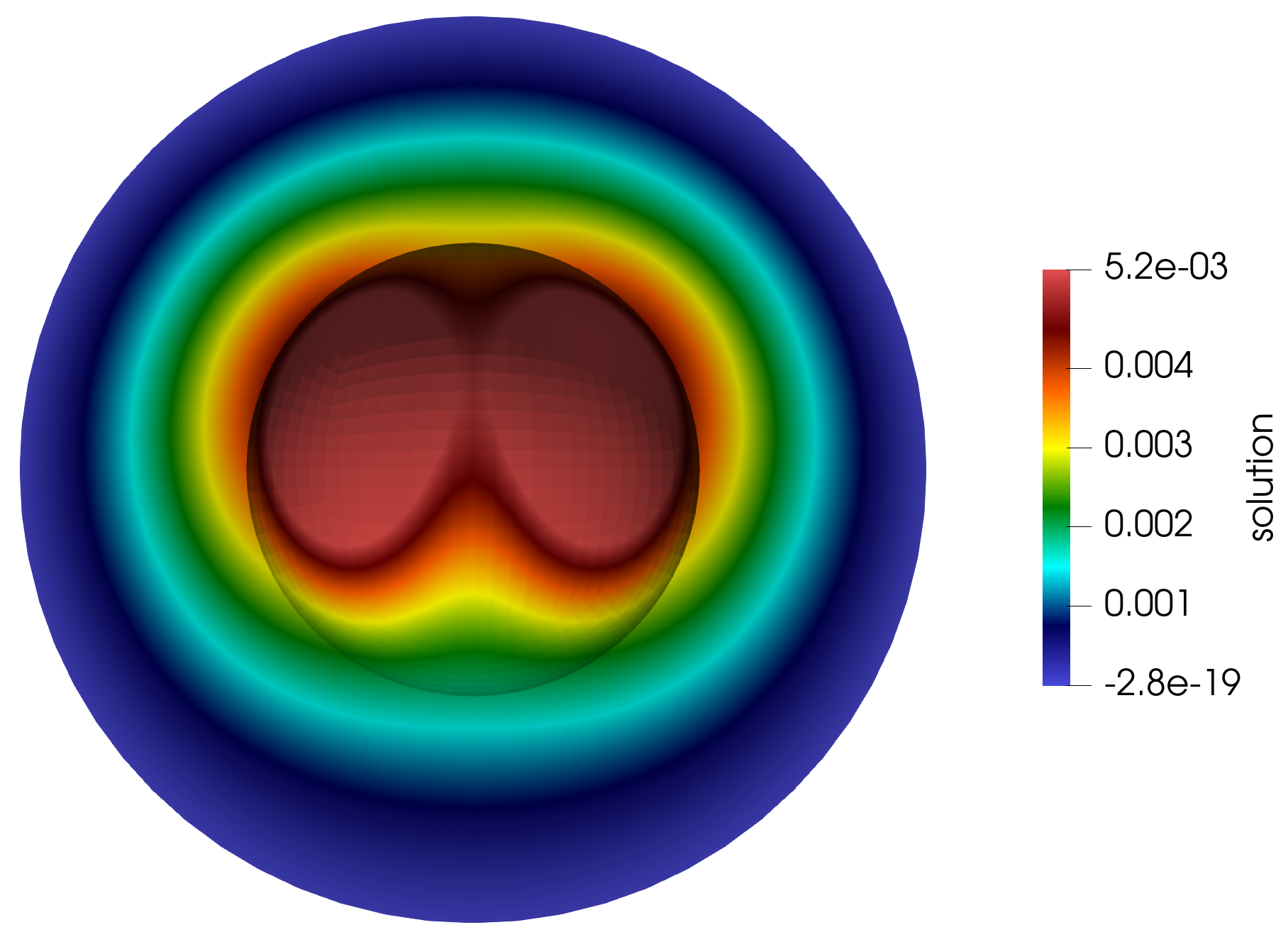}
	$\quad$
	\includegraphics[width=0.3\textwidth]{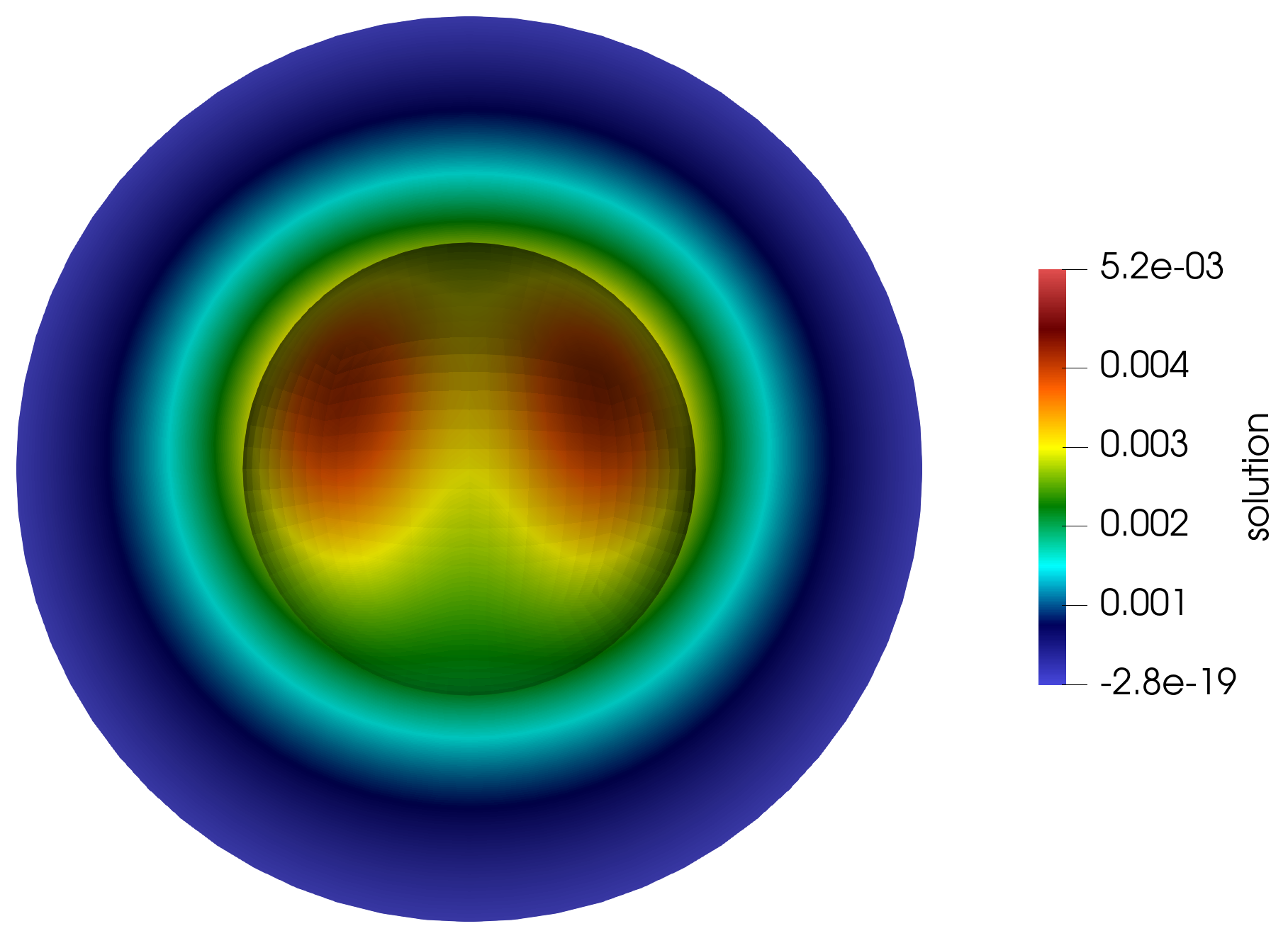}
	\\
	\includegraphics[width=0.3\textwidth]{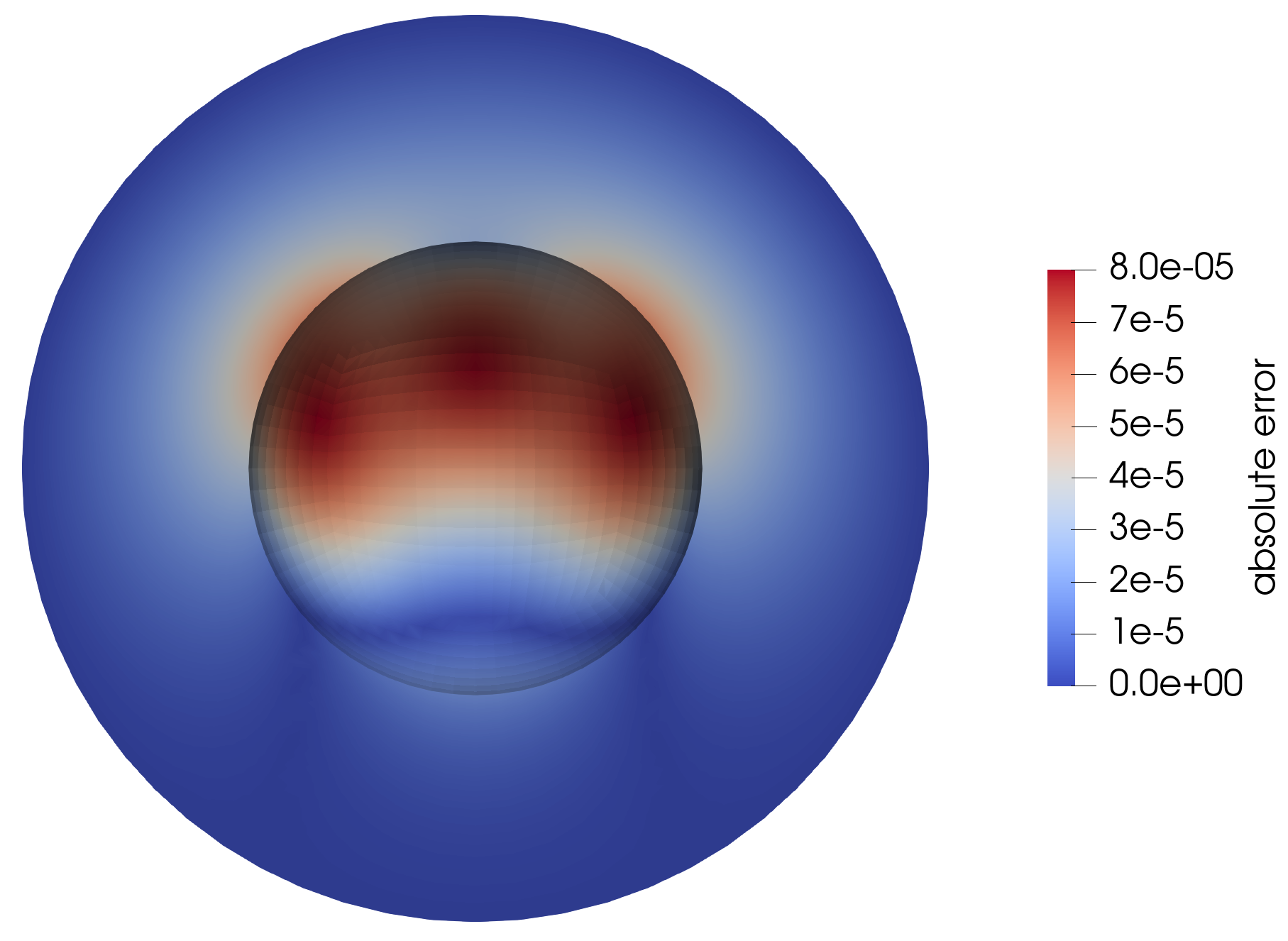}
	$\quad$
	\includegraphics[width=0.3\textwidth]{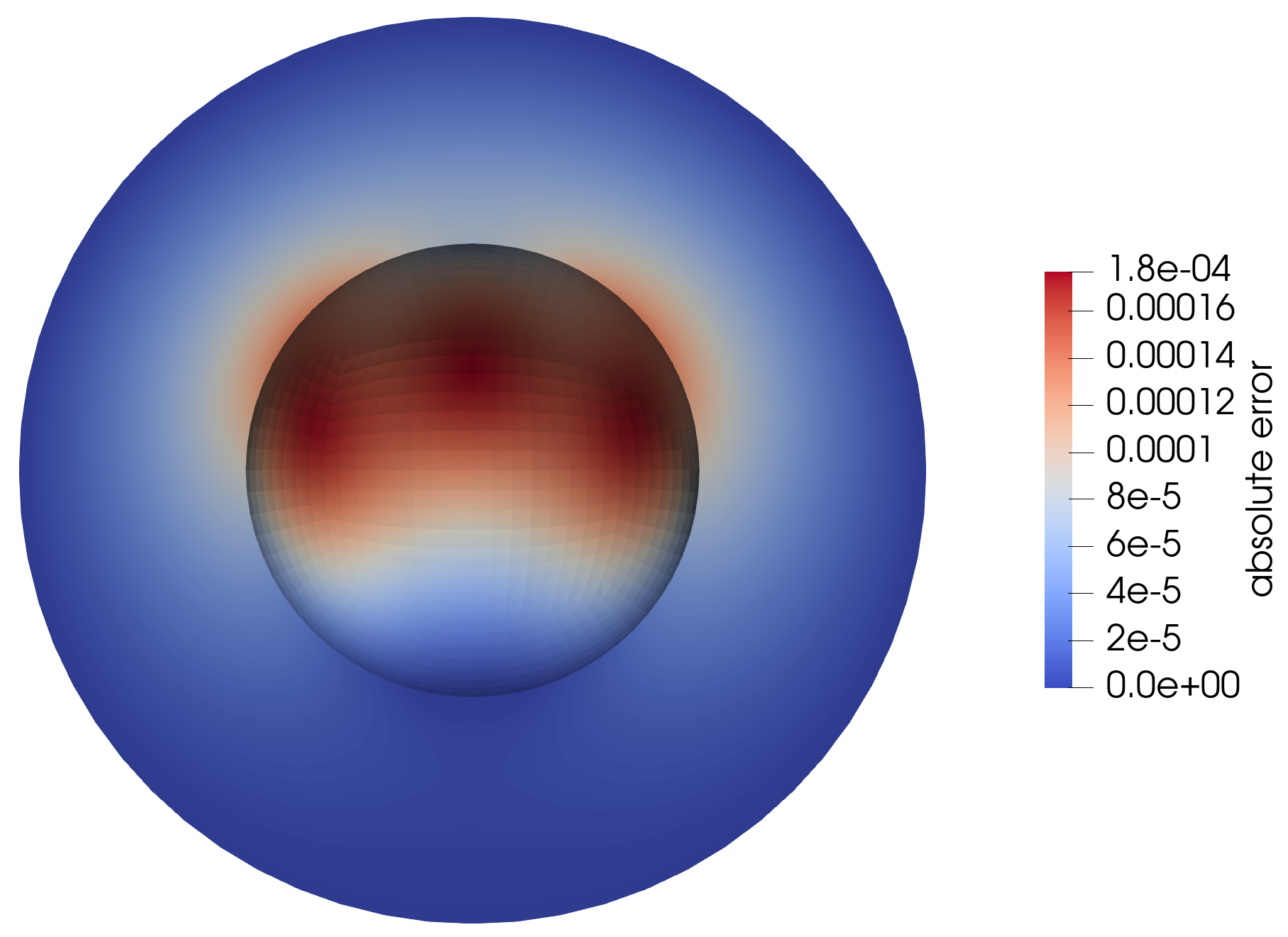}
	$\quad$
	\includegraphics[width=0.3\textwidth]{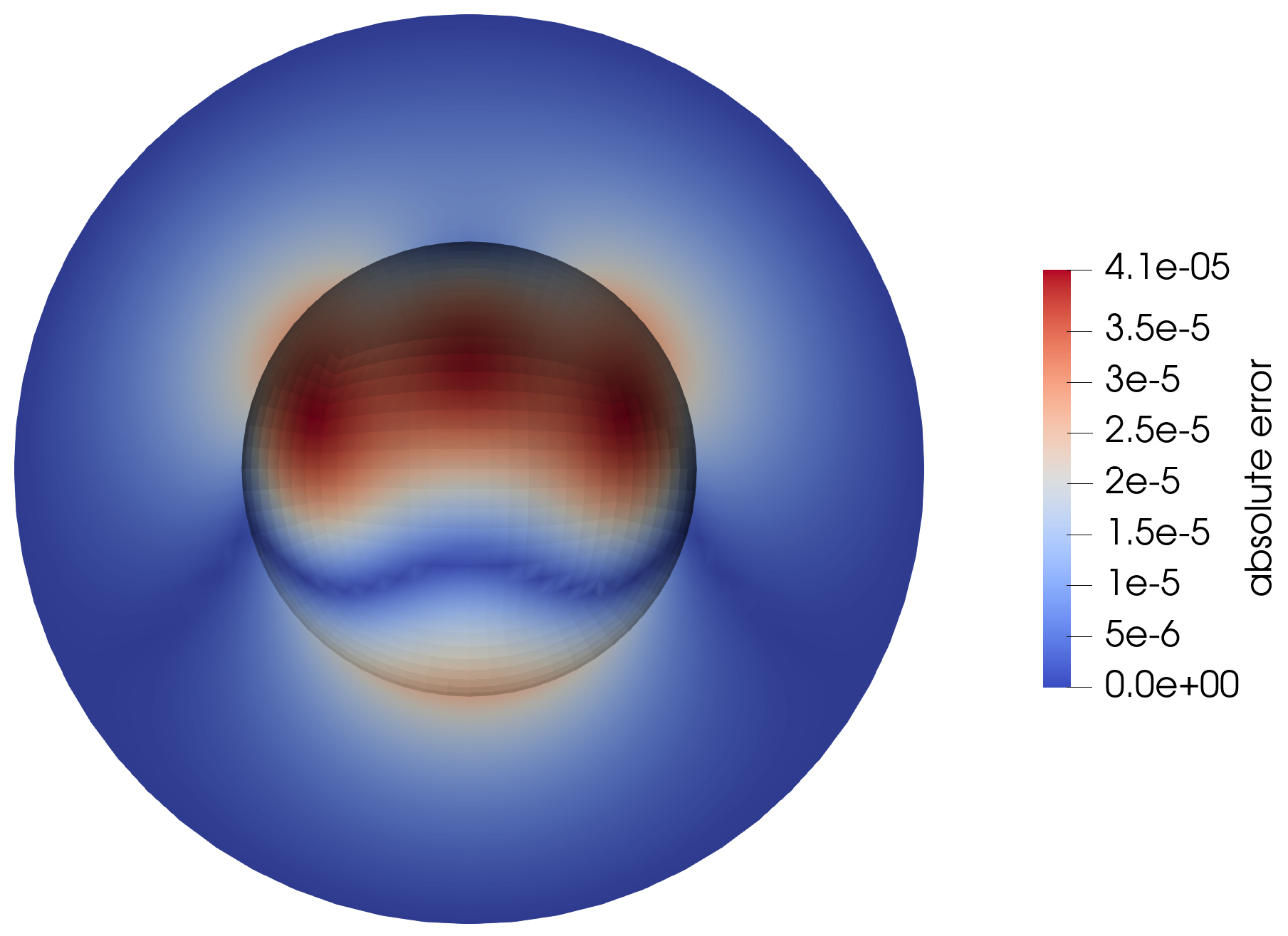}
	\caption{\emph{Test\#1.} Master solution FOM (top), ROM (center) solutions, and absolute error (bottom) for three different vectors of testing parameters.}
	\label{fig:snapshots_laplace_slave}
\end{figure} 


%

As highlighted in Section \ref{sec:Two_way_models}, Remark \ref{rem:complete_ROM}, for each set of parameters we solve twice the high--fidelity FOM (once with a conforming mesh on $\Omega$ designed starting from the mesh size in $\Omega_1$, and the other time with another conforming mesh on $\Omega$ designed starting from the mesh size in $\Omega_2$) and once the ROM. Computational relative errors are evaluated considering the difference between the FOM and ROM solutions in the $H^1(\Omega_i)$ norm, normalized with the $H^1(\Omega_i)$ norm of the FOM solutions, \emph{i.e.}
\begin{equation}
\label{eq:error_2norm}
\frac{\| \mathbf{u}_{FOM}(\boldsymbol{\mu}) - \mathbf{u}_{ROM}(\boldsymbol{\mu})\|_{H^1(\Omega_i)}}{\| \mathbf{u}_{FOM}(\boldsymbol{\mu}) \|_{H^1(\Omega_i)}}, \quad i=1,2.
\end{equation}

\begin{remark} 
The FOM solutions are computed on conforming meshes, whereas the ROM solutions are computed on non-conforming ones. Consequently, FOM and ROM solutions inherently differ. This discrepancy, which should be related to the distinction between the FOM computed on conforming meshes (as conducted in this work) and non-conforming ones, should become apparent when performing a comprehensive convergence study of the RB method proposed here. However, this topic exceeds the scope of this work and will be the focus of future research endeavors.
\end{remark}

\begin{figure}[h!]
	\centering
    \includegraphics[width=0.42\textwidth]{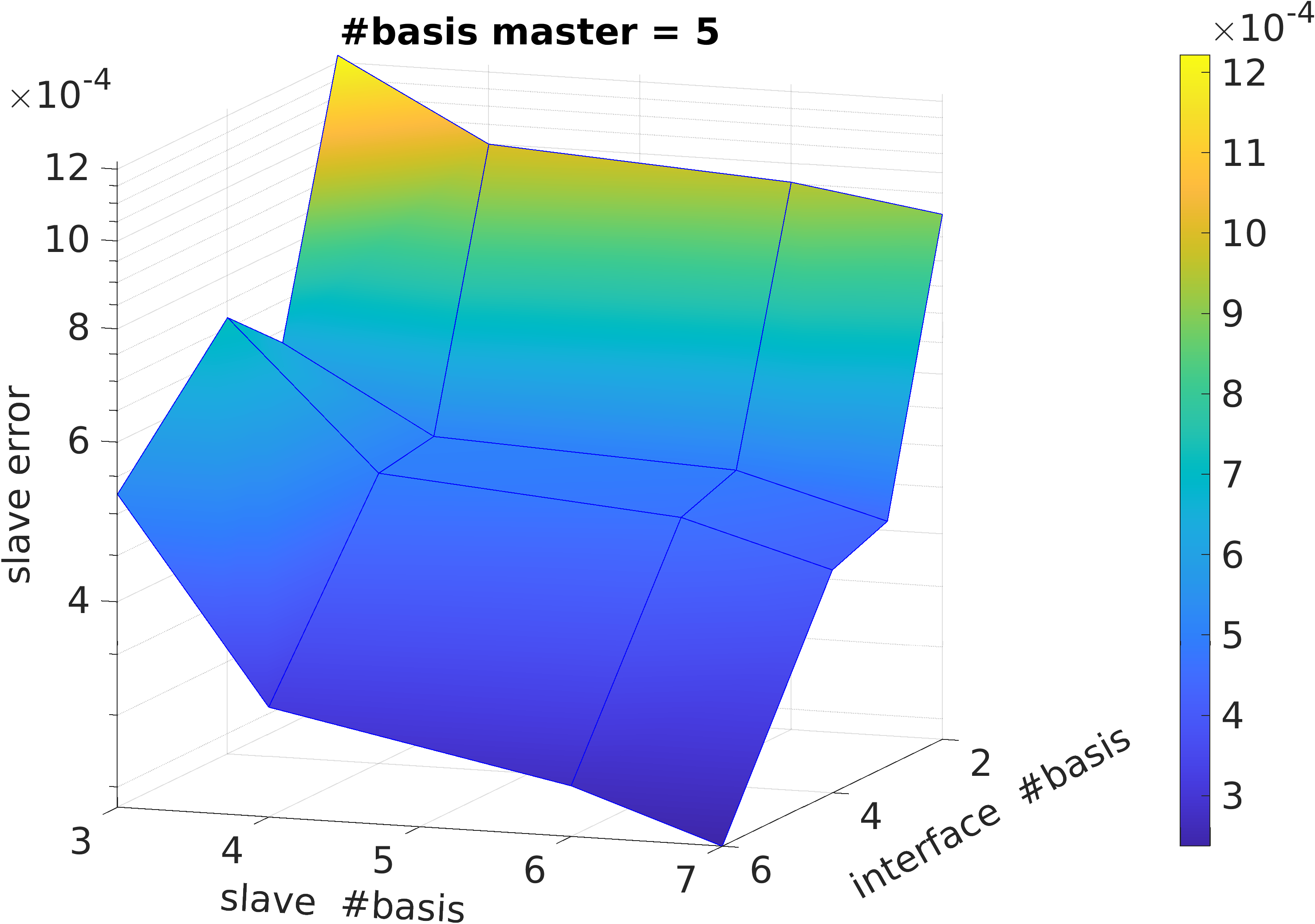}
    \quad
	\includegraphics[width=0.42\textwidth]{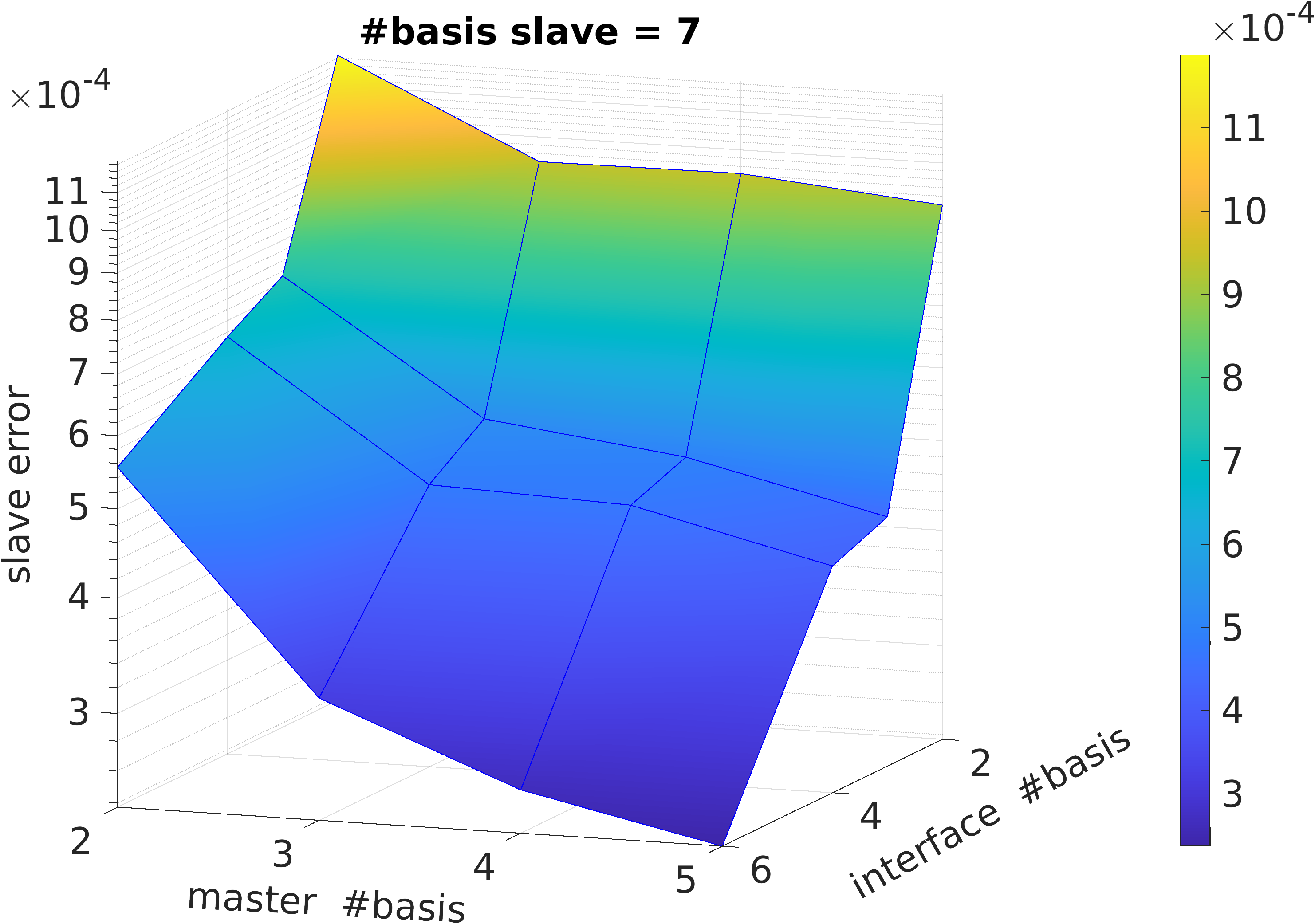}
	\\
	\bigskip
	\includegraphics[width=0.42\textwidth]{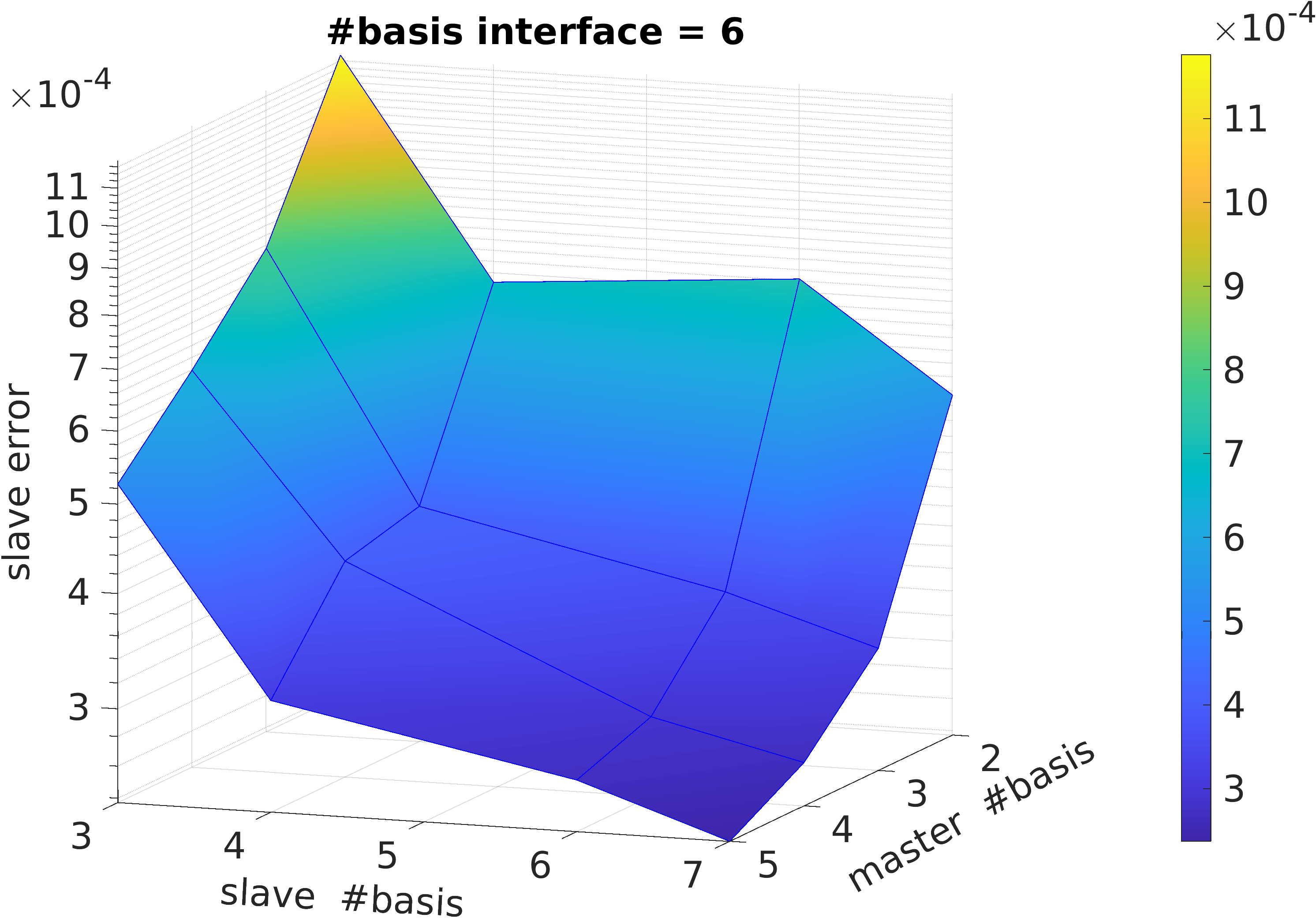}
    \quad
	\includegraphics[width=0.42\textwidth]{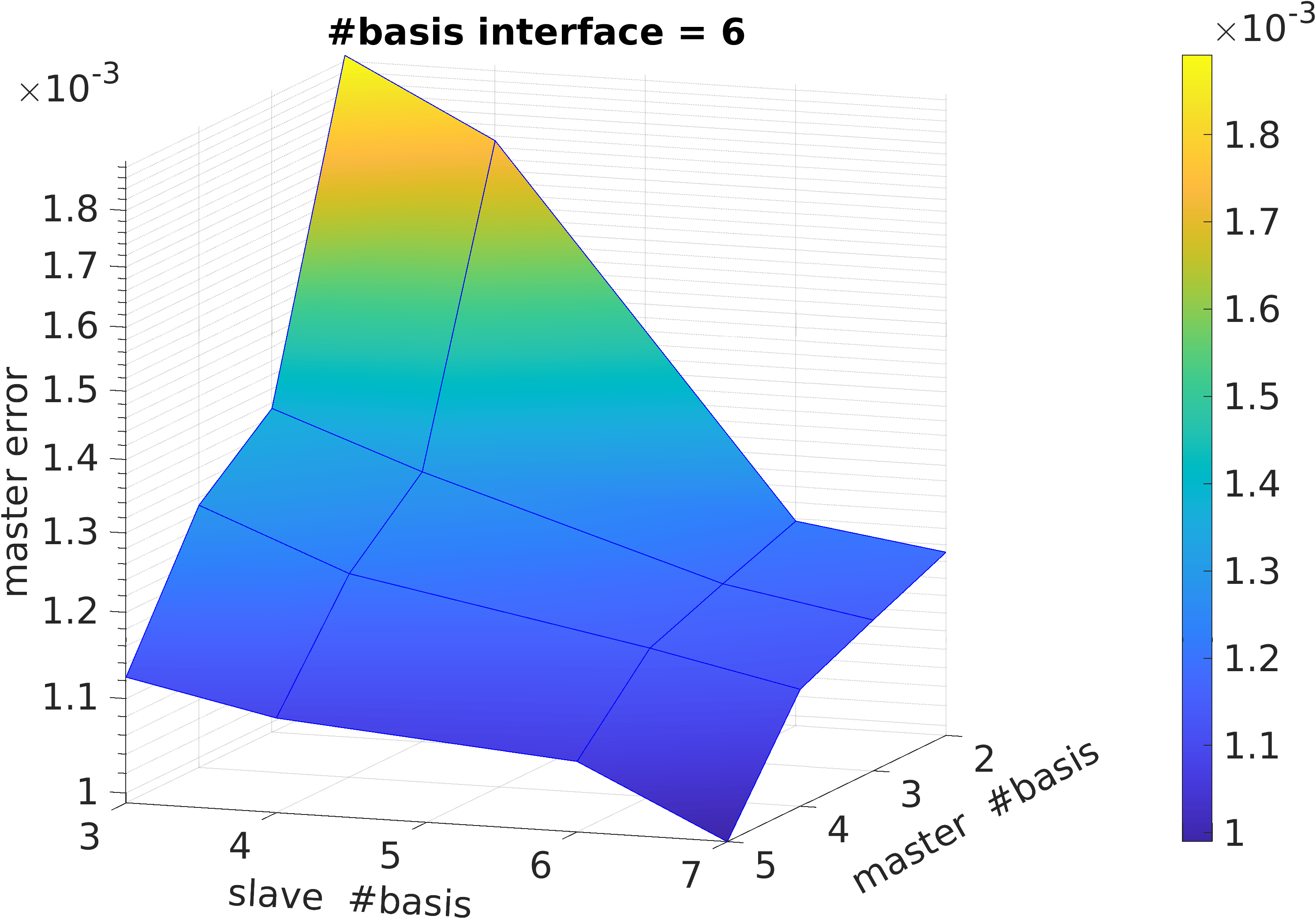}
	\caption{\emph{Test\#1.} $H^1(\Omega_i)$ mean relative error ($z$--axis) over the solution for $N_\text{test} = 20$ different instances of the parameters between the FOM and ROM solutions varying the number of basis functions used to represent the slave and the master solution $n_1$ and $n_2$, and the interface data $M_1$ and $M_2$ ($x$-- and $y$--axis). On the top row, we fix the number of basis functions of the master problem to 5 (on the left) and to 7 for the slave problem (on the right), while on the bottom we fix the number of basis functions equal to 6 for the interface data representation.}
	\label{fig:Laplace_error}
\end{figure}

\begin{figure}[h!]
	\centering
    \includegraphics[width=0.45\textwidth]{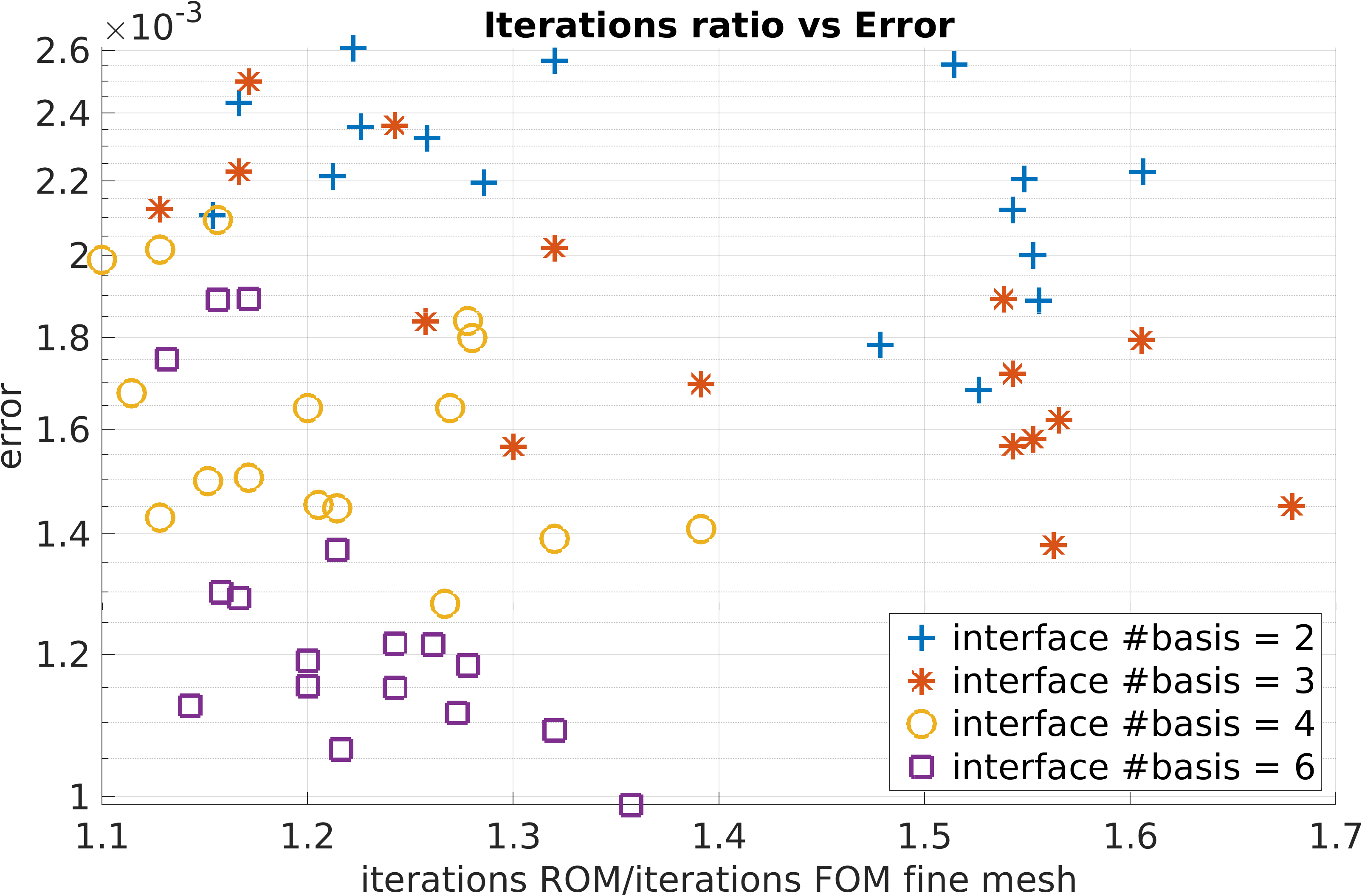}
    \quad
	\includegraphics[width=0.45\textwidth]{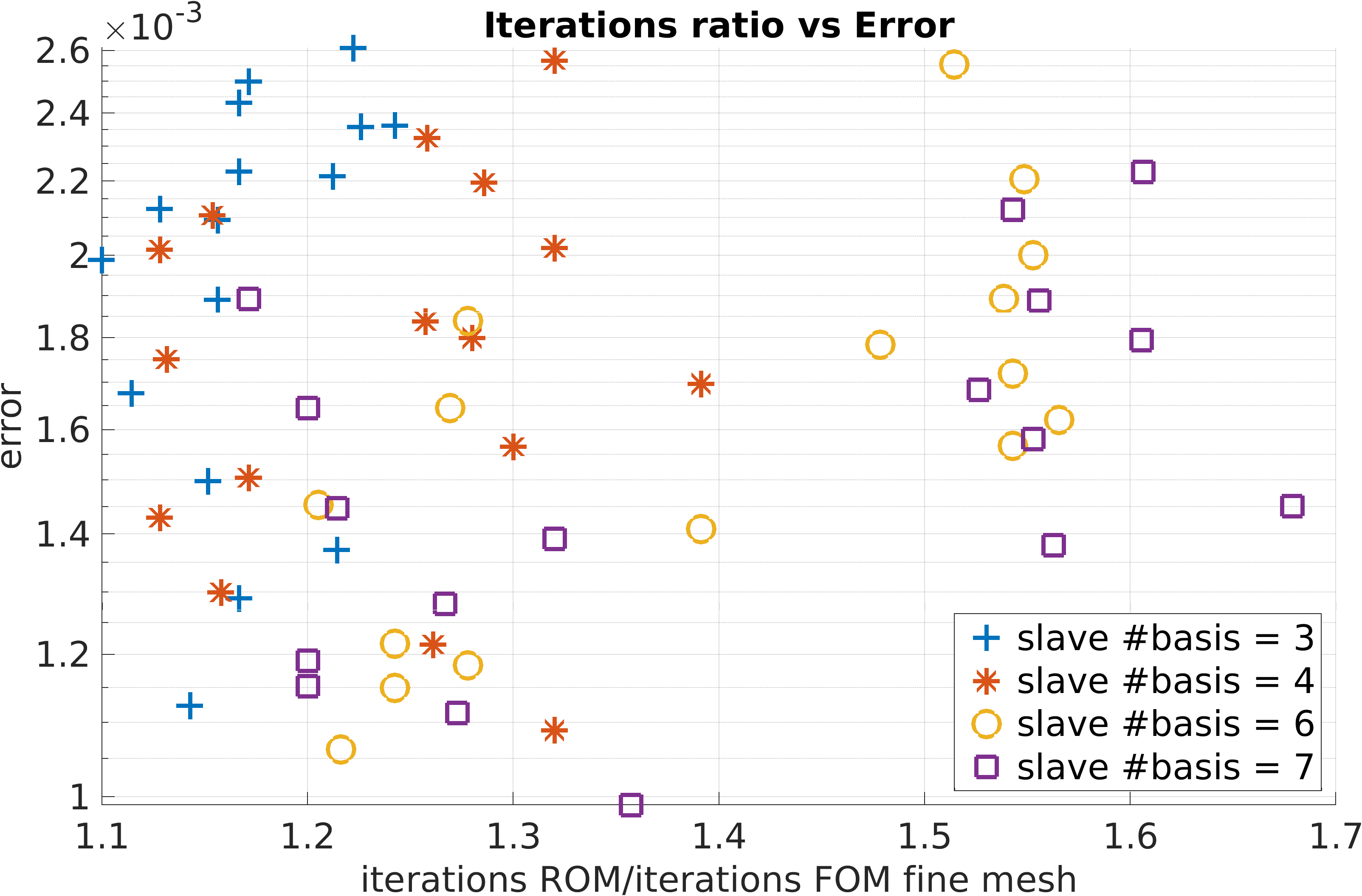}
	\caption{\emph{Test\#1.} Ratio between the number of iterations obtained with ROM and FOM schemes versus the achieved approximation error depending on the number of basis functions employed for the interface data (left) and the slave solution (right).}
	\label{fig:Laplace_ratio_vs_error}
\end{figure}

Fig. \ref{fig:snapshots_laplace_master} and \ref{fig:snapshots_laplace_slave} show three FOM and ROM solutions and the absolute error on the slave and the master domains, respectively, while in Fig. \ref{fig:Laplace_error} we plot the $H^1(\Omega_i)$--norm relative errors for different dimensions of the ROM, \emph{i.e.} for different values of the parameters $n_1$, $n_2$, $M_1$ and $M_2$ defined in Sections \ref{sec:reduced_order_formulation} and \ref{sec:interface_reduction}. After some tests, we chose to consider the same number of basis functions to get a fixed accuracy of the solution for the Dirichlet and Neumann data \emph{i.e.} $M_1 = M_2$, while we treat independently the number of basis functions $n_1$ and $n_2$ used for the master and slave reduction. Then, Fig. \ref{fig:Laplace_error} shows an increase, even if small, in the solution accuracy when the number of basis functions of one of the reduced quantities involved in the procedure is increased, as expected from RB theory. Major effect on the approximation accuracy of either master and slave solutions can be observed when $M_1$ and $M_2$ are varied, \emph{i.e.} depending on the level of accuracy required on the interface data approximation (see also Fig.\ref{fig:Laplace_ratio_vs_error}, left). A more complete analysis of the effect of varying the number of basis functions for each involved reduced quantity can be found in Appendix \ref{appendix:analysis_case_test_1}.

\begin{figure}[h!]
	\centering
    \includegraphics[width=0.45\textwidth]{Iteration_ration_vs_basis_interface_fine.png}
    \quad
	\includegraphics[width=0.45\textwidth]{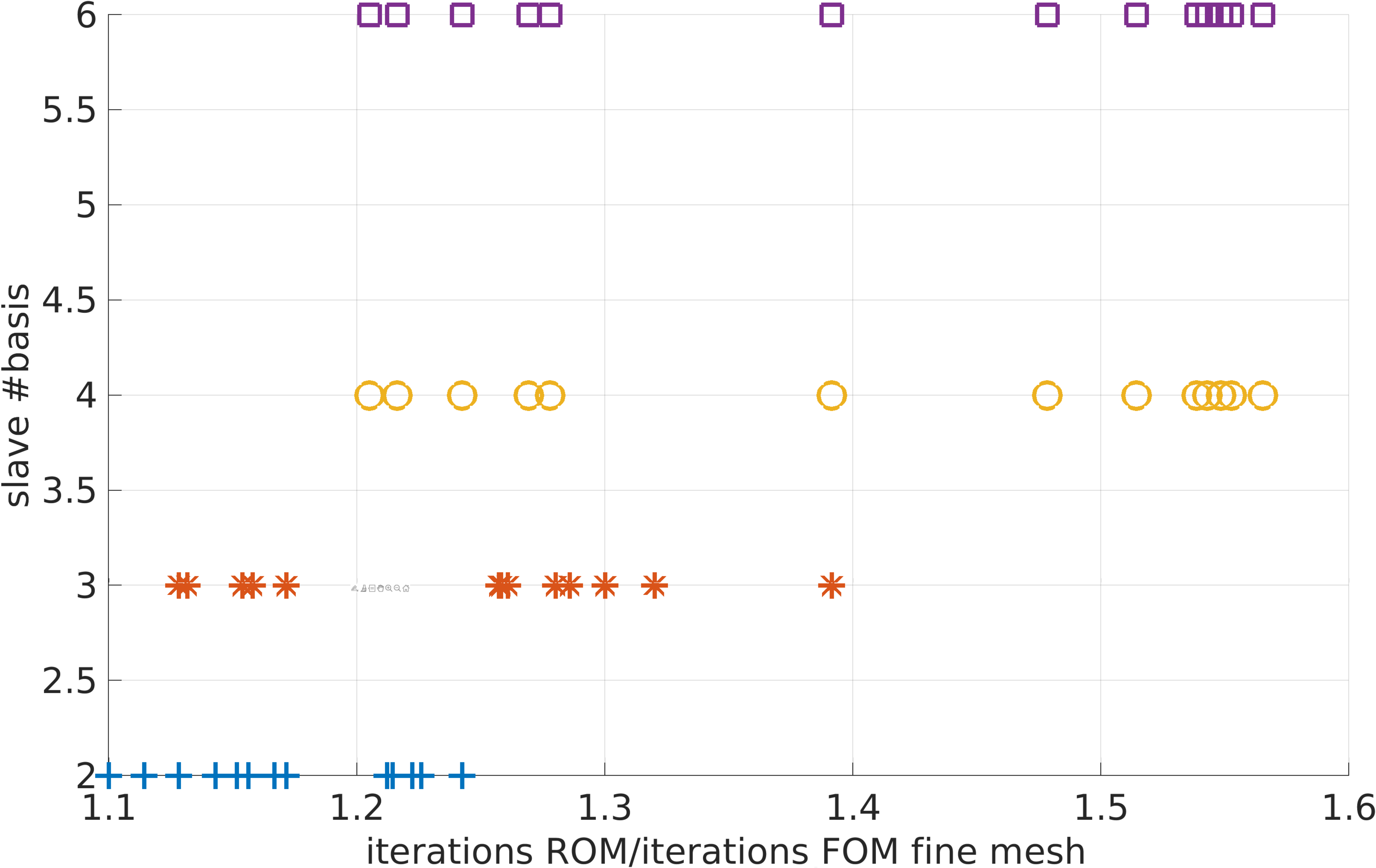}
	\\
	\bigskip
	\includegraphics[width=0.45\textwidth]{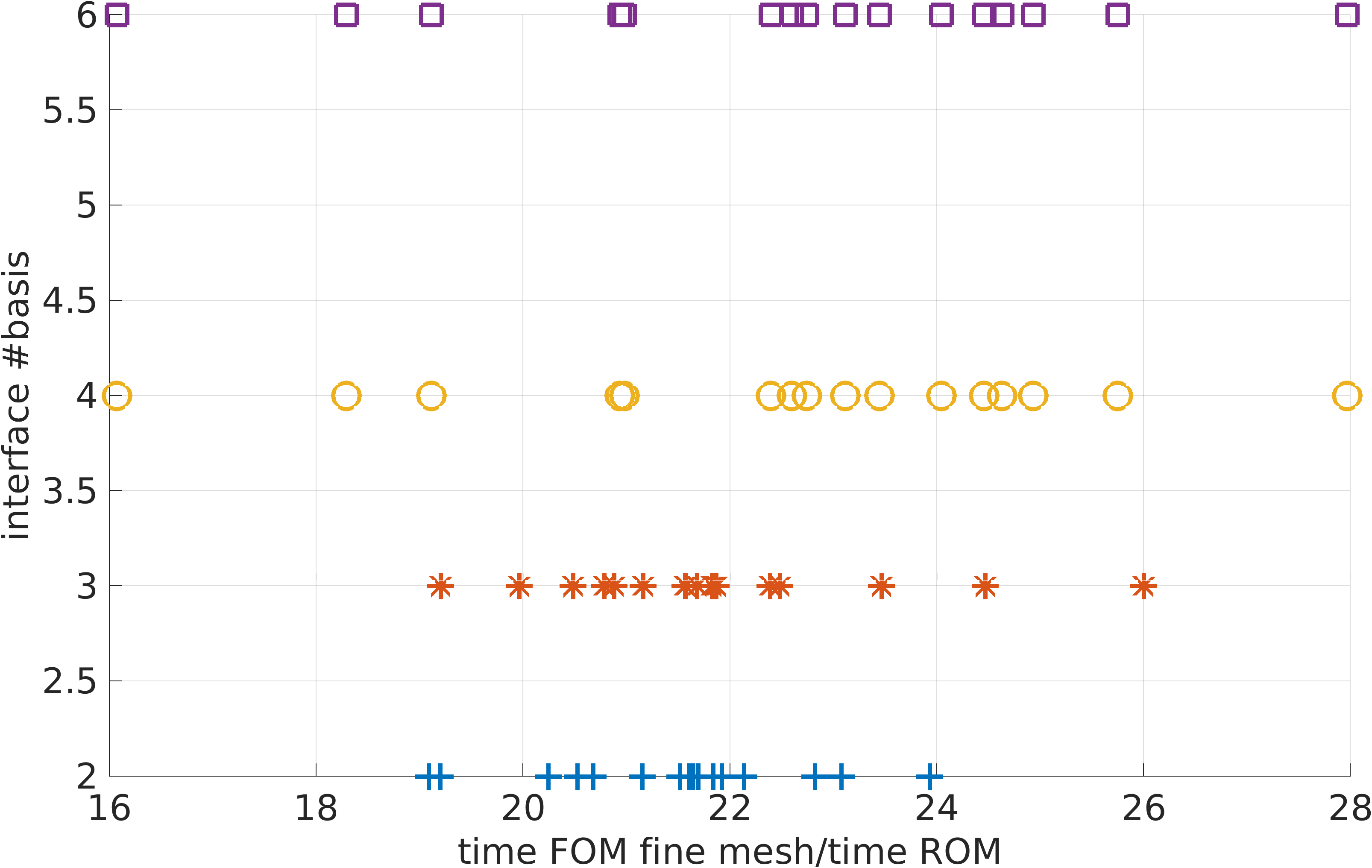}
    \quad
	\includegraphics[width=0.45\textwidth]{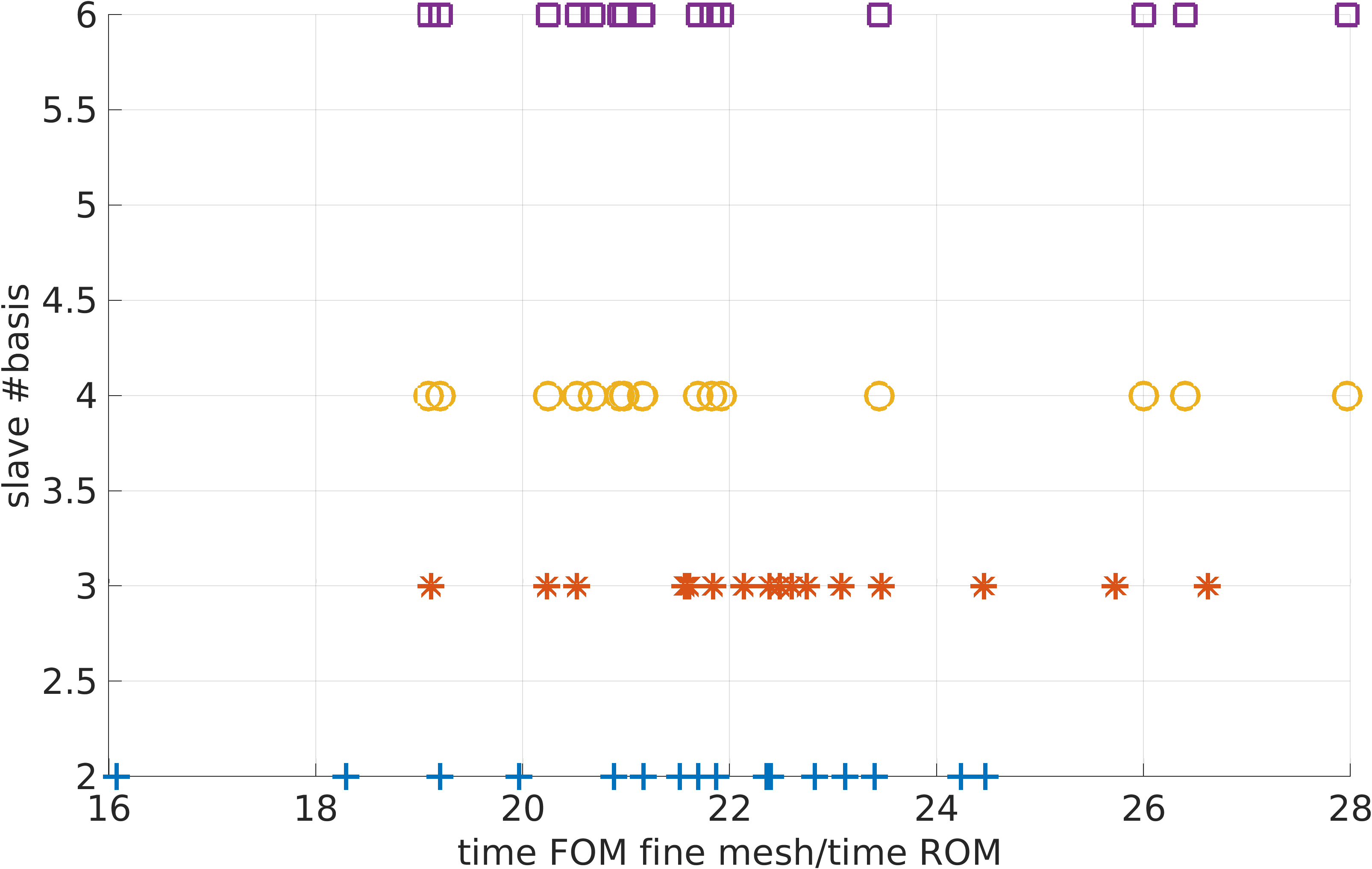}
	\caption{\emph{Test\#1.} Top row: the ratio between the number of iterations obtained with ROM and FOM schemes versus the number of basis functions employed to approximate the interface data (left) and the slave solution (right). Bottom row: the ratio between the FOM and ROM computational time versus the number of basis functions employed to approximate the interface data (left) and the slave solution (right). The FOM simulation refers to the finer discretization.}
	\label{fig:Laplace_ratio_vs_basis}
\end{figure}

Table \ref{Tab:laplace_time} compares the number of DoFs and basis functions employed in each FOM and ROM simulation. On average, a FOM solution is found after 33 iterations for the coarser discretization and 32 iterations for the finer one, while a ROM solution is found after 42 iterations. Since, however, the number of iterations to reach convergence of the interface solutions depends on the model parameter instance, we investigate the overall iterations trends of the ROM model versus the FOM ones by computing the average ratio between the number of iterations required by the ROM, and the number of iterations required by the corresponding FOM, with either coarse and fine discretizations. Fig. \ref{fig:Laplace_ratio_vs_basis} (top row) shows that increasing the number of basis functions to approximate the interface data, the number of iterations to reach the solution convergence decreases, whereas increasing the number of basis functions for the slave solutions, the number of iterations increases. The iterations ratio is overall not influenced by the master hyper--parameters. The same behavior can be also observed by plotting the iterations ratio against the approximation error in Fig. \ref{fig:Laplace_ratio_vs_error}.

\begin{remark}
    The computational expense of the ROM should be compared with that of a FOM utilizing non-conforming meshes in both domains. It is important to note that the FOM computations are performed on conforming meshes -- with the finest discretization employed in the slave domain and the coarsest discretization in the master domain of the ROM simulations, respectively. Moreover, the FOM calculations do not incorporate interface interpolation methods. As a result, comparing the computational times of FOM and ROM simulations would not yield a fair assessment. However, the appropriate FOM simulation for comparison with the ROM should possess a computational expense that falls between the two FOM cases under consideration.
\end{remark}

The same behaviour can be observed in terms of computational costs (see the bottom row of Fig. \ref{fig:Laplace_ratio_vs_basis}), where we show that the speed up gained by employing the ROM instead of the FOM scheme increases when $M_1$ and $M_2$ increase, whereas it is overall the same independently of $n_1$ and $n_2$. Of the complete ROM simulation, about 20\% of the computational cost is required to set up and assemble the parameter--independent matrices and vectors, whereas the remaining 80\% is devoted to the ROM iterations.  7\% of the costs of each ROM iteration depends instead on the extraction and application of Dirichlet interface data, 33\% on extraction, computation, and application of the Neumann interface data, 32\% on the assembling and solution of slave ROM, 25\% on the assembling and solution of the master ROM, and 3\% on the computation of the interface residual. Specifically, for an instance of the parameters, with a prescribed accuracy of the solution of $10^{-5}$, the computational costs of solving problem \eqref{Eq:test_case_1} with a coarse (fine) discretization of both subdomains is of about $3.39s$ ($37.16s$), whereas the ROM simulation requires an online time of about $1.55s$, and a total time of $1h$ and $45min$ for training the model. The ROM is thus about 2.5 times faster than the FOM employing the coarser discretization, corresponding to a CPU time reduction of about 61\%, while we are able to achieve a speedup of about 24 times compared to the FOM employing the finer discretization, corresponding to a reduction of 96\% of the computational costs. 



\begin{table}[h!]
	\centering 
	\begin{tabular}{c|ccccc}
		\toprule
		&&Master solution &Slave solution &Master interface &Slave interface\\
		\hline &&&&&\\[-2ex]
		FOM -- coarse mesh &\#DoFs &3474 &3474 &386 &386 \\
		FOM -- fine mesh &\#DoFs &26146 &26146 &1538 &1538\\
		\hline &&&&&\\ [-2ex]
		\multirow{2}{*}{ROM} &\#Basis &7 &5 &6 &6  \\
        &\#DoFs &26146 &3474 &1538 &386\\
		\bottomrule	
	\end{tabular}
	\caption{\emph{Test\#1.} High fidelity and reduced order model dimensions of subdomains and interface discretization, as well as the number of basis functions required to achieve an approximation error of $10^{-5}$.}
	\label{Tab:laplace_time}
\end{table}	

\subsection{Test\#2. Steady case: diffusion reaction equation with parametrized sources}
\label{Subsect:steady_case_source}
\begin{figure}[h!]
	\centering
	\hspace{-40pt} $\alpha = 6.63, ~\beta = 5.28$ \qquad  \qquad  \qquad	$\alpha = 2.13, ~\beta = 5.73$ \qquad \qquad \qquad	$\alpha = 3.03, ~\beta = 2.58$ \\
    \hspace{-40pt} $\gamma_1 = 3.38, ~\gamma_2 = 1.88$ \qquad  \qquad  \qquad	$\gamma_1 = 8.63, ~\gamma_2 = 1.13$ \qquad \qquad \qquad	$\gamma_1 = 5.63, ~\gamma_2 = 12.38$ \\
	\vspace{7pt}
	\includegraphics[width=0.3\textwidth]{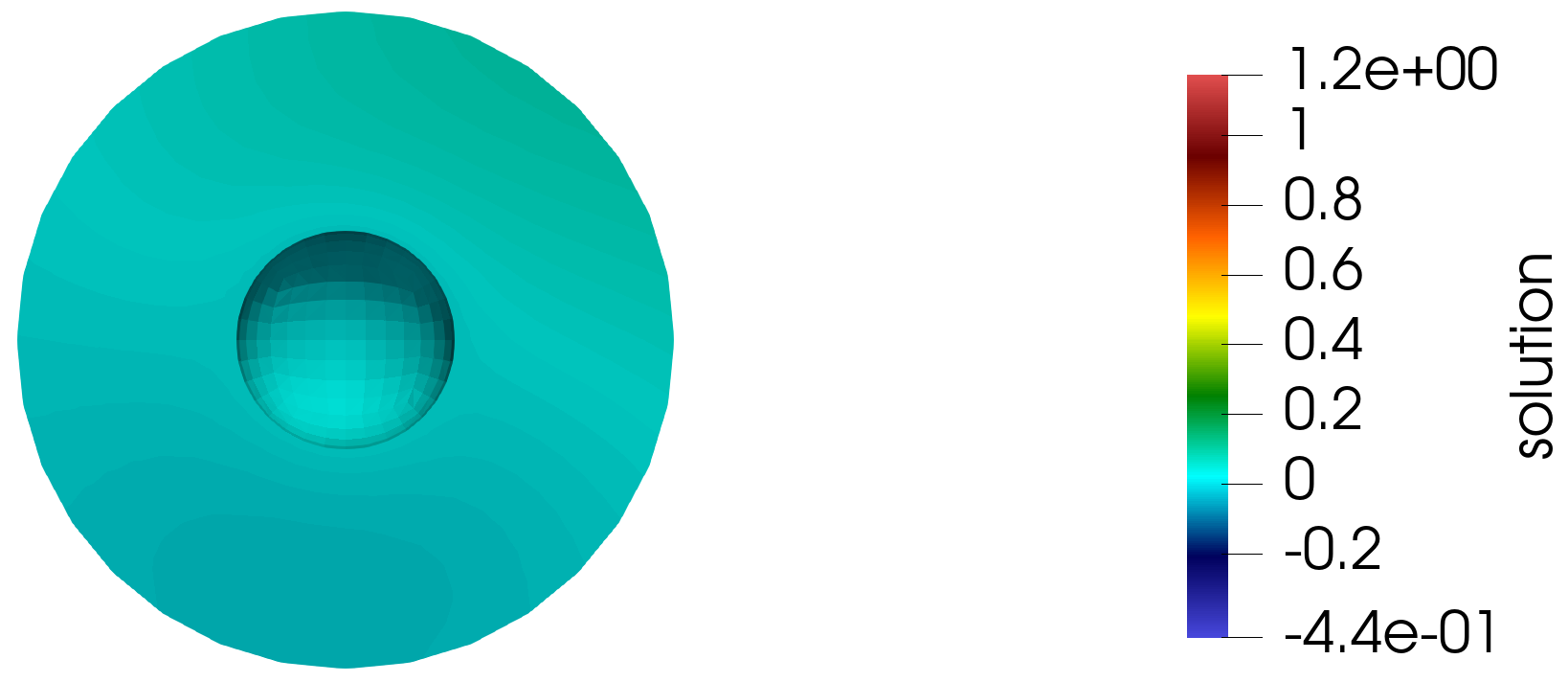}	
	\quad
	\vspace{7pt}
	\includegraphics[width=0.3\textwidth]{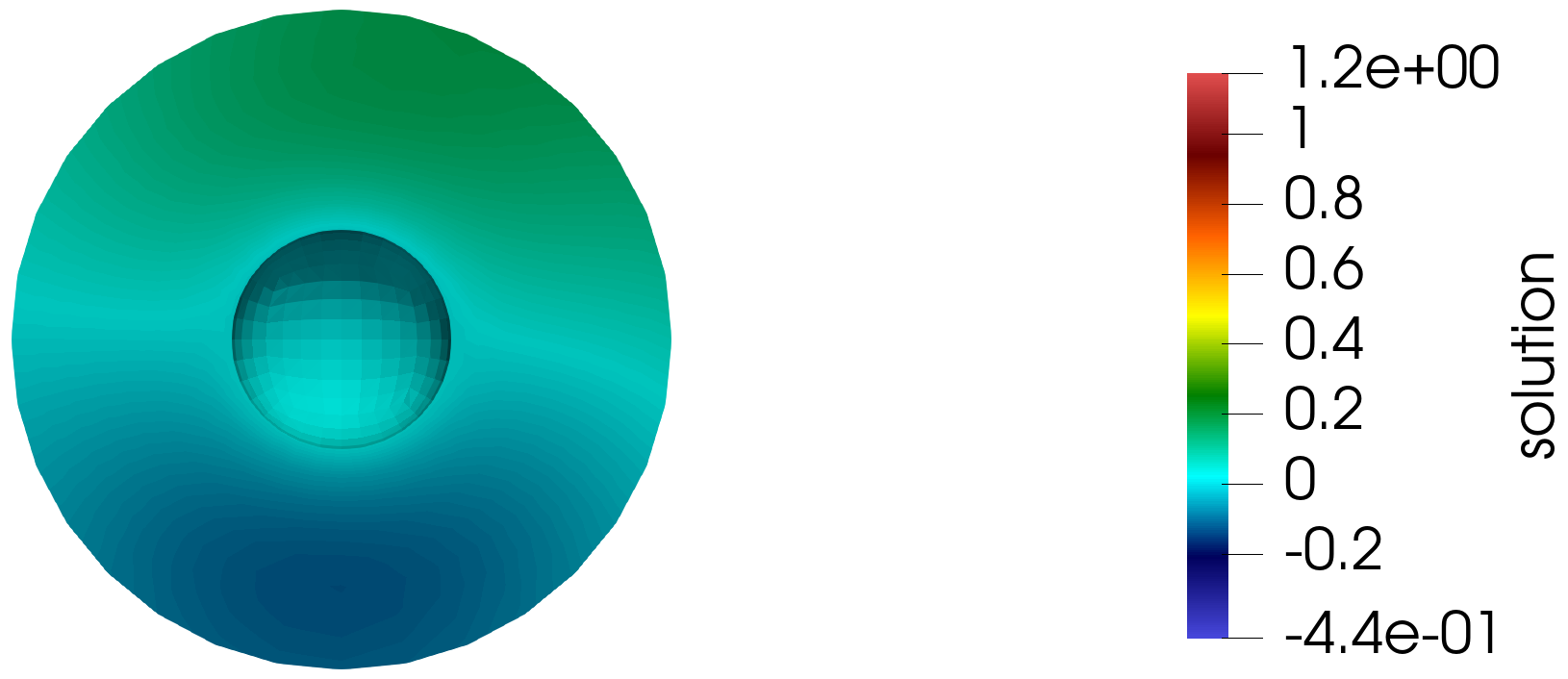}	
	\quad
	\vspace{7pt}
	\includegraphics[width=0.3\textwidth]{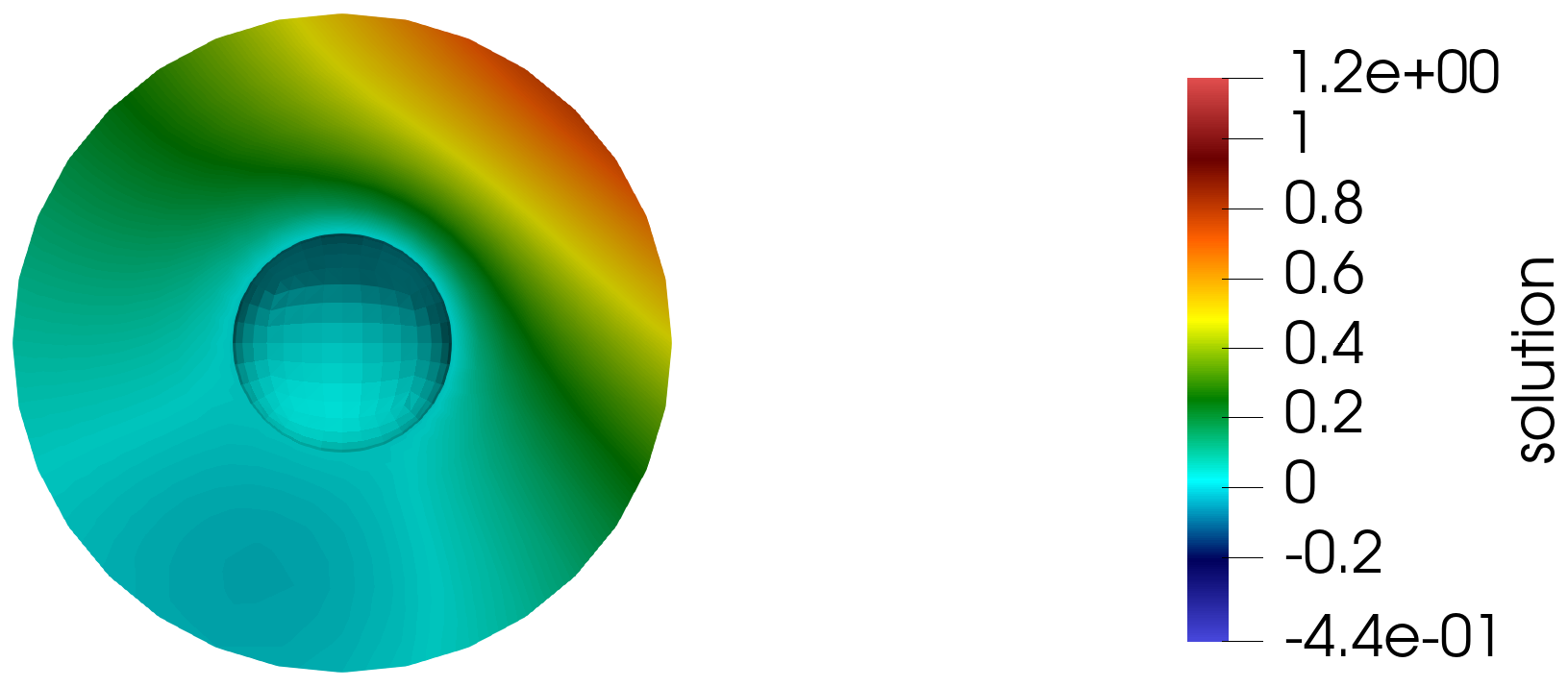}
	  \\
	\includegraphics[width=0.3\textwidth]{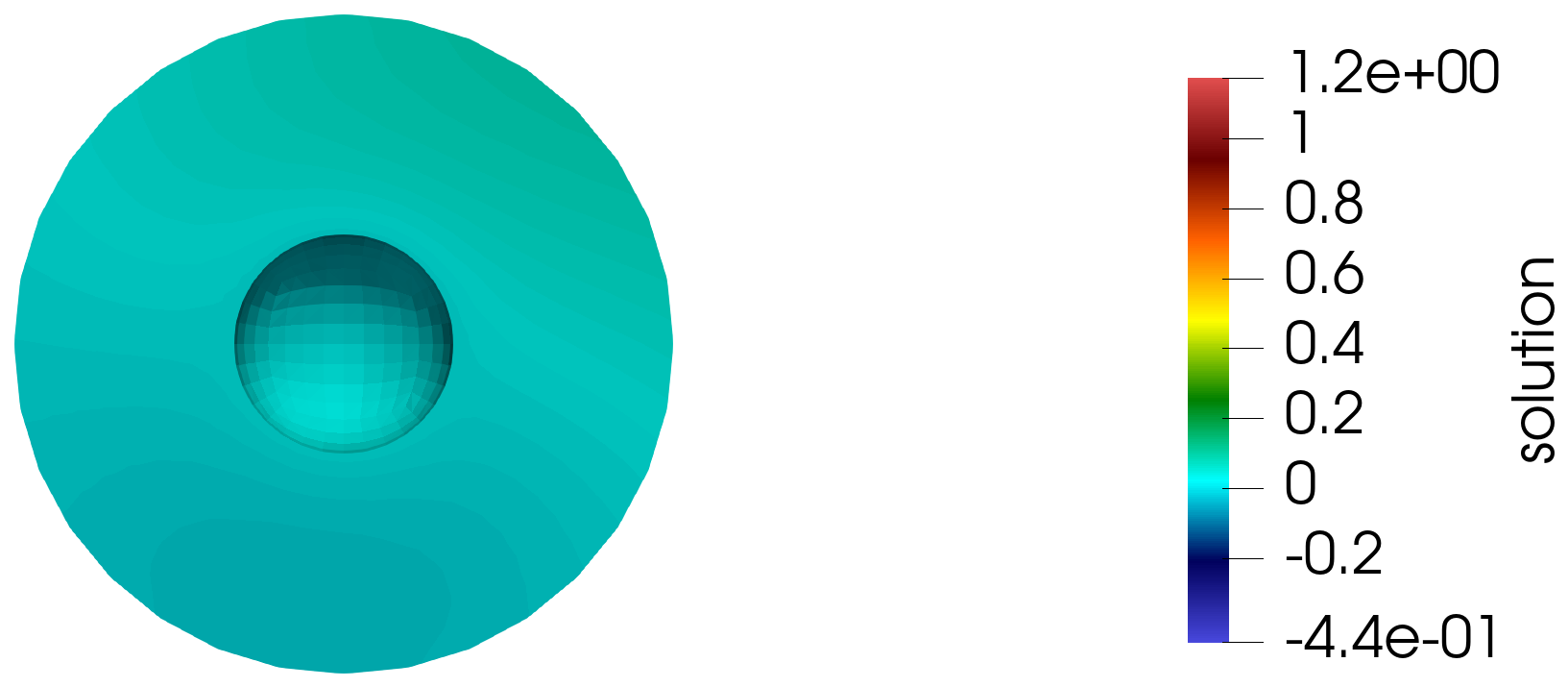}
	\quad
	\includegraphics[width=0.3\textwidth]{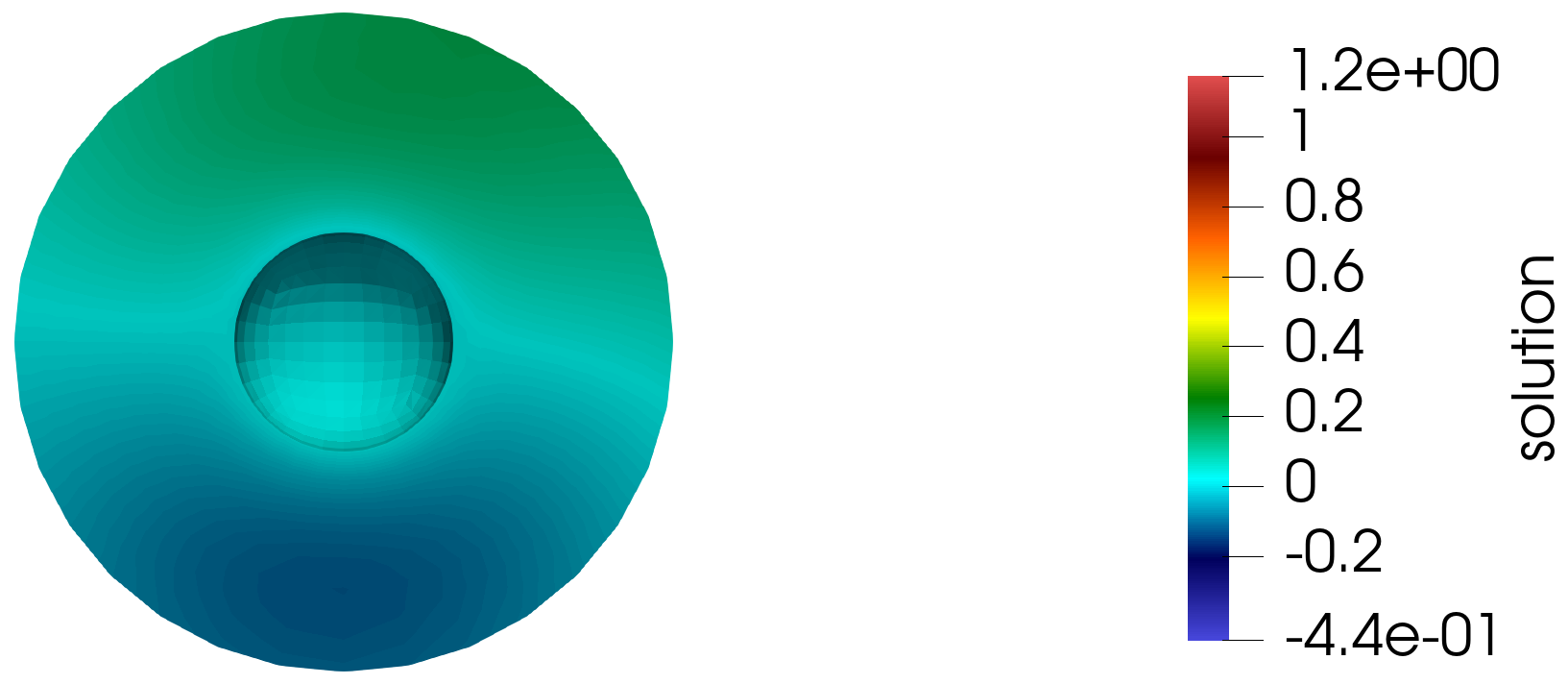}
	\quad
	\includegraphics[width=0.3\textwidth]{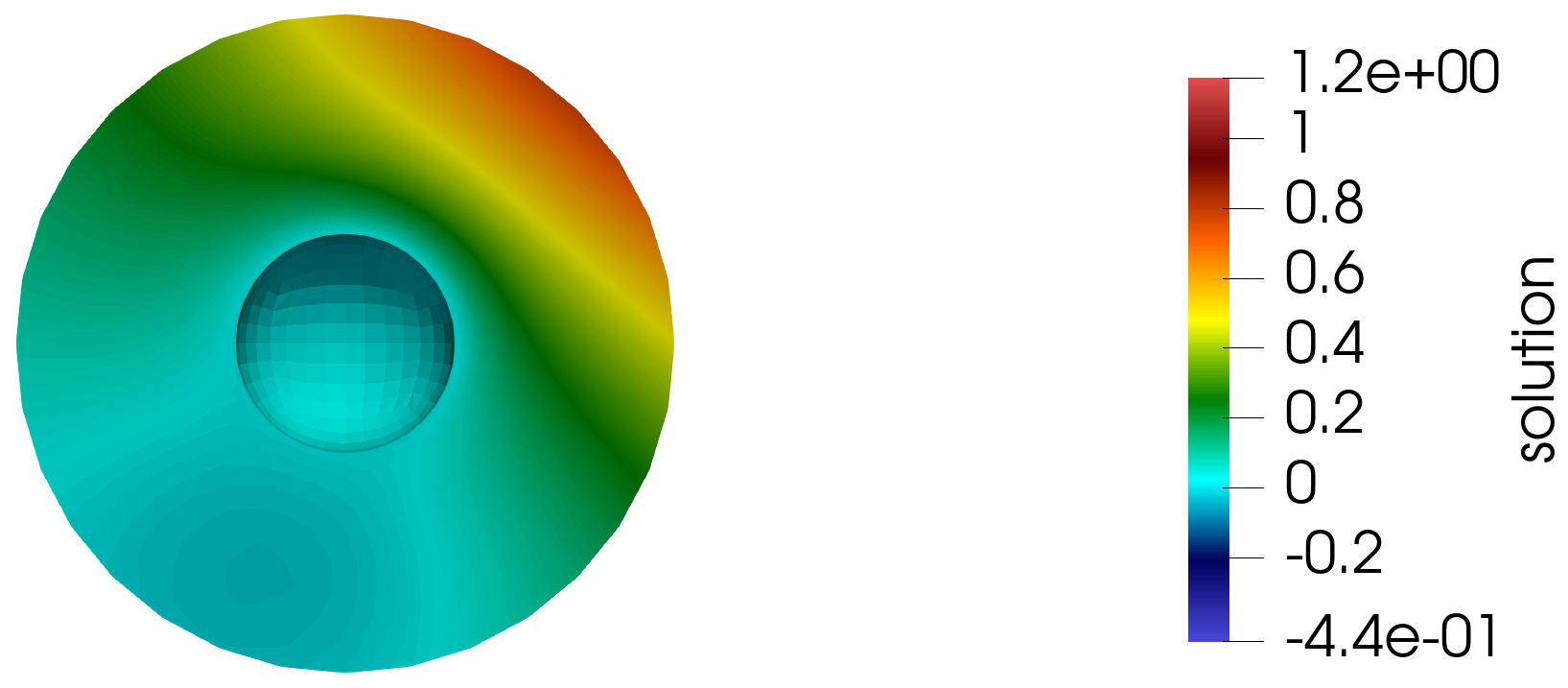}
	\\
	\vspace{12pt}
	\includegraphics[width=0.3\textwidth]{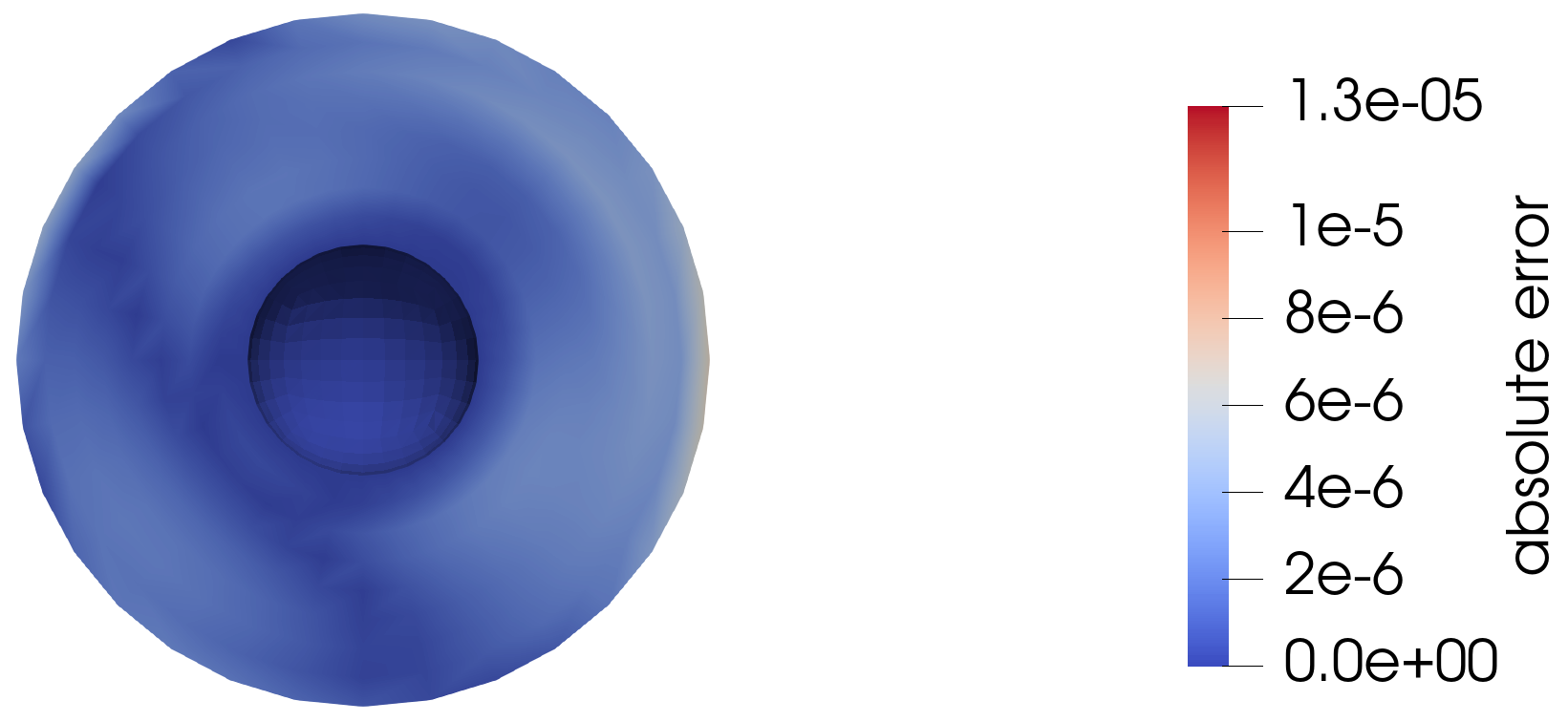}	
	\quad
	\includegraphics[width=0.3\textwidth]{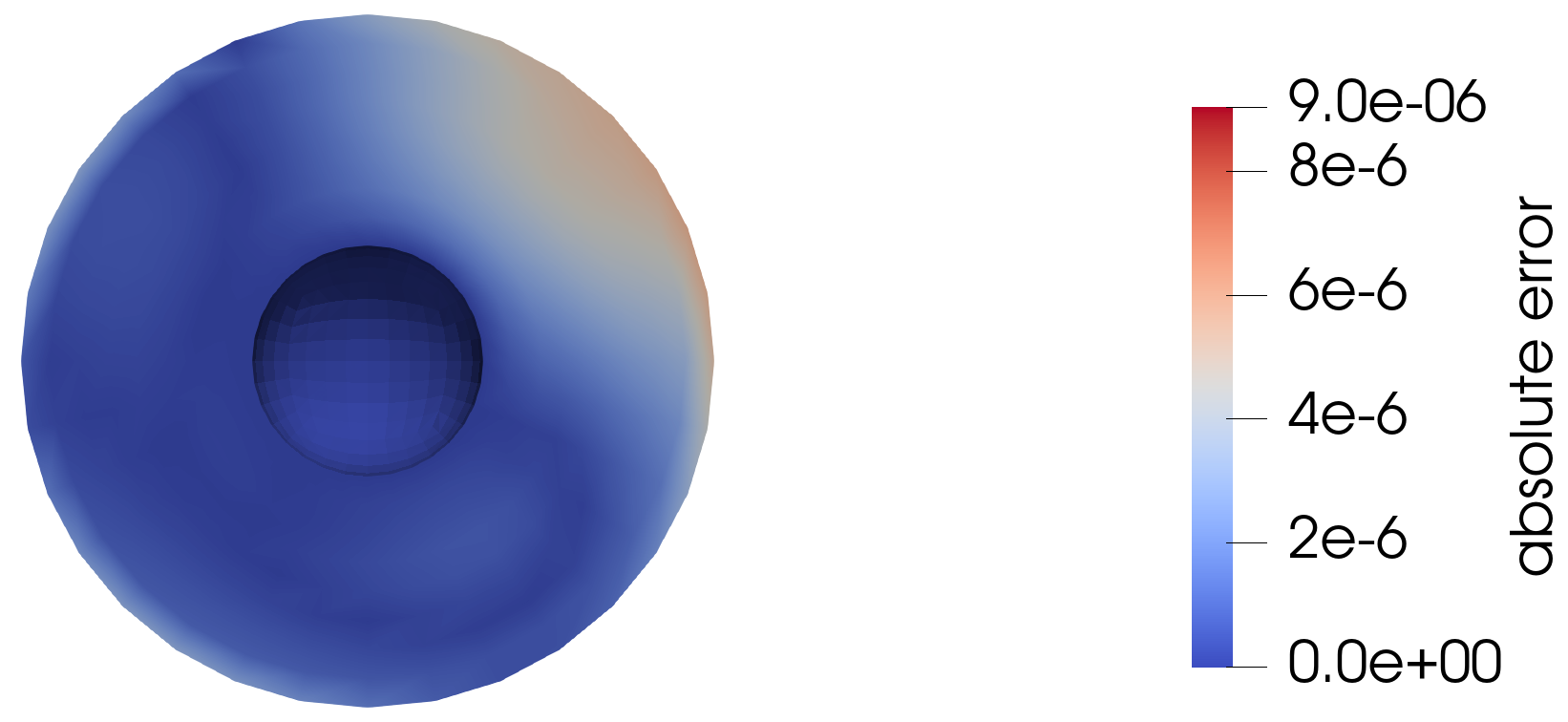}
	\quad
	\includegraphics[width=0.3\textwidth]{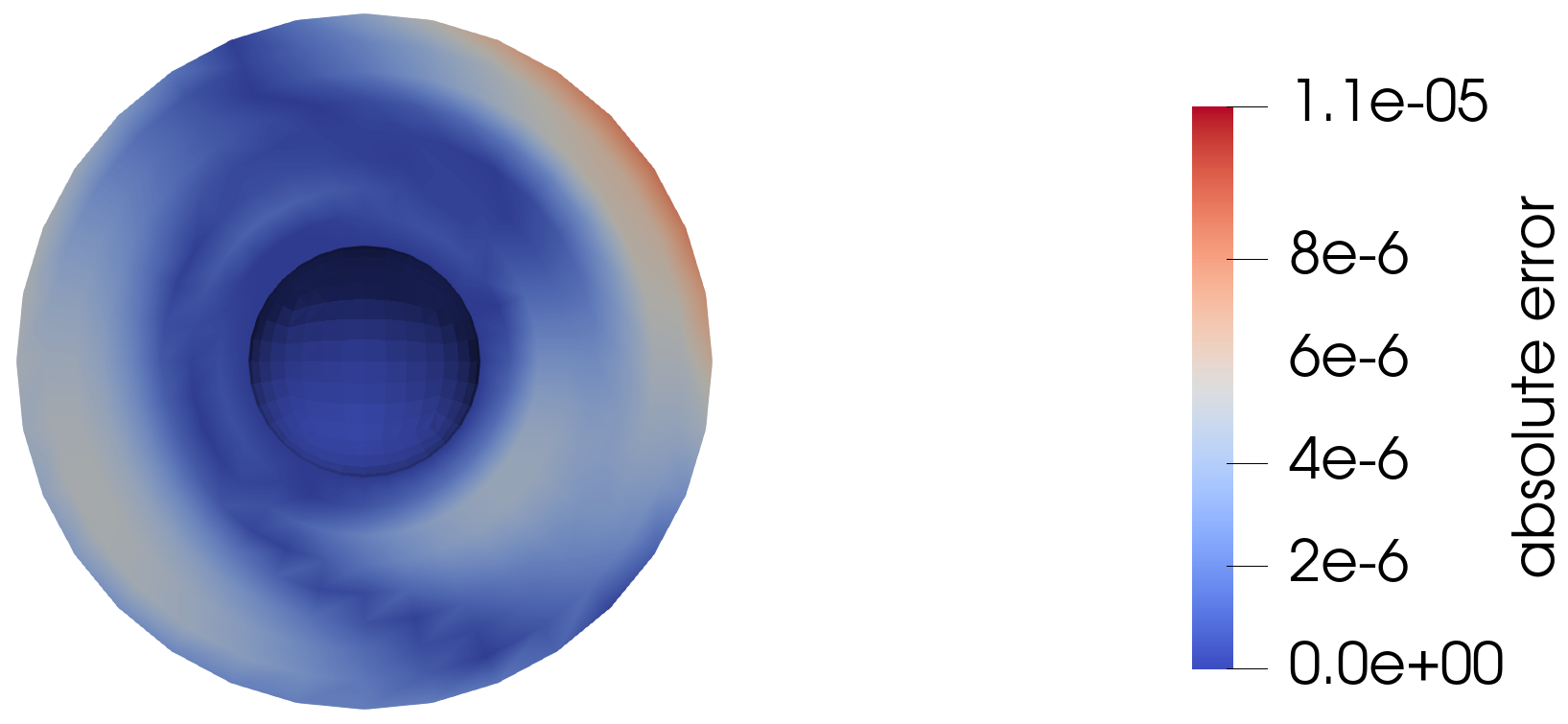}
	\caption{\emph{Test\#2.} Slave solution FOM (top), ROM (center) solutions, and absolute error (bottom) for three different vectors of testing parameters.}
	\label{fig:snapshots_laplace_slave_source}
\end{figure}

\begin{figure}[h!]
	\centering
	\hspace{-40pt} $\alpha = 6.63, ~\beta = 5.28$ \qquad  \qquad  \qquad	$\alpha = 2.13, ~\beta = 5.73$ \qquad \qquad \qquad	$\alpha = 3.03, ~\beta = 2.58$ \\
  \hspace{-40pt} $\gamma_1 = 3.38, ~\gamma_2 = 1.88$ \qquad  \qquad  \qquad	$\gamma_1 = 8.63, ~\gamma_2 = 1.13$ \qquad \qquad \qquad	$\gamma_1 = 5.63, ~\gamma_2 = 12.38$ \\
	\vspace{7pt}
	\includegraphics[width=0.3\textwidth]{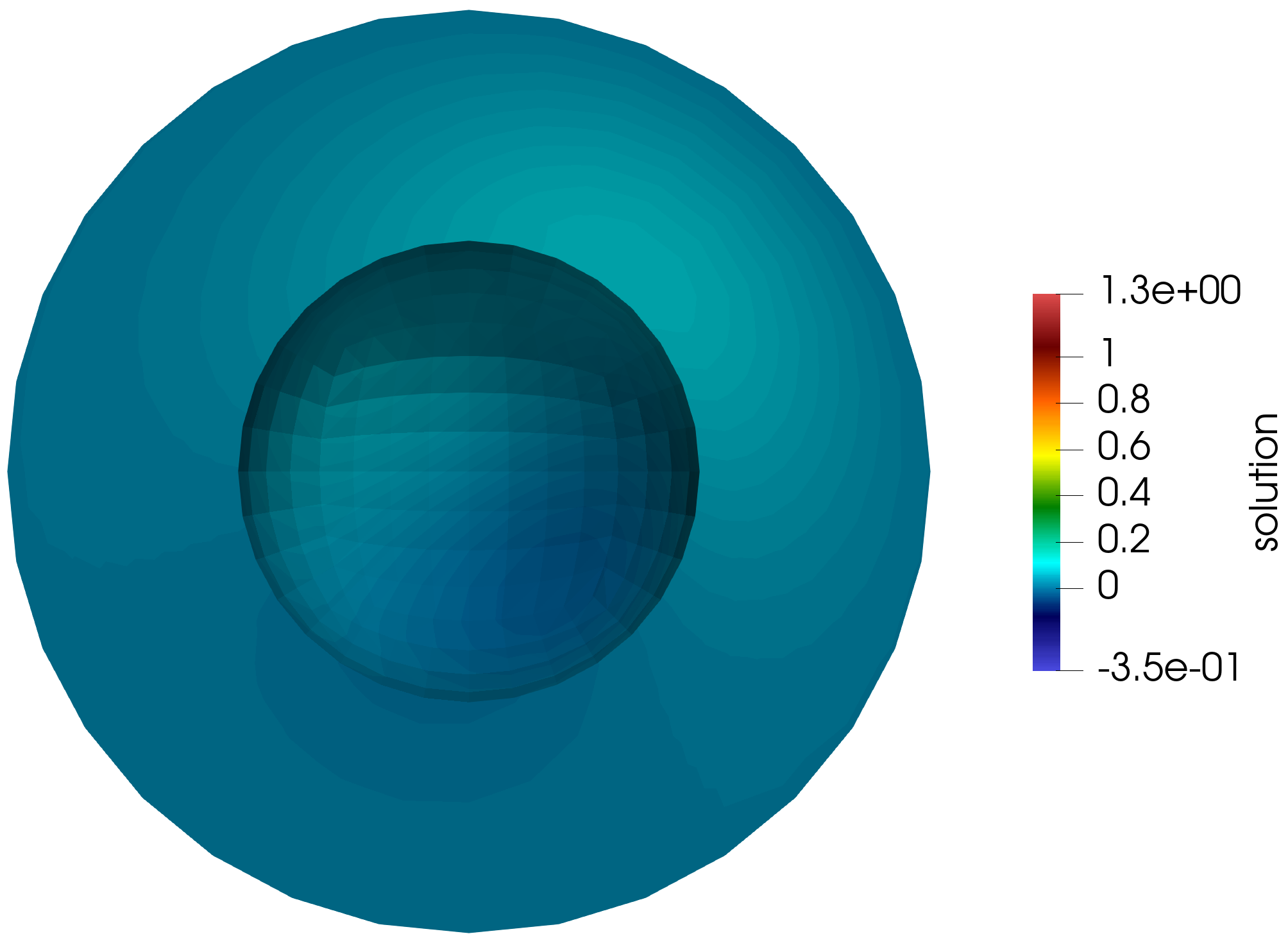}
	\quad
	\vspace{7pt}
	\includegraphics[width=0.3\textwidth]{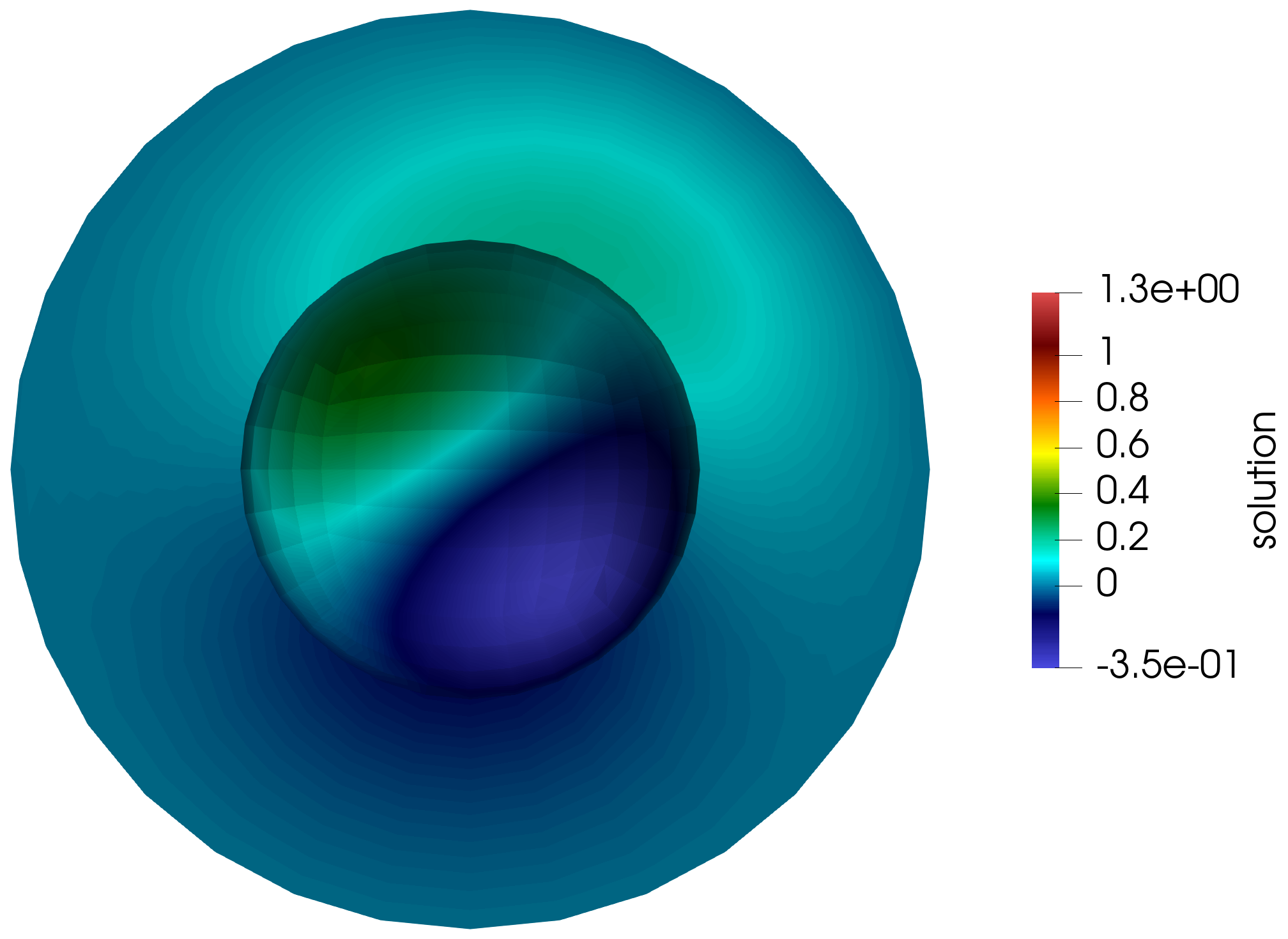}%
	\quad
	\vspace{7pt}
	\includegraphics[width=0.3\textwidth]{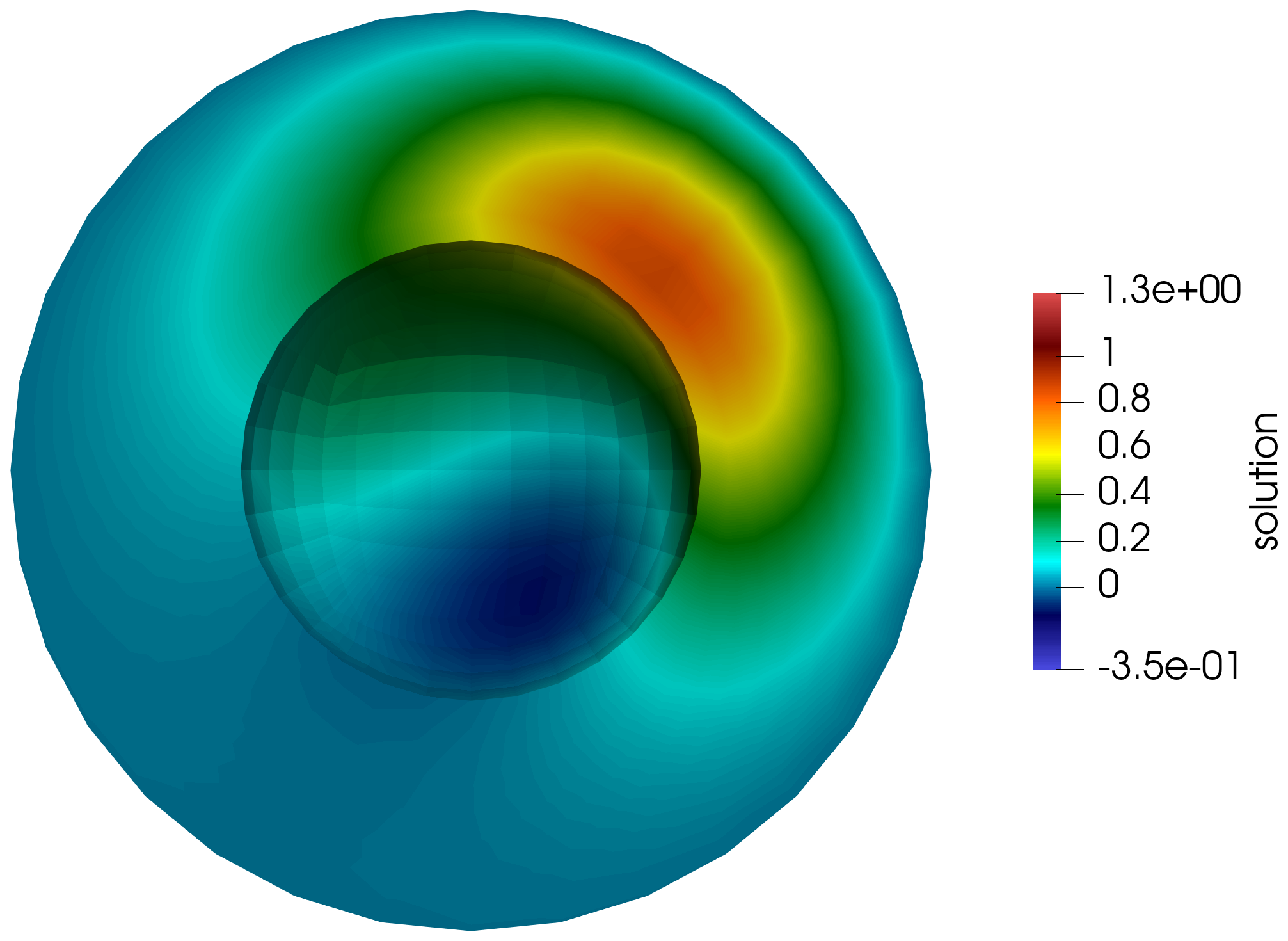}
	\\
	\includegraphics[width=0.3\textwidth]{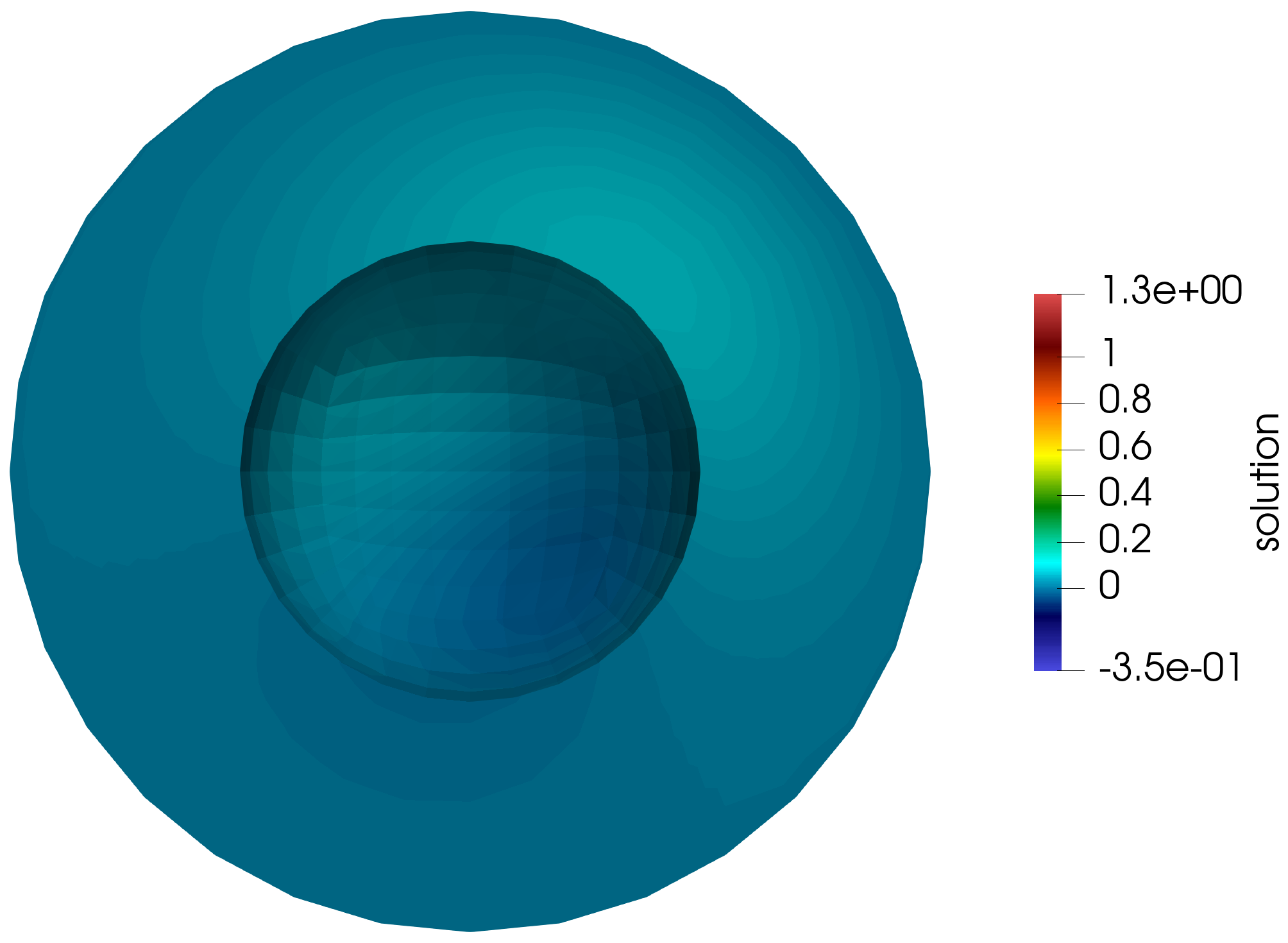}
	\quad
	\includegraphics[width=0.3\textwidth]{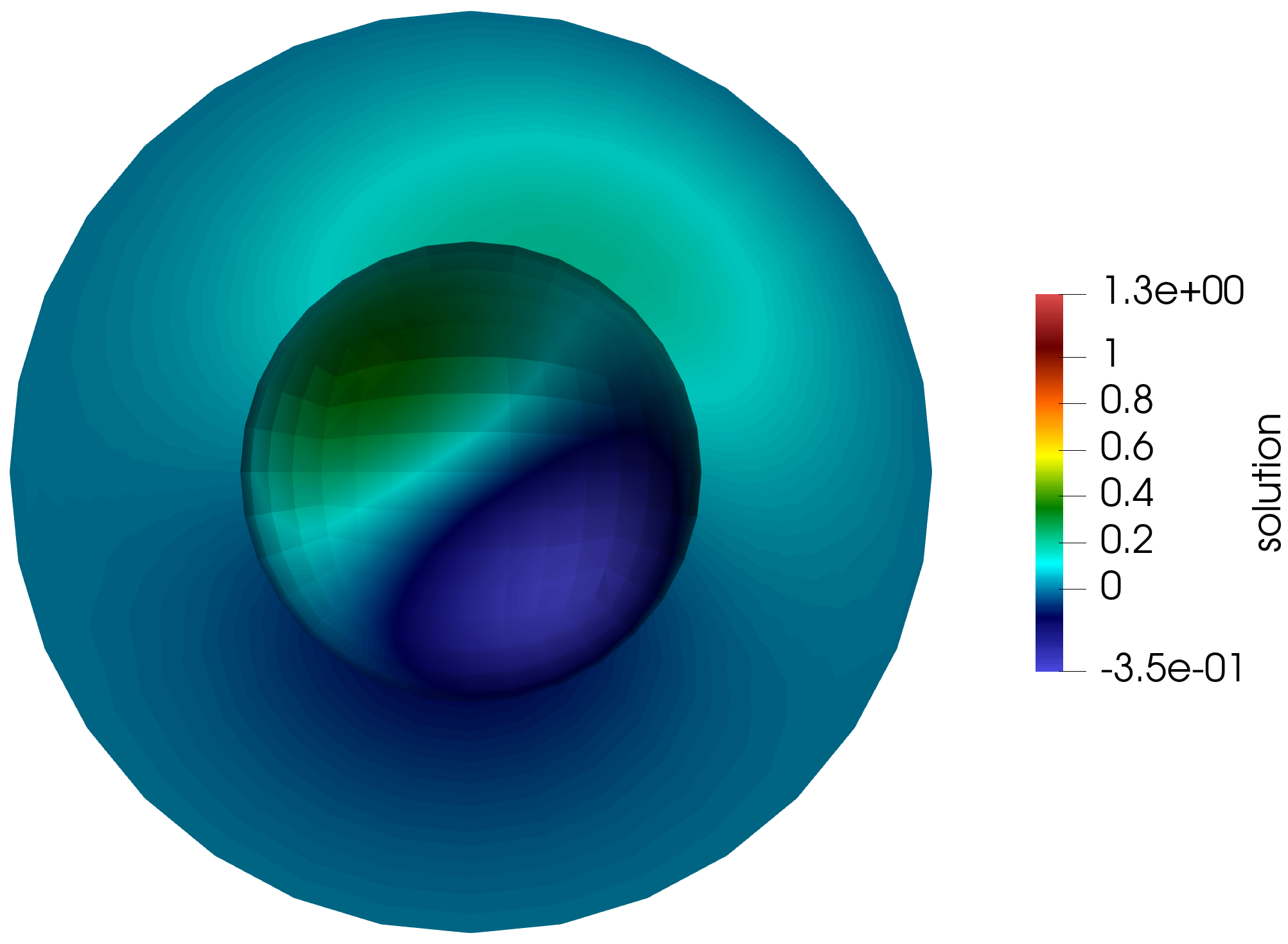}
	\quad
	\includegraphics[width=0.3\textwidth]{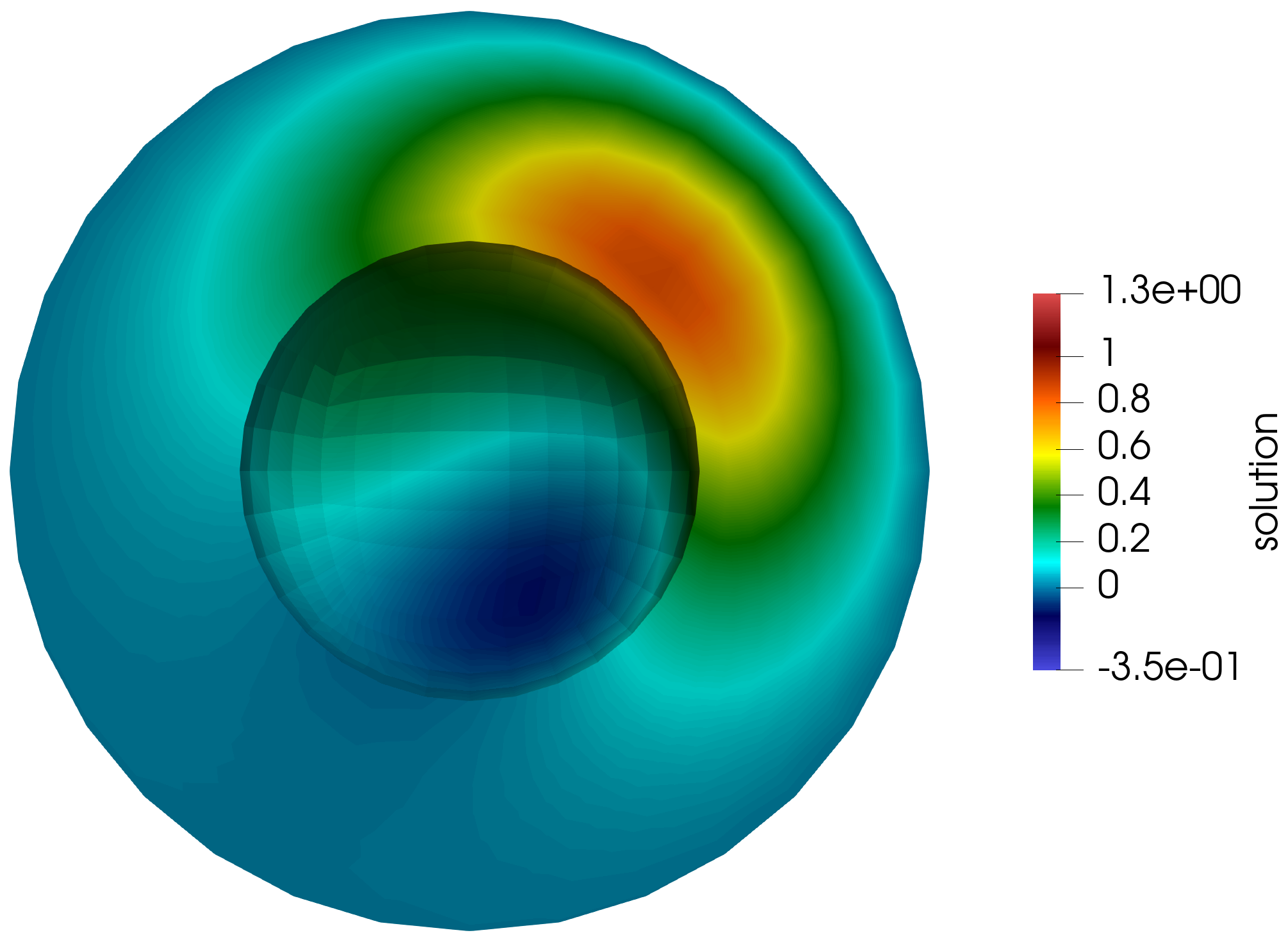}
	\\
	\includegraphics[width=0.3\textwidth]{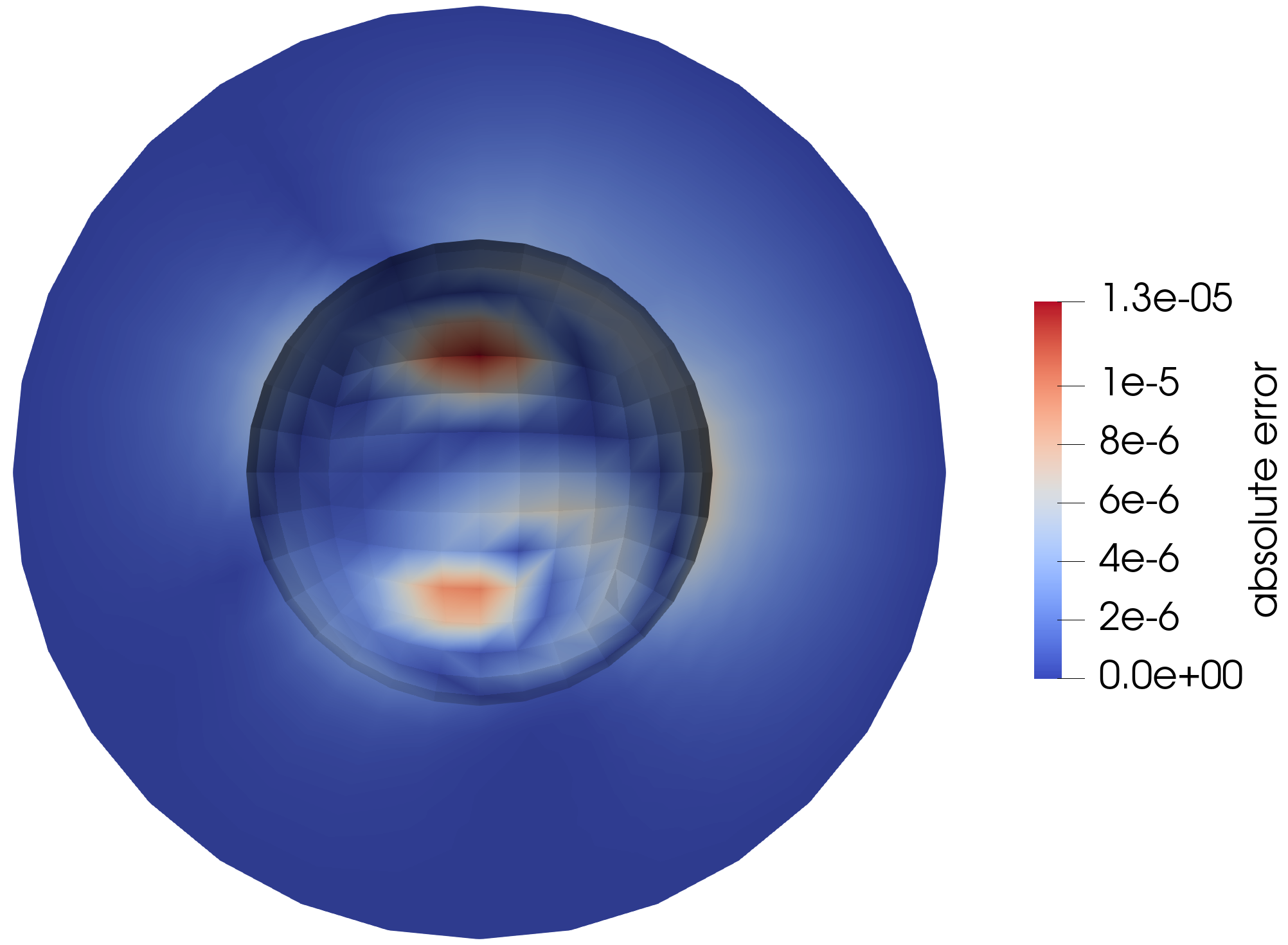}
	\quad
	\includegraphics[width=0.3\textwidth]{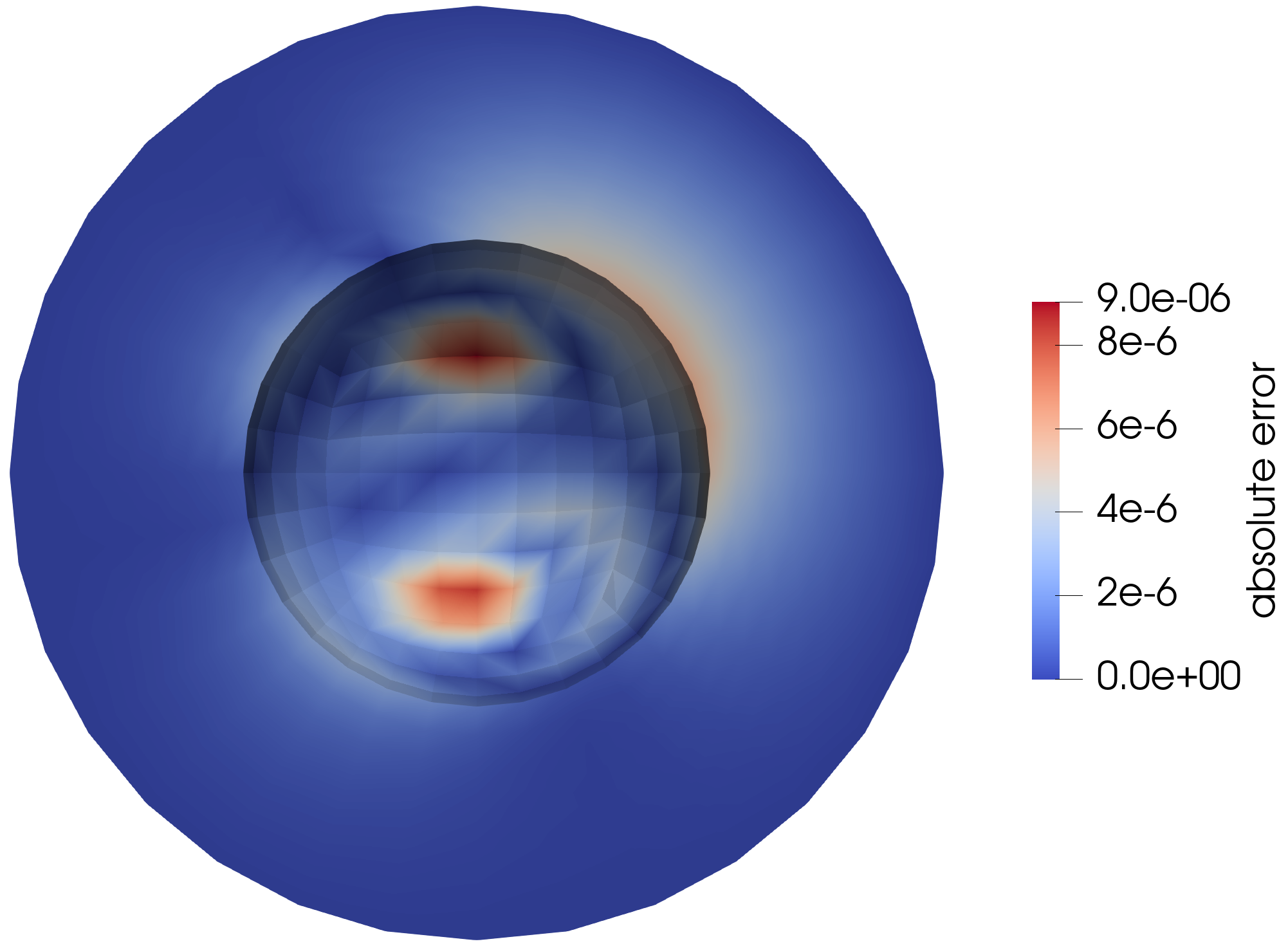}
	\quad
	\includegraphics[width=0.3\textwidth]{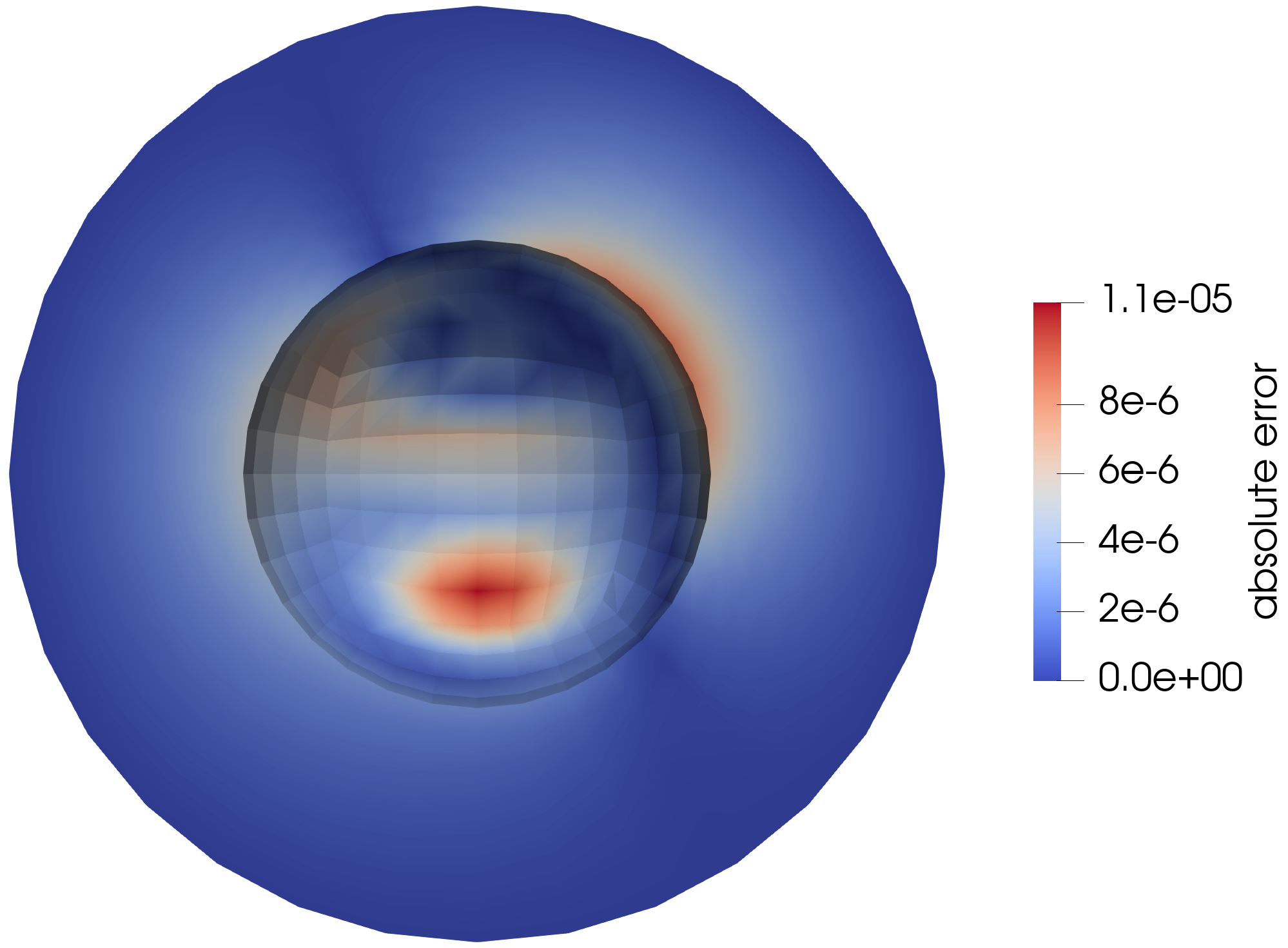}
	\caption{\emph{Test\#2.} Master solution FOM (top), ROM (center) solutions, and absolute error (bottom) for three different vectors of testing parameters.}
	\label{fig:snapshots_laplace_master_source}
\end{figure}

Employing the same modeling framework of Subsection \ref{Subsect:steady_case}, we increase the level of complexity of the model by modifying the right--hand side function $f$ as
\begin{equation}
\label{Eq:source_term_case_3}
f(x,y,z) = \begin{cases}
\gamma_1 \left ( \sin\left(  \frac{\pi}{2} x^2 z\right) + xy \right) &\text{in }\Omega_1, \\
\displaystyle \gamma_2 \text{e}^{-\frac{(x-1)^2+(y-1)^2+(z-1)^2}{2}} &\text{in }\Omega_2,
\end{cases}
\end{equation}
where $\gamma_1$ and $\gamma_2$ represent two varying parameters, defined in the interval $[0,15]$. This $f$ thus represents two different physical sources in $\Omega_1$ and $\Omega_2$, that increase the interface solution mismatch, as well as the gradient of the interface solution, if compared to the previous test case.
Model solutions for different instances of $\boldsymbol{\mu} = [\alpha, \beta,\gamma_1,\gamma_2]$ on either slave or master subdomains can be found in Fig. \ref{fig:snapshots_laplace_slave_source} and \ref{fig:snapshots_laplace_master_source}. 


\begin{figure}[h!]
	\centering
    \includegraphics[width=0.42\textwidth]{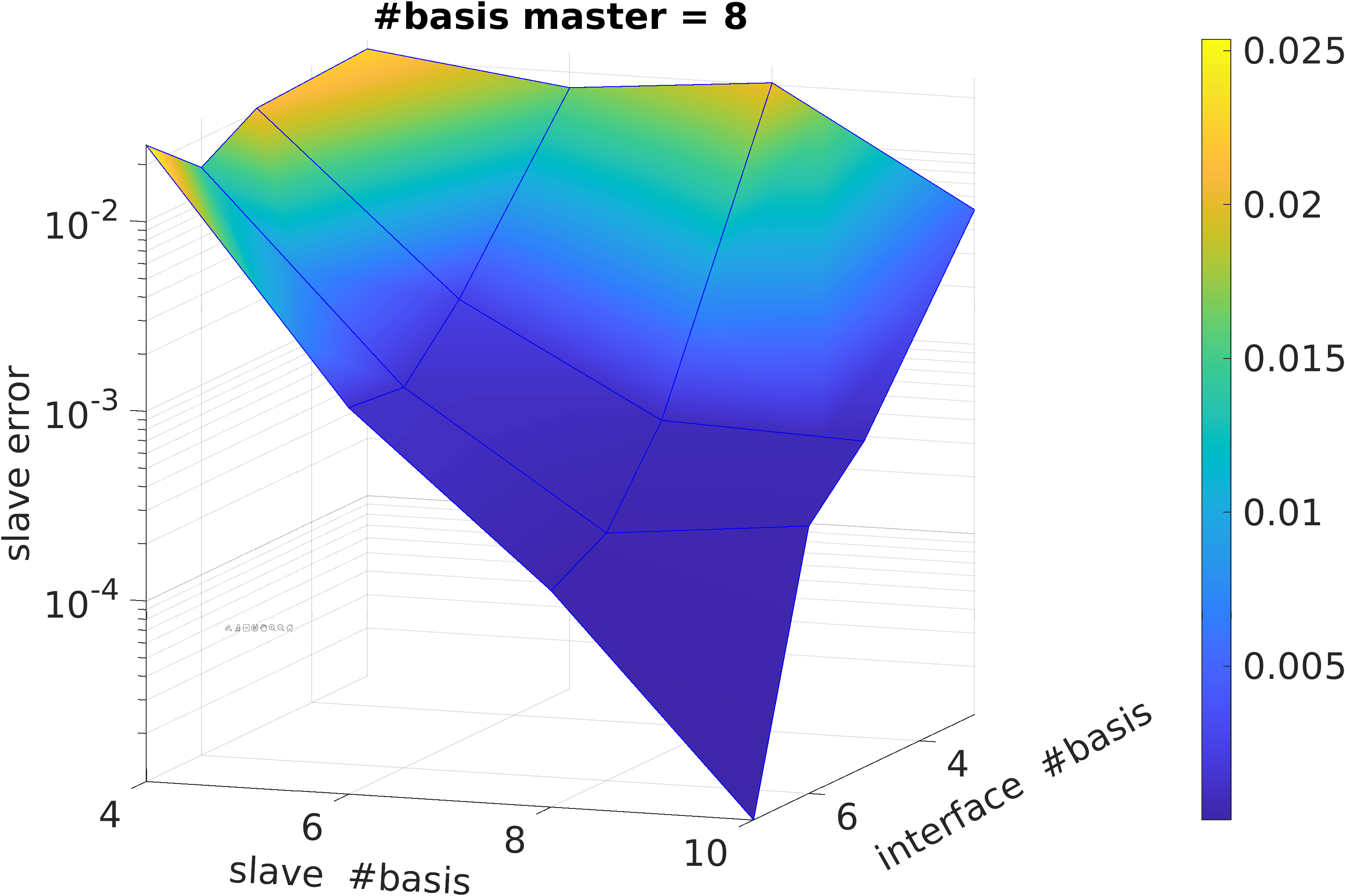}
    \quad
	\includegraphics[width=0.42\textwidth]{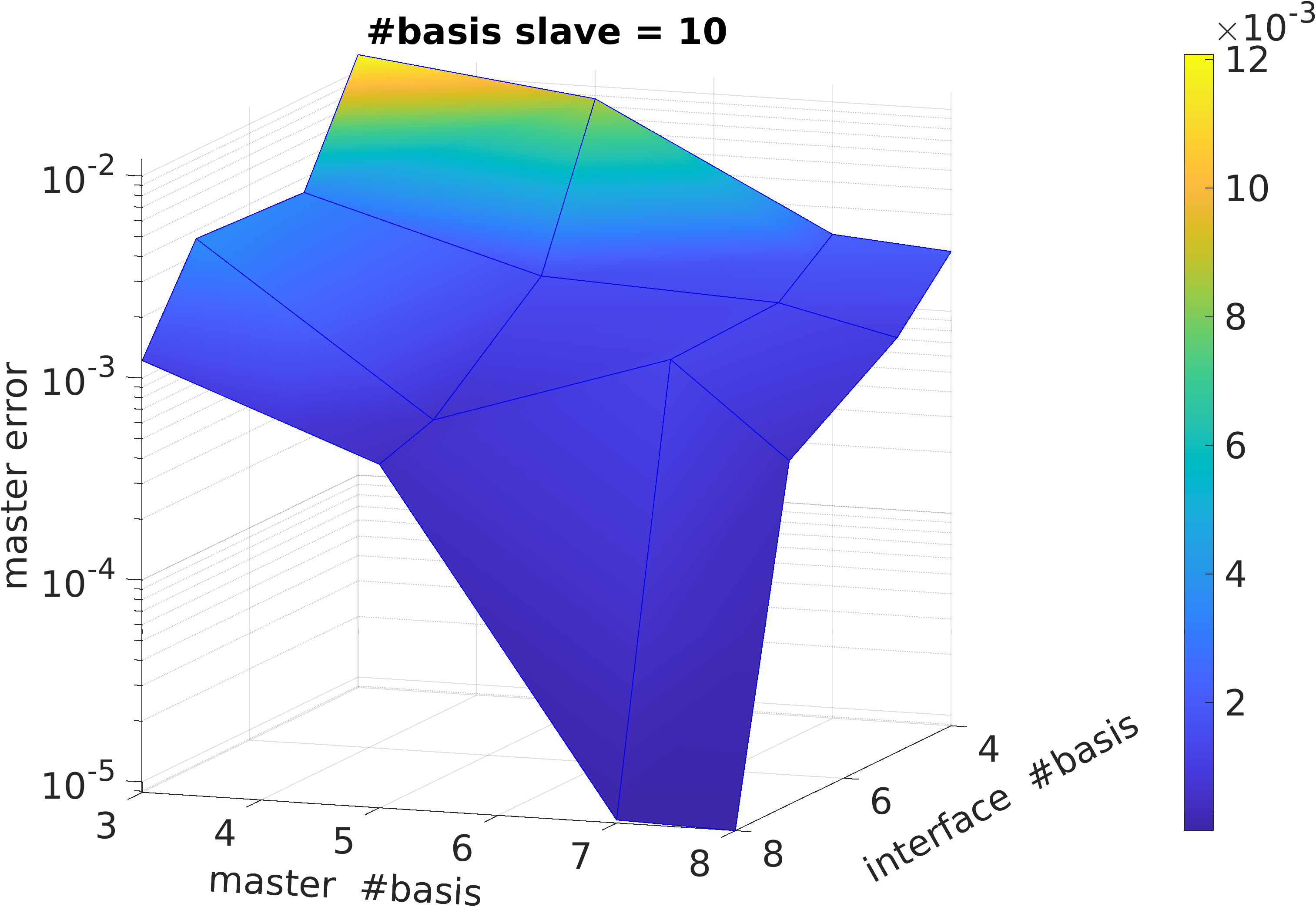}
	\\
	\bigskip
	\includegraphics[width=0.42\textwidth]{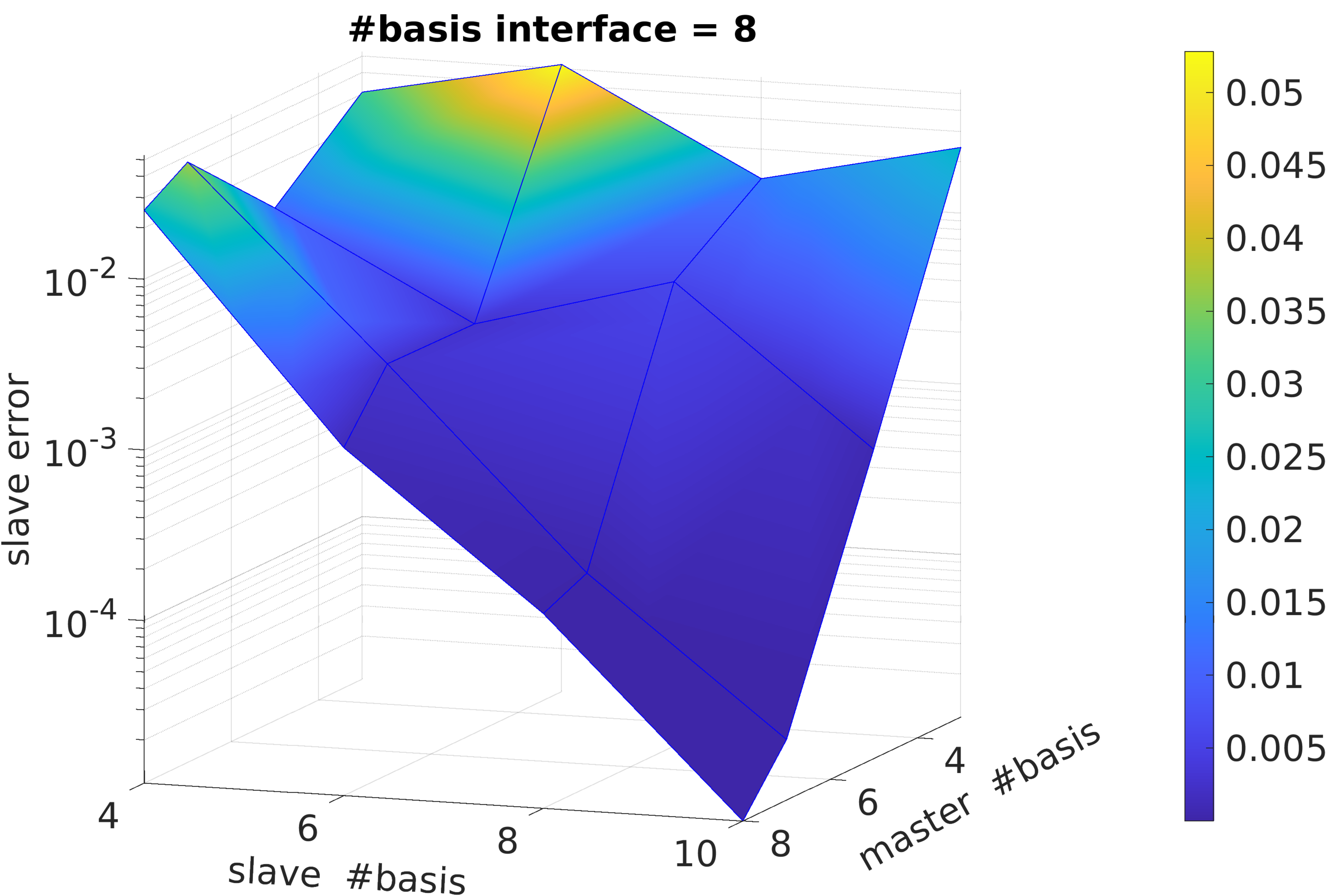}
    \quad
	\includegraphics[width=0.42\textwidth]{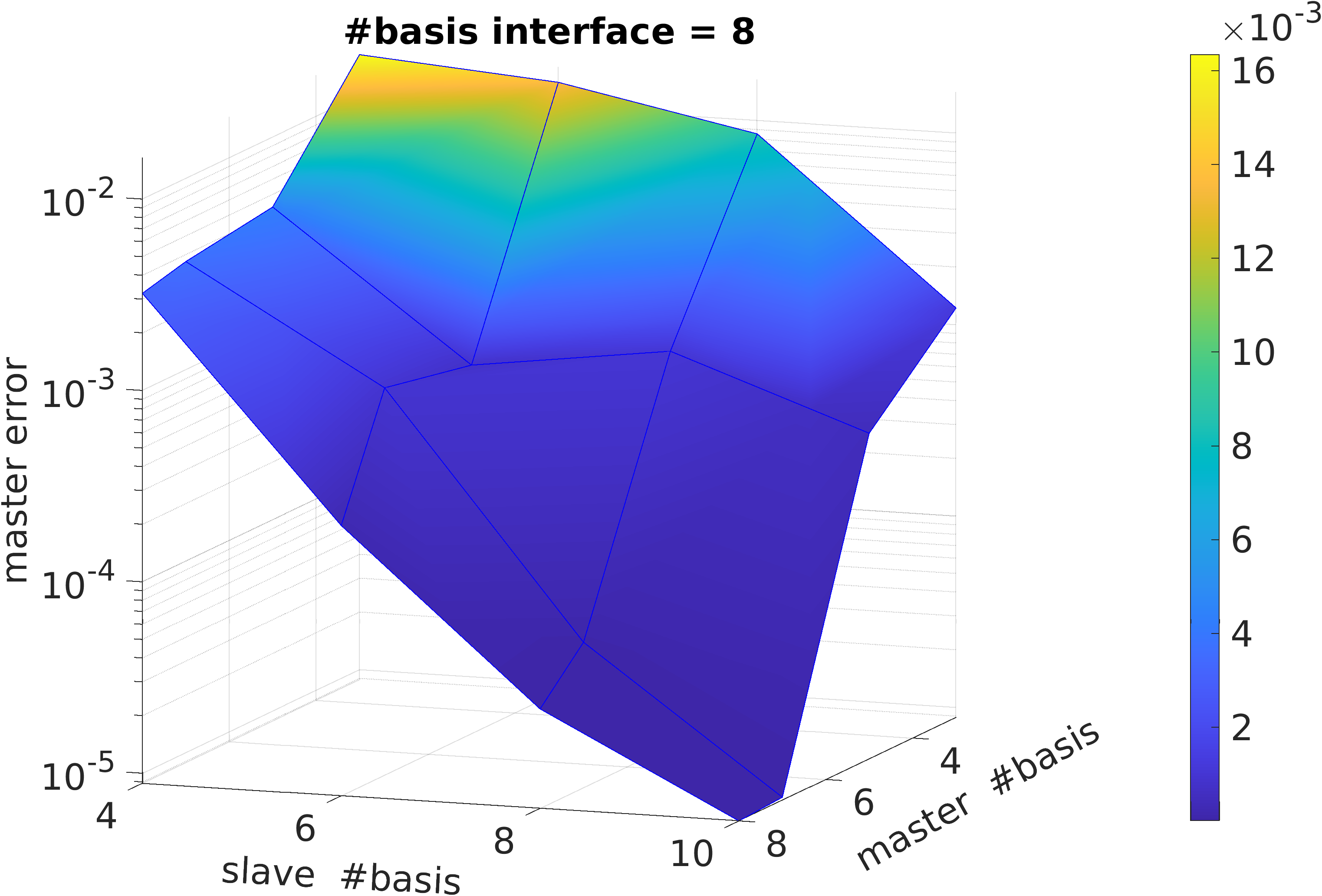}
	\caption{\emph{Test\#2.} $H^1(\Omega_i)$ mean relative error ($z$--axis) over the solution for $N_\text{test} = 20$ different instances of the parameters between the FOM and ROM solutions varying the number of basis functions used to represent the slave and the master solution $n_1$ and $n_2$, and the interface data $M_1$ and $M_2$ ($x$-- and $y$--axis). On the top row, we fix the number of basis functions of the master problem to 8 (on the left) and to 10 for the slave problem (on the right), while on the bottom we fix to 8 the number of basis function for the Dirichlet and Neumann interface data representation.}
	\label{fig:Laplace_error_source}
\end{figure}

\begin{figure}[h!]
	\centering
    \includegraphics[width=0.45\textwidth]{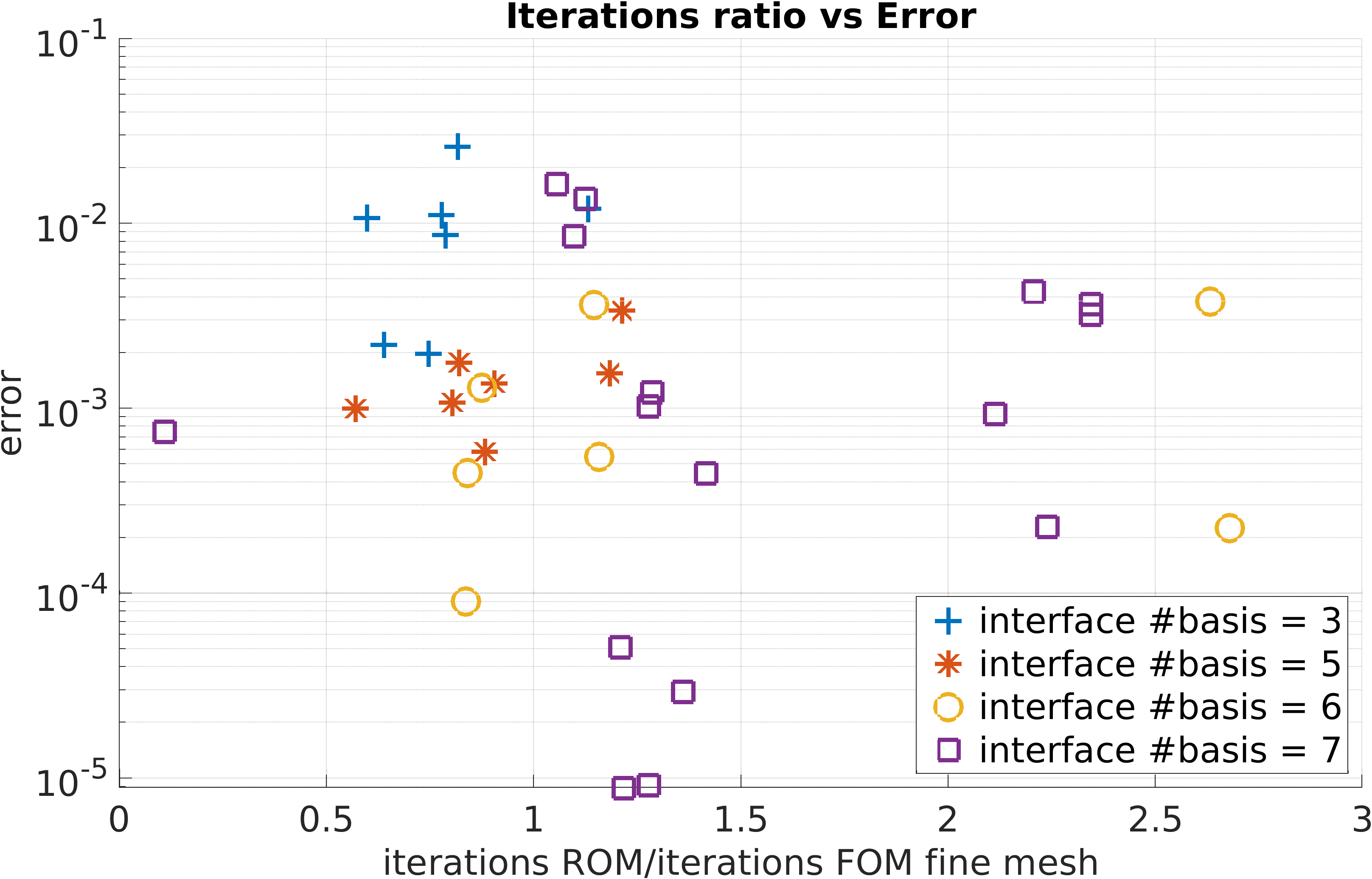}
    \quad
	\includegraphics[width=0.45\textwidth]{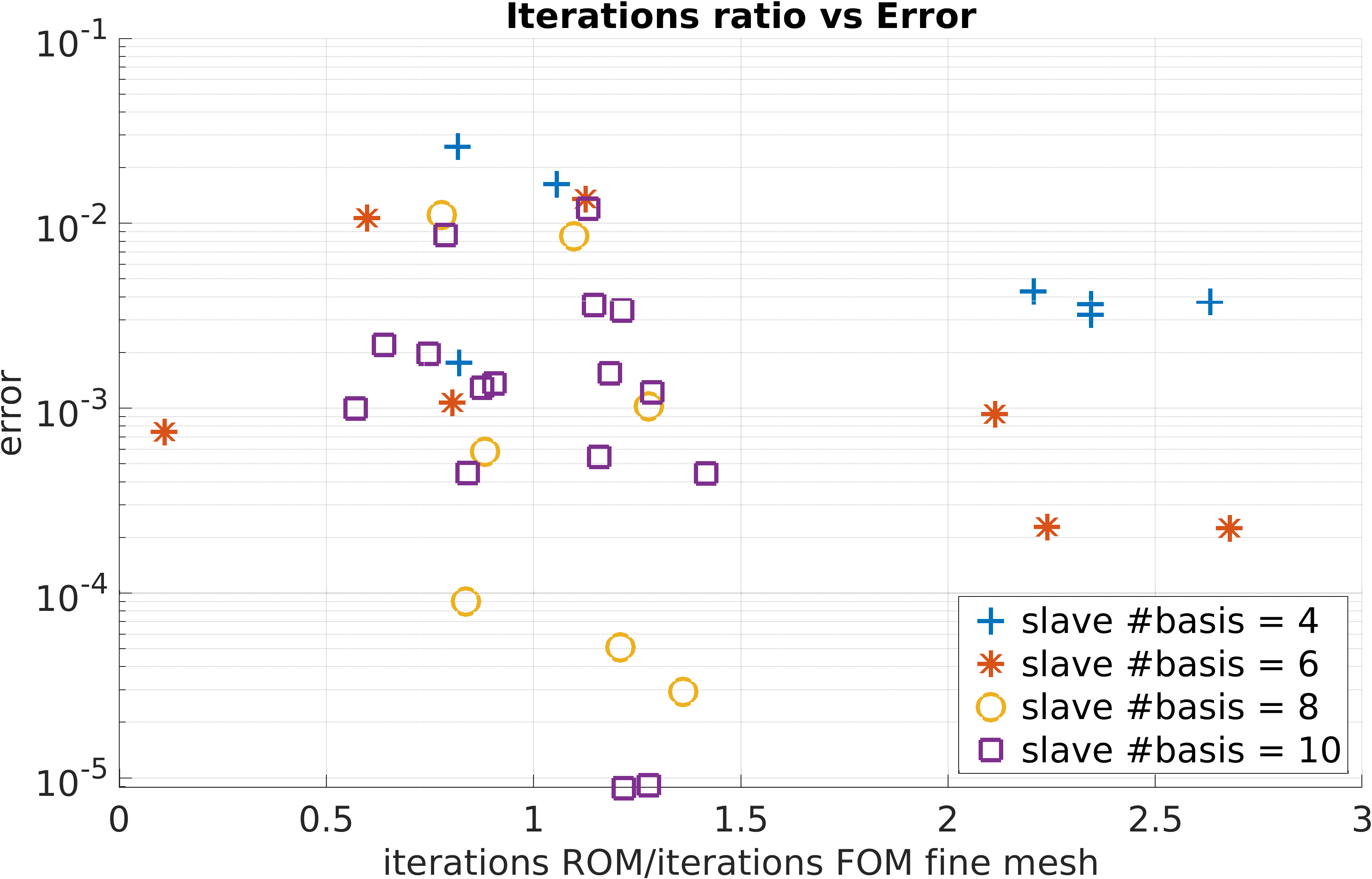}
	\caption{\emph{Test\#2.} Ratio between the number of iterations obtained with ROM and FOM schemes versus the achieved approximation error depending on the number of basis functions employed for the interface data (left) and the slave solution (right).}
	\label{fig:Laplace_source_ratio_vs_error}
\end{figure}

Employing an LHS to randomly select the parameter values, we define $N_{\text{train}} = 250$ to train the ROM, and $N_{\text{test}} = 20$ to test the online procedure. Fig. \ref{fig:Laplace_error_source} shows the variation of the approximation error as a function of the number of basis functions employed to compute one of the three reduced quantities involved. Fixing two of the reduced order hyper--parameters (between $n_1$, $n_2$, or $M_1 = M_2$), the decrease of the approximation error follows the growth of the third hyper--parameter involved, as expected from the ROM theory. However, differently from the test case in Subsection \ref{Subsect:steady_case}, this test shows a similar dependency on the accuracy of either the slave solution, master solution, or interface data, as can be observed from Fig. \ref{fig:Laplace_source_ratio_vs_error}.

\begin{figure}[h!]
	\centering
    \includegraphics[width=0.45\textwidth]{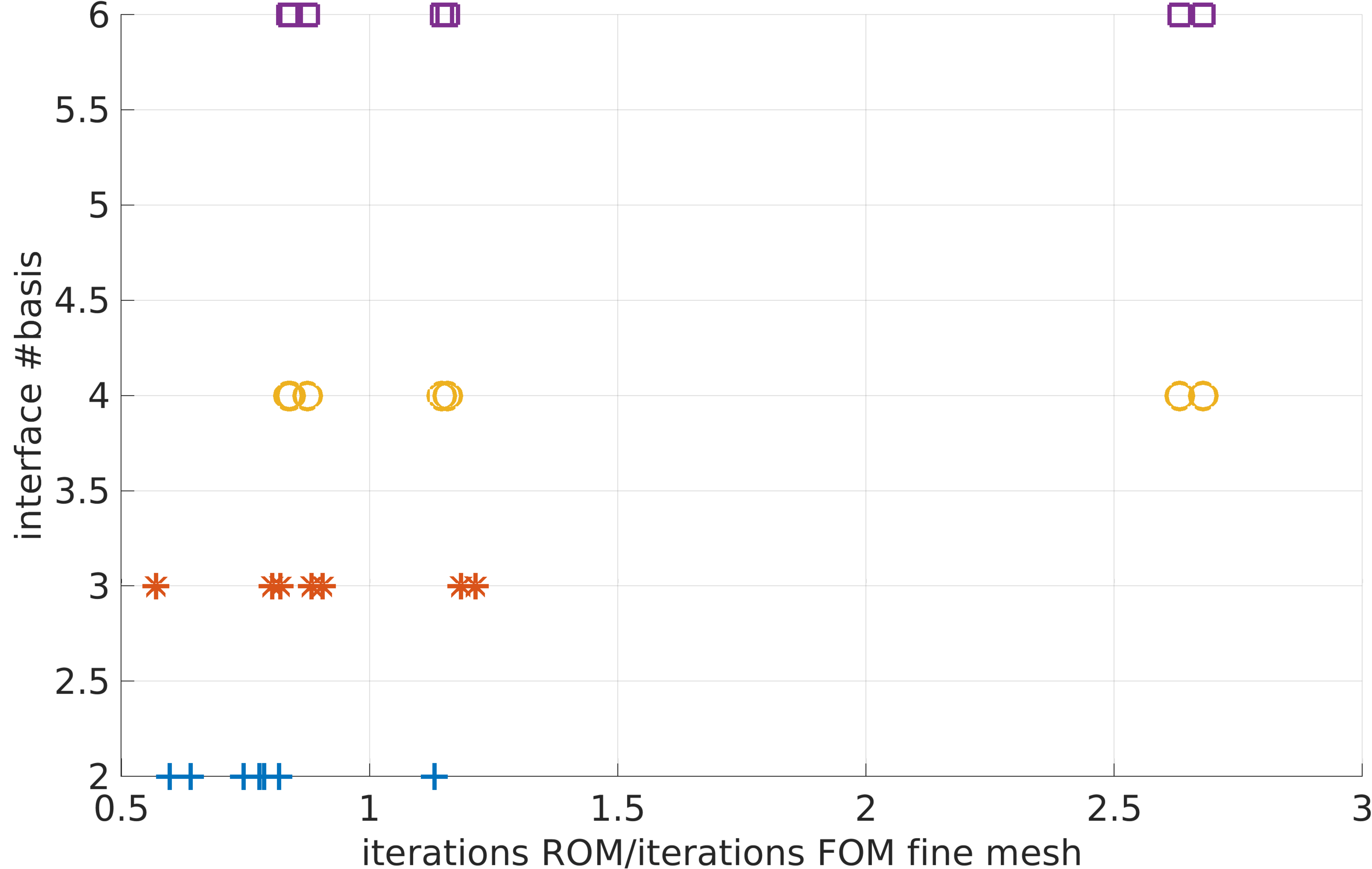}
    \quad
	\includegraphics[width=0.45\textwidth]{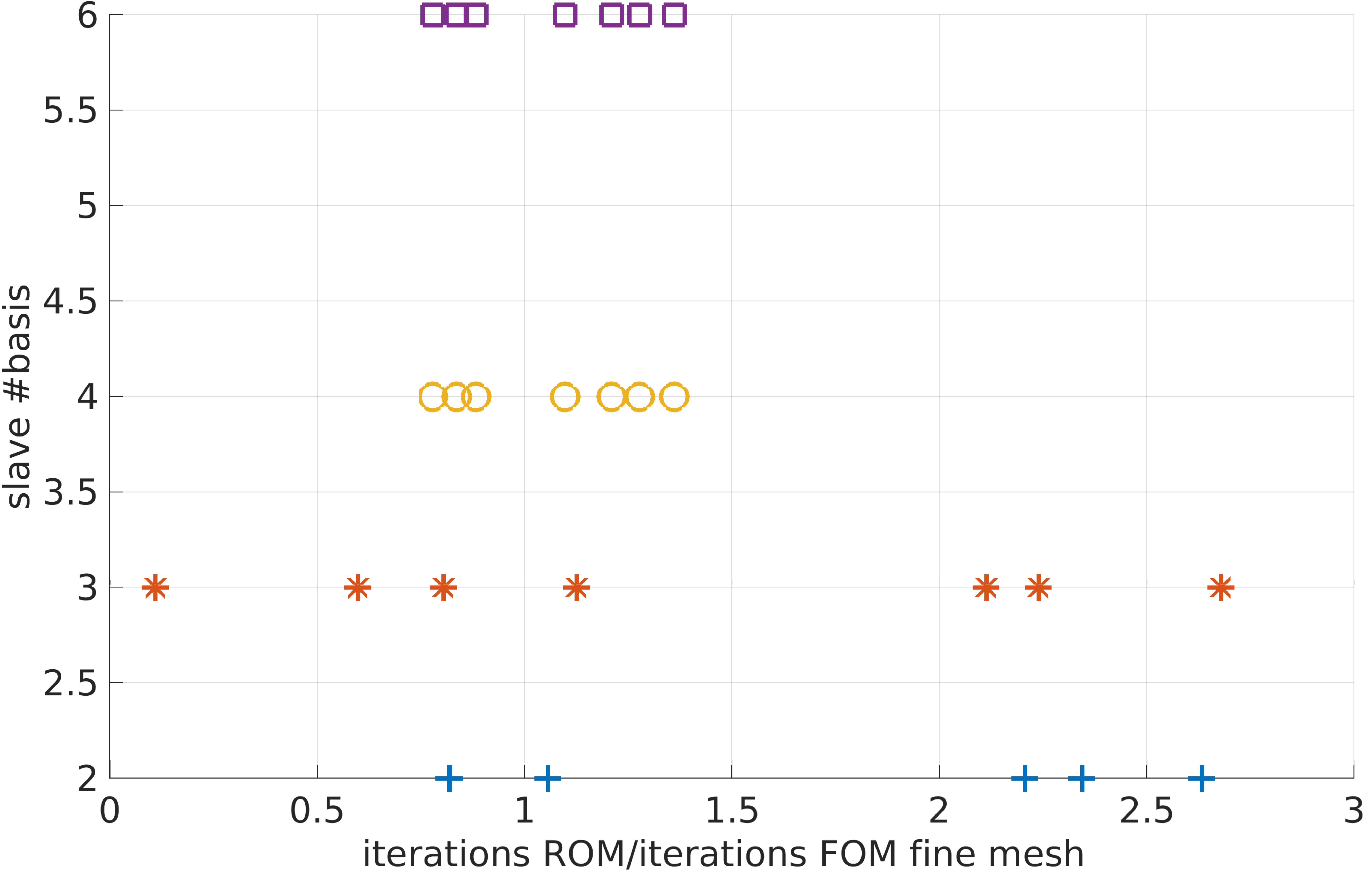}
	\\
	\bigskip
	\includegraphics[width=0.45\textwidth]{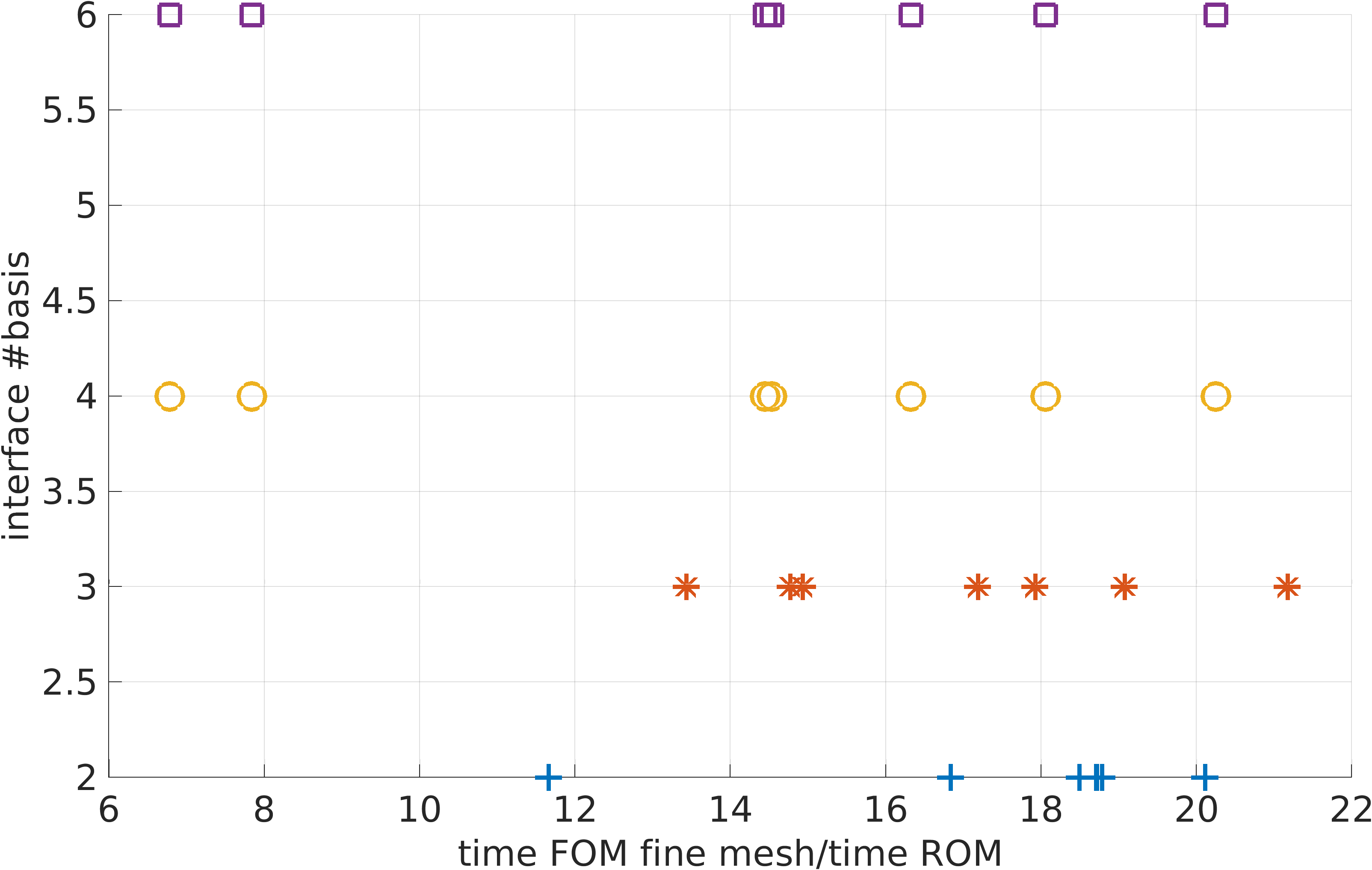}
    \quad
	\includegraphics[width=0.45\textwidth]{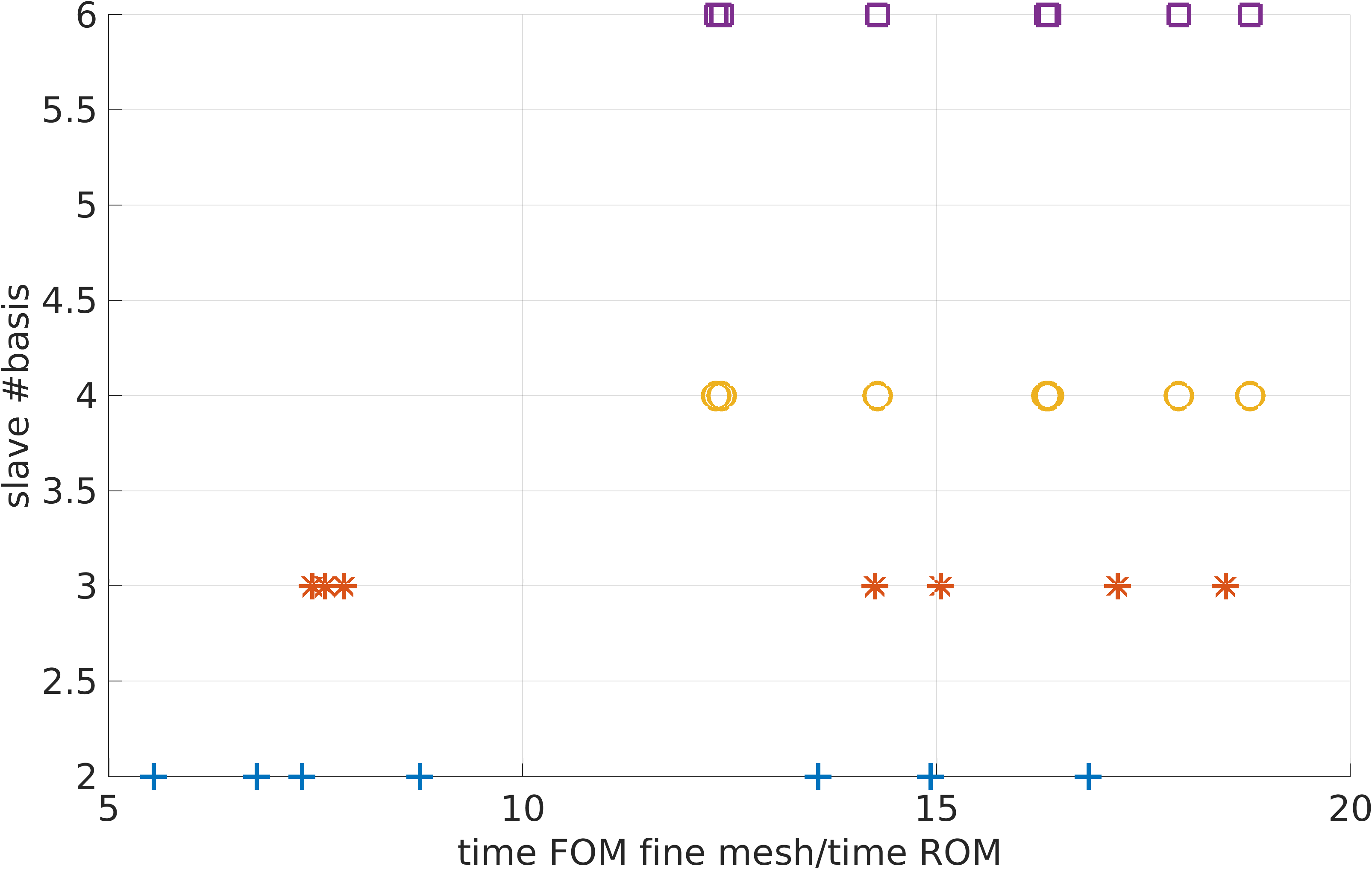}
	\caption{\emph{Test\#2.} Top row: the ratio between the number of iterations obtained with ROM and FOM schemes versus the number of basis functions employed to approximate the interface data (left) and the slave solution (right). Bottom row: the ratio between the FOM and ROM computational time versus the number of basis functions employed to approximate the interface data (left) and the slave solution (right). The FOM simulation referred to fine discretization.}
	\label{fig:Laplace_source_ratio_vs_basis}
\end{figure}

To solve the model with the coarse and the fine discretization in both $\Omega_1$ and $\Omega_2$, on average 34 and 33 iterations are required, respectively, whereas the ROM algorithm needs 45 iterations to reach interface convergence (see Table \ref{Tab:laplace_time} for the number of DoFs employed in each subdomain). The iterations trend of the ROM vs the FOM solutions (see Fig. \ref{fig:Laplace_source_ratio_vs_basis}) indicates that, increasing $M_1$ and $M_2$, the number of iterations required by the ROM algorithm increases, which can be explained by the complex interface solution to be approximated. By increasing $n_1$, instead, the cost of the iterations decreases. The same behaviour can be observed also in Fig. \ref{fig:Laplace_source_ratio_vs_error} by comparing the approximation errors and the iterations ratio, as well as looking at the computational time as a function of the ROM hyper--parameters (see Fig. \ref{fig:Laplace_source_ratio_vs_basis}, bottom row) where, by increasing $M_1$ and $M_2$, the ROM speed up compared to the FOM decreases, while the overall time gained by employing the ROM is overall the same independently on $n_1$ and $n_2$. Imposing a tolerance of $10^{-5}$ to approximate all slave and master solutions, and interface data, about $18.1s$ are needed by the ROM algorithm to compute the solution, $26.8s$ for the FOM computation with the coarse discretization, and $225.33s$ for the FOM computation with the fine discretization. Therefore, the CPU time is reduced by about 34\% compared to the coarse discretization, and 92\% compared to the fine discretization.

\subsection{Test\#3. Unsteady case: time--dependent heat equation}
\label{sub:heat_problem}
\begin{figure}[h!]
	\centering
	\includegraphics[width=0.23\textwidth]{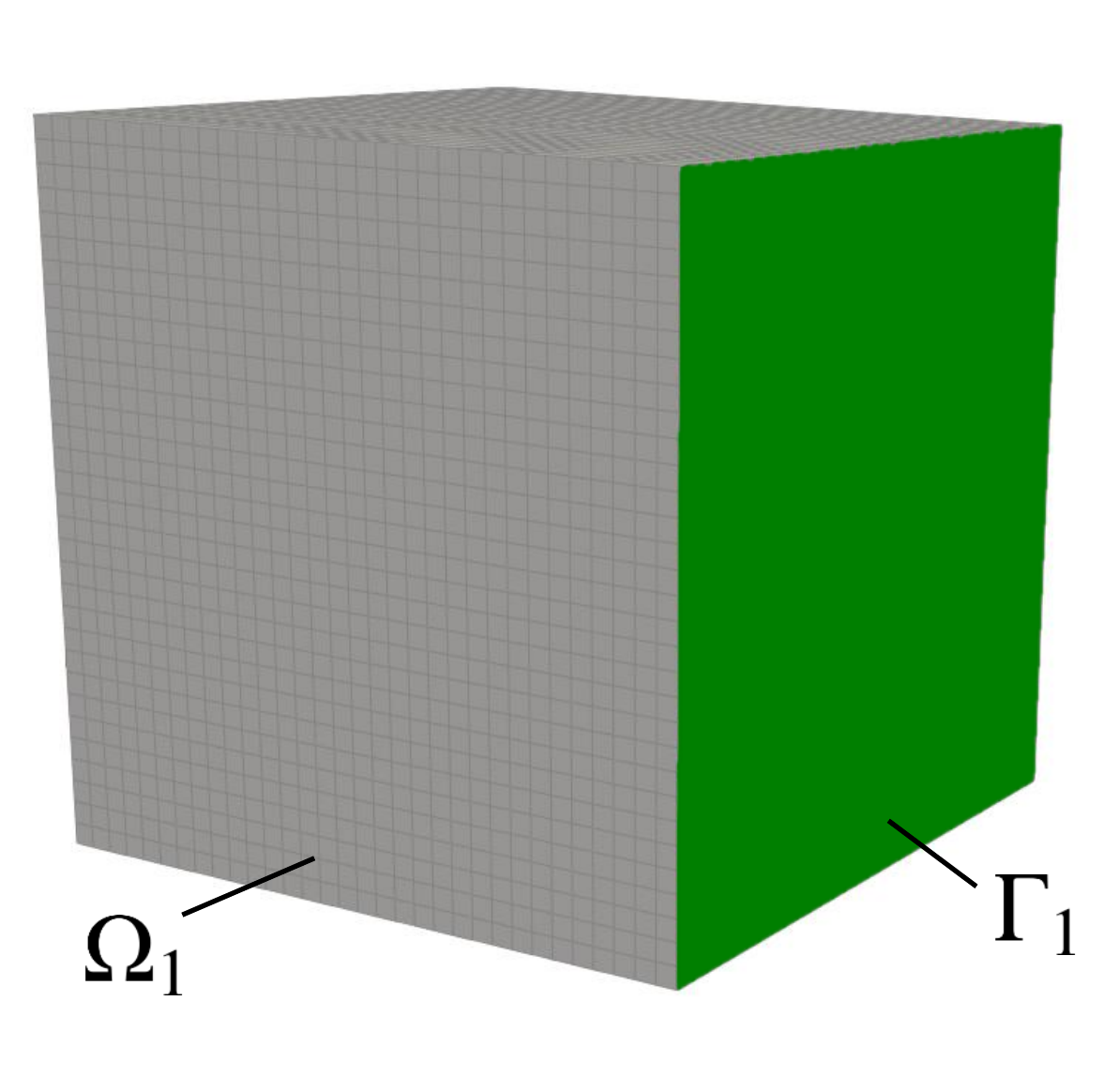}
	\qquad\quad
	\includegraphics[width=0.23\textwidth]{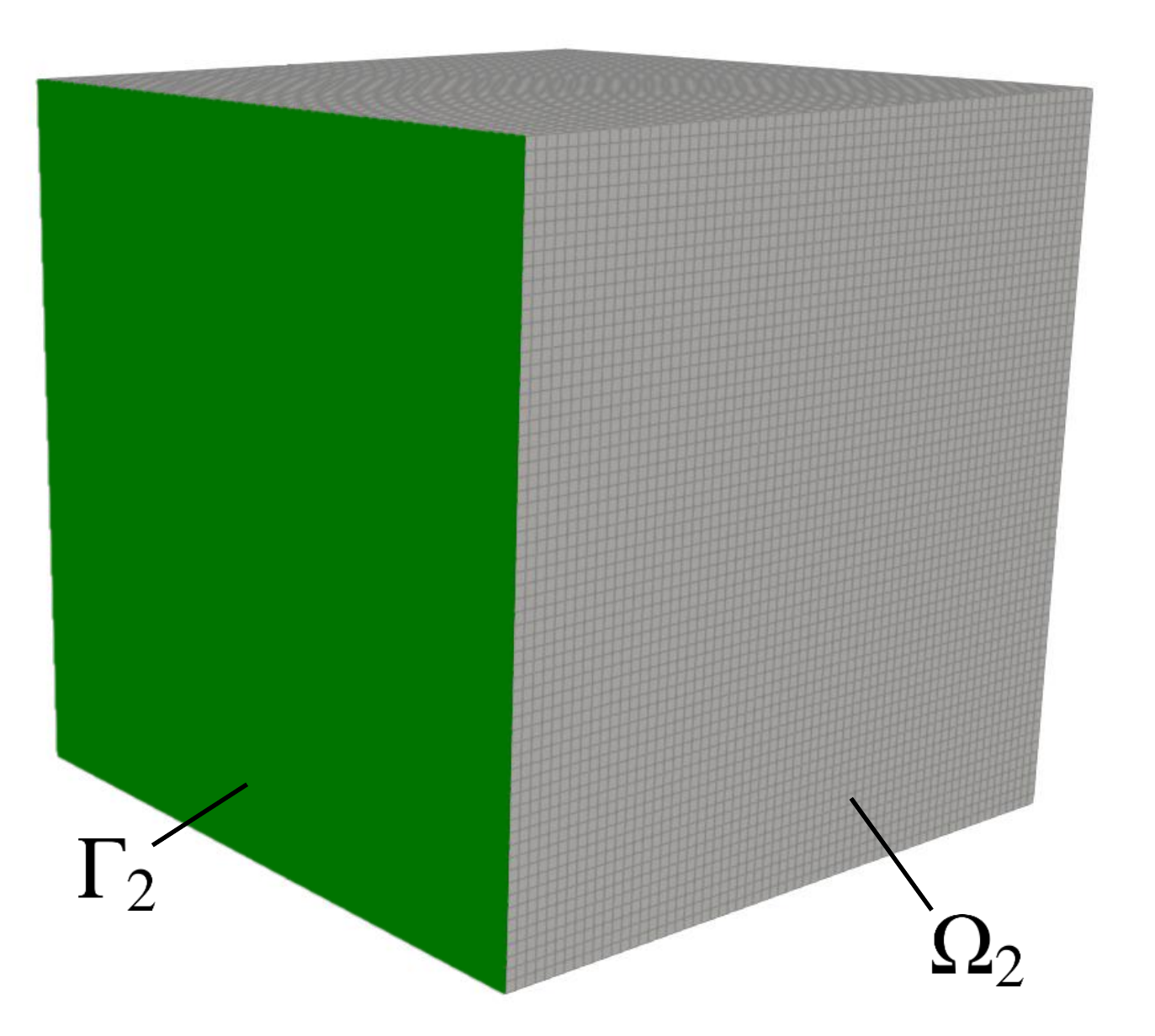}
	\qquad\quad
	\includegraphics[width=0.32\textwidth]{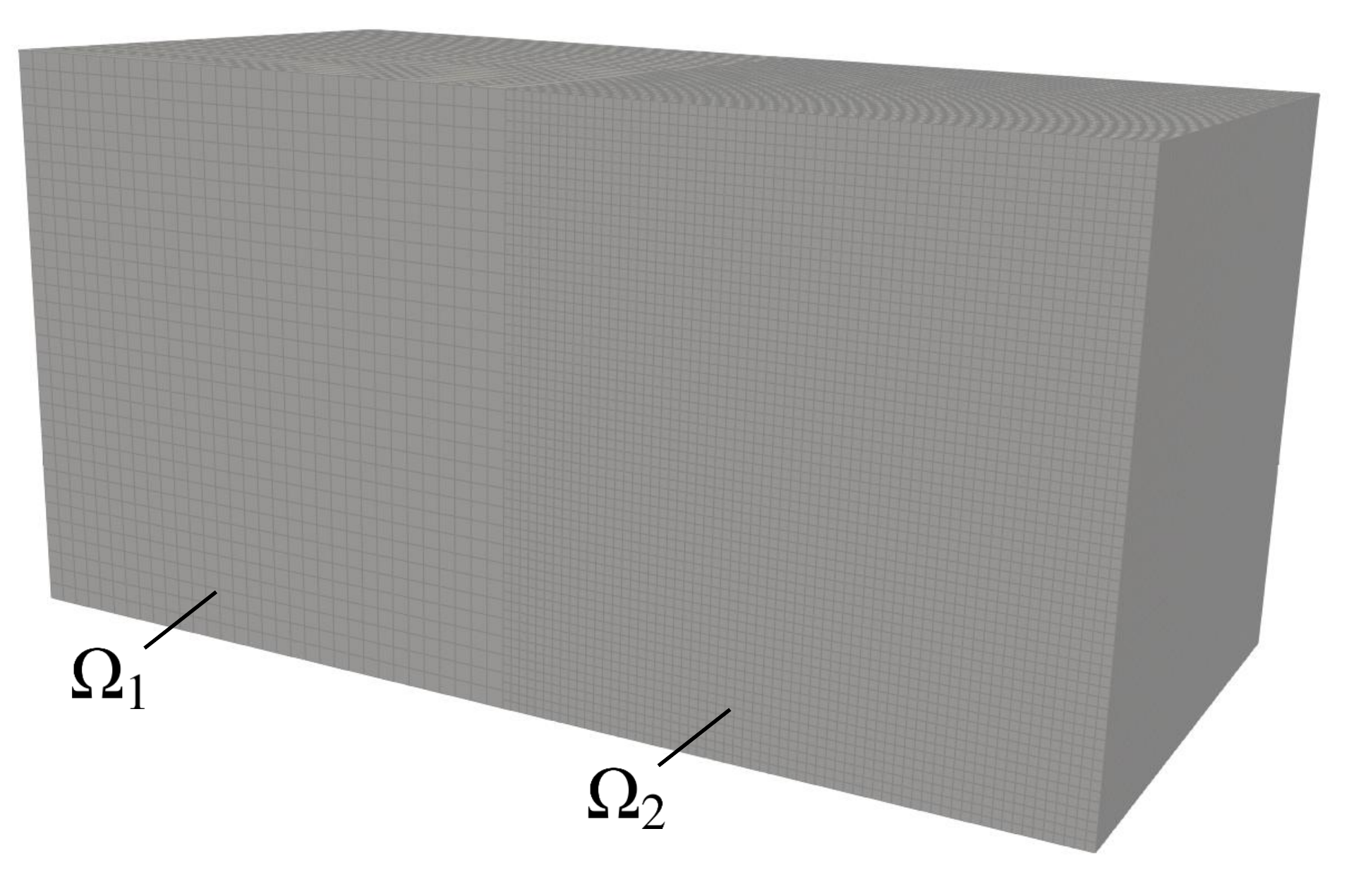}
	\caption{\emph{Test\#3.} $\Omega_1$ (left) and $\Omega_2$ (center) and $\Omega$ (right). In green the interface $\Gamma$.}
	\label{Fig:Heat_domains}
\end{figure}

We now apply the proposed technique to a time--dependent problem. In particular, we consider the following initial--boundary value problem for the heat equation
\begin{equation}
\label{eq:heat_eq_complete}
\begin{cases}
\displaystyle \frac{\partial u}{\partial t} - \nabla \cdot (\alpha \nabla u) = f &\text{in }\Omega\times (0,T)\\
\alpha \nabla u \cdot \mathbf{n} = 0 &\text{on }\partial \Omega\times (0,T) \\
u = 0 &\text{in } \Omega\times \{t = 0\} \\
\end{cases}
\end{equation}
being $\Omega = (-0.5,1.5) \times (-0.5,0.5) \times (-0.5,0.5)$  (see Figure \ref{Fig:Heat_domains}) and 
\[ f = \begin{cases}
1 &\text{if } x < 0 \land 0.2 < t < 0.5\\
0 &\text{otherwise}. 
\end{cases}\]
The problem is parametrized through the coefficient $\alpha$ and -- indirectly -- through the time variable $t$. 

To apply the reduced technique presented in this work, we split $\Omega$ into two cubes with a common plane $x = 0.5$ that represents the interface $\Gamma$ (see Fig. \ref{Fig:Heat_domains}). The application of the Dirichlet--Neumann iterative scheme leads therefore to the two following sub--problems: given $\lambda_2^0$ on $\Gamma_2$, for each $k\geq 0$, solve until convergence of the interface solution
\begin{equation}
\label{Eq:heat_sub_1}
\begin{cases}
\displaystyle \frac{\partial u_1^{k+1}}{\partial t} - \nabla \cdot (\alpha \nabla u_1^{k+1}) = f &\text{in }\Omega_1 \times (0,T)\\
\alpha \nabla u_1^{k+1} \cdot \mathbf{n}_1 = 0 &\text{on }\partial \Omega_{1,N}\times (0,T)\\
u_1^{k+1} = \lambda_2^k &\text{on }\Gamma_1 \times (0,T)\\
u_1^{k+1} = 0 &\text{in }\Omega_1 \times \{t = 0 \}
\end{cases}
\end{equation}
and
\begin{equation}
\label{Eq:heat_sub_2}
\begin{cases}
\displaystyle \frac{\partial u_2^{k+1}}{\partial t} - \nabla \cdot (\alpha \nabla u_2^{k+1}) = f &\text{in }\Omega_2 \times (0,T)\\
\alpha \nabla u_2^{k+1} \cdot \mathbf{n}_2 = 0 &\text{on }\partial \Omega_{2,N}\times (0,T)\\
\alpha \nabla u_2^{k+1} \cdot \mathbf{n}_2 = - \alpha \nabla u_1^{k+1} \cdot \mathbf{n}_1 &\text{on }\Gamma_2 \times (0,T)\\
u_2^{k+1} = 0 &\text{in }\Omega_2 \times \{t = 0 \},
\end{cases}
\end{equation}
while
$$ \lambda_2^{k+1} = \omega u_{2_{|\Gamma_2}}^k + (1-\omega)\lambda_2^k.$$

As for the steady test case, a fixed point acceleration strategy with parameter $\omega = 0.25$ is applied to accelerate the convergence of the Dirichlet--Neumann scheme for both FOM and ROM computations.

Two different discretizations are employed, corresponding to a FOM dimension equal to $N_1 = N_2 = 35937$ for the first discretization, and $N_1 = N_2 = 274625$ for the second discretization. The ROM is then solved with $N_2 = 274625$ DoFs for the slave domain and $N_1 = 35937$ DoFs for the master domain. We define $t \in [0,1]s$, and we use the first order backward Euler scheme with $\Delta t = 10^{-2}s$ to handle time discretization. Moreover, the diffusivity coefficient $\alpha$ is varied between 0.5 and 5, employing also, in this case, the LHS strategy to select the snapshots to be computed for the ROM training. In particular, we select $N_{S} = N_t N_{\text{train}}$ snapshots to train the ROM, corresponding to $N_{\text{train}} = $10 complete simulations in time, being $N_t = 100$ the number of time steps in a model simulation. Then, $N_{\text{test}} = 2$ simulations of 100 time--steps each are considered to test the ROM. The tolerance to stop the Dirichlet--Neumann iteration is $10^{-10}$, which corresponds to 18 iterations to find the FOM solution with the coarser discretization, 19 iterations for the FOM solution with the finer discretizations, and 17 iterations to find the ROM solution for each time step.

\begin{figure}[h!]
	\centering
	\hspace{-40pt} $t = 0.2s$\qquad \qquad \qquad \qquad \qquad \qquad$t = 0.4s$ \qquad \qquad \qquad\qquad\qquad \qquad $t = 0.60s$\\
	\vspace{3pt}
	\includegraphics[width=0.3\textwidth]{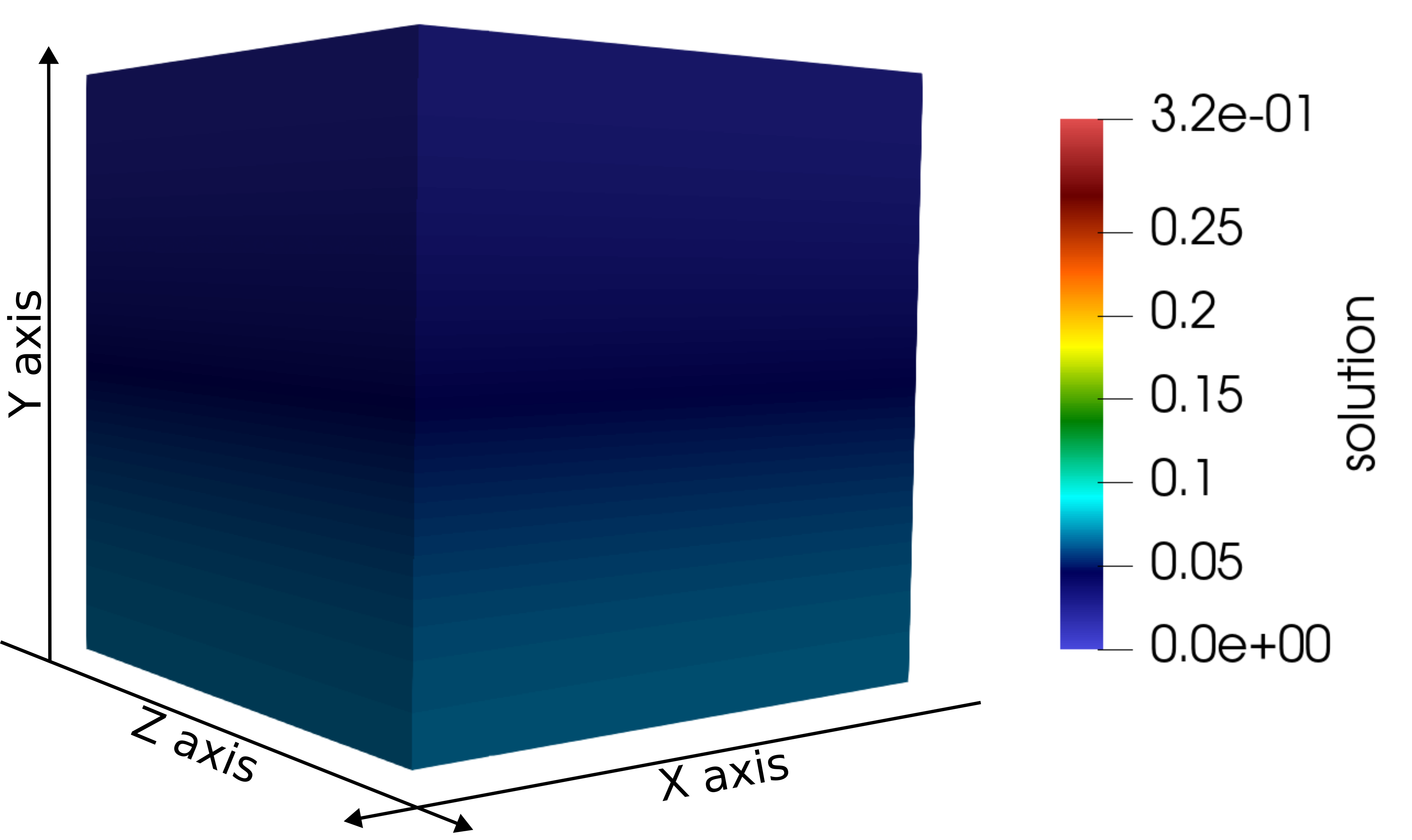}
	$\quad$
	\includegraphics[width=0.3\textwidth]{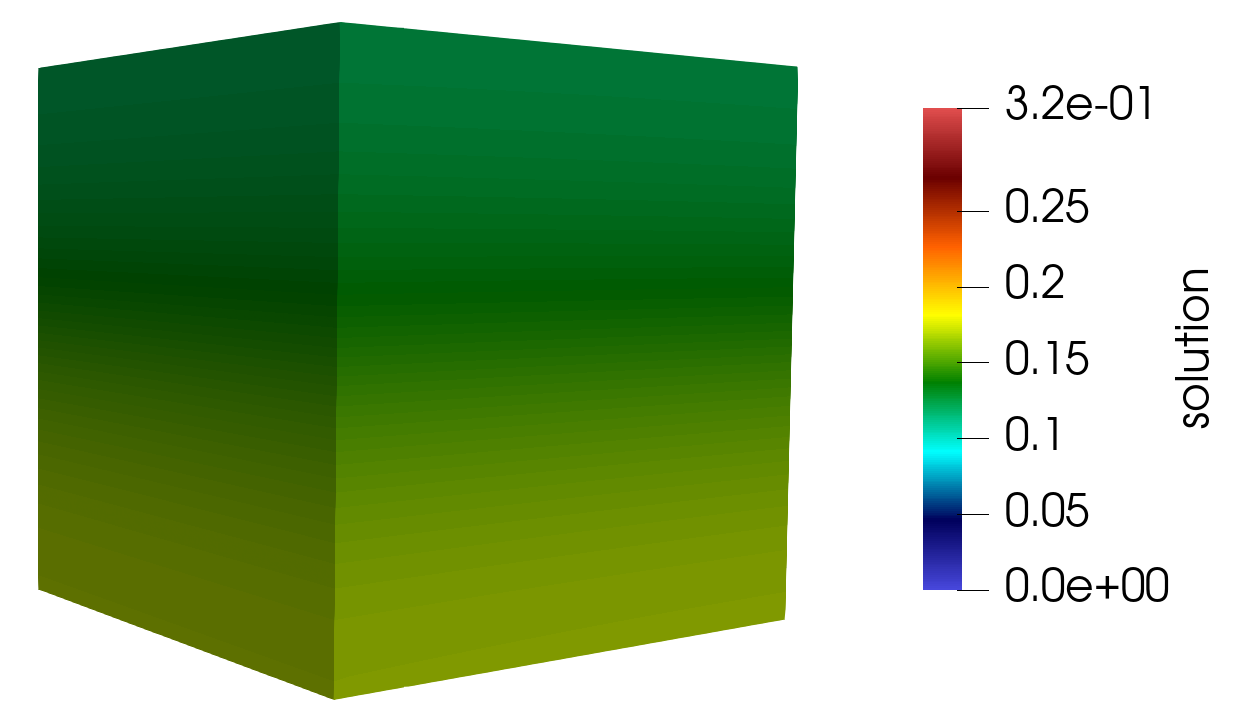}
	$\quad$
	\includegraphics[width=0.3\textwidth]{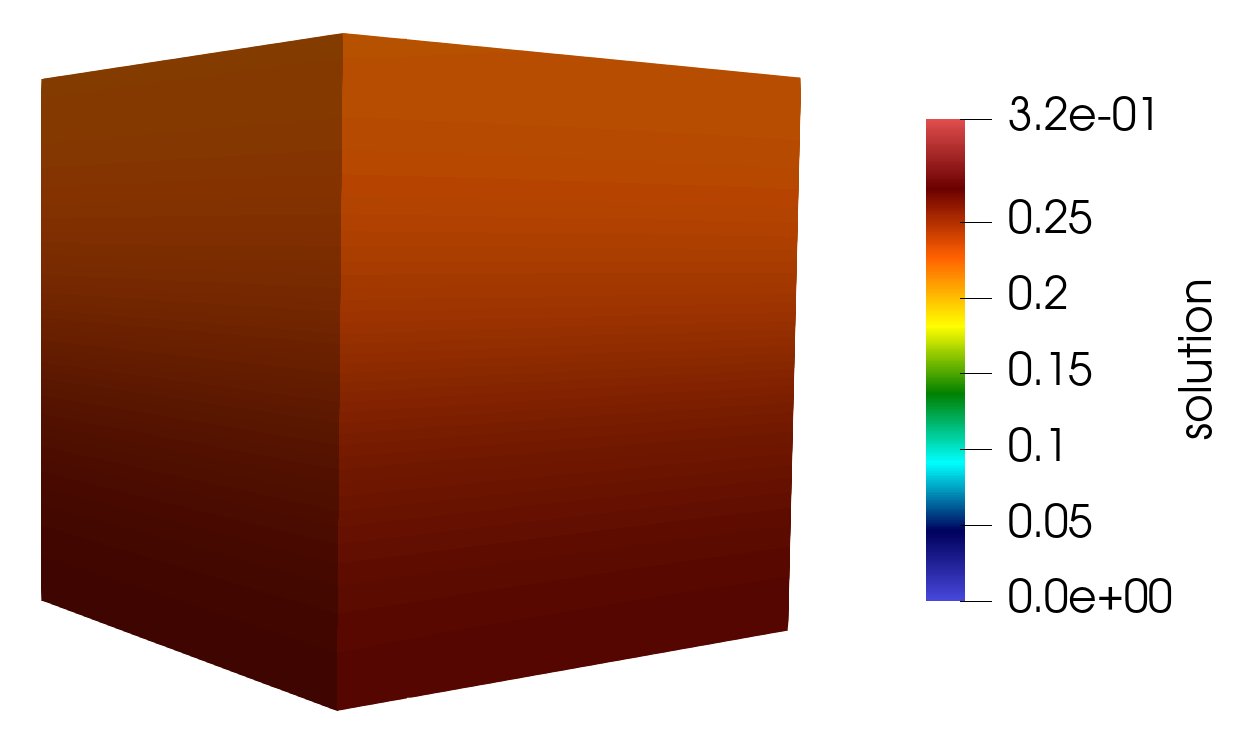}
	 \\
	\includegraphics[width=0.3\textwidth]{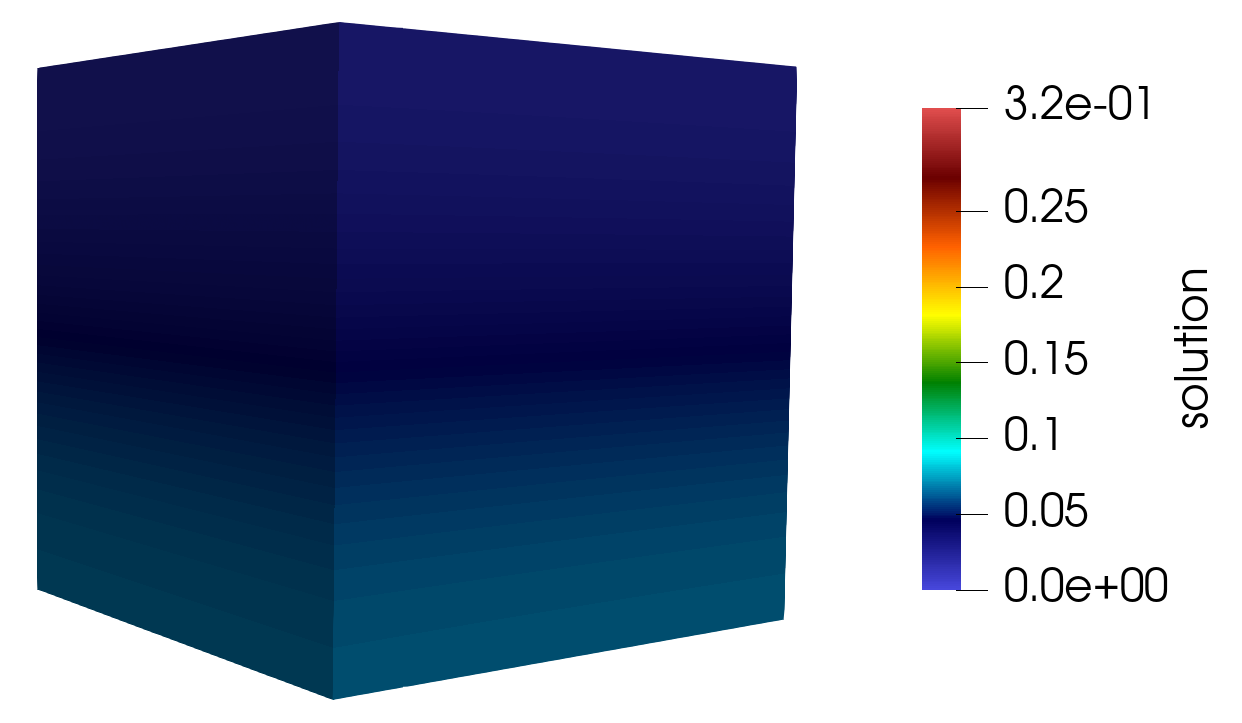}
	$\quad$
	\includegraphics[width=0.3\textwidth]{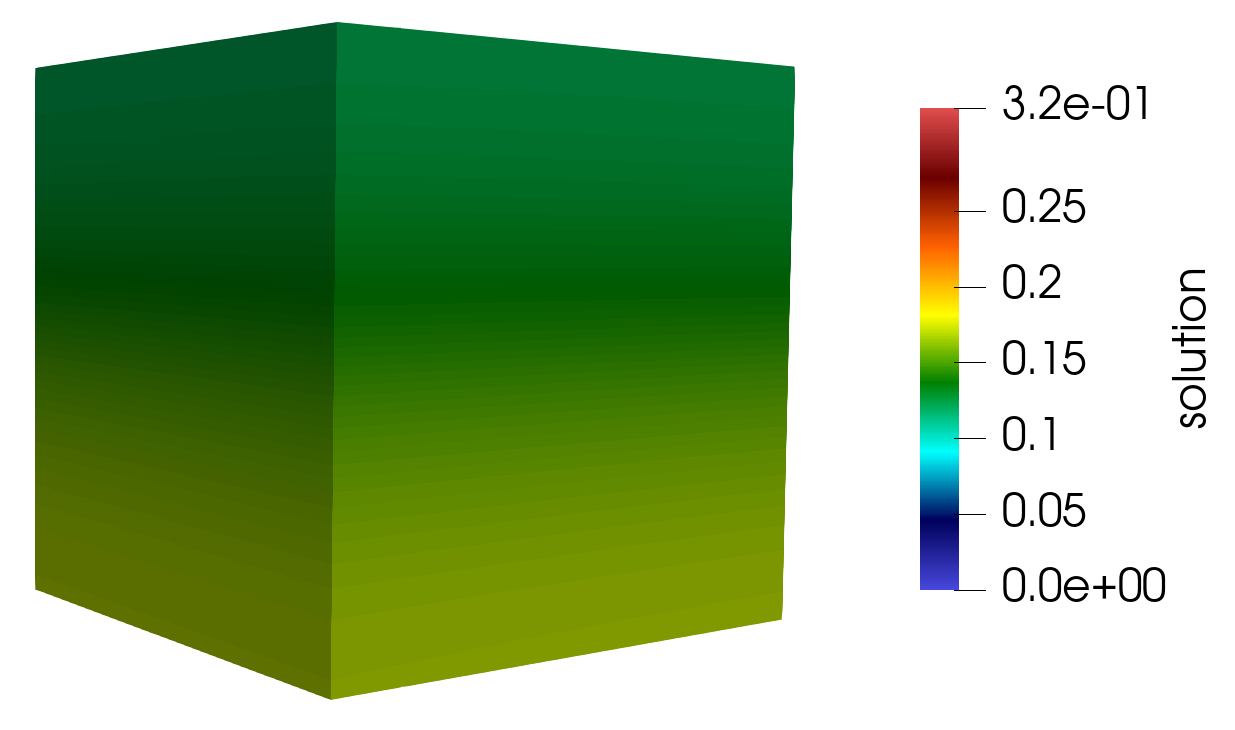}
	$\quad$
	\includegraphics[width=0.3\textwidth]{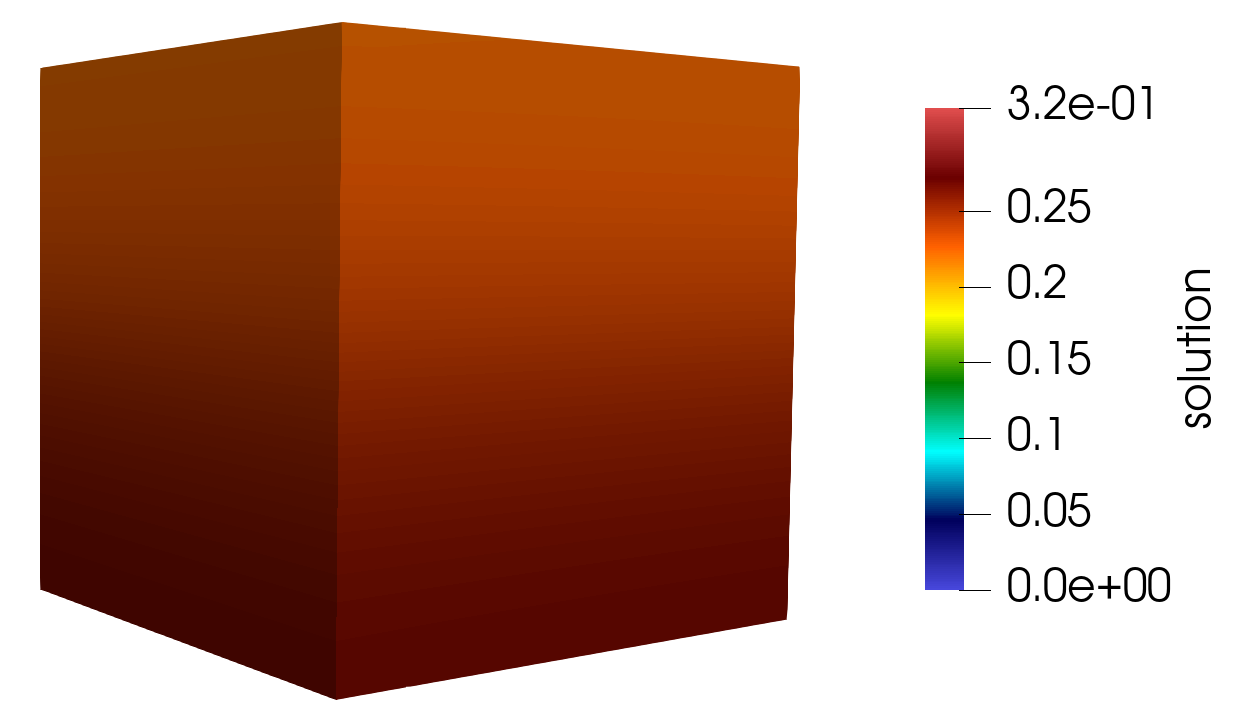}
	\\
	\vspace{2pt}
	\includegraphics[width=0.3\textwidth]{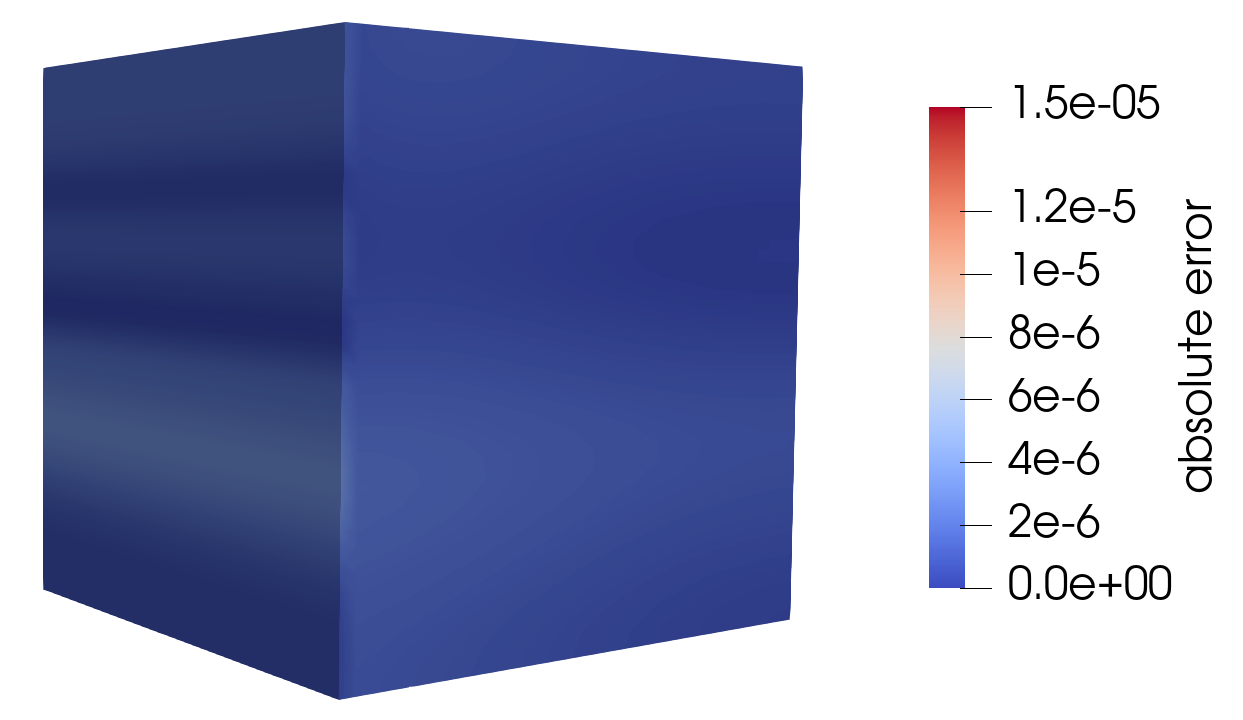}
	$\quad$
	\includegraphics[width=0.3\textwidth]{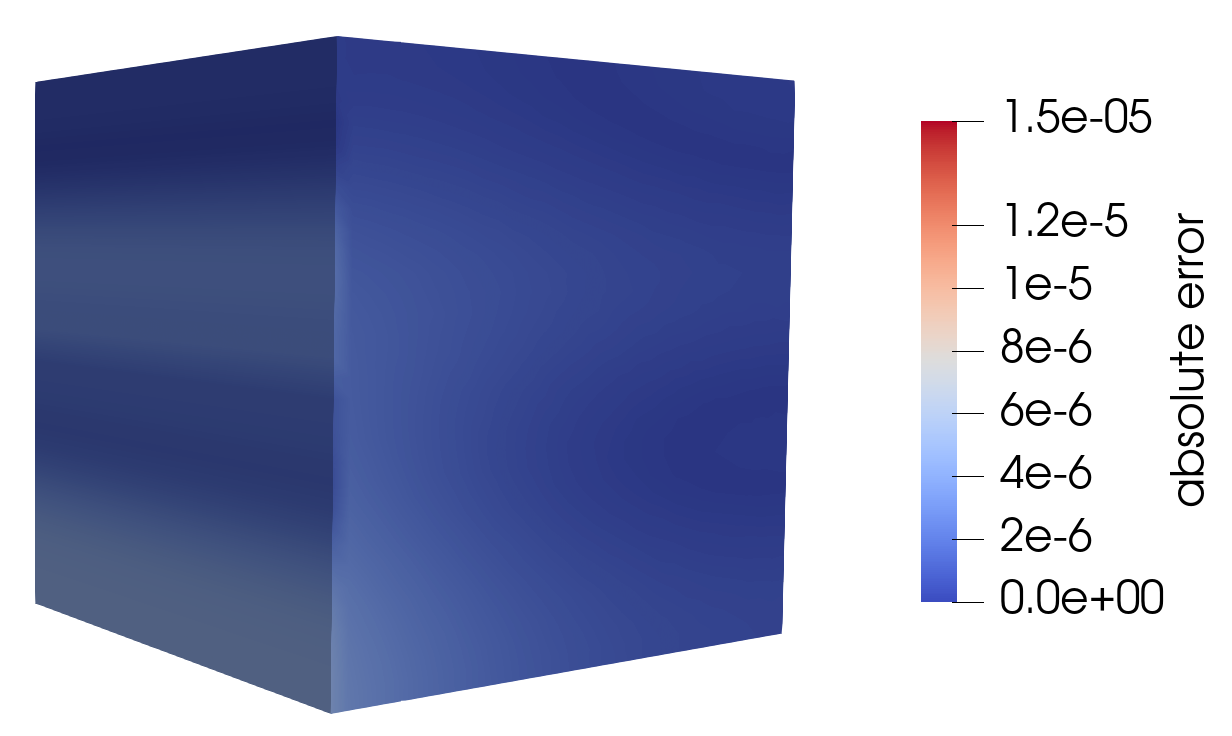}
	$\quad$
	\includegraphics[width=0.3\textwidth]{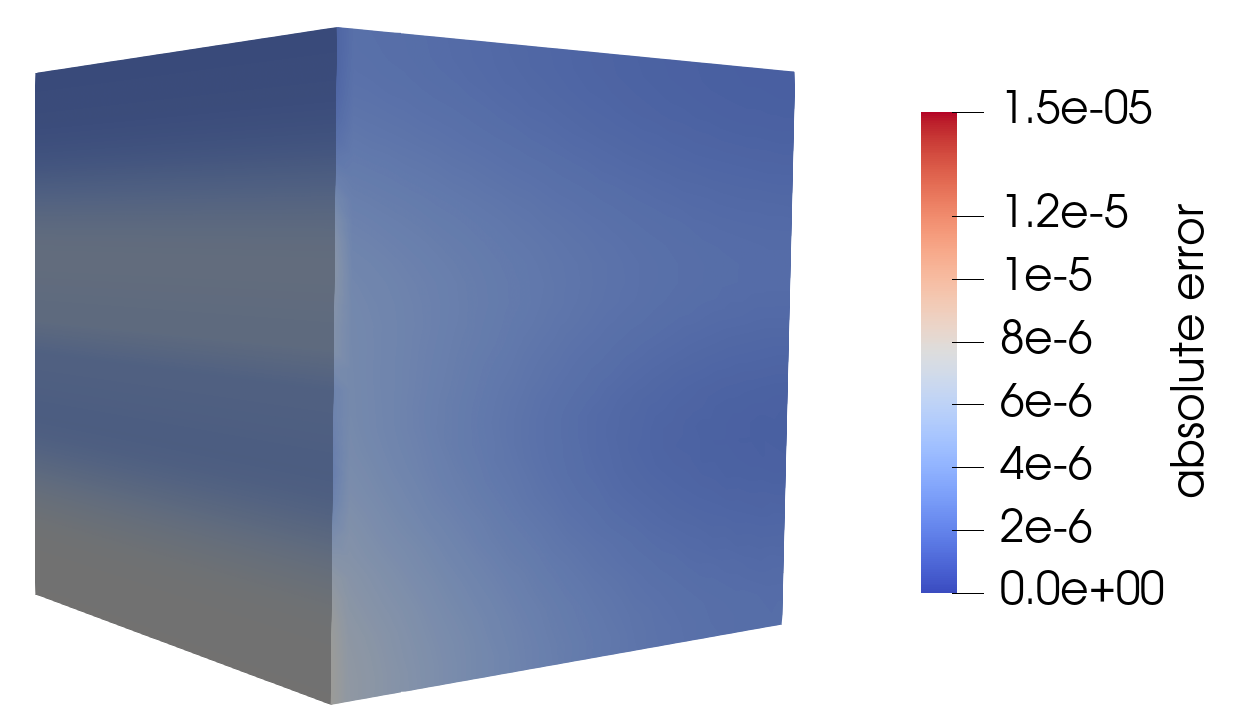}
	\caption{\emph{Test\#3.} Slave solution FOM (top), ROM (center) solutions and absolute error (bottom) for $\alpha=2.75$ and three different time instants. }
	\label{fig:snapshots_heat_master}
\end{figure} 

\begin{figure}[h!]
	\centering
	\hspace{-40pt} $t = 0.2s$\qquad \qquad \qquad \qquad \qquad \qquad$t = 0.4s$ \qquad \qquad \qquad\qquad\qquad \qquad $t = 0.6s$\\
	\vspace{3pt}
	\includegraphics[width=0.3\textwidth]{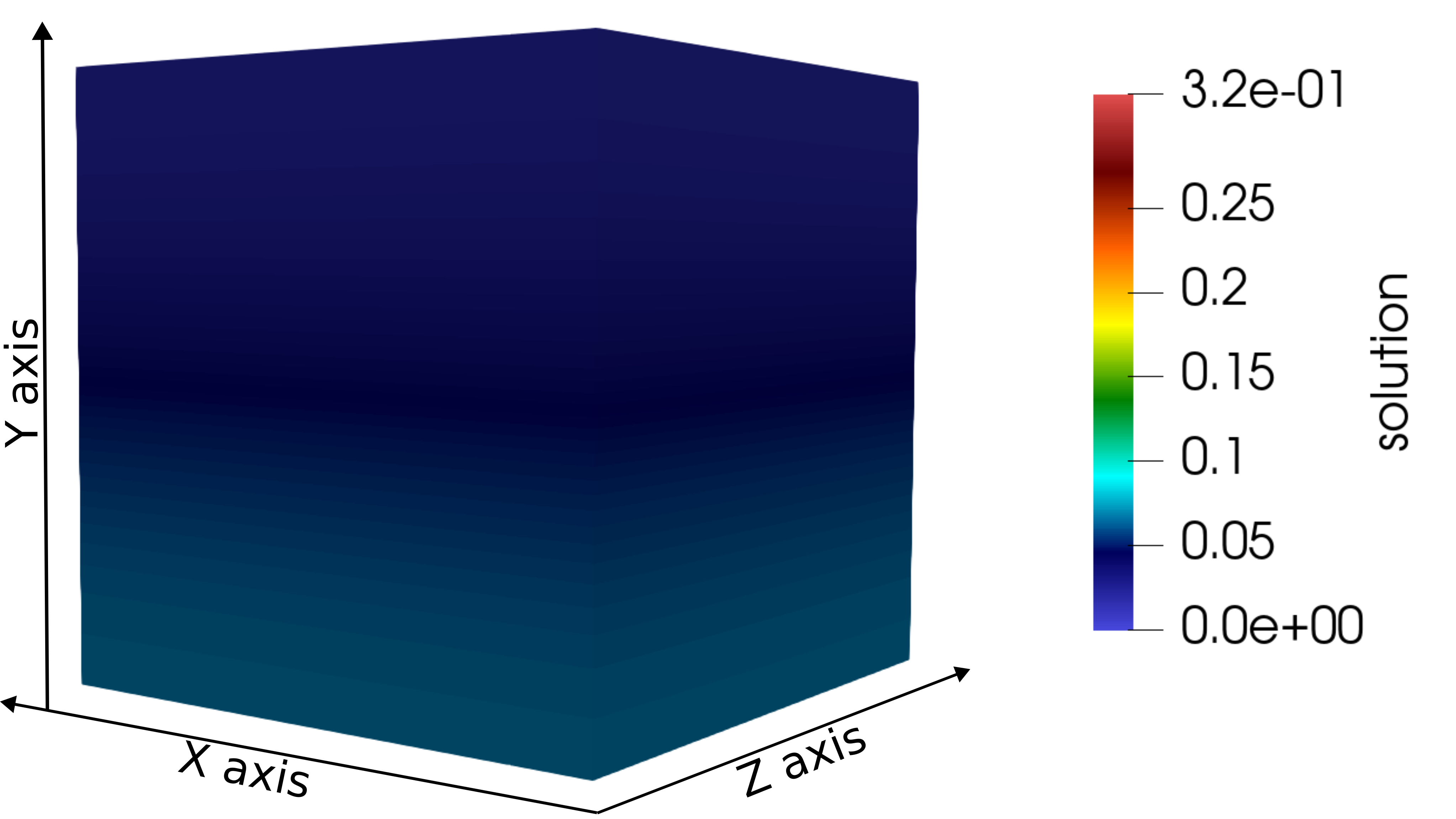}
	$\quad$
	\includegraphics[width=0.3\textwidth]{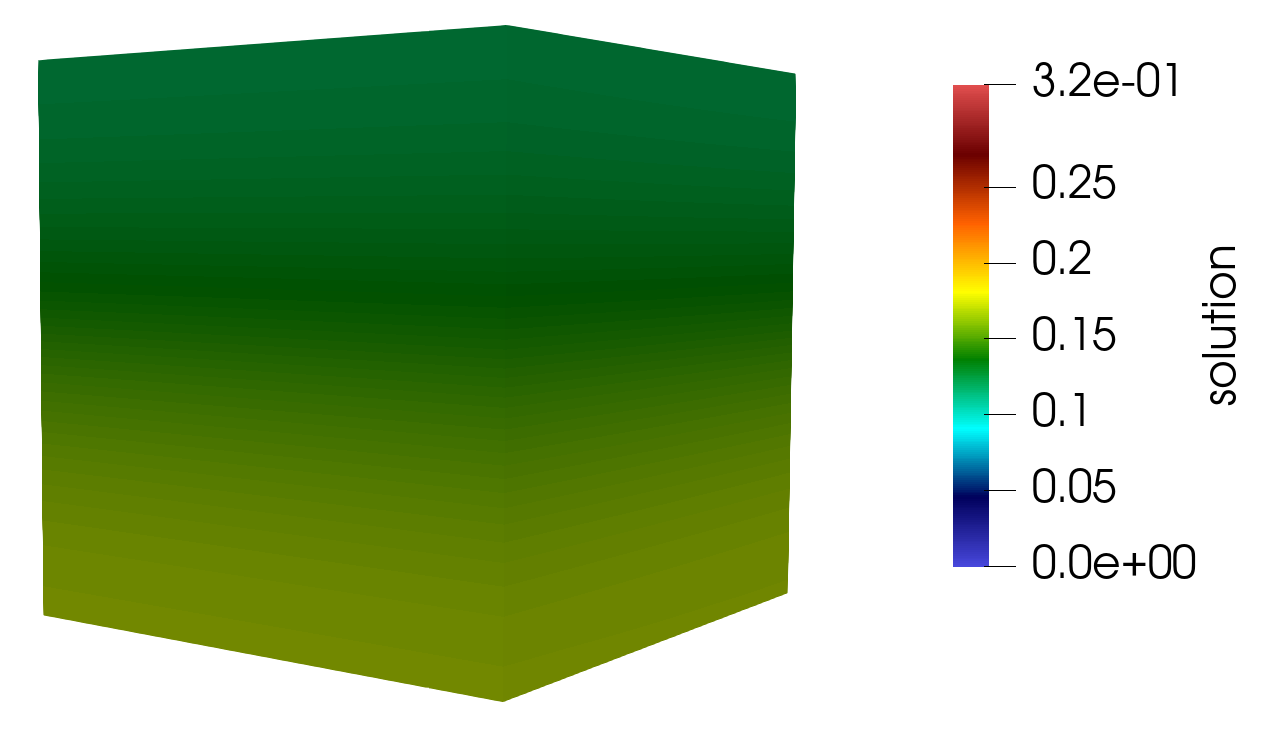}
	$\quad$
	\includegraphics[width=0.3\textwidth]{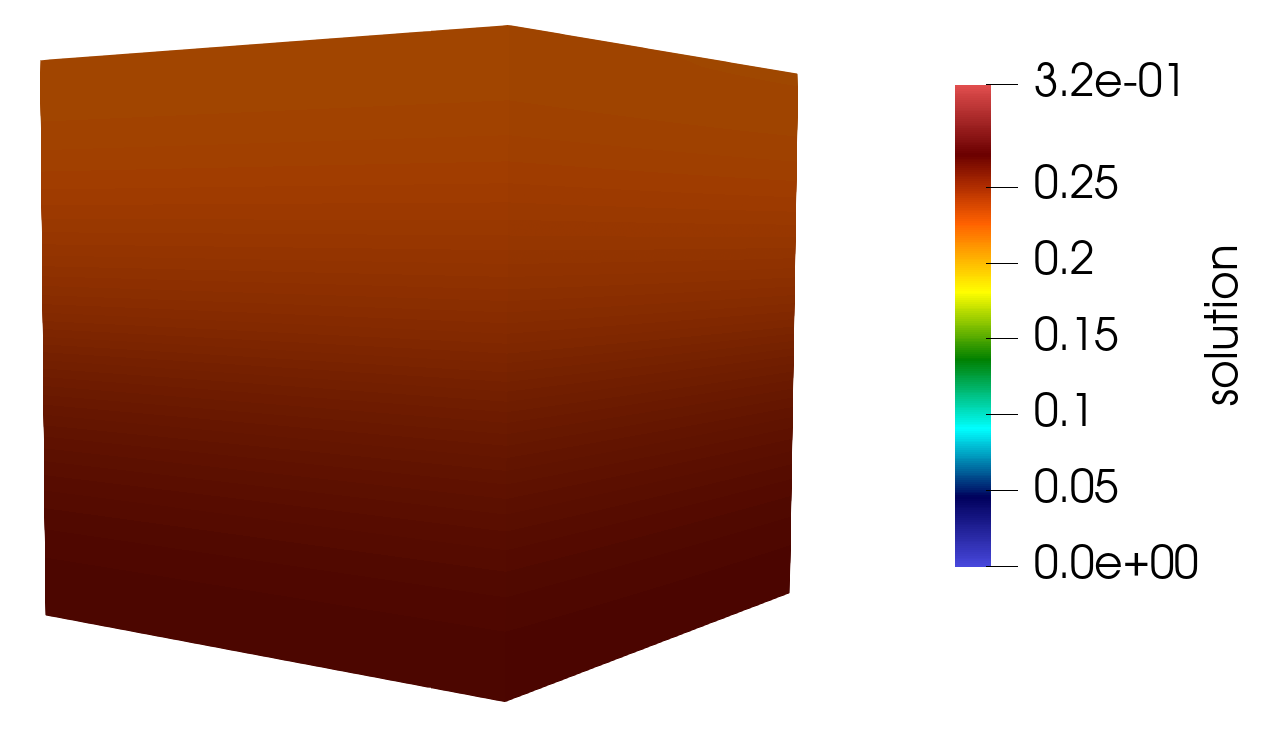}
	 \\
	\includegraphics[width=0.3\textwidth]{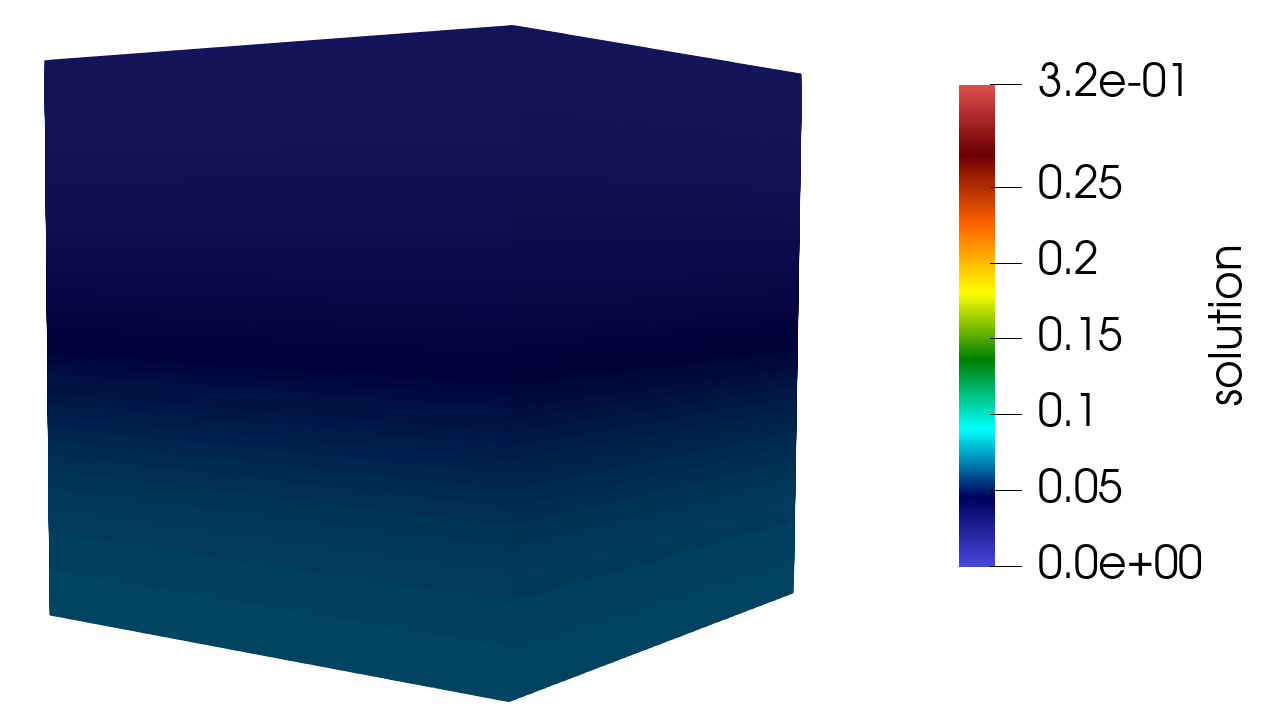}
	$\quad$
	\includegraphics[width=0.3\textwidth]{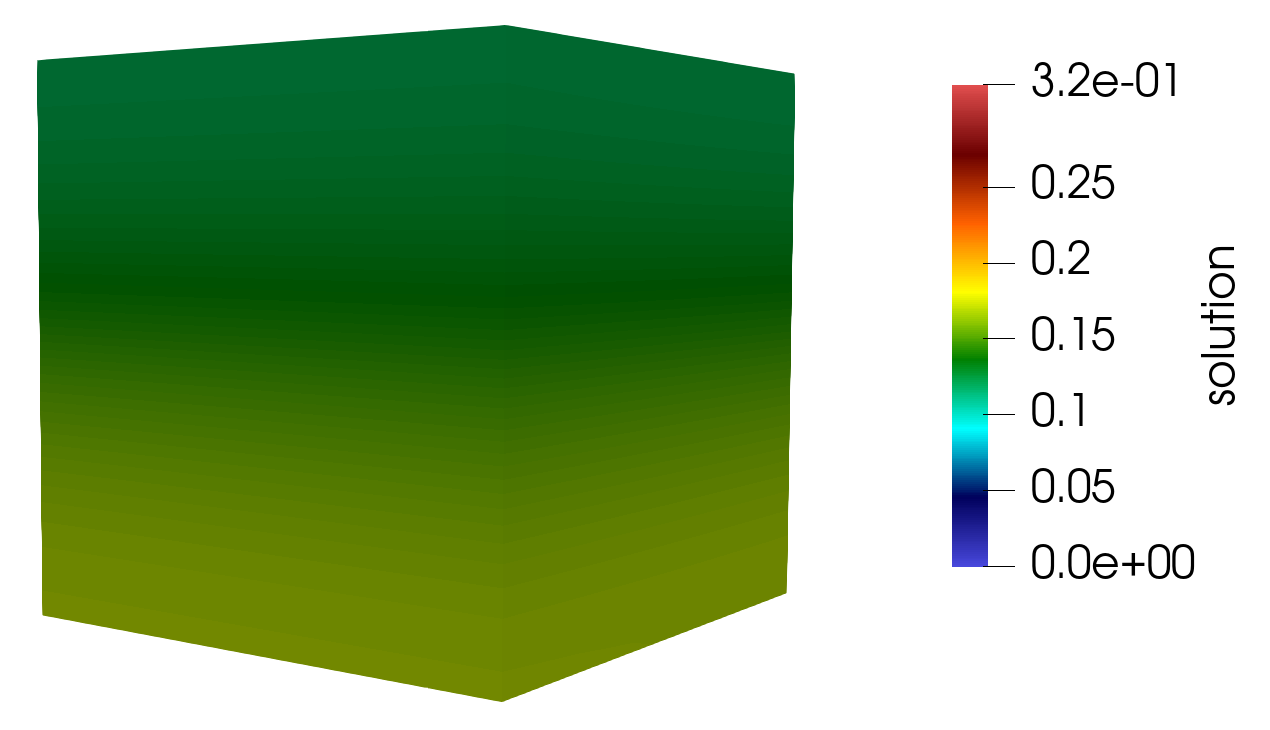}
	$\quad$
	\includegraphics[width=0.3\textwidth]{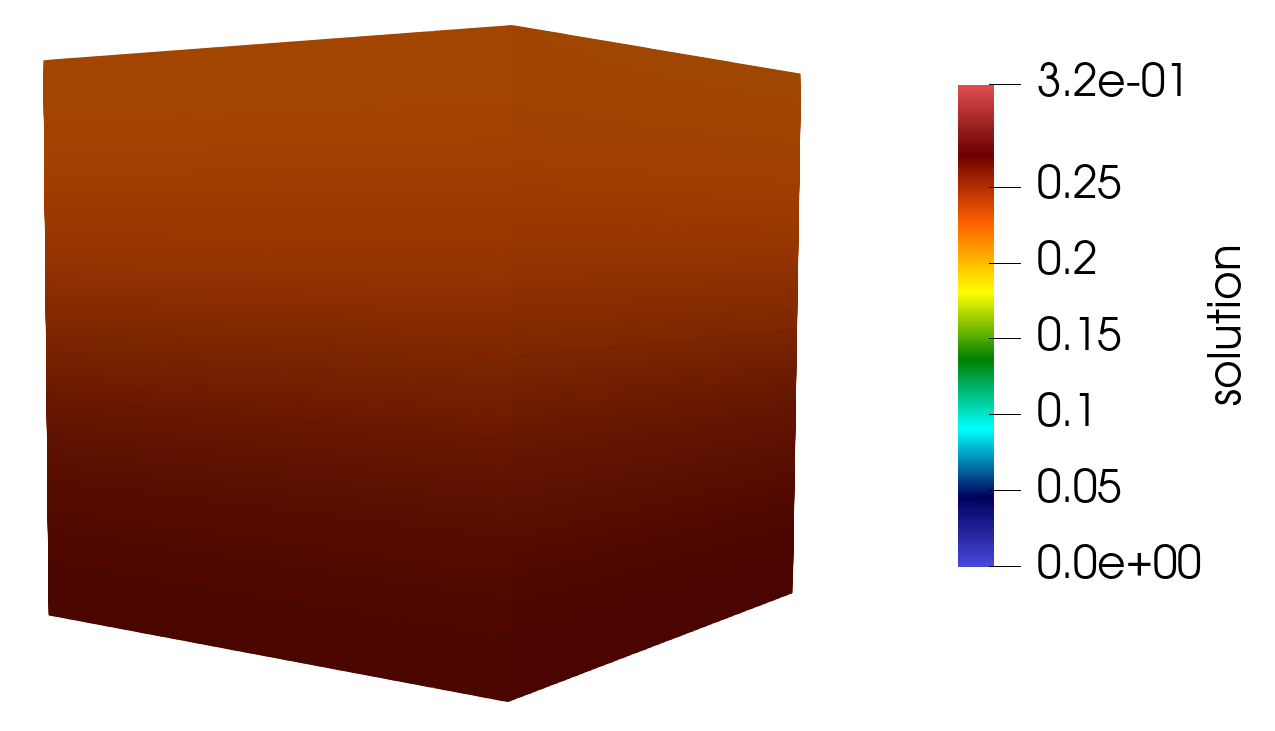}
	\\
	\vspace{3pt}
	\includegraphics[width=0.3\textwidth]{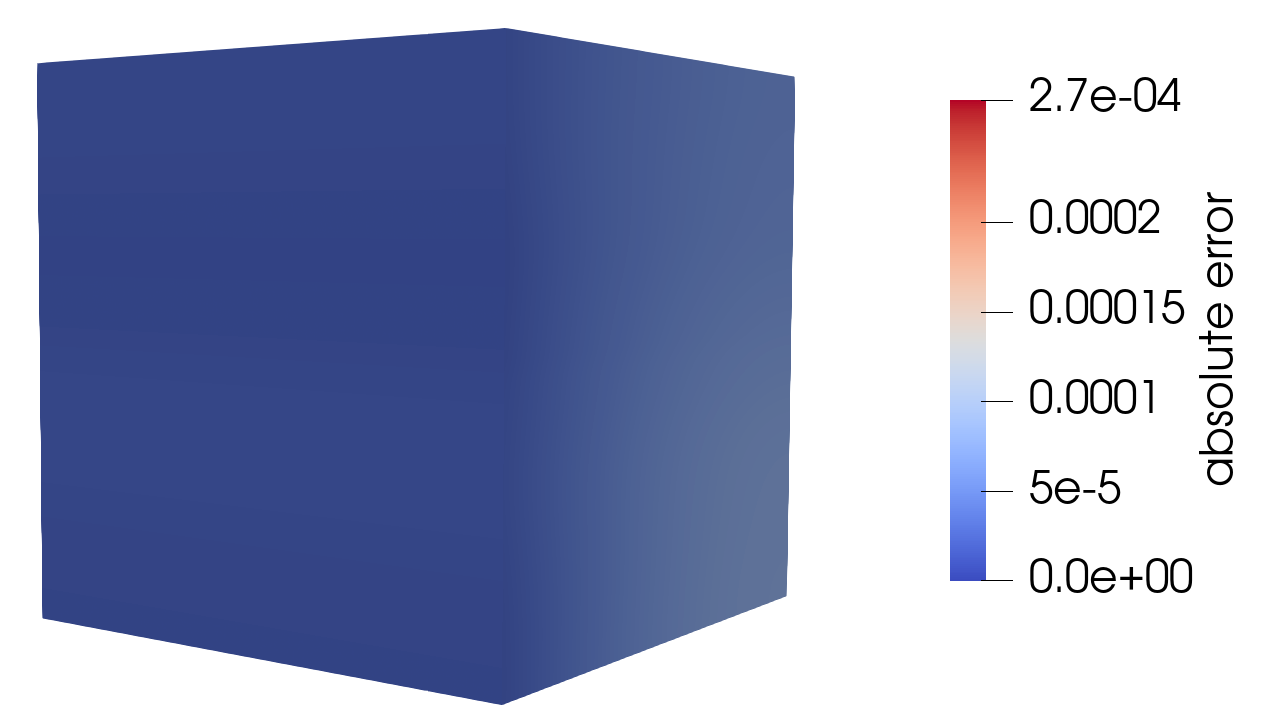}
	$\quad$
	\includegraphics[width=0.3\textwidth]{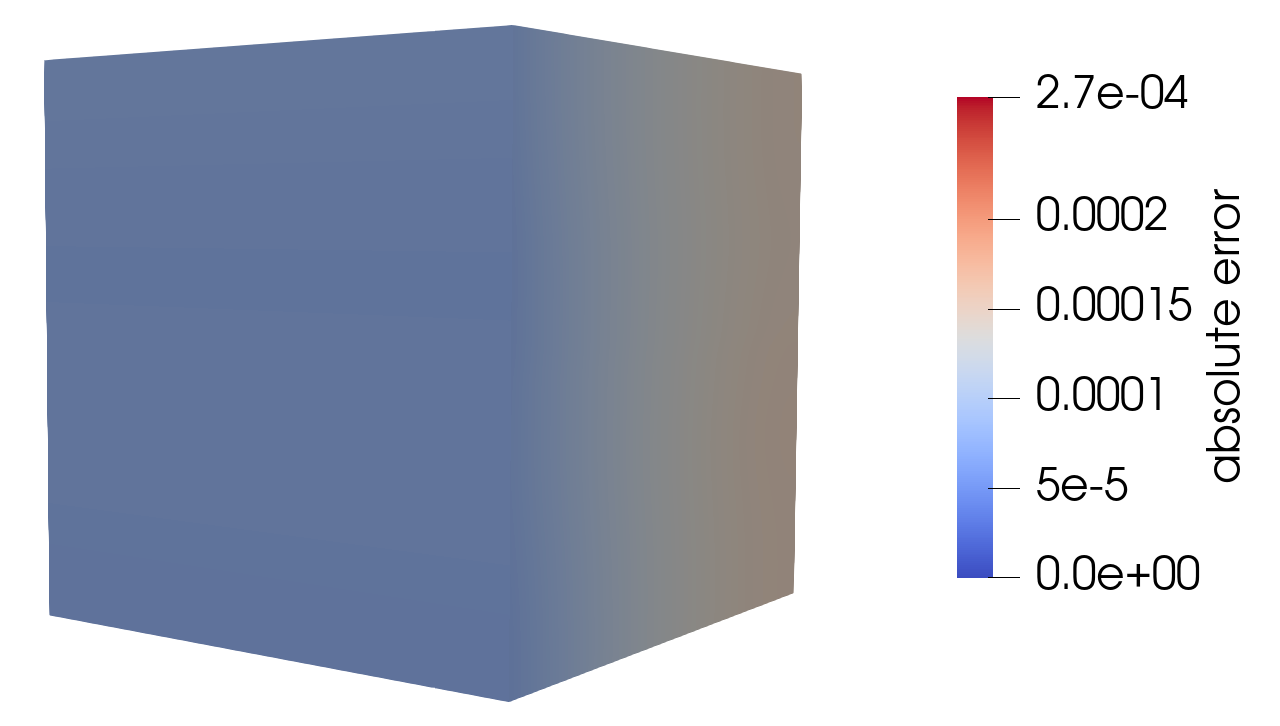}
	$\quad$
	\includegraphics[width=0.3\textwidth]{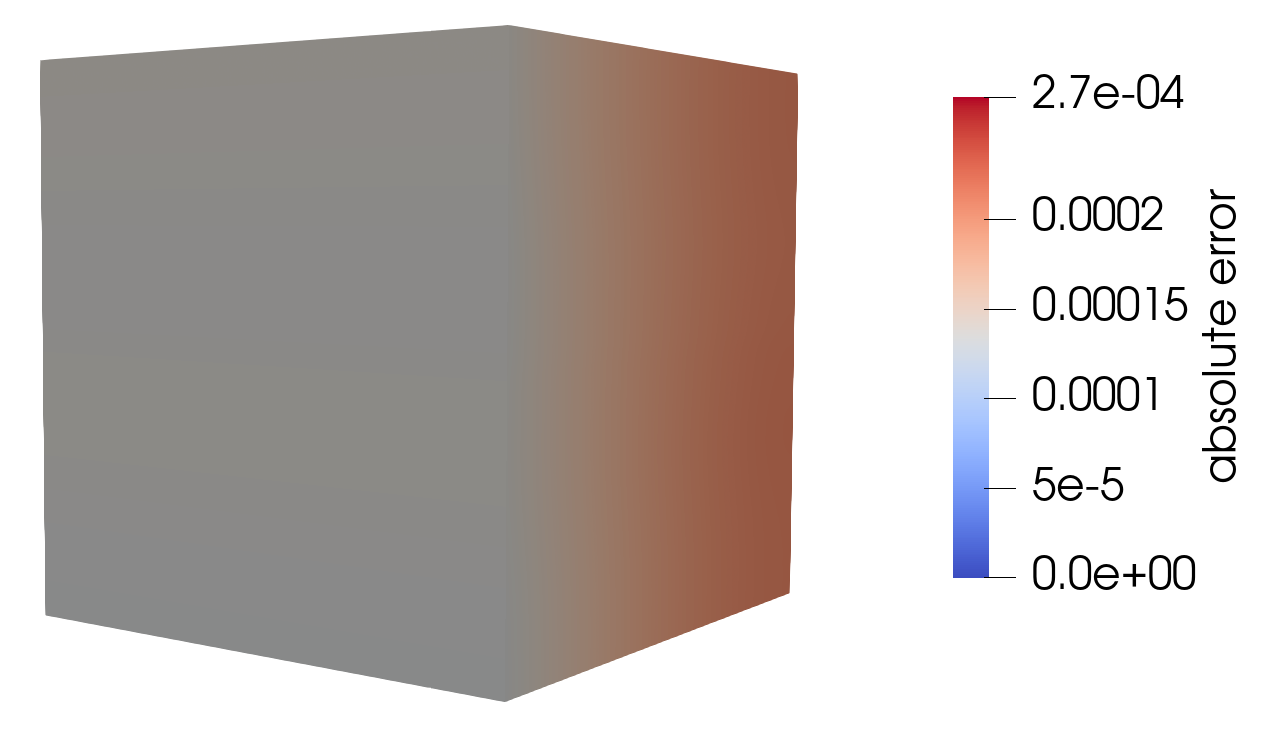}
	\caption{\emph{Test\#3.} Master solution FOM (top), ROM (center) solutions and absolute error (bottom) for $\alpha= 2.75$ and three different time instants. }
	\label{fig:snapshots_heat_slave}
\end{figure} 

\begin{figure}[h!]
	\centering
	\includegraphics[width=0.4\textwidth]{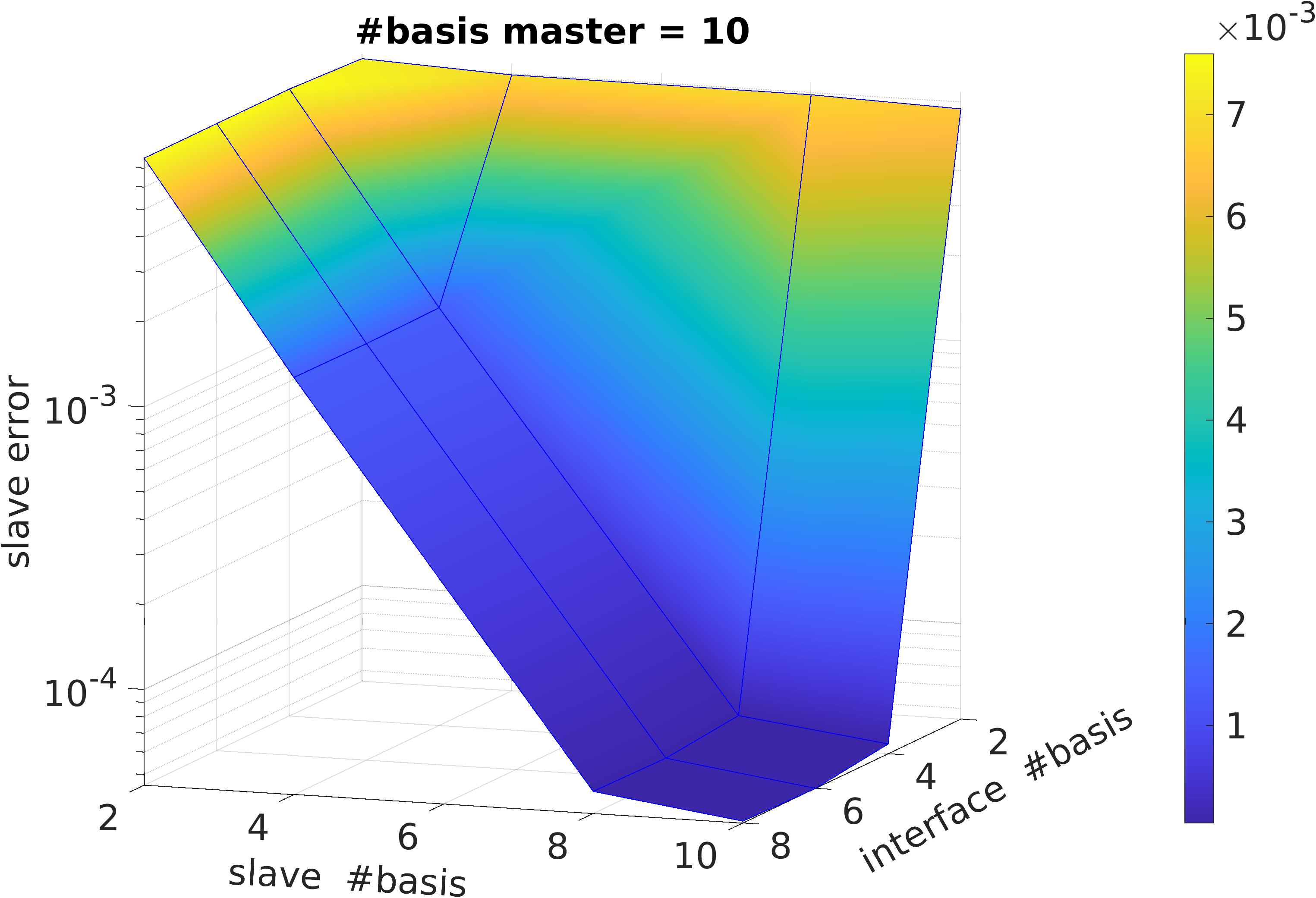}
	\includegraphics[width=0.4\textwidth]{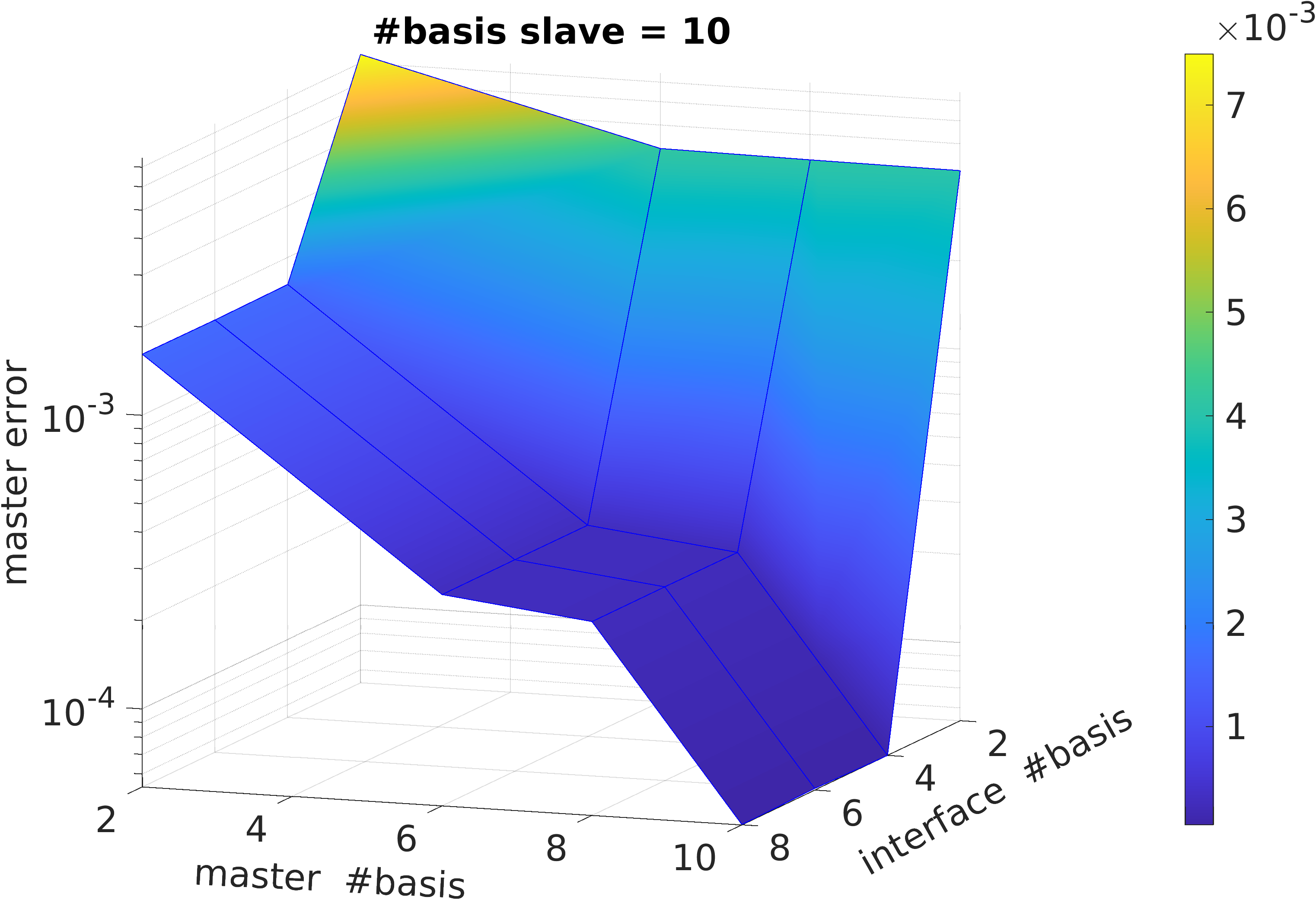}
	\\
	\bigskip
	\includegraphics[width=0.4\textwidth]{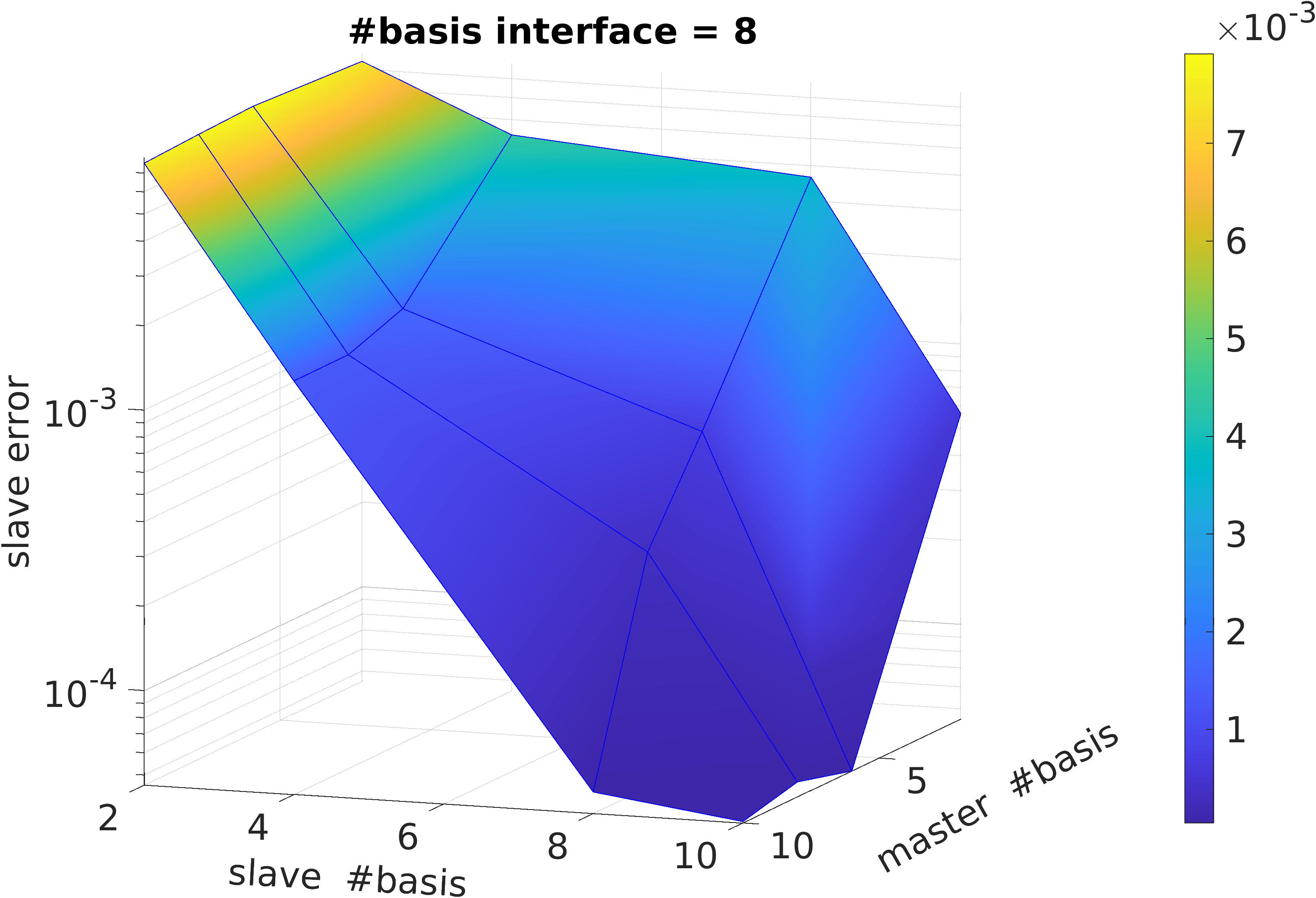}
	\includegraphics[width=0.4\textwidth]{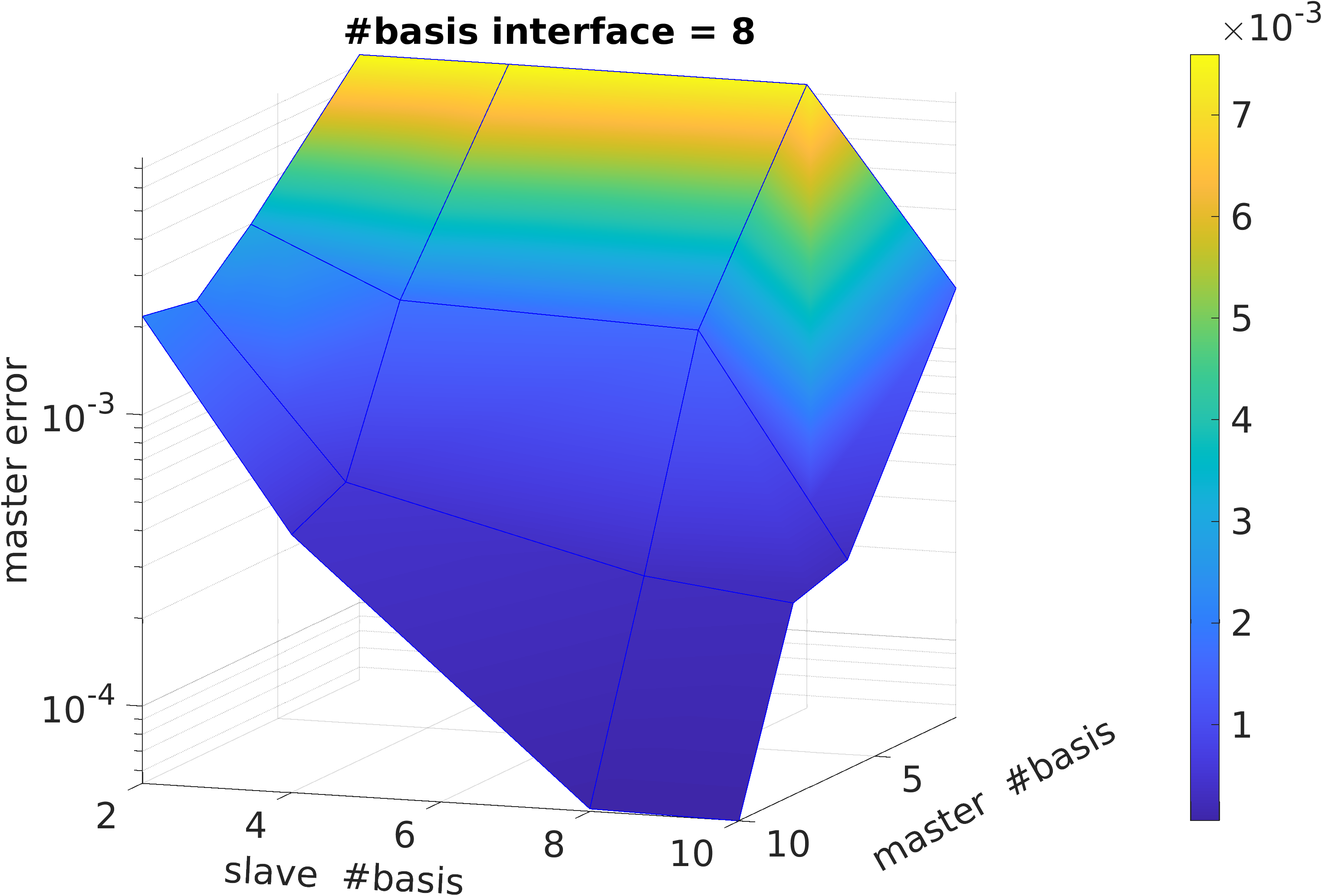}
	\caption{\emph{Test\#3.} $H^1(\Omega_i)$ mean relative error ($z$--axis) over the solution for $N_t N_\text{test} = 20$ different instances of the parameters between the FOM and ROM solutions varying the number of basis functions used to represent the slave and the master solution, and the interface data ($x$-- and $y$--axis). On the top row, we fix the number of basis functions of the master problem to 10 (on the left) and to 10 for the slave problem (on the right), while on the bottom we fix 8 the number of basis functions for the interface data representation.}
	\label{fig:Heat_error}
\end{figure}

\begin{figure}[h!]
	\centering
    \includegraphics[width=0.45\textwidth]{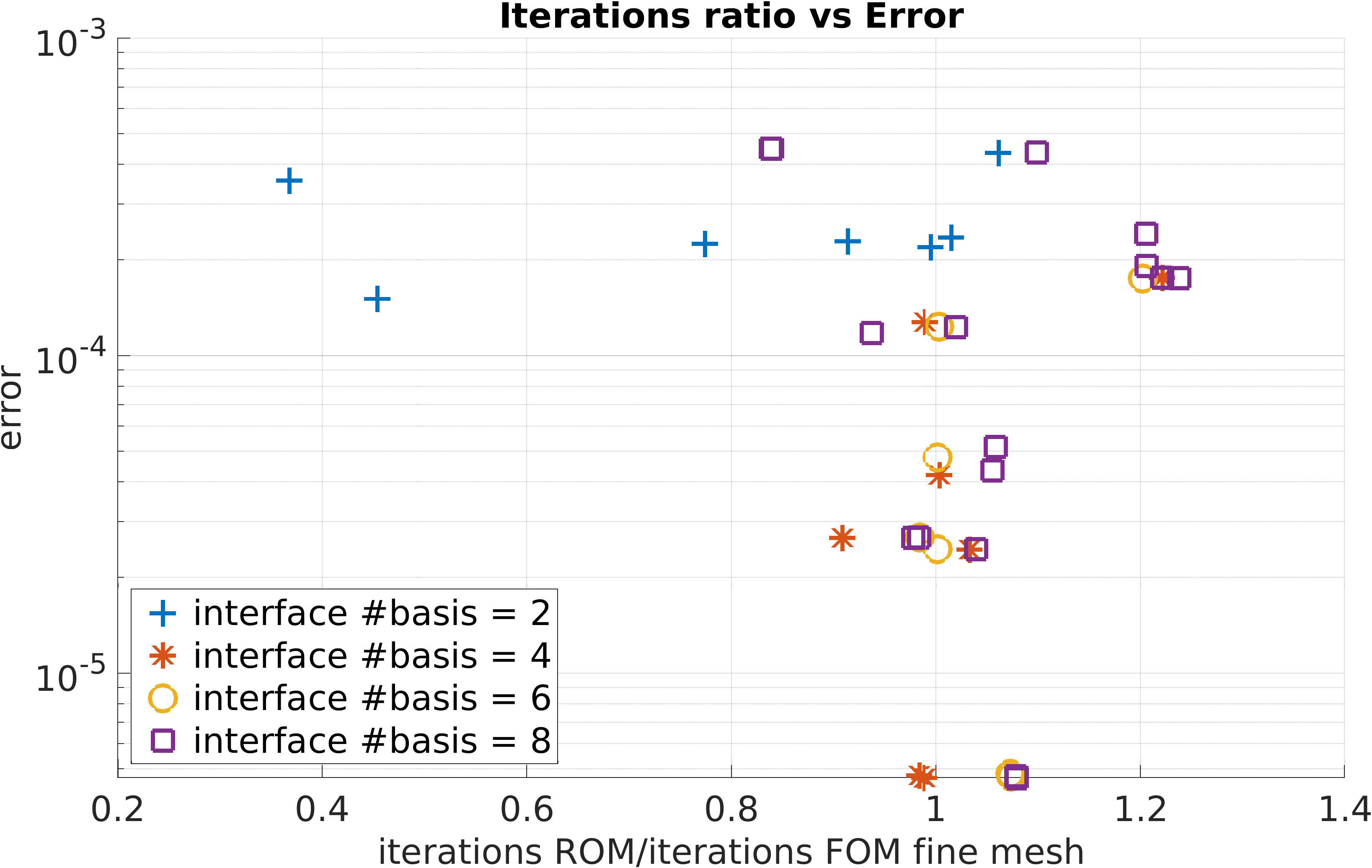}
    \quad
	\includegraphics[width=0.45\textwidth]{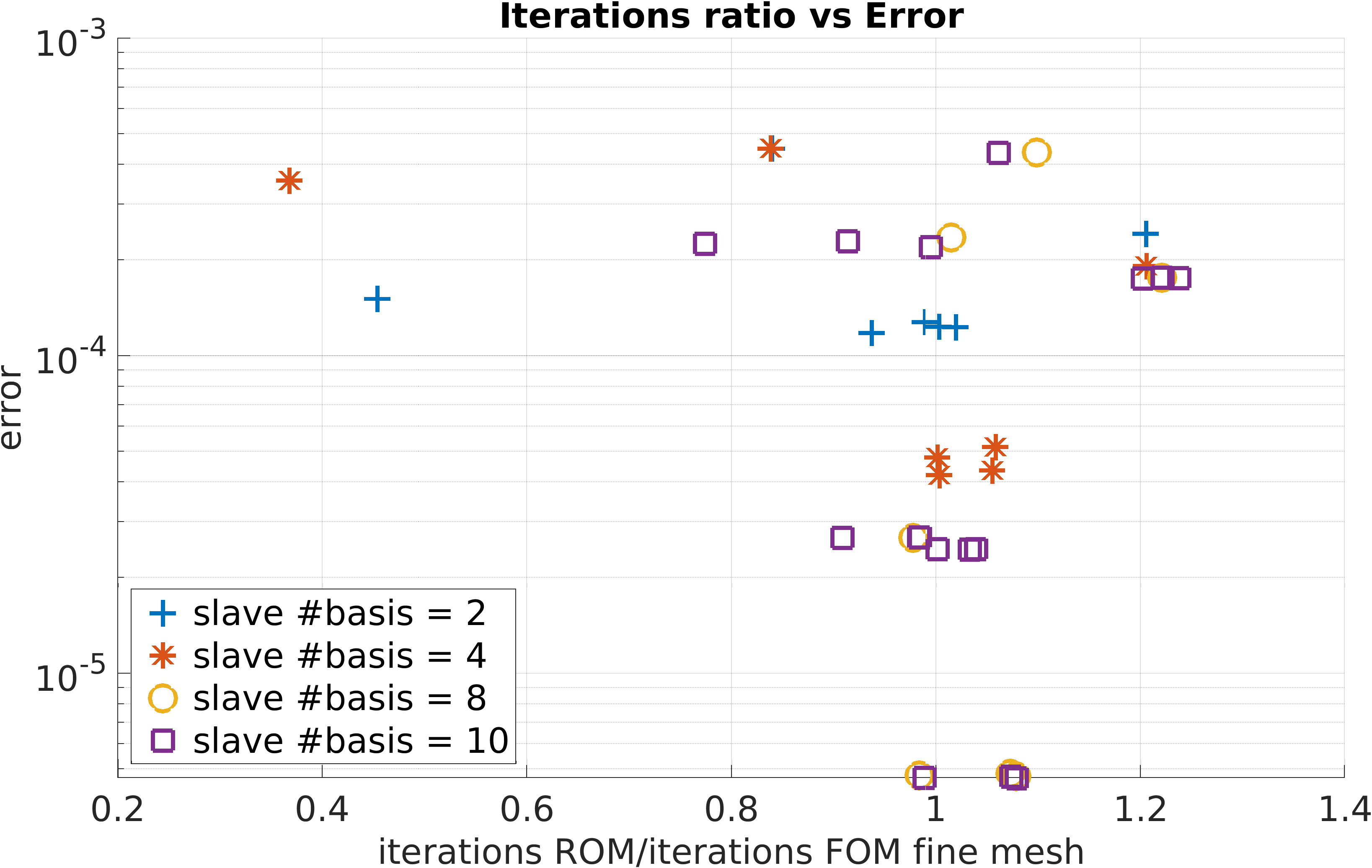}
	\caption{\emph{Test\#3.} Ratio between the number of iterations obtained with ROM and FOM schemes versus the achieved approximation error depending on the number of basis functions employed for the interface data (left) and the slave solution (right).}
	\label{fig:Heat_ratio_vs_error}
\end{figure}

We compute the $H^1(\Omega_i)$ relative error \eqref{eq:error_2norm} between the ROM and the FOM solutions with the coarser mesh (see Table \ref{Tab:heat_time}) for the slave problem, and between the ROM and the FOM solutions with the coarser mesh (see Table \ref{Tab:heat_time}) for the master problem. In Figs. \ref{fig:snapshots_heat_master} and \ref{fig:snapshots_heat_slave} we draw the ROM and FOM solutions and the absolute errors for three different time instants and some selected values of $\alpha$, while in Fig. \ref{fig:Heat_error} we plot the $H^1(\Omega_i)$ relative error for different ROM dimensions. Note that whenever we fix the number of basis functions used to approximate the master solution, we vary the number of basis functions for both the slave solution and the interface data. Similarly, when we fix the number of basis functions for either the slave solution or the Dirichlet--Neumann interface data, we vary the number of basis functions for the other sets. Both the slave and the master approximation errors show the expected behaviour: they decrease whenever we increase the employed number of basis functions. Different from the steady test cases (see Subsections \ref{Subsect:steady_case} and \ref{Subsect:steady_case_source}), the approximations error is here equally dependent on the hyper--parameters $n_1$, $n_2$ and $M_1$ and $M_2$, as can also be observed in Fig. \ref{fig:Heat_ratio_vs_error}.

\begin{figure}[h!]
	\centering
    \includegraphics[width=0.45\textwidth]{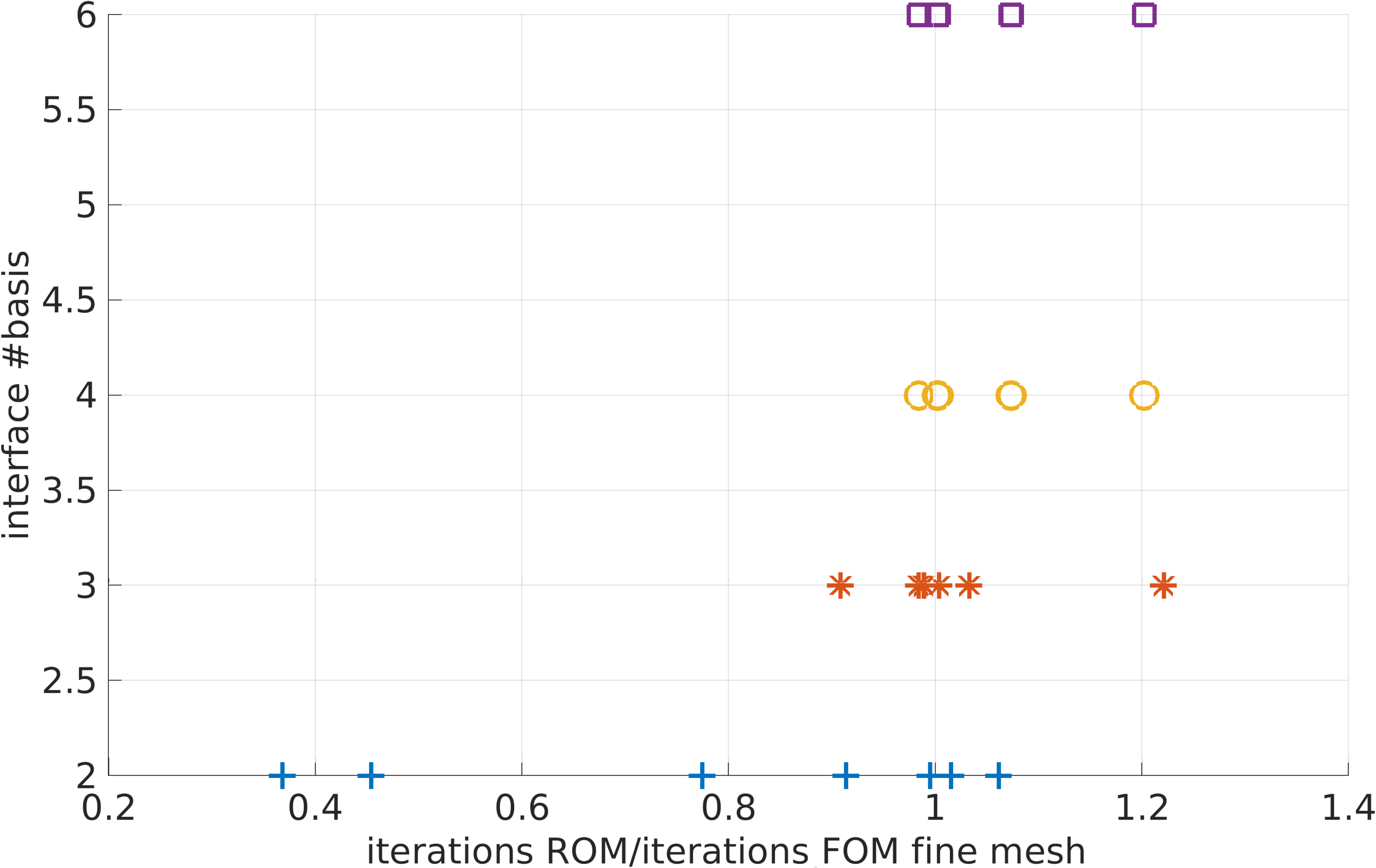}
    \quad
	\includegraphics[width=0.45\textwidth]{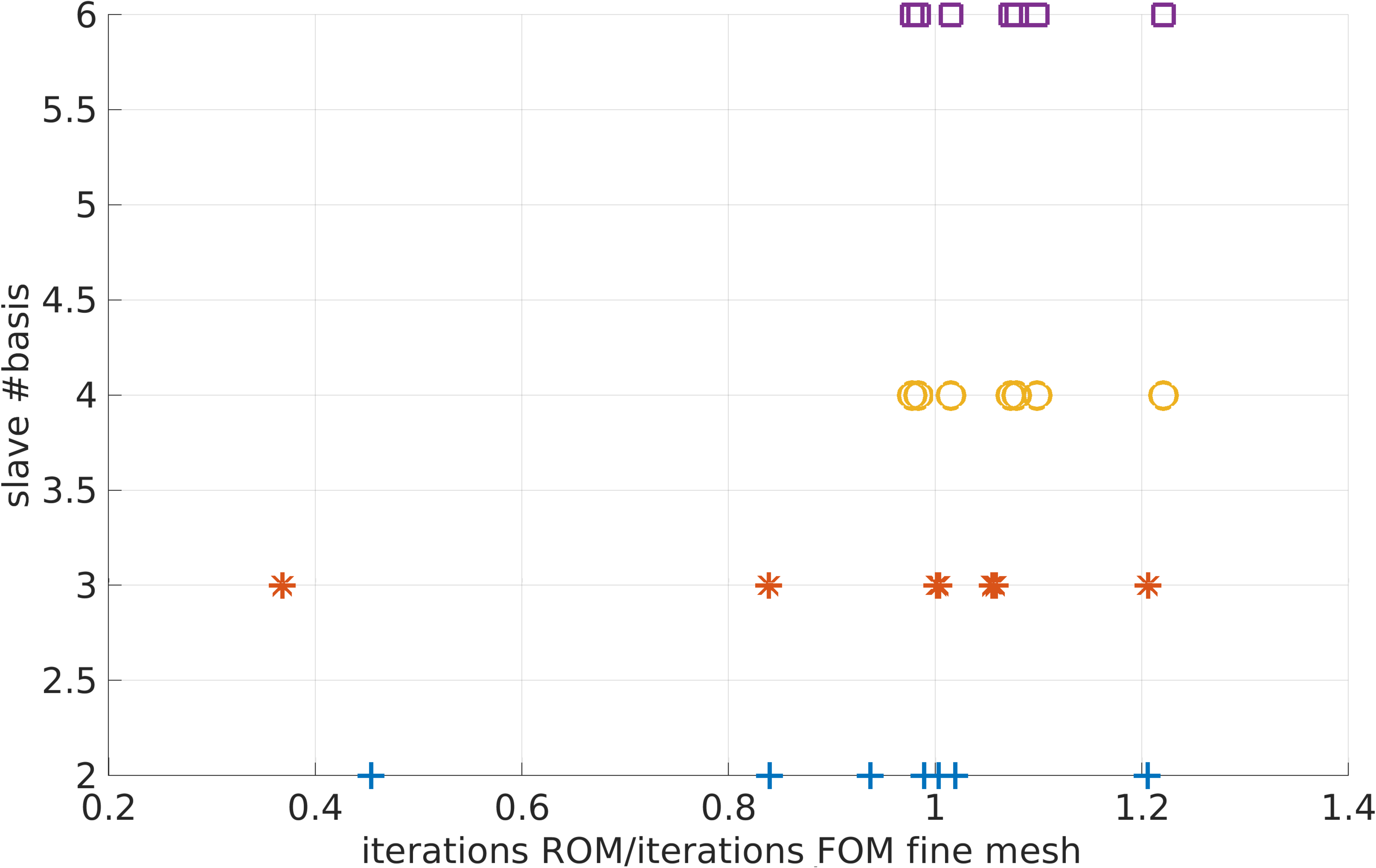}
	\\
	\bigskip
	\includegraphics[width=0.45\textwidth]{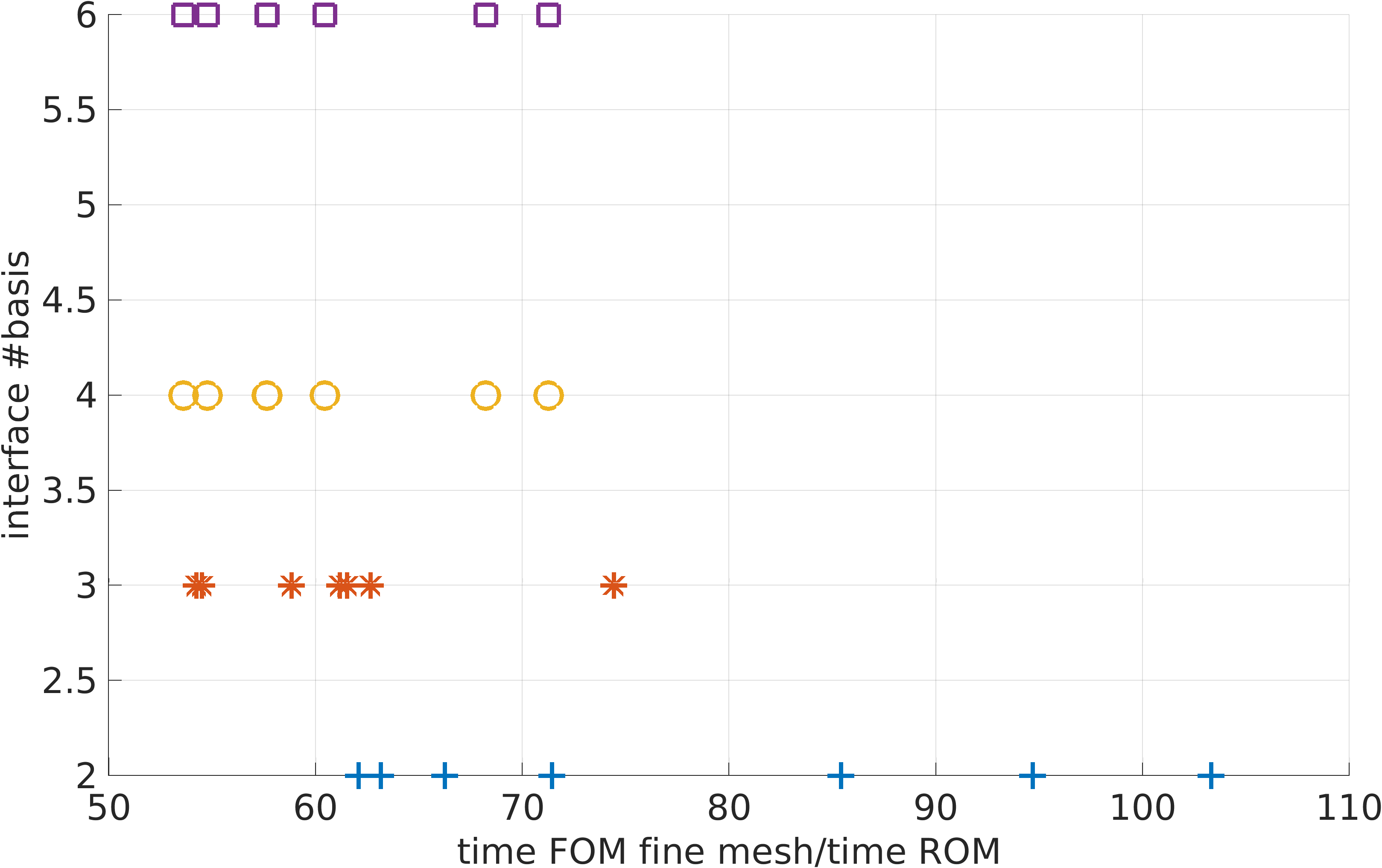}
    \quad
	\includegraphics[width=0.45\textwidth]{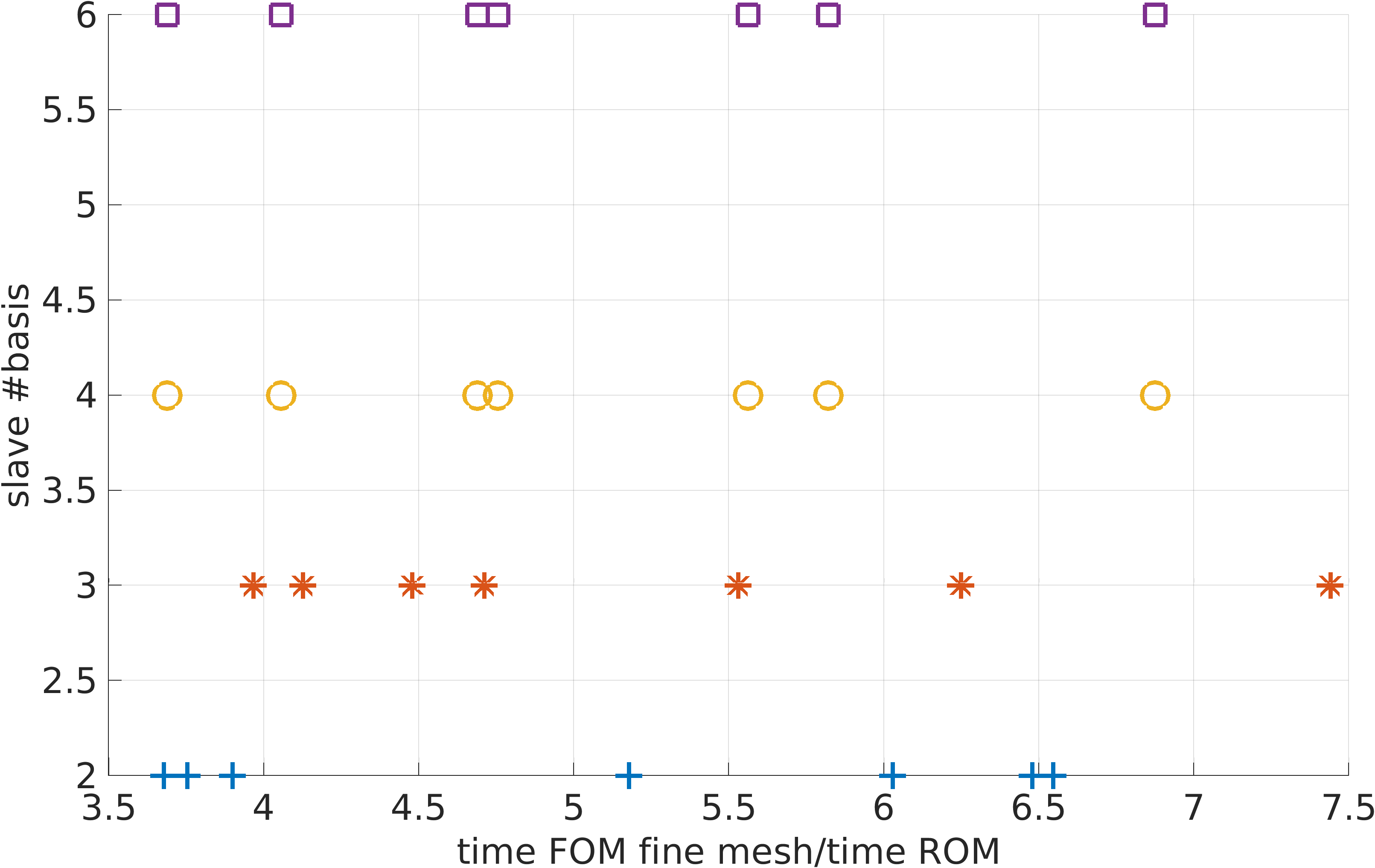}
	\caption{\emph{Test\#3.} Top row: the ratio between the number of iterations obtained with ROM and FOM schemes versus the number of basis functions employed to approximate the interface data (left) and the slave solution (right). Bottom row: the ratio between the FOM and ROM computational time versus the number of basis functions employed to approximate the interface data (left) and the slave solution (right). The FOM simulation referred to fine discretization.}
	\label{fig:Heat_ratio_vs_basis}
\end{figure}

Finally, in Table \ref{Tab:heat_time} we compare the number of DoFs employed to solve the FOM and ROM. The time average iterations varying the ROM interface hyper--parameters or the $n_1$ can be found in Fig. \ref{fig:Heat_ratio_vs_basis} (top row). Such graphs show that increasing both $M_1$ (and $M_2$) and $n_1$, the number of iterations required by the ROM to reach convergence of the interface solutions increases. A corresponding decrease in speed up can be seen in the bottom row of Fig. \ref{fig:Heat_ratio_vs_basis}. When the approximation tolerance of $10^{-5}$ is selected for each reduced quantity with the coarsest (finest) discretization, the FOM solution requires $37min$ and $18s$ ($7h$ and $52 min$) of computations, whereas only $7min$ and $45s$ are needed to solve the ROM. As before, our reduced technique shows high computational efficiency, effectively achieving a speedup of 5 times with respect to the FOM solution with the coarsest discretization and of about 60.9 times compared to the finest discretization, corresponding to a reduction of about 80\% and 98\% of the computational costs, respectively. These results refer to a complete simulation of 100 time--steps, including some repeated tasks such as the assembling of the right--hand side, according to the time discretization scheme implemented. 

\begin{table}[h!]
	\centering 
	\begin{tabular}{c|ccccc}
		\toprule
		&&Master solution &Slave solution &Master interface &Slave interface\\
		\hline &&&&&\\ [-2ex]
		FOM -- coarse mesh &\#DoFs &35937 &35937 &1089 &1089 \\
		FOM -- fine mesh &\#DoFs &274625 &274625 &4225 &4225  \\
		\hline &&&&&\\ [-2ex]
		\multirow{2}{*}{ROM} &\#Basis &10 &10 &10 &10 \\
        &\#DoFs &274625 &35937 &4225 &1089\\
		\bottomrule	
	\end{tabular}
	\caption{\emph{Test\#3.} High fidelity and reduced order model dimensions of subdomains and interface discretization, as well as the number of basis functions required to achieve an approximation error of $10^{-5}$.}
	\label{Tab:heat_time}
\end{table}

\section{Conclusion}
\label{sec:conclusion}
In this work, we have introduced a reduced order modeling technique based on RB methods to decrease the computational costs entailed by the solution of two--way coupled problems employing Dirichlet--Neumann iterative schemes. The modularity of the procedure ensures the efficiency of the ROM through different treatments of the master, slave, and interface data reduction. Indeed, the master and the slave modes can be reduced using appropriate RB strategies, while the Dirichlet and the Neumann interface data can be handled through the DEIM inside the INTERNODES method, highlighting the possibility of using such reduced scheme to transfer the interface data between the conforming and non--conforming interface grids. The proposed algorithm can also be applied, in principle, to more general multi--physics problems that are solved through Dirichlet--Neumann domain--decomposition iterative schemes. 

The numerical test cases show that our ROM is very cheap in the online phase, outperforming the online CPU time of the FOM when either fine or coarse conforming meshes in the two domains are considered. Indeed, a decrease of 98\% of the CPU time can be achieved in the case of a time--dependent parabolic PDE problem (see Subsection \ref{sub:heat_problem}) with respect to a finite element FOM built over a fine discretization, and of 85\% with respect to a finite element FOM built over a coarser one. An analysis of the iteration cost as a function of the reduced order hyper--parameters is also carried out.


This paper represents (to our knowledge) the first attempt toward the achievement of a fully--reduced RB interpolation--based numerical scheme for the two--way coupled model. 

Future developments concern the use of the INTERNODES method on the FOM offline phase and the use of interpolation operators in the ROM online phase (see Remark \ref{rem:interporetation_DEIM_interface}) that are more accurate than those used in this paper. We are confident in this way to improve the accuracy as well as the efficiency (of the offline phase) of our approach. Finally, future work will be the derivation of an error estimate for the proposed method.

Considering the results proposed in the last section, we expect to be able to apply the presented method to more complex and relevant physics--based coupled problems, \emph{e.g.} cardiac electrophysiological and fluid--structure interaction problems. 
Moreover, having been able to correctly reduce Dirichlet and Neumann interface data, the final extension of this work would be the application of our technique to Robin--like interface conditions, too. All these aspects will be the focus of future papers.

\bibliography{Ref}

\section*{Statements and Declarations}
This research has been funded partly by the European Research Council (ERC)
under the European Union's Horizon 2020 research and innovation programme
(grant agreement No. 740132, iHEART ``An Integrated Heart Model
for the simulation of the cardiac function'', P.I. Prof. A. Quarteroni)
and partly by the Italian Ministry of University and Research (MIUR) within the PRIN
(Research projects of relevant national interest 2017
``Modeling the heart across the scales: from cardiac cells to the whole organ''
Grant Registration number 2017AXL54F).\\
This manuscript has no associated data.  

\begin{appendices}
\section{A detailed error and computational costs analysis for the steady test case}
\label{appendix:analysis_case_test_1}
We report here a detailed analysis of the model performances for the test case \emph{Test\#1.} described in Subsection \ref{Subsect:steady_case}, including a study of the effect of the ROM hyper--parameters on the approximation error, number of iterations before solution convergence, and overall computational costs. A simplified analysis for the other test cases, which confirm the results reported here, is given in Subsections \ref{sub:heat_problem}.

\begin{figure}[h!]
	\centering
	\includegraphics[width=0.24\textwidth]{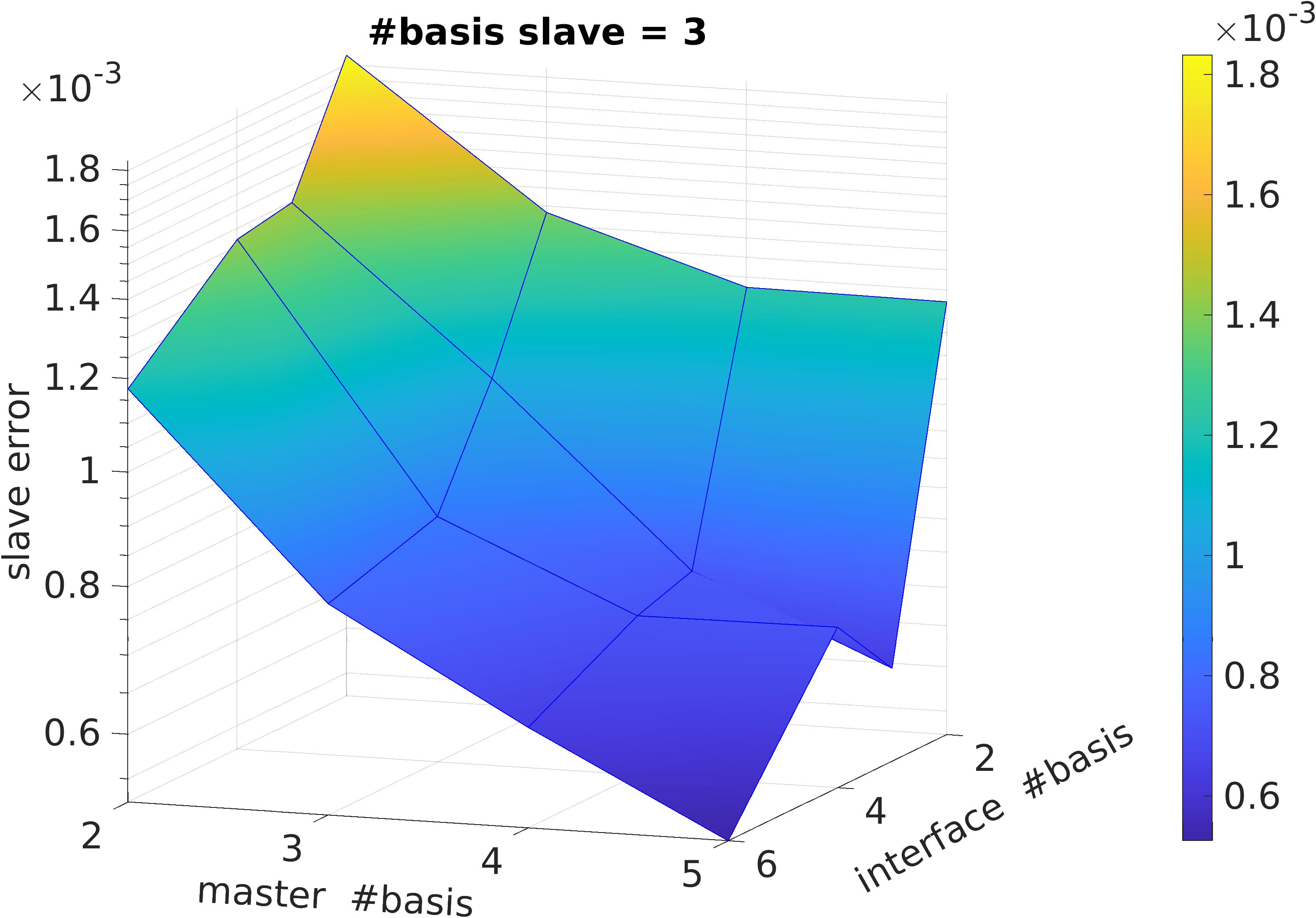}
	\includegraphics[width=0.24\textwidth]{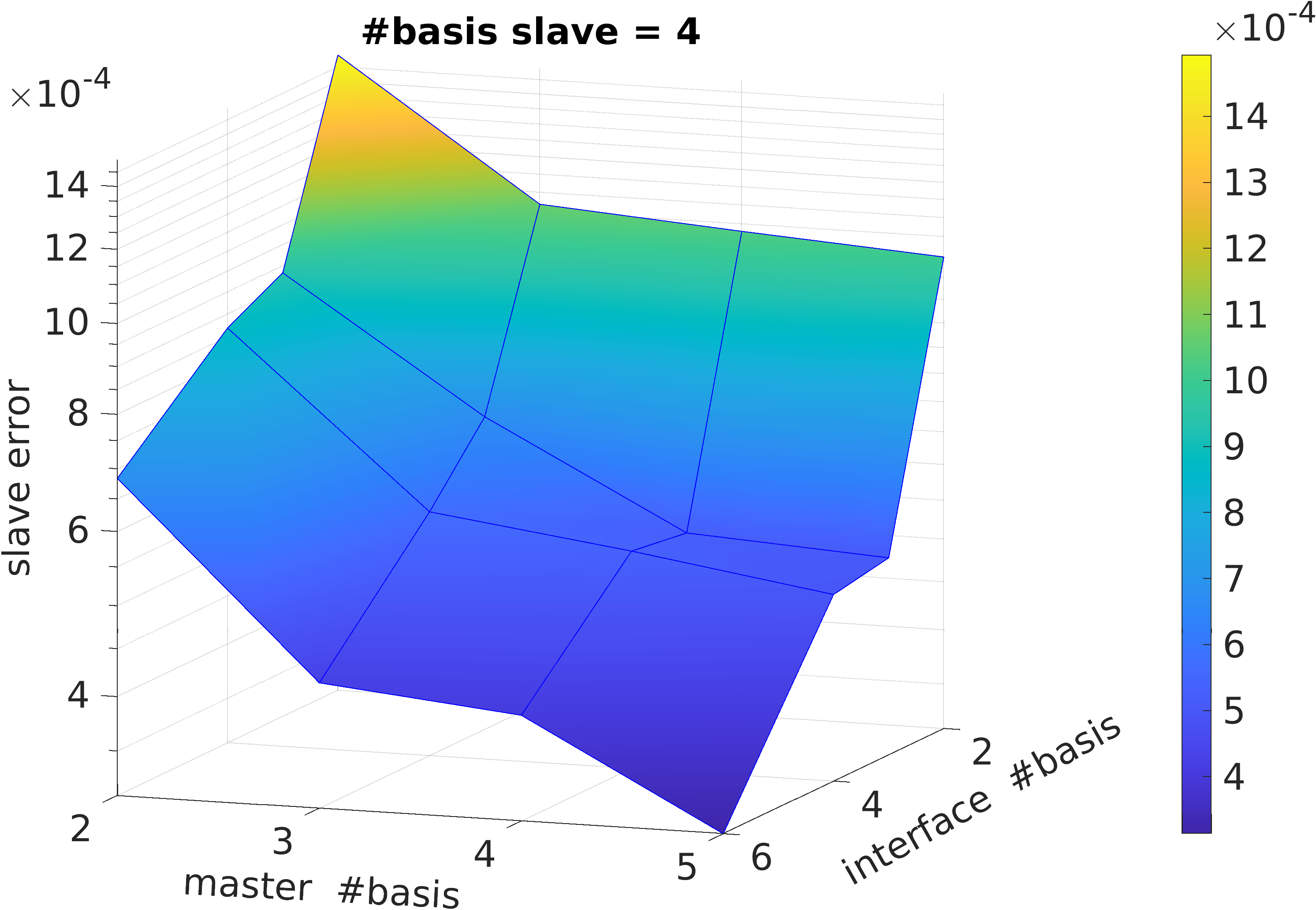}
	\includegraphics[width=0.24\textwidth]{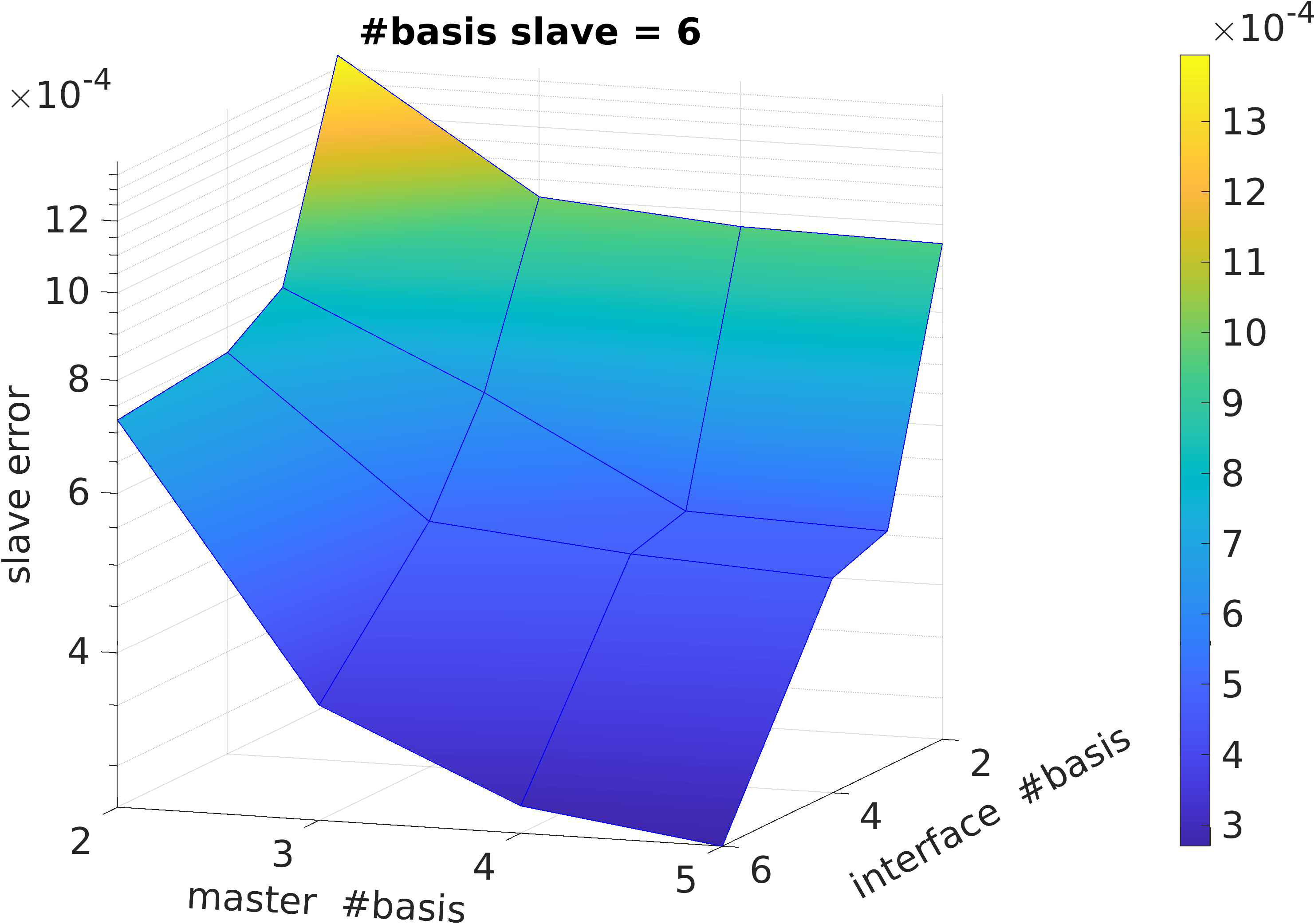}
	\includegraphics[width=0.24\textwidth]{slave_error_slave_7.png}\\
    \bigskip
     \includegraphics[width=0.24\textwidth]{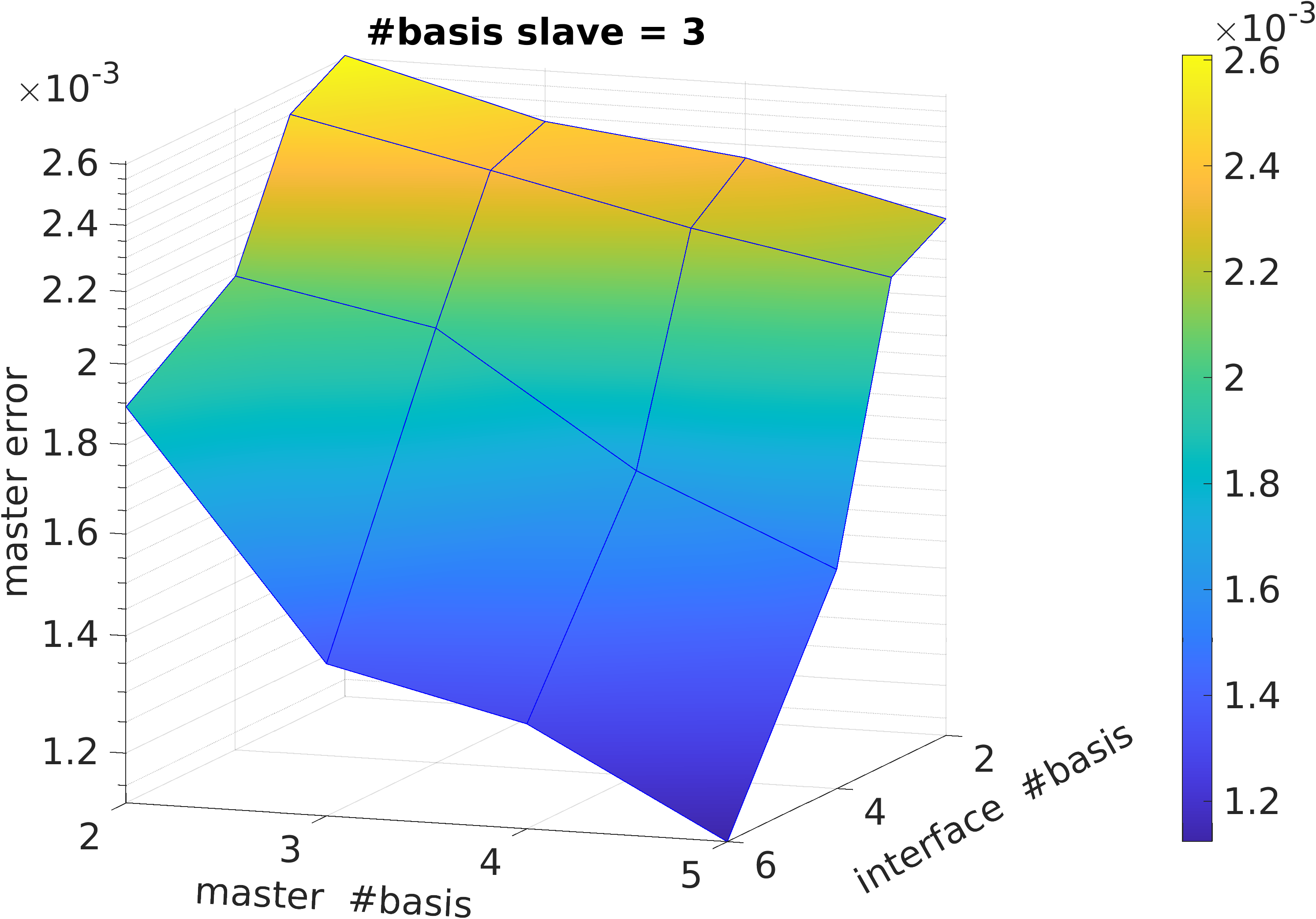}
	\includegraphics[width=0.24\textwidth]{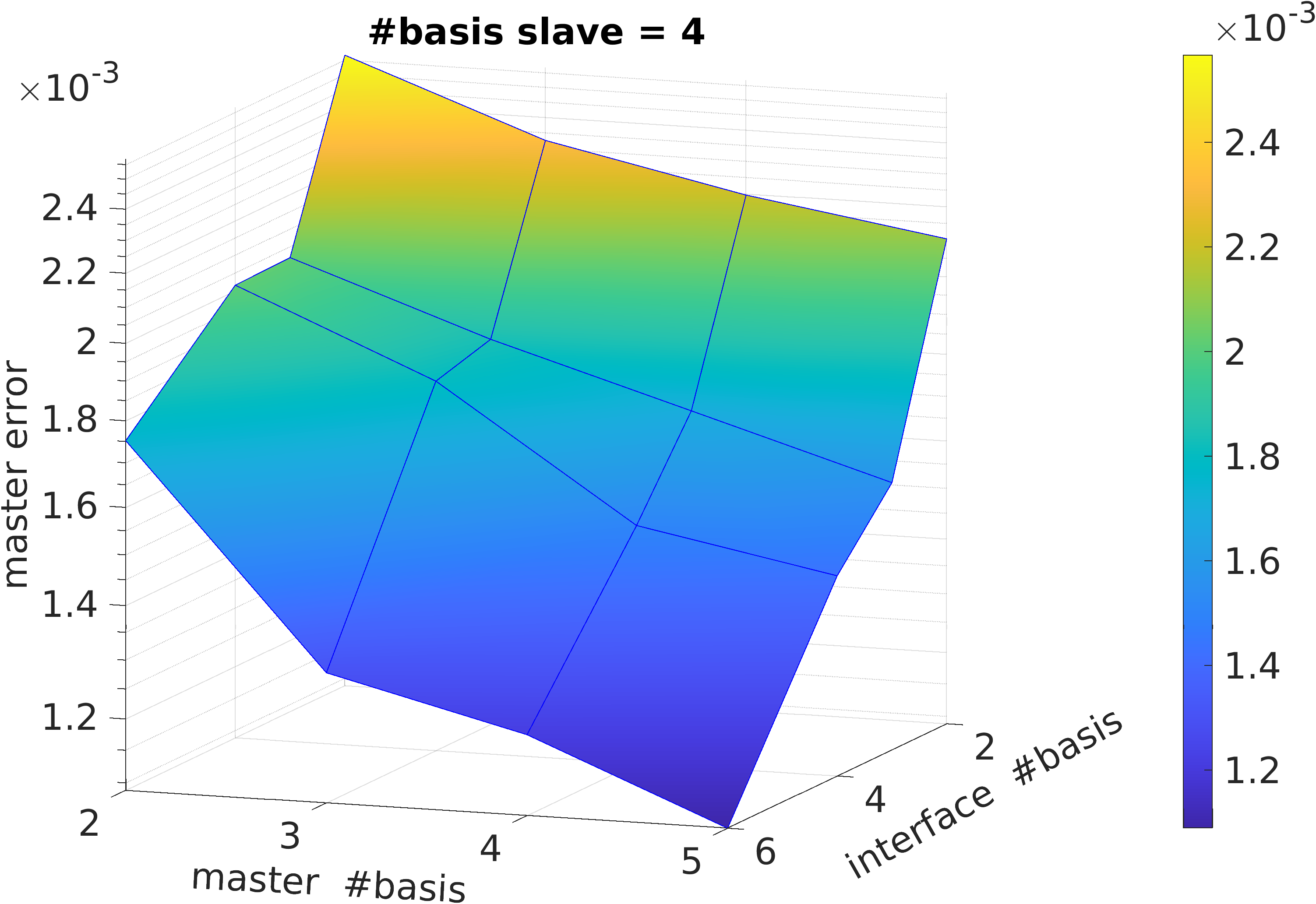}
	\includegraphics[width=0.24\textwidth]{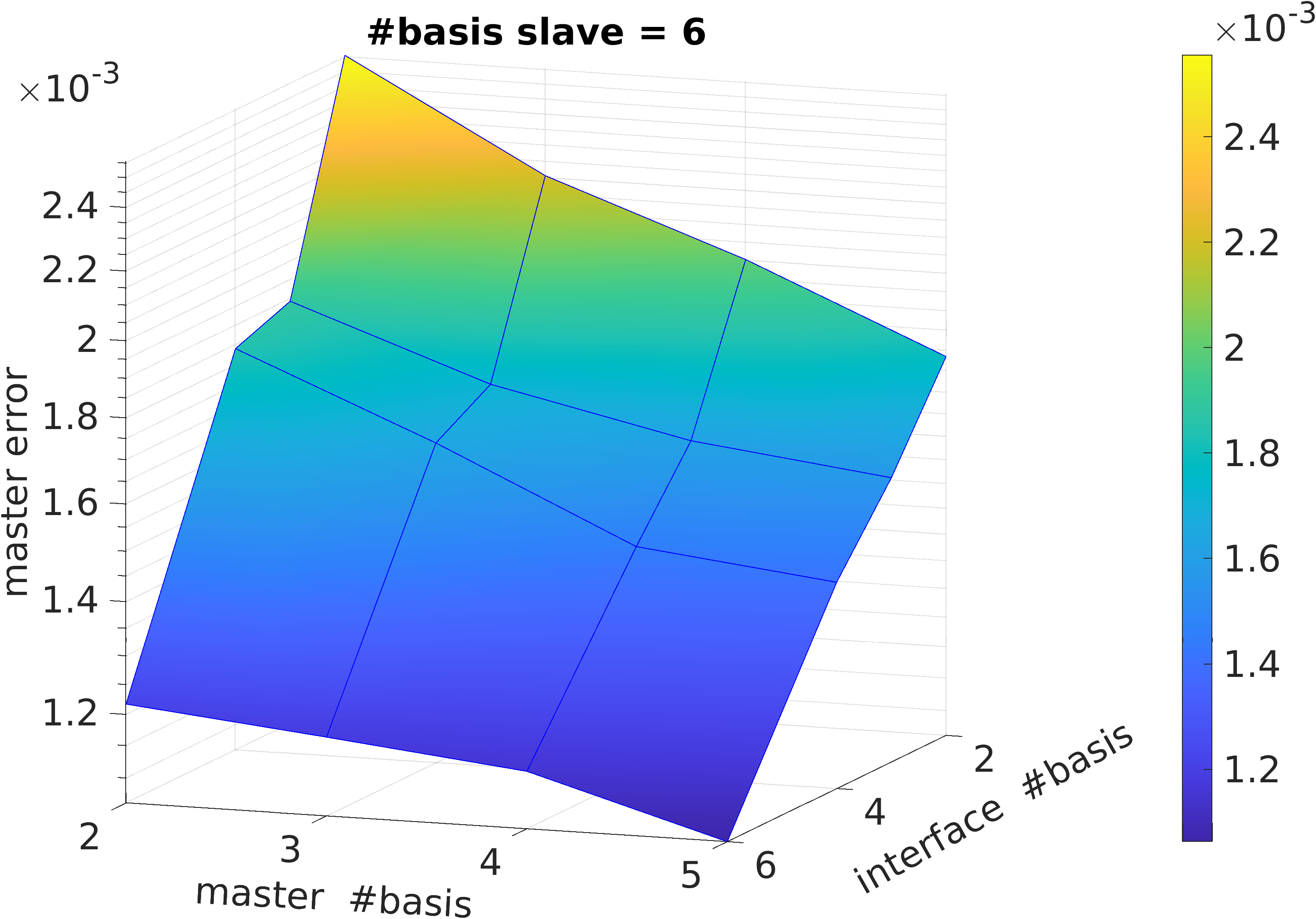}
	\includegraphics[width=0.24\textwidth]{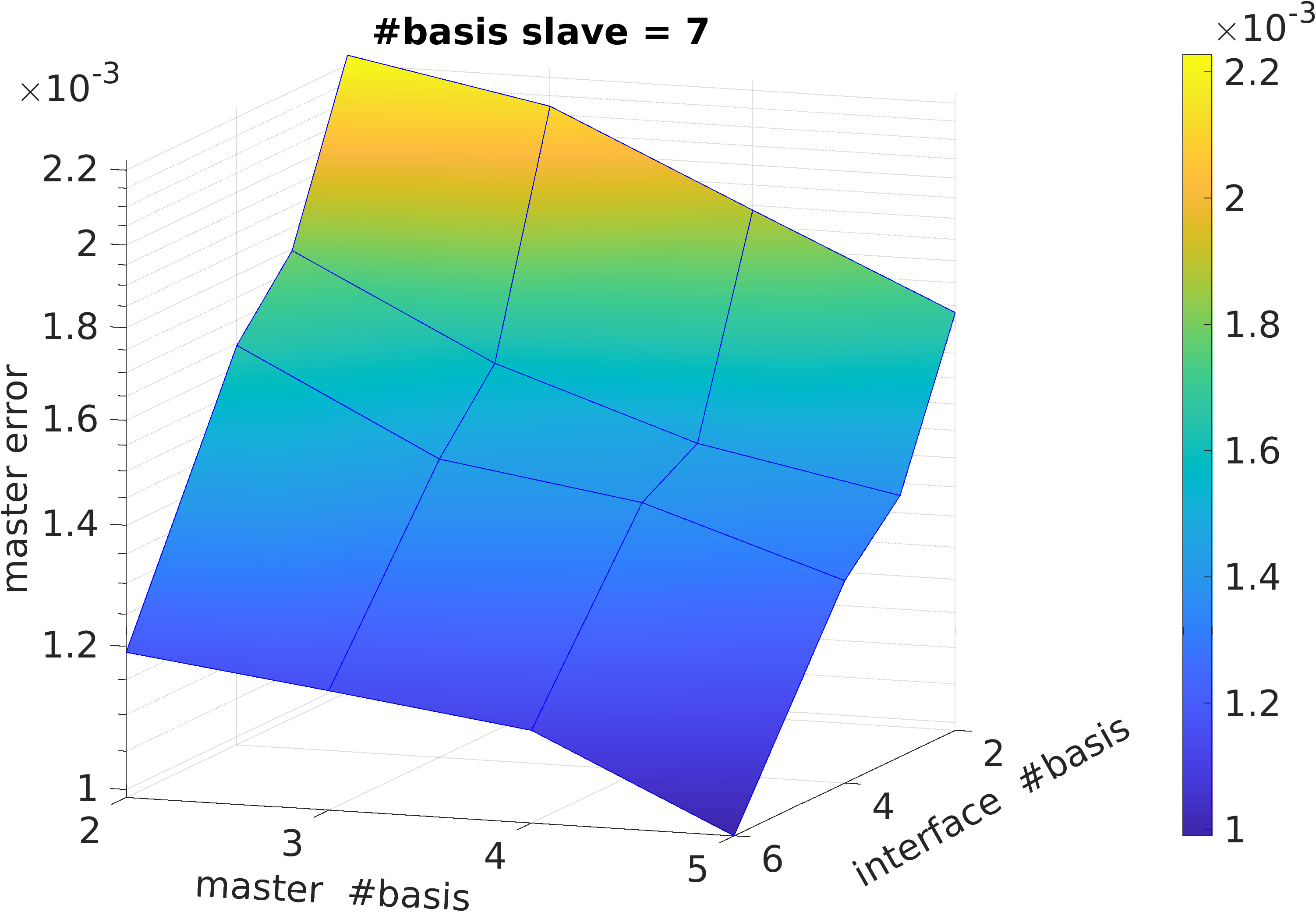}
	\caption{\emph{Test\#1.} $H^1(\Omega_i)$ mean relative error ($z$--axis) over the slave (top row) and master (bottom row) solution for $N_\text{test}=20$ different instances of the parameters between the FOM and ROM solutions varying the number of basis functions used to represent the master solution and $n_2$, and the interface data $M_1$ and $M_2$ ($x$-- and $y$--axis), and fixing the number of basis functions employed to approximate the slave solution.}
	\label{fig:fixed_slave_basis}
\end{figure}

\begin{figure}[t!]
	\centering
	\includegraphics[width=0.24\textwidth]{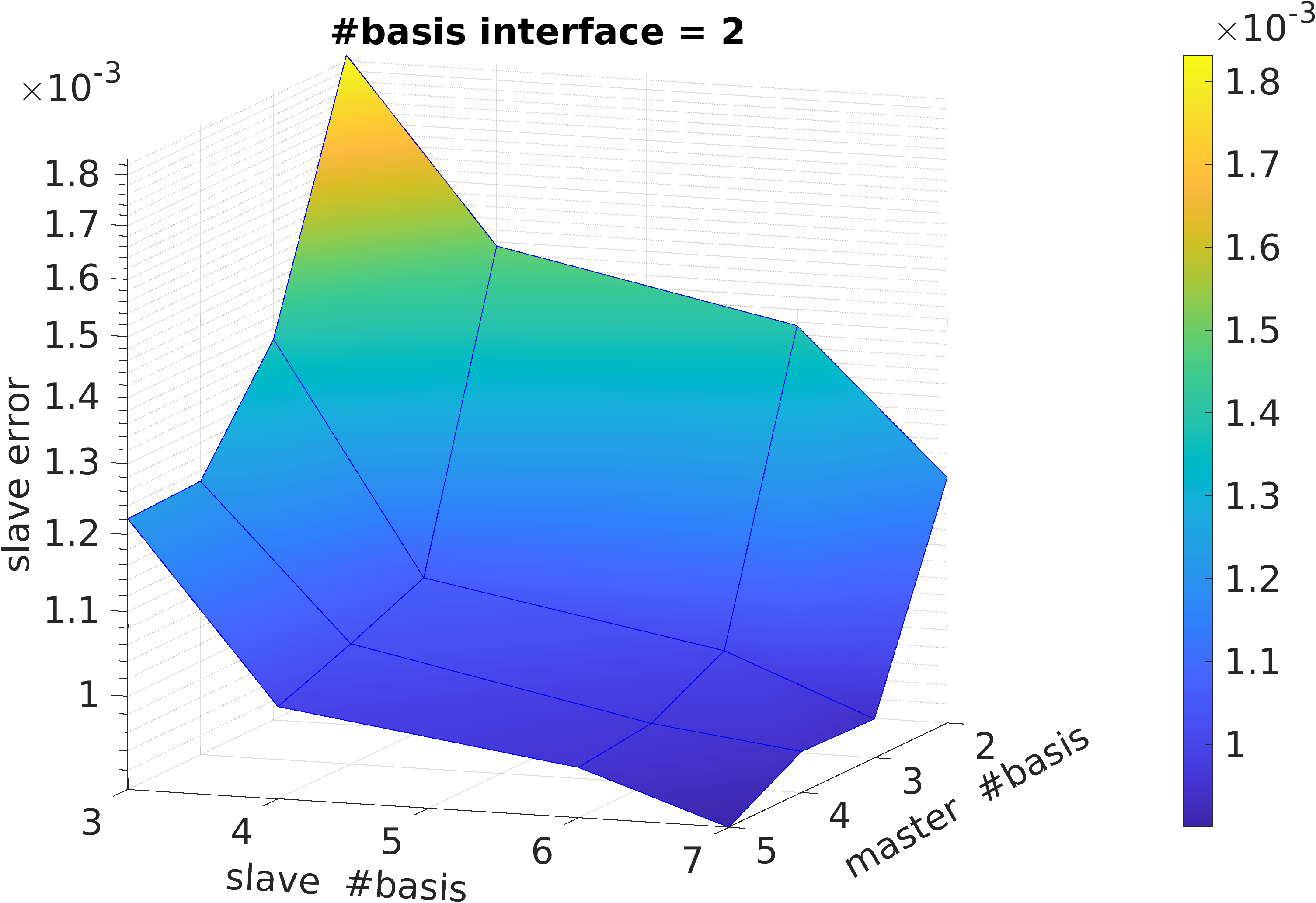}
	\includegraphics[width=0.24\textwidth]{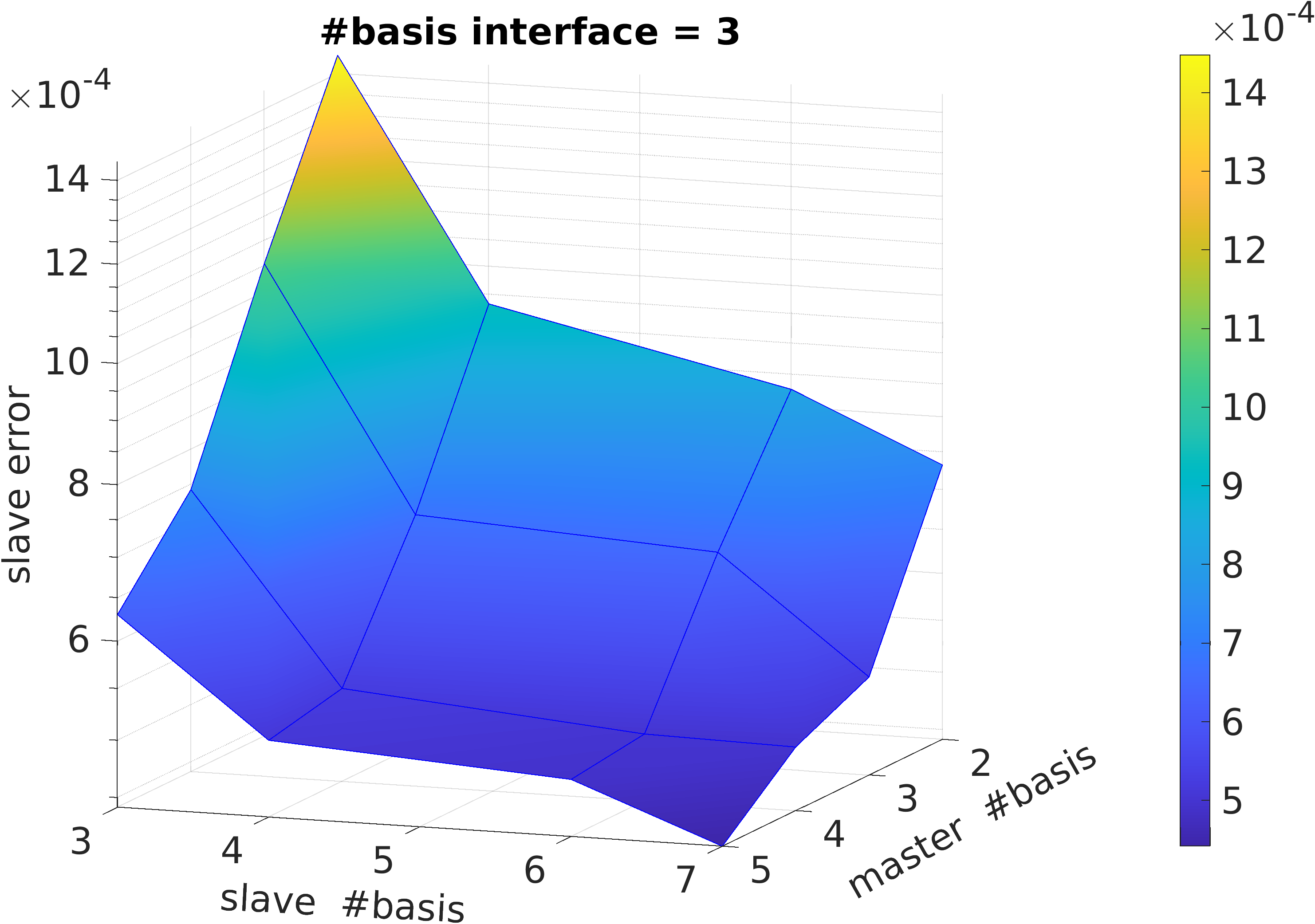}
	\includegraphics[width=0.24\textwidth]{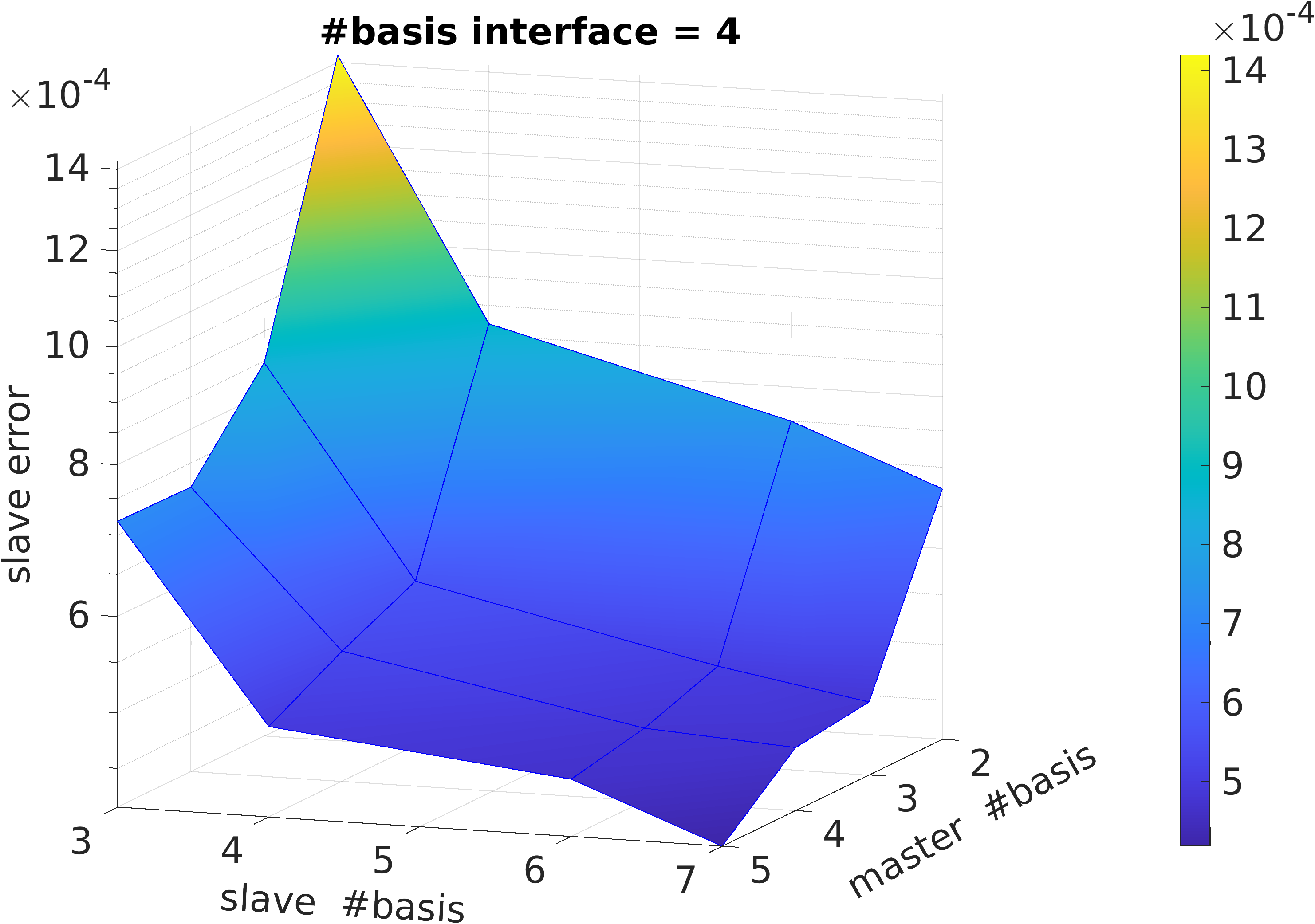}
	\includegraphics[width=0.24\textwidth]{slave_error_interface_6.png}\\
    \bigskip
     \includegraphics[width=0.24\textwidth]{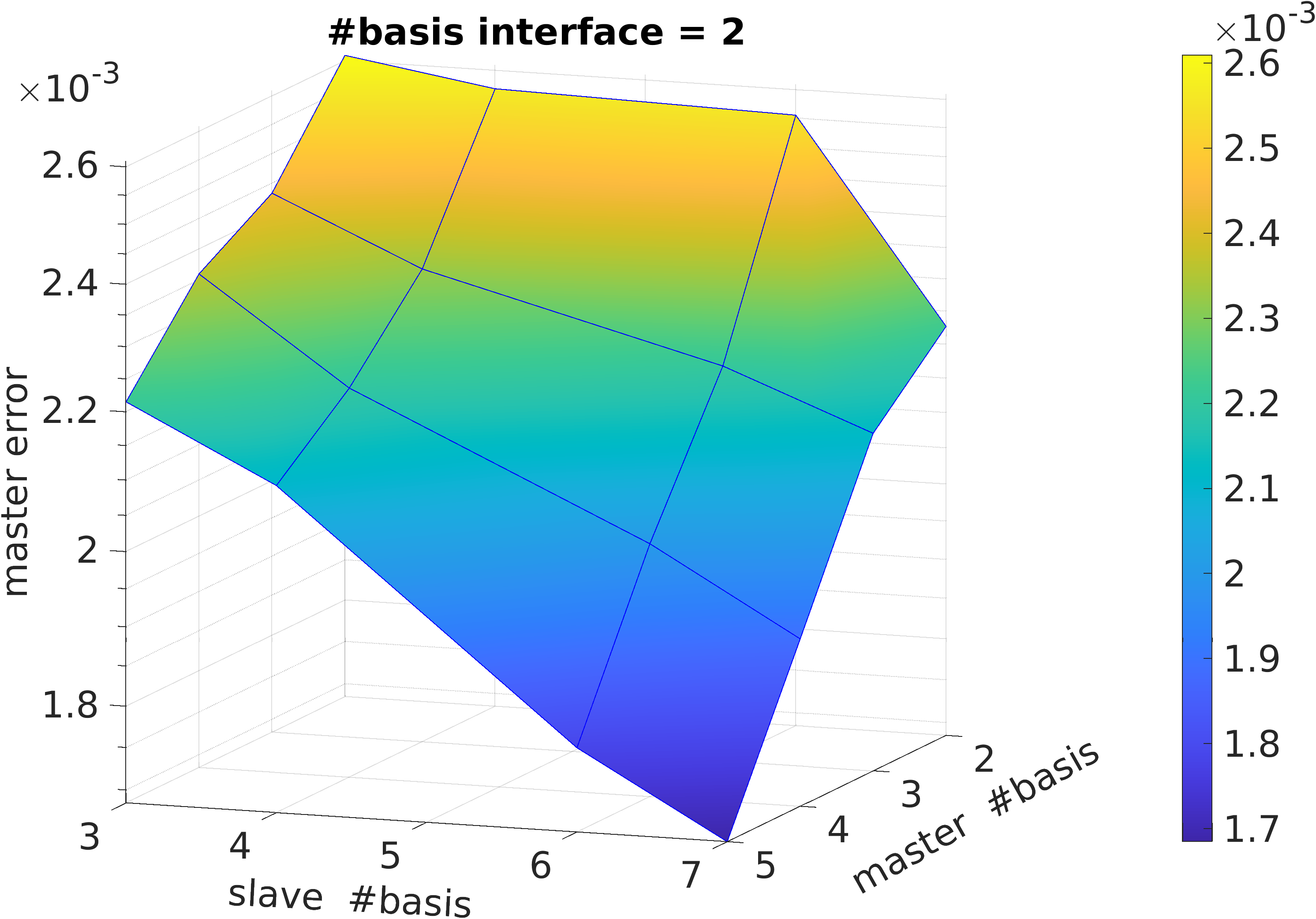}
	\includegraphics[width=0.24\textwidth]{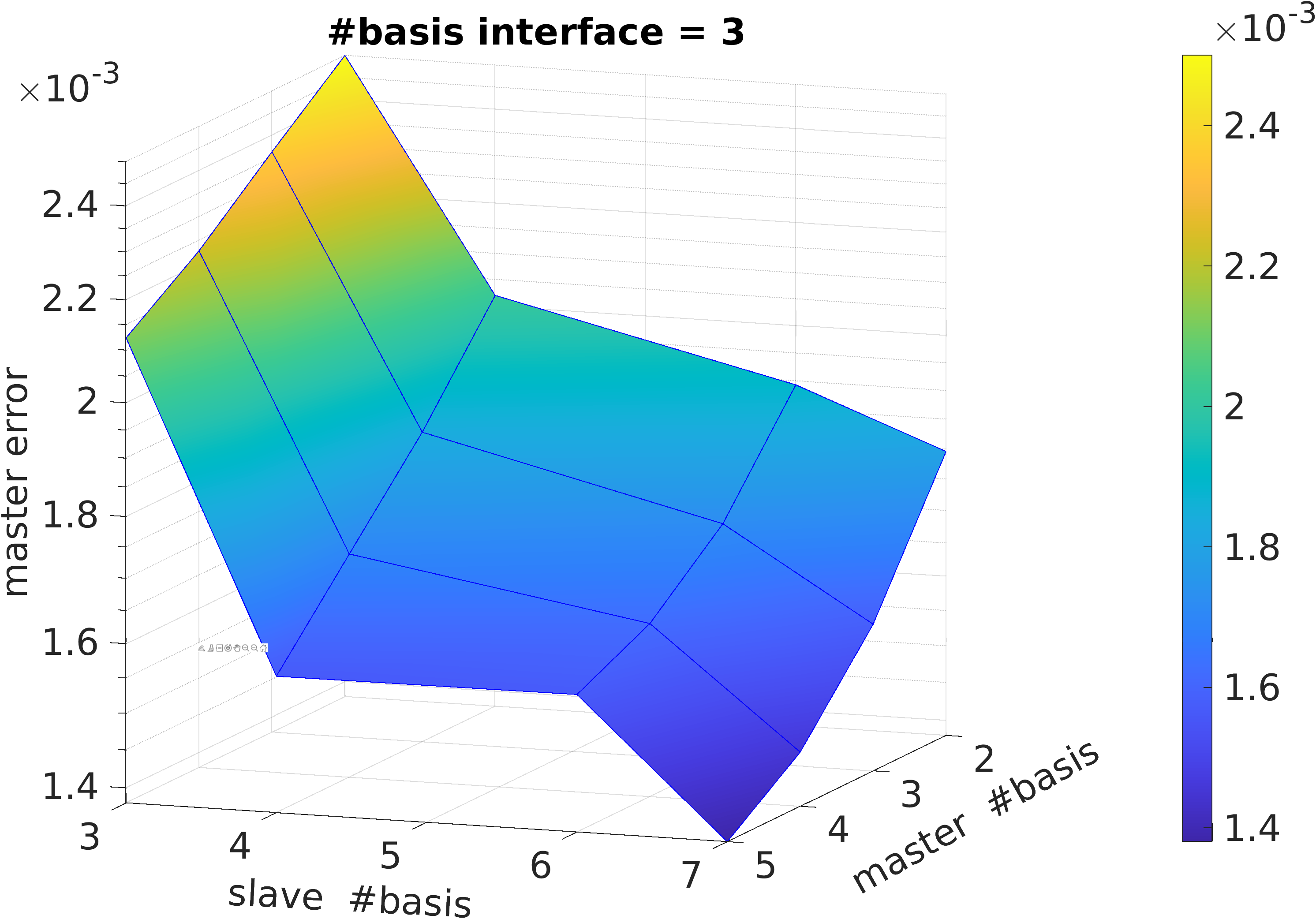}
	\includegraphics[width=0.24\textwidth]{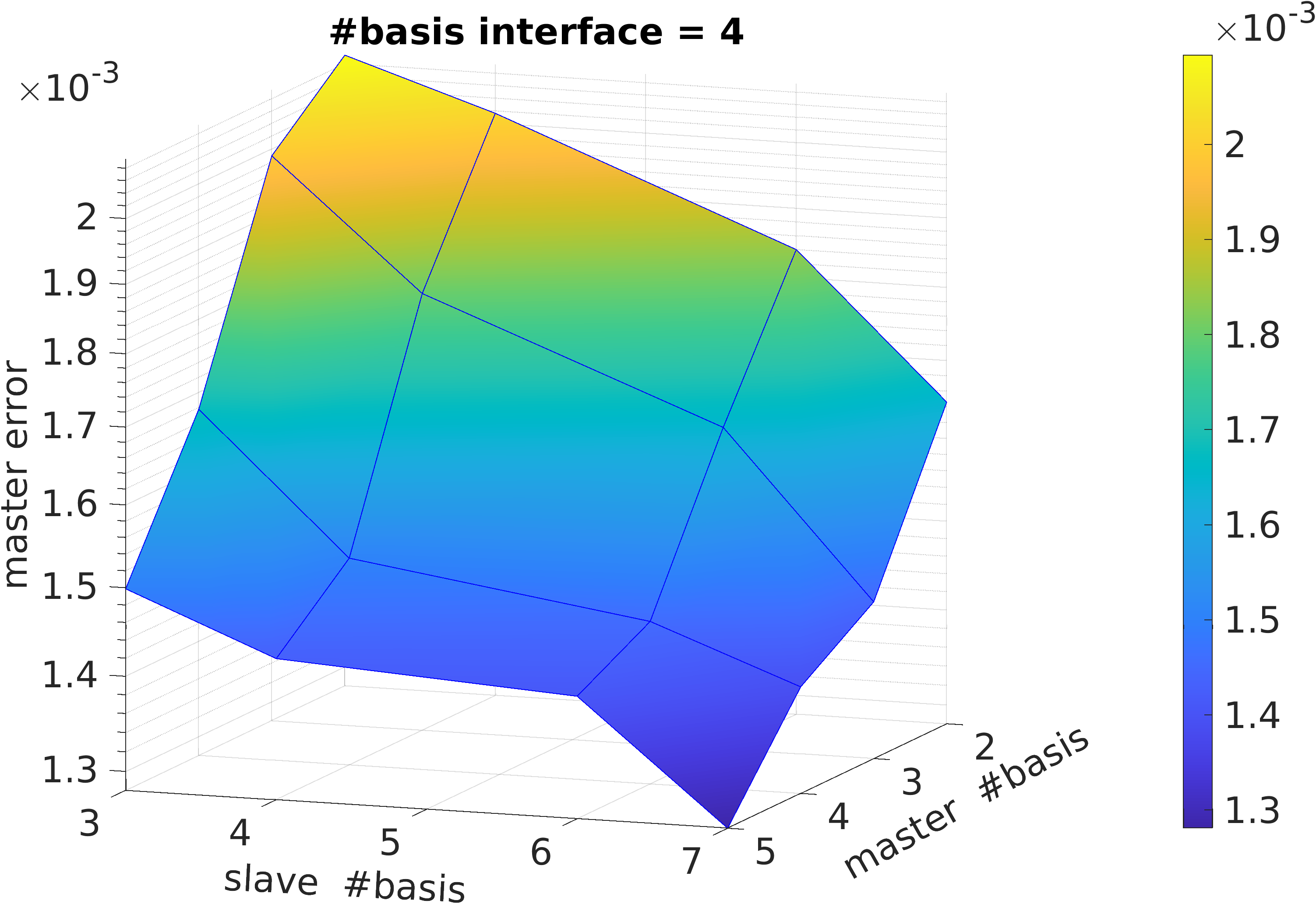}
	\includegraphics[width=0.24\textwidth]{master_error_interface_6.png}
	\caption{\emph{Test\#1.} $H^1(\Omega_i)$ mean relative error ($z$--axis) over the slave (top row) and master (bottom row) solution for $N_\text{test}=20$ different instances of the parameters between the FOM and ROM solutions varying the number of basis functions used to represent the slave and the master solution $n_1$ and $n_2$ ($x$-- and $y$--axis), and fixing the number of basis functions employed to approximate interface data $M_1$ and $M_2$.}
	\label{fig:fixed_interface_basis}
\end{figure}

\begin{figure}[t!]
	\centering
	\includegraphics[width=0.24\textwidth]{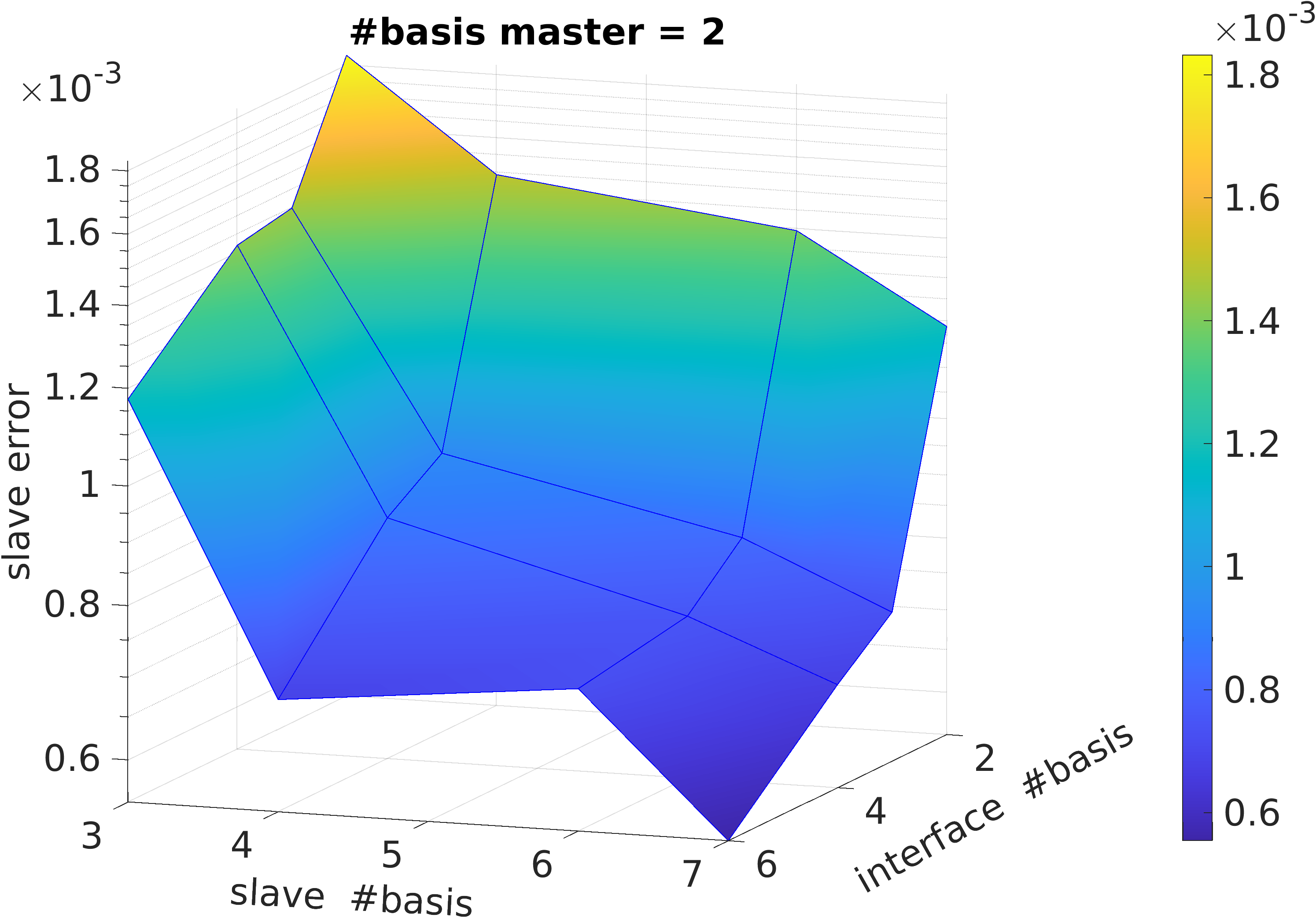}
	\includegraphics[width=0.24\textwidth]{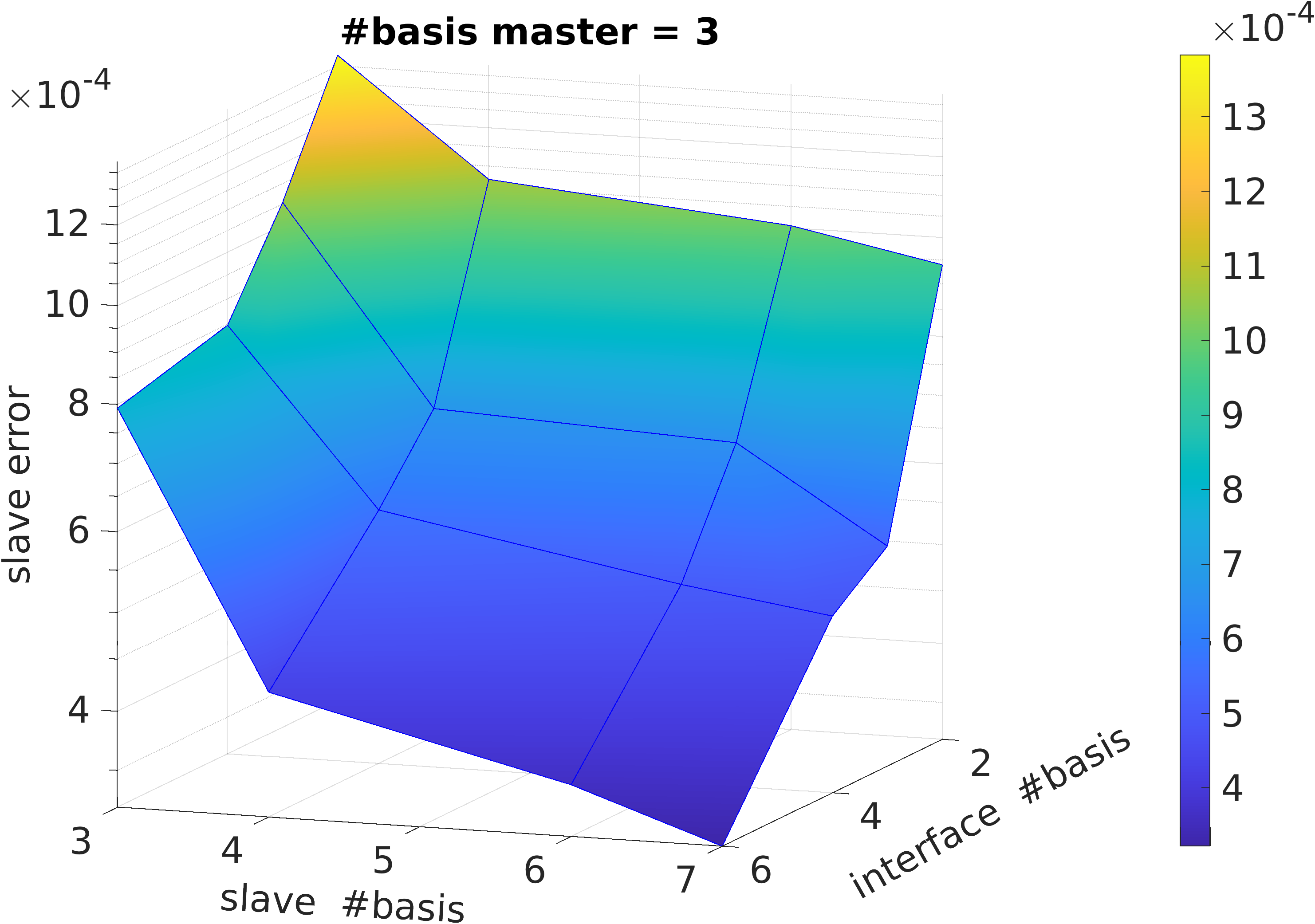}
	\includegraphics[width=0.24\textwidth]{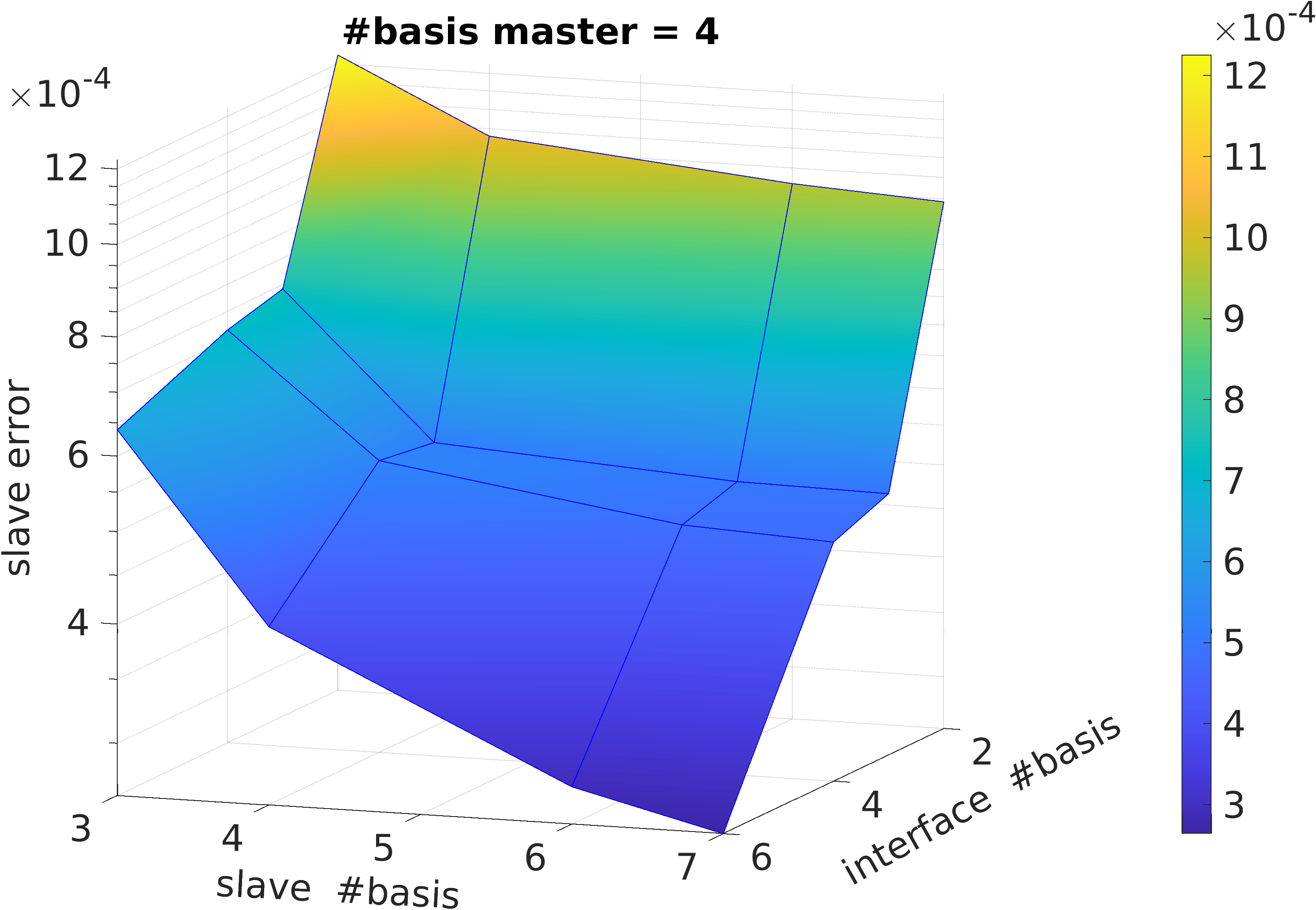}
	\includegraphics[width=0.24\textwidth]{slave_error_master_5.png}\\
    \bigskip
     \includegraphics[width=0.24\textwidth]{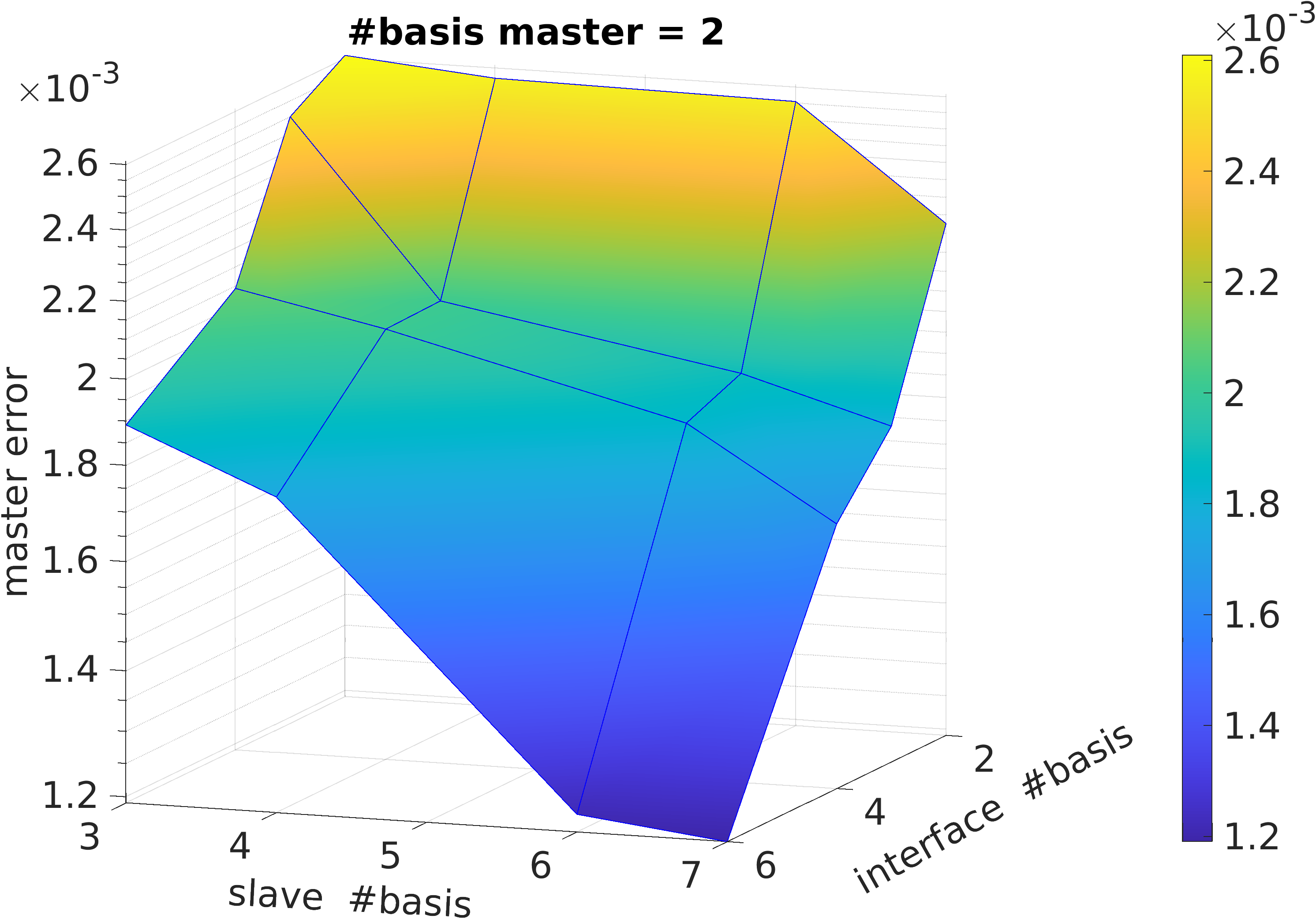}
	\includegraphics[width=0.24\textwidth]{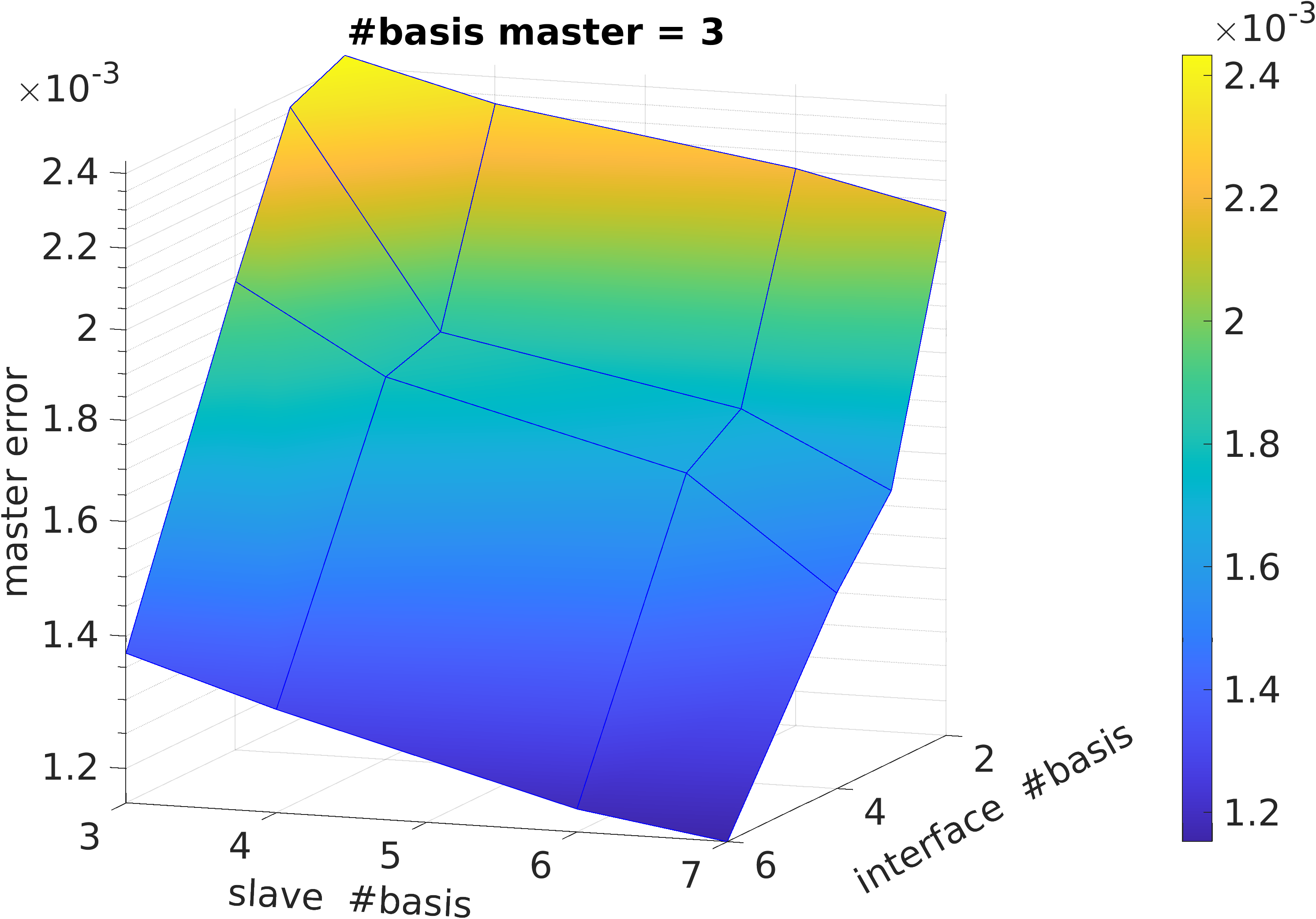}
	\includegraphics[width=0.24\textwidth]{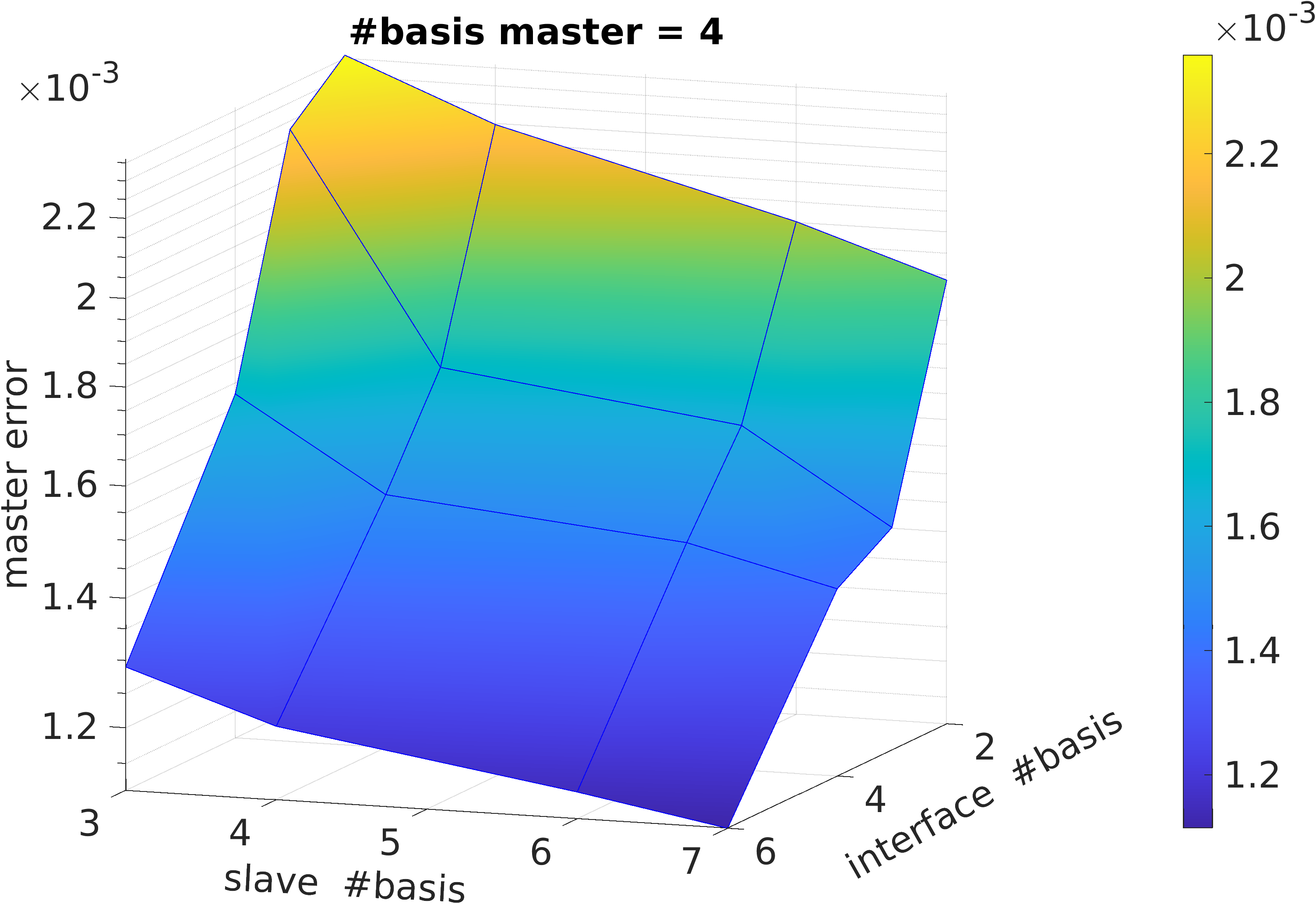}
	\includegraphics[width=0.24\textwidth]{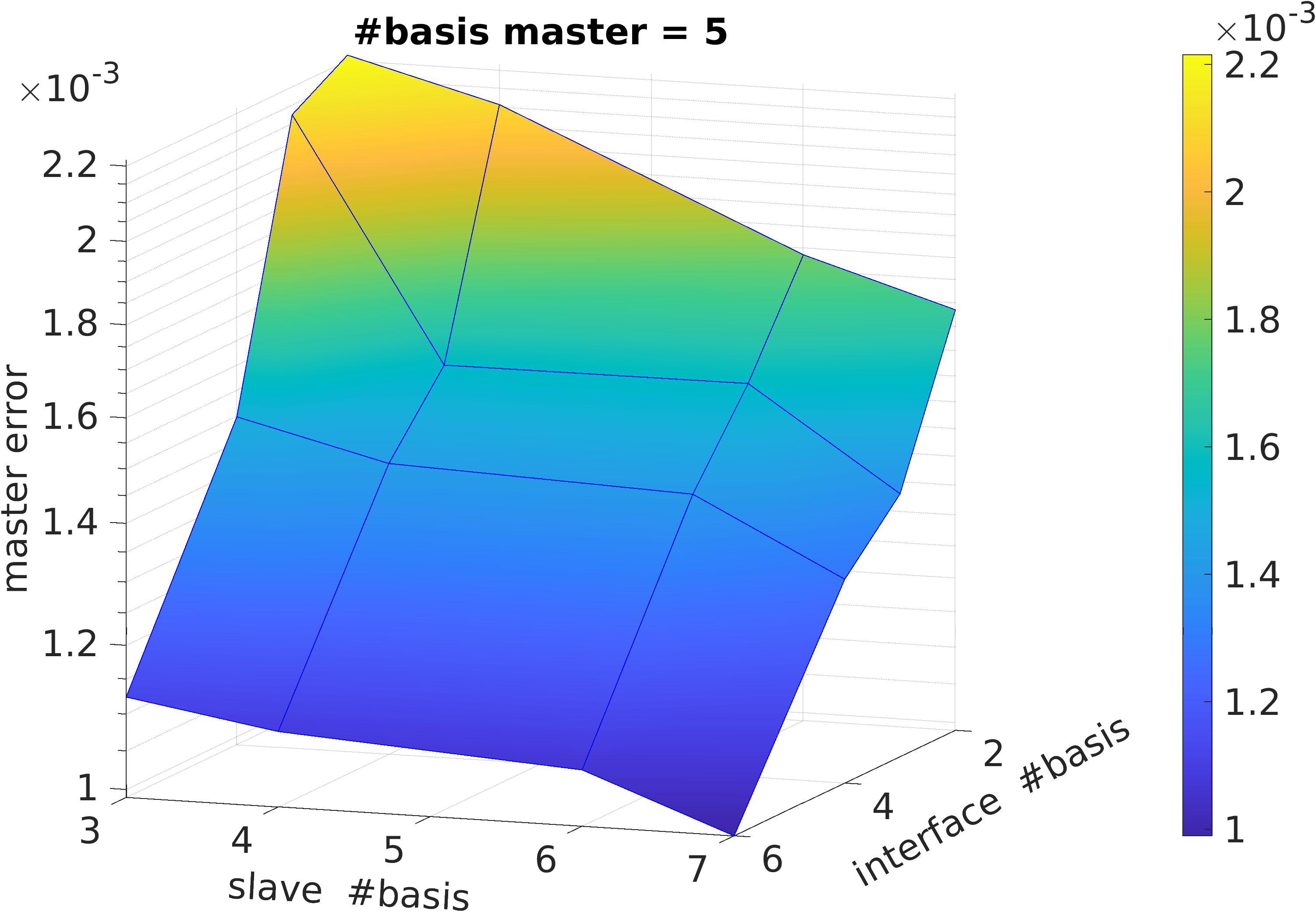}
	\caption{\emph{Test\#1.} $H^1(\Omega_i)$ mean relative error ($z$--axis) over the slave (top row) and master (bottom row) solution for $N_\text{test}=20$ different instances of the parameters between the FOM and ROM solutions varying the number of basis functions used to represent the slave solution $n_1$ and the interface data $M_1$ and $M_2$  ($x$-- and $y$--axis), and fixing the number of basis functions employed to approximate the master solution.}
	\label{fig:fixed_master_basis}
\end{figure}

Relative approximation errors employing the $H^1(\Omega_i)$ norm, $i = 1,2$, for both slave and master solutions are computed according to \eqref{eq:error_2norm}. Figs. \ref{fig:fixed_slave_basis}--\ref{fig:fixed_master_basis} depict such errors as functions of the reduced order hyper--parameters, \emph{i.e.} of the variations of the number of basis functions chosen to approximate the slave solution, the master solution, and interface data. Both problem solutions are shown to depend on all the reduced quantities involved, and a decreasing approximation error can be achieved every time that the number of basis functions for one between the reduced slave solution, master solution, or interface data, is increased. The major influence on the approximation of both master and slave solution comes from the reduction operated on in the interface data: it can be observed, indeed, a higher reduction of the approximation error when $M_1$ and $M_2$ are increased, whereas the error decrease is much slower when either $n_1$ or $n_2$ are increased. A more expensive and careful approximation of the interface data is therefore to be preferred over an expensive approximation of the slave and master solutions.

\begin{figure}[t!]
	\centering
    \includegraphics[width=0.45\textwidth]{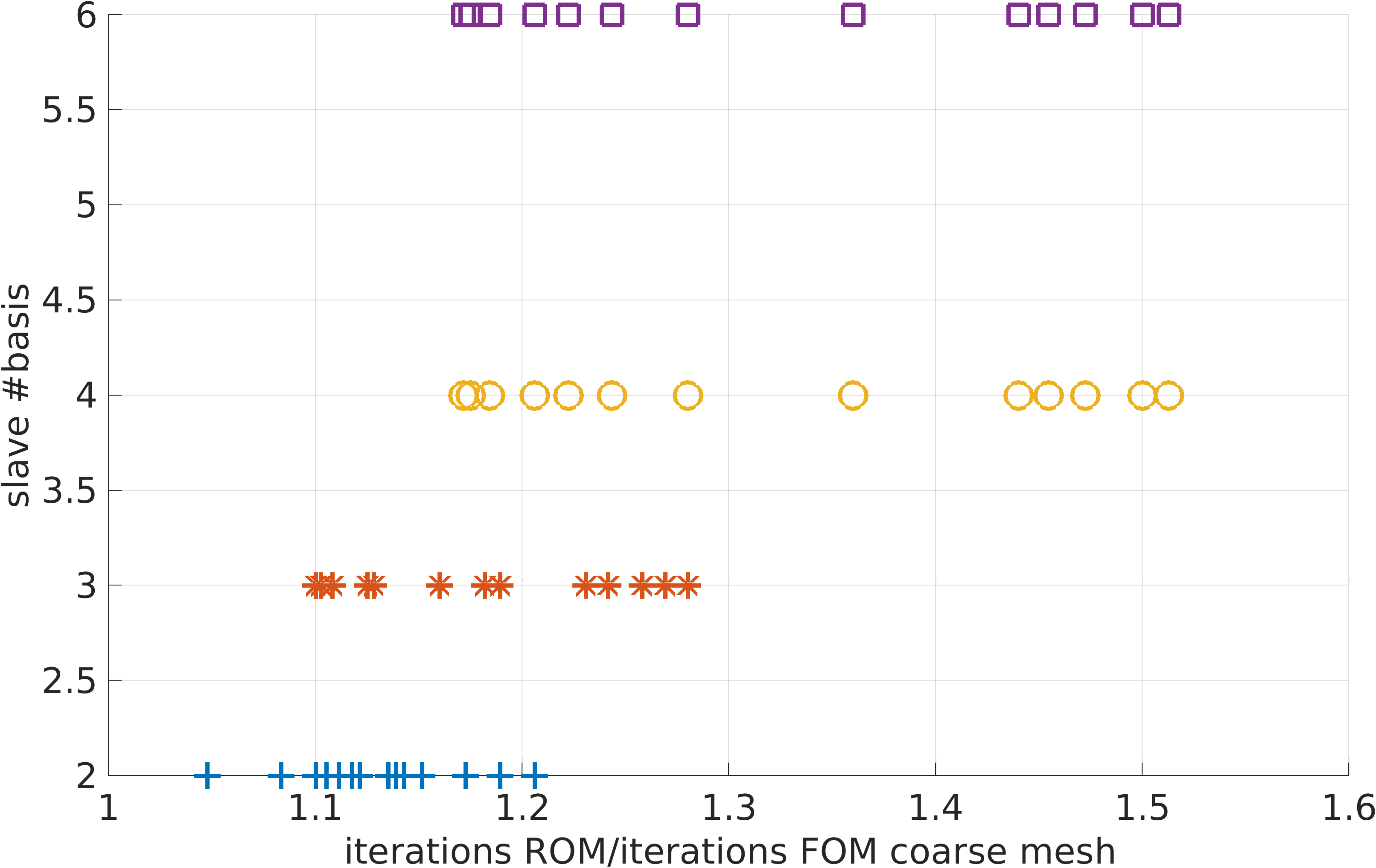}
    \quad
	\includegraphics[width=0.45\textwidth]{Iteration_ration_vs_basis_slave_fine.png}\\
    \bigskip
    \includegraphics[width=0.45\textwidth]{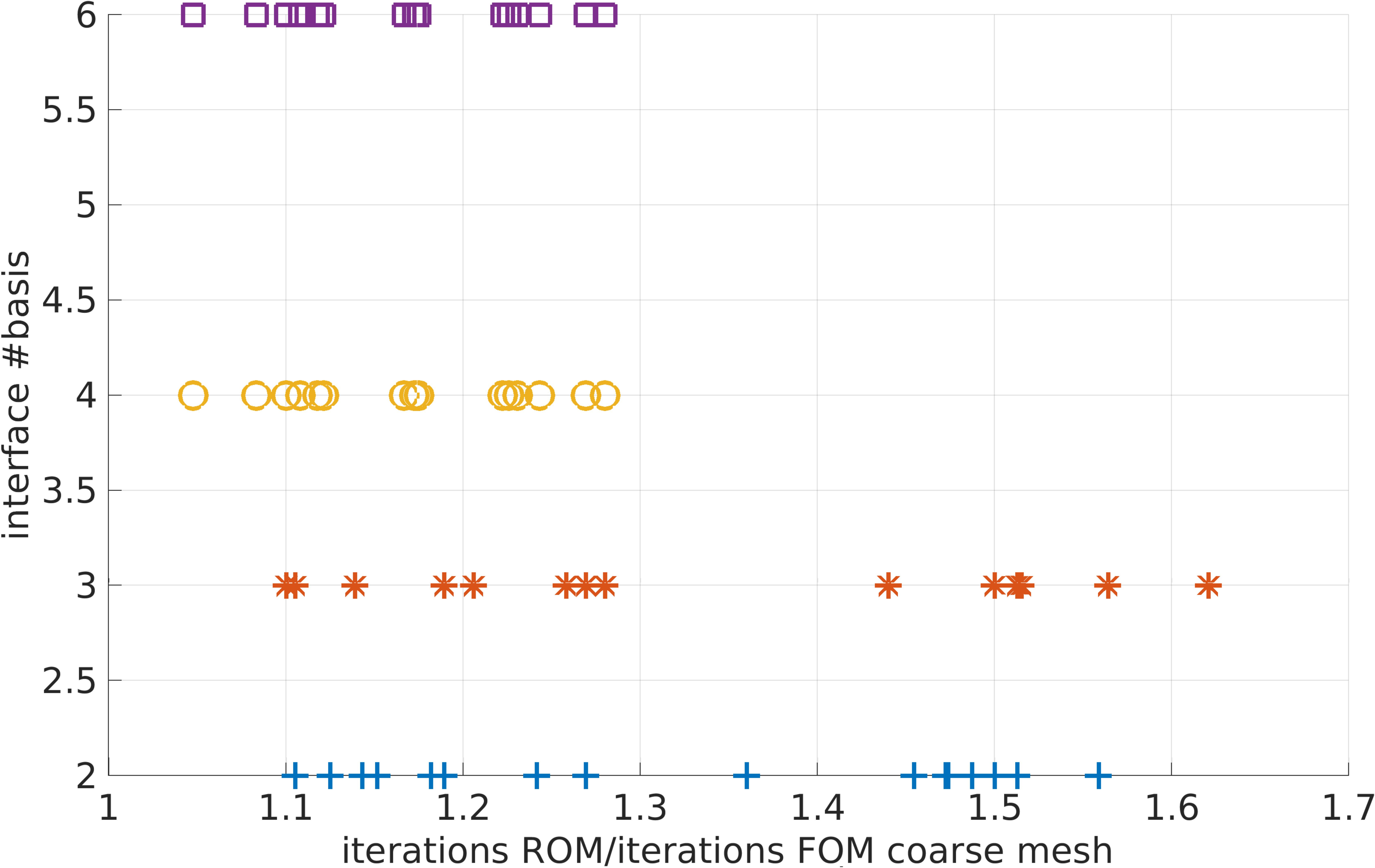}
    \quad
	\includegraphics[width=0.45\textwidth]{Iteration_ration_vs_basis_interface_fine.png}
	\\
	\bigskip
	\includegraphics[width=0.45\textwidth]{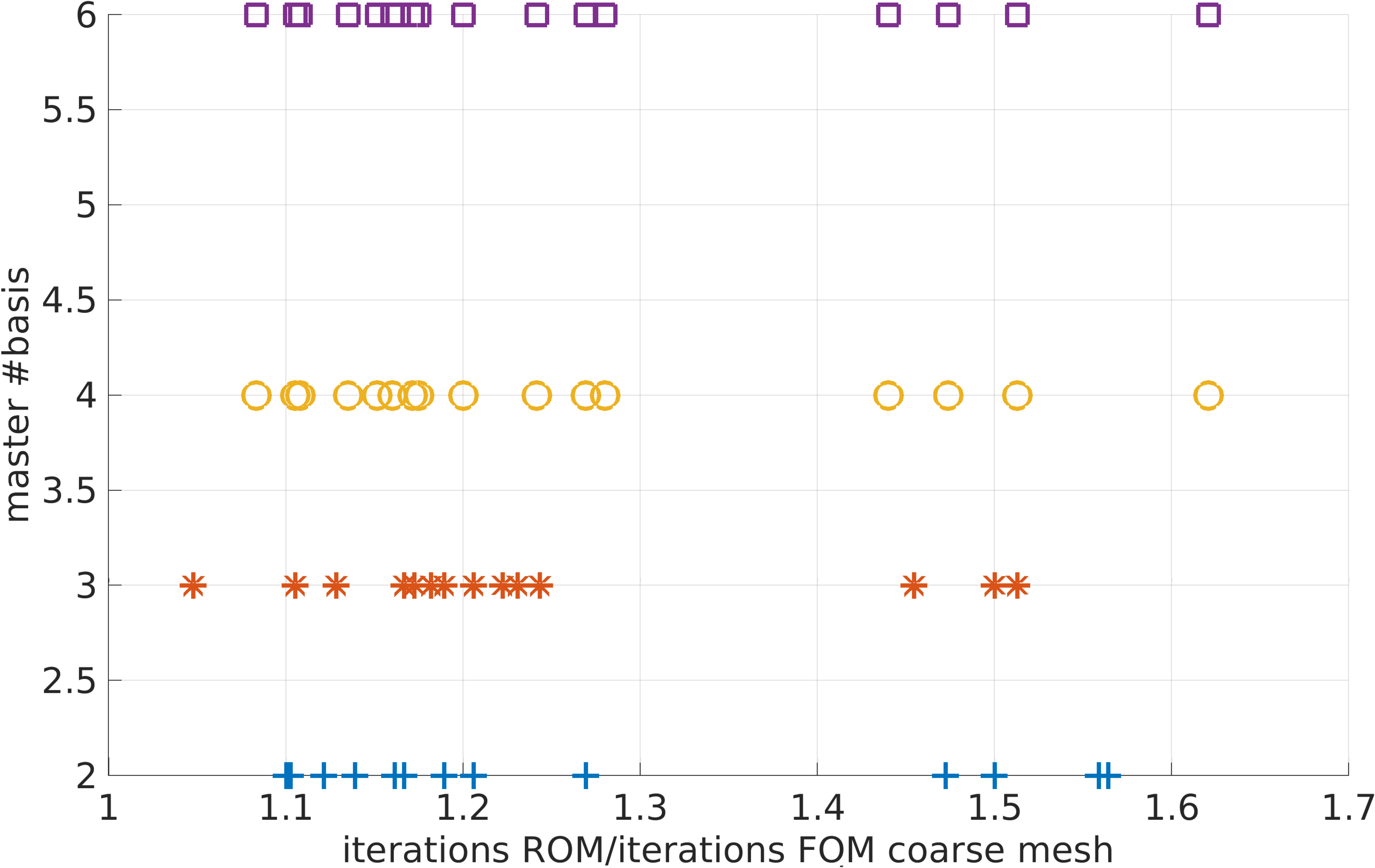}
    \quad
	\includegraphics[width=0.45\textwidth]{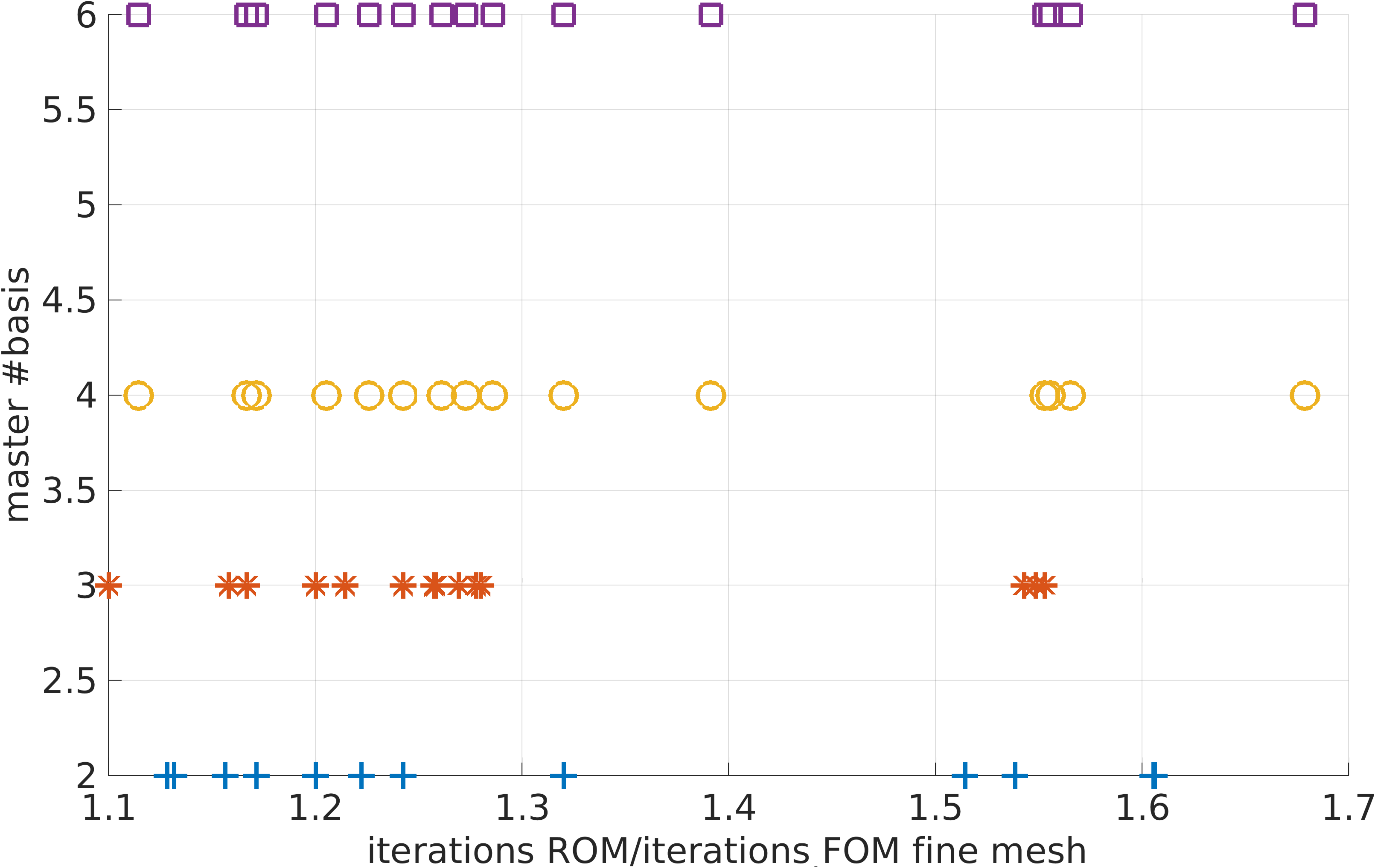}
	\caption{\emph{Test\#1.} Ratio between the number of iterations obtained with ROM and FOM schemes versus the number of basis functions employed to approximate the slave solutions (first row), the interface data (second row), or the master solution (third row), either employing the coarse discretization (left) or the fine discretization (right) for the FOM computation.}
	\label{fig:Laplace_ratio_vs_basis_full}
\end{figure}

\begin{figure}[t!]
	\centering
    \includegraphics[width=0.45\textwidth]{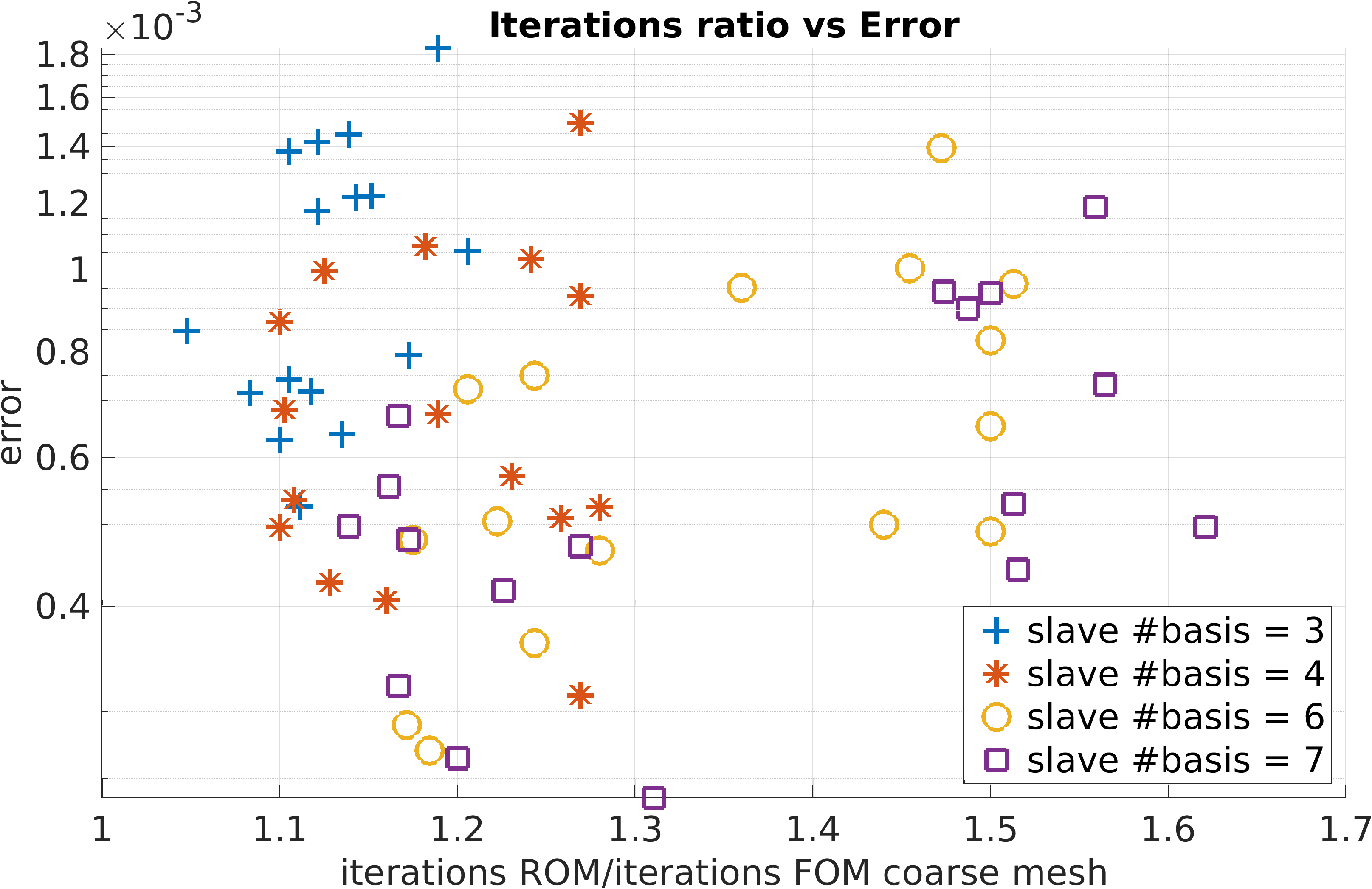}
    \quad
	\includegraphics[width=0.45\textwidth]{Iteration_ration_vs_error_slave_fine.png}\\
    \bigskip
    \includegraphics[width=0.45\textwidth]{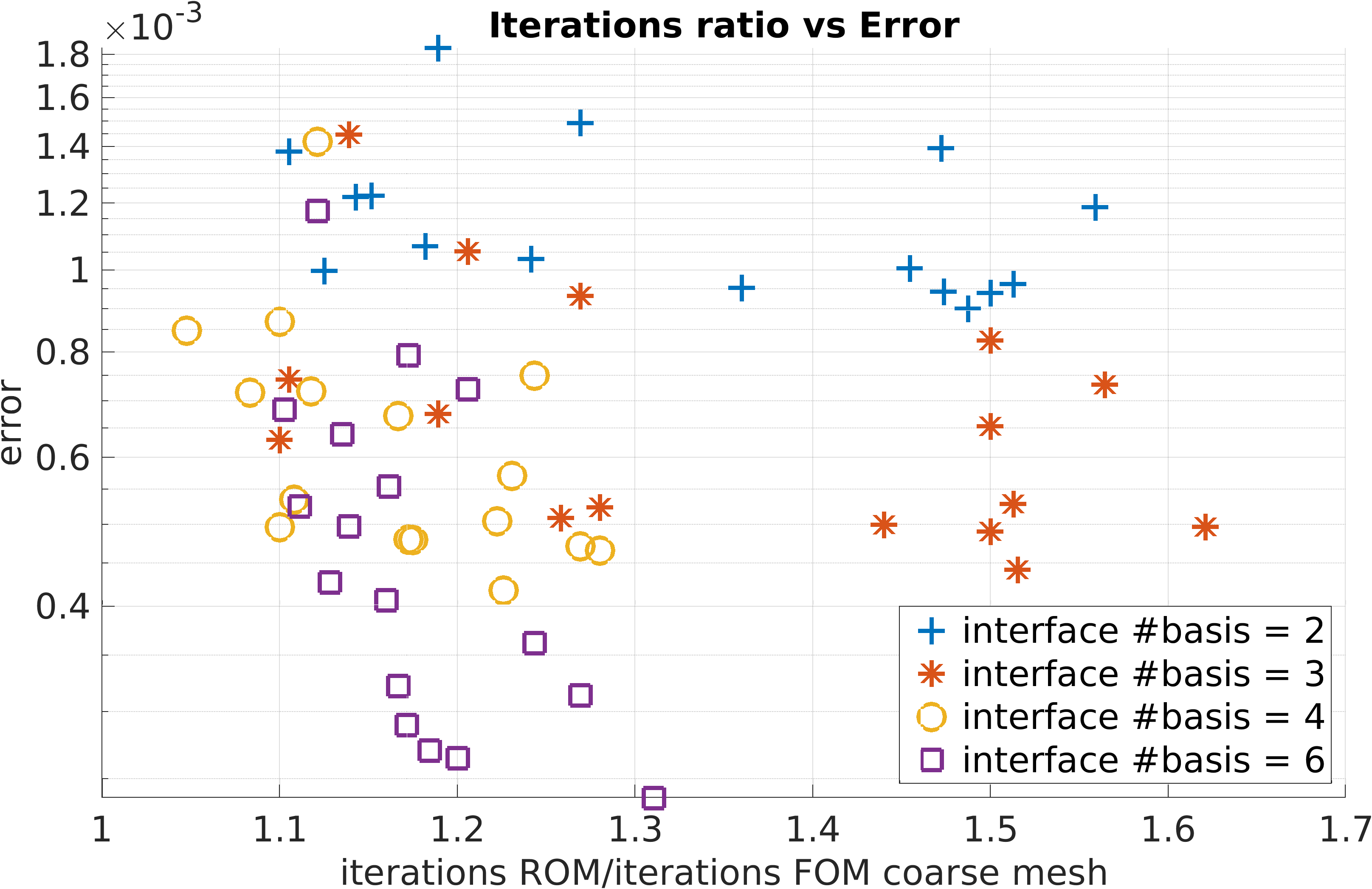}
    \quad
	\includegraphics[width=0.45\textwidth]{Iteration_ration_vs_error_interface_fine.png}
	\\
	\bigskip
	\includegraphics[width=0.45\textwidth]{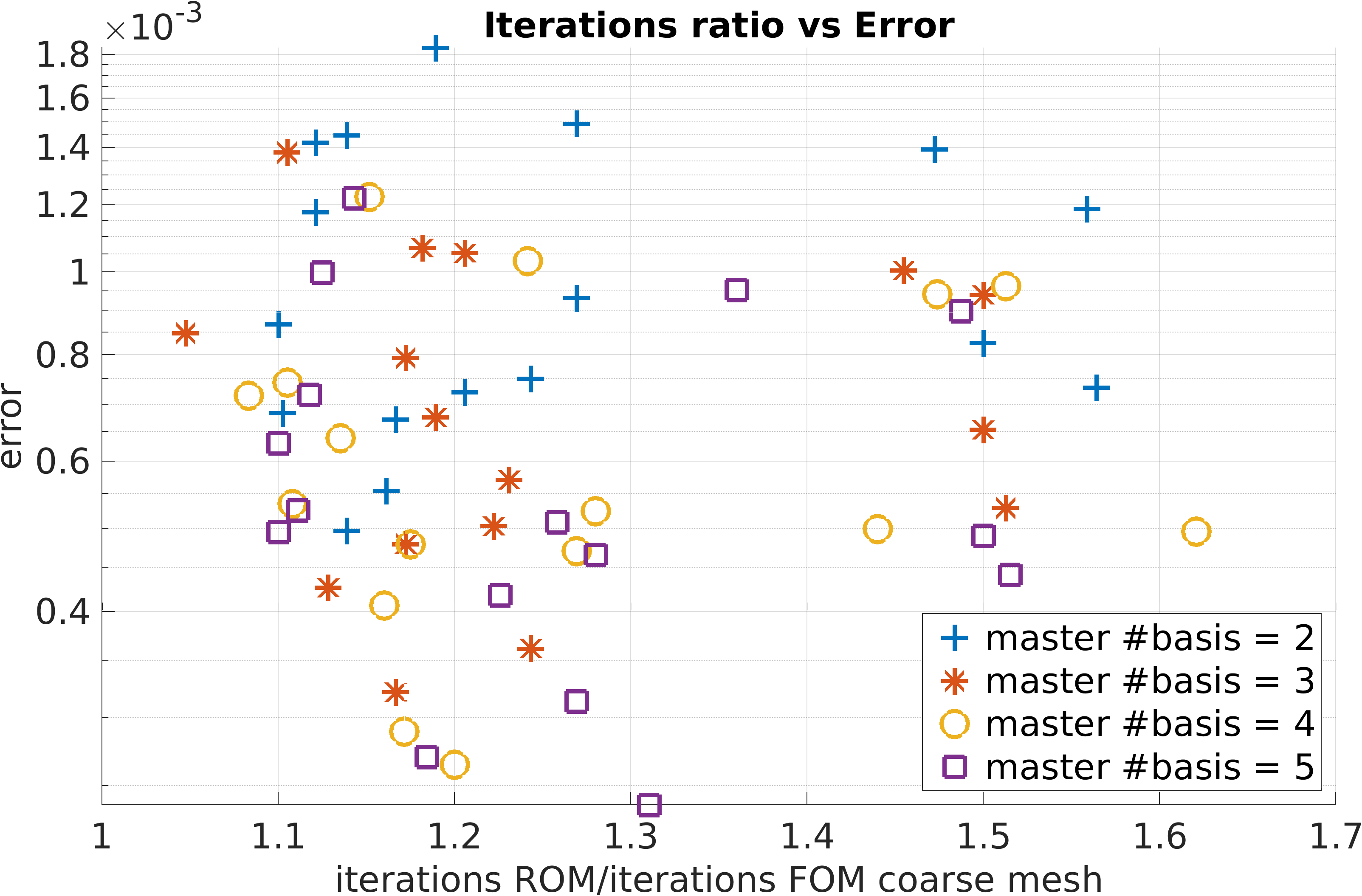}
    \quad
	\includegraphics[width=0.45\textwidth]{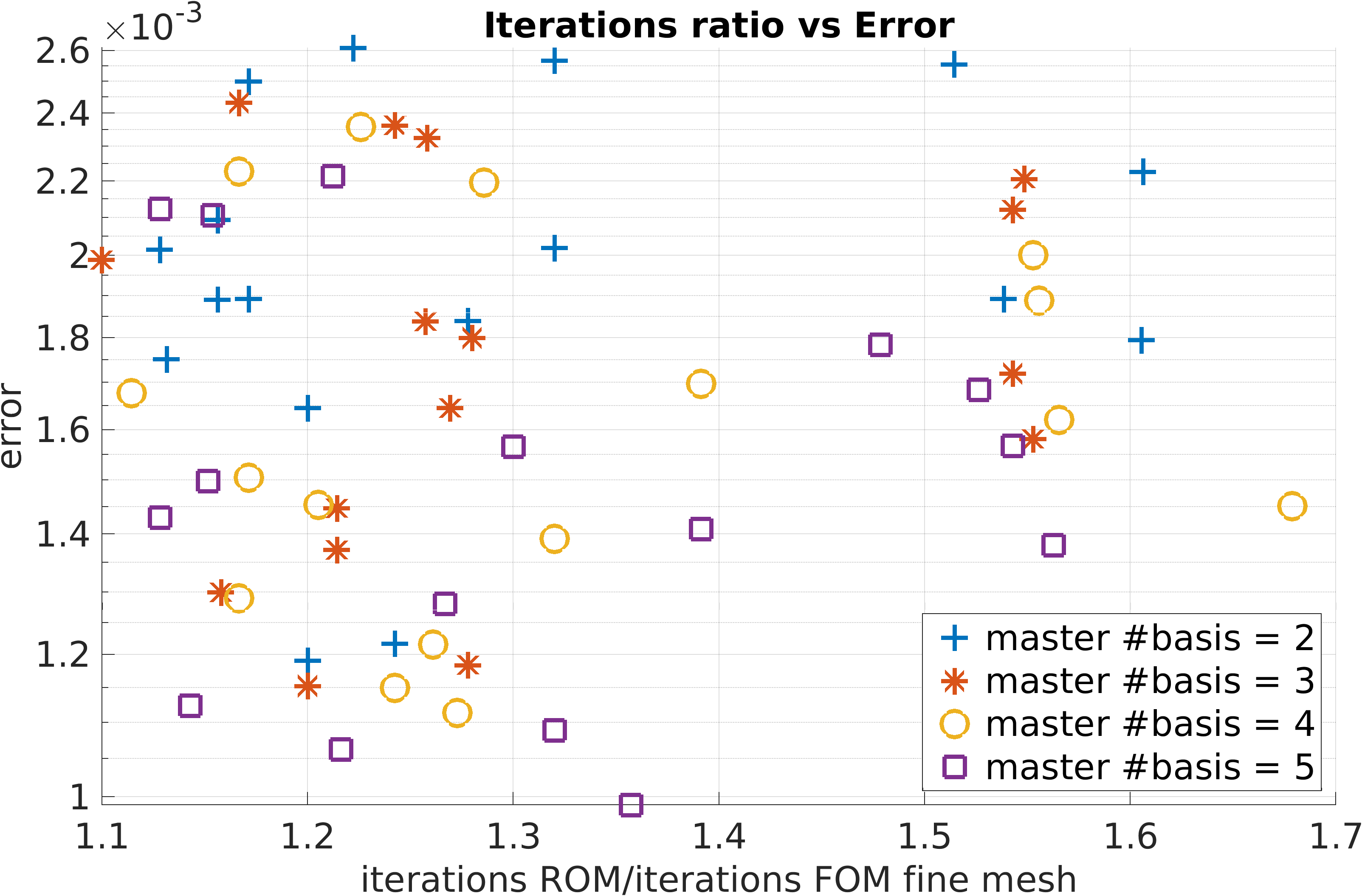}
	\caption{\emph{Test\#1.} Iterations ratio versus approximation error
 depending on the number of basis functions employed to approximate the slave solutions (first row), the interface data (second row), or the master solution (third row), either employing the coarse discretization (left) or the fine discretization (right) for the FOM computation.}
	\label{fig:Laplace_ratio_vs_error_full}
\end{figure}

The computational costs of the ROM are investigated by means of the number of iterations needed by the scheme to reach the solution interface convergence, as well as by the effective CPU time. Since both quantities depend on parameter instance, we compute the ratio between the iterations number of the ROM and FOM computations, as well as the ratio between the FOM and ROM computational times. The last ratio is able to describe the speed--up achieved by employing the ROM. Fig. \ref{fig:Laplace_ratio_vs_basis_full} shows the variation of such iterations ratio depending on $n_1$, $n_2$, $M_1$ and $M_2$, comparing the results obtained with both the coarse and the fine discretization. The graphs show a dependency of the iterations number from the number of basis functions employed to compute the slave solution and interface data. Specifically, a greater approximation accuracy (and a number of basis functions) of the slave solution increases the number of iterations required to achieve the solution convergence, whereas a higher approximation accuracy of the interface data decreases the number of iterations. The same effect can be observed also in Fig. \ref{fig:Laplace_ratio_vs_error_full}, where we compare the iterations ratio with the approximation error. Varying the number of basis functions to compute the master solution has instead a minor effect on the number of iterations.

\begin{figure}[t!]
	\centering
    \includegraphics[width=0.45\textwidth]{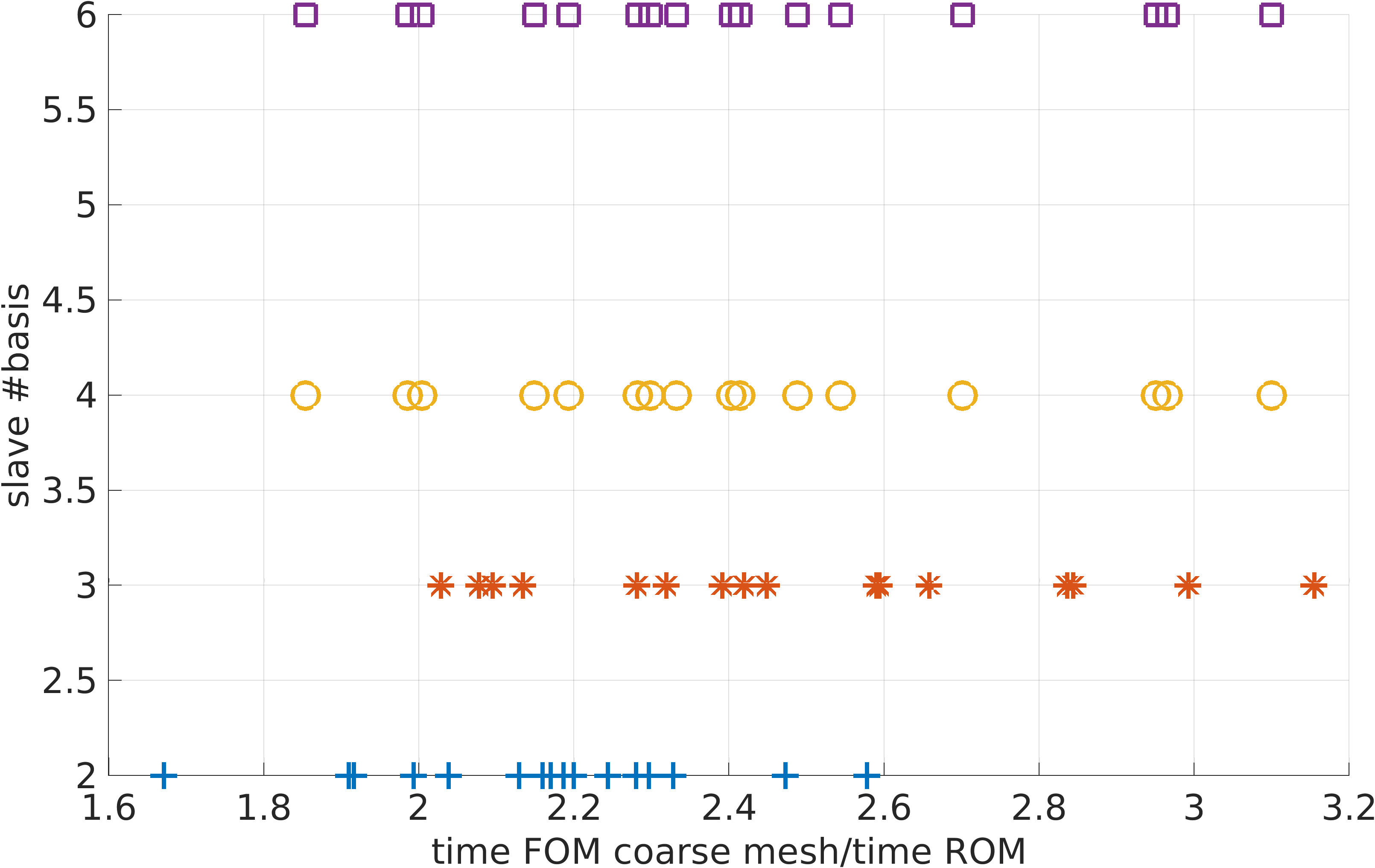}
    \quad
	\includegraphics[width=0.45\textwidth]{Time_ration_vs_basis_slave_fine.png}\\
    \bigskip
    \includegraphics[width=0.45\textwidth]{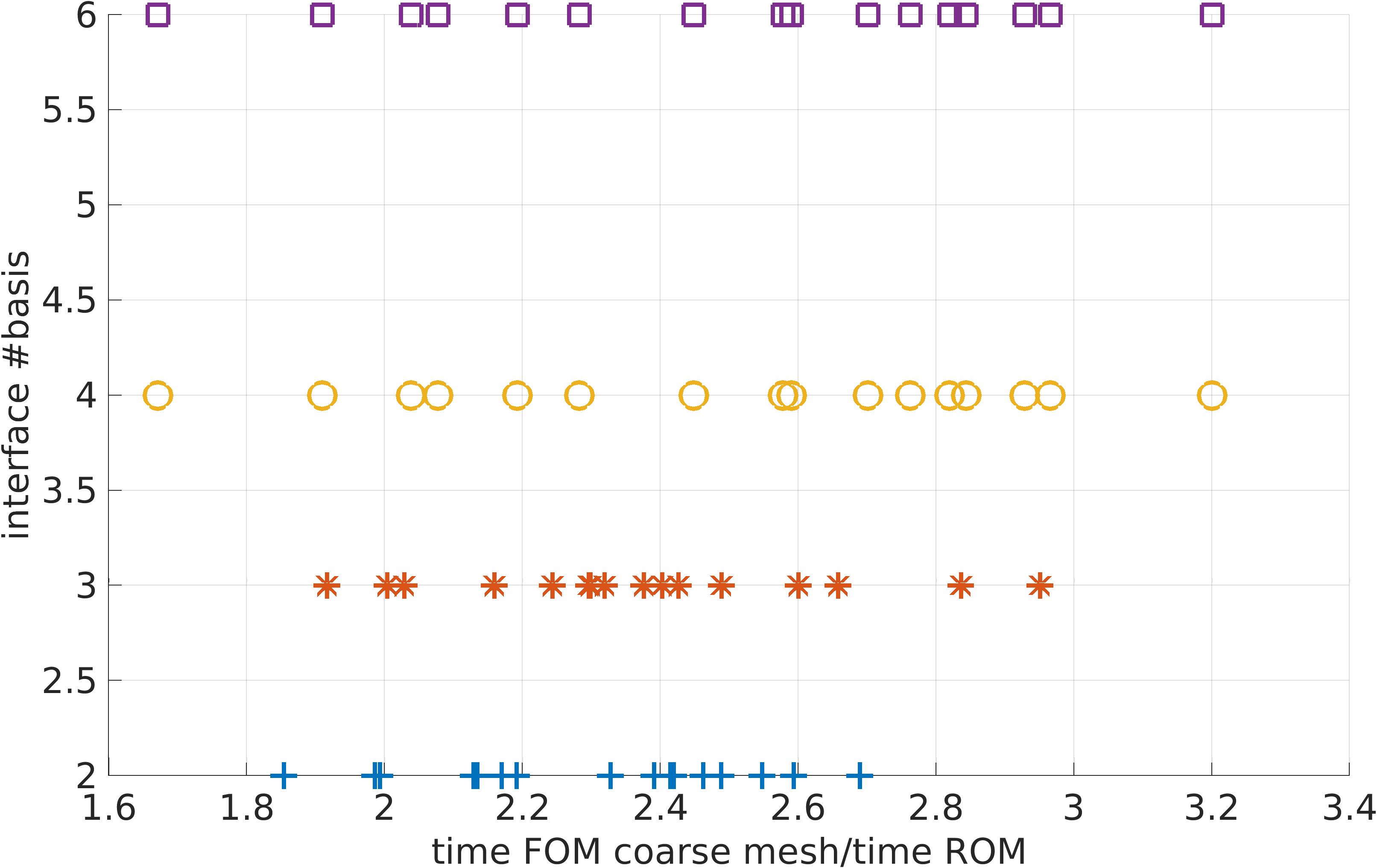}
    \quad
	\includegraphics[width=0.45\textwidth]{Time_ration_vs_basis_interface_fine.png}
	\\
	\bigskip
	\includegraphics[width=0.45\textwidth]{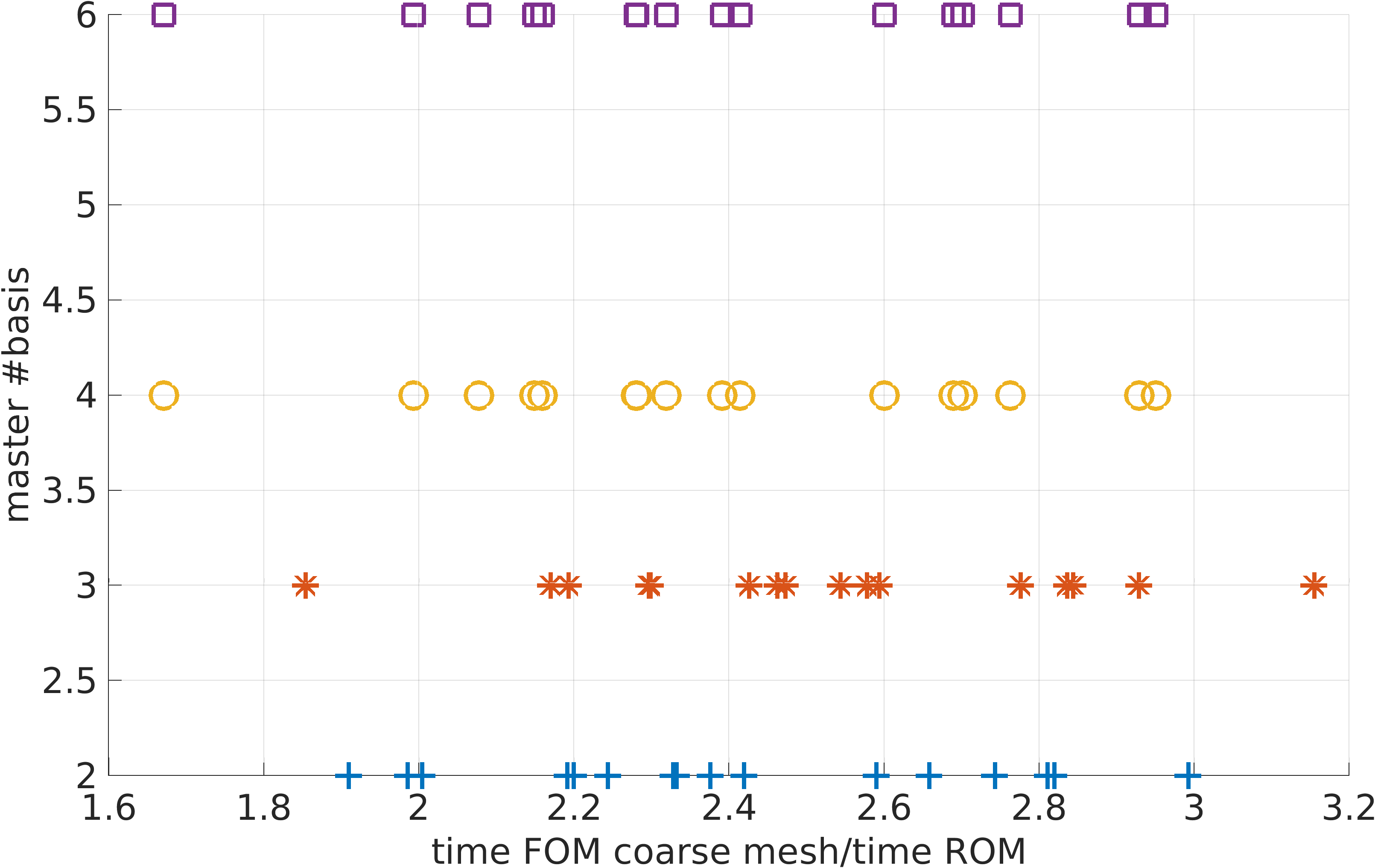}
    \quad
	\includegraphics[width=0.45\textwidth]{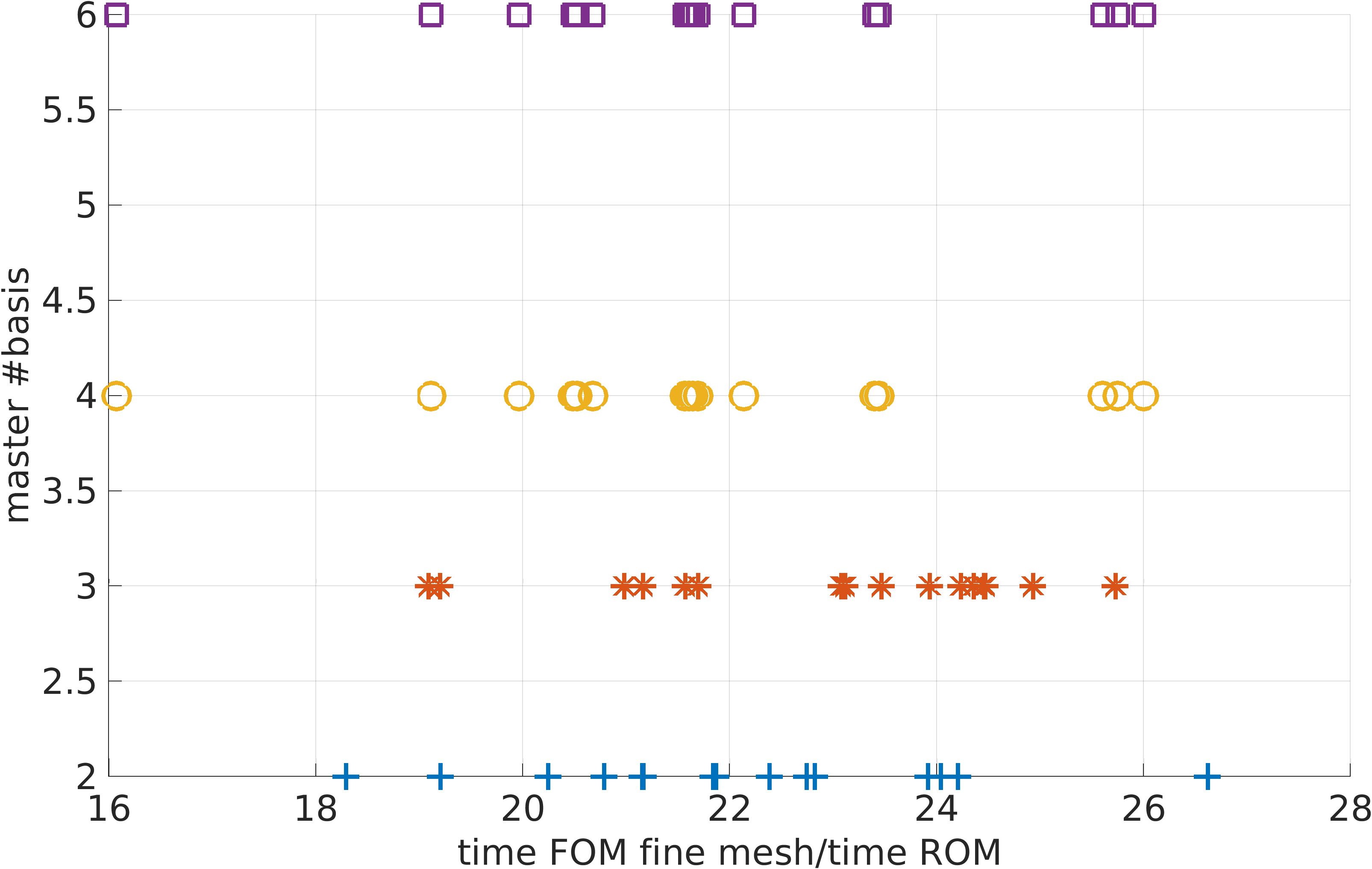}
	\caption{\emph{Test\#1.} Ratio between the CPU time of FOM and ROM versus the number of basis functions employed to approximate the slave solutions (first row), the interface data (second row), or the master solution (third row), either employing the coarse discretization (left) or the fine discretization (right) for the FOM computation.}
	\label{fig:Laplace_time_vs_basis_full}
\end{figure}

We report the CPU time ratio compared to the basis functions number in Fig. \ref{fig:Laplace_time_vs_basis_full}, for both coarse and fine FOM discretizations, showing that independently of the number of basis functions imposed for any reduced quantity, we achieve very similar computational speed up.

\section{Preliminary results employing Radial Basis Functions to interpolate data across the domain interfaces}
\label{App_B}
In this appendix, we present a very preliminary investigation of the effects of utilizing a more advanced interpolation technique compared to the nearest neighbor approach outlined and utilized throughout the rest of the paper, to interpolate in the ROM algorithm the Dirichlet--Neumann data across the non-conforming interfaces. Specifically, in accordance with Remark \ref{Rem:RBF_formula} of Section \ref{sec:interface_reduction}, we have implemented Radial Basis Function (RBF) interpolation to interpolate the model solution between $\Gamma_2$ and $\Gamma_1$ (i.e., the Dirichlet interface data), as well as the residual vector between $\Gamma_1$ and $\Gamma_2$ (i.e., the Neumann interface data) at the subdomain interfaces. We subsequently apply this RBF scheme to \emph{Test\#2} in Subsection \ref{Subsect:steady_case_source}, and a partial application to the test case \emph{Test\#3} in Subsection \ref{sub:heat_problem}.

\subsection{Test\#2 - steady case: diffusion reaction equation with parametrized source}
\label{App_sub:test_case_2}
In this subsection we compute a comprehensive set of approximated ROM solutions of problem \eqref{Eq:test_case_1} with source term \eqref{Eq:source_term_case_3}, by employing varying numbers of basis functions for the master, slave, and interface data, and employing the RBF method to interpolate reduced order Dirichlet and Neumann data across the domain interface. This approach allows us to achieve different levels of approximation accuracy for each quantity depending on the parameters choice, reflecting the methodology outlined in Section \ref{sec:numerical_results}. 

For the sake of fairness, we employing the same parameters -- encompassing both physical and reduced basis parameters (such as $N_{\text{train}}$ and $N_{\text{test}}$  dimensions) - of \emph{Test\#2}, presenting the same error analysis procedure as depicted in Fig. \ref{fig:Laplace_error_source}.  

\begin{figure}[!h]
	\centering
	\includegraphics[width=0.4\textwidth]{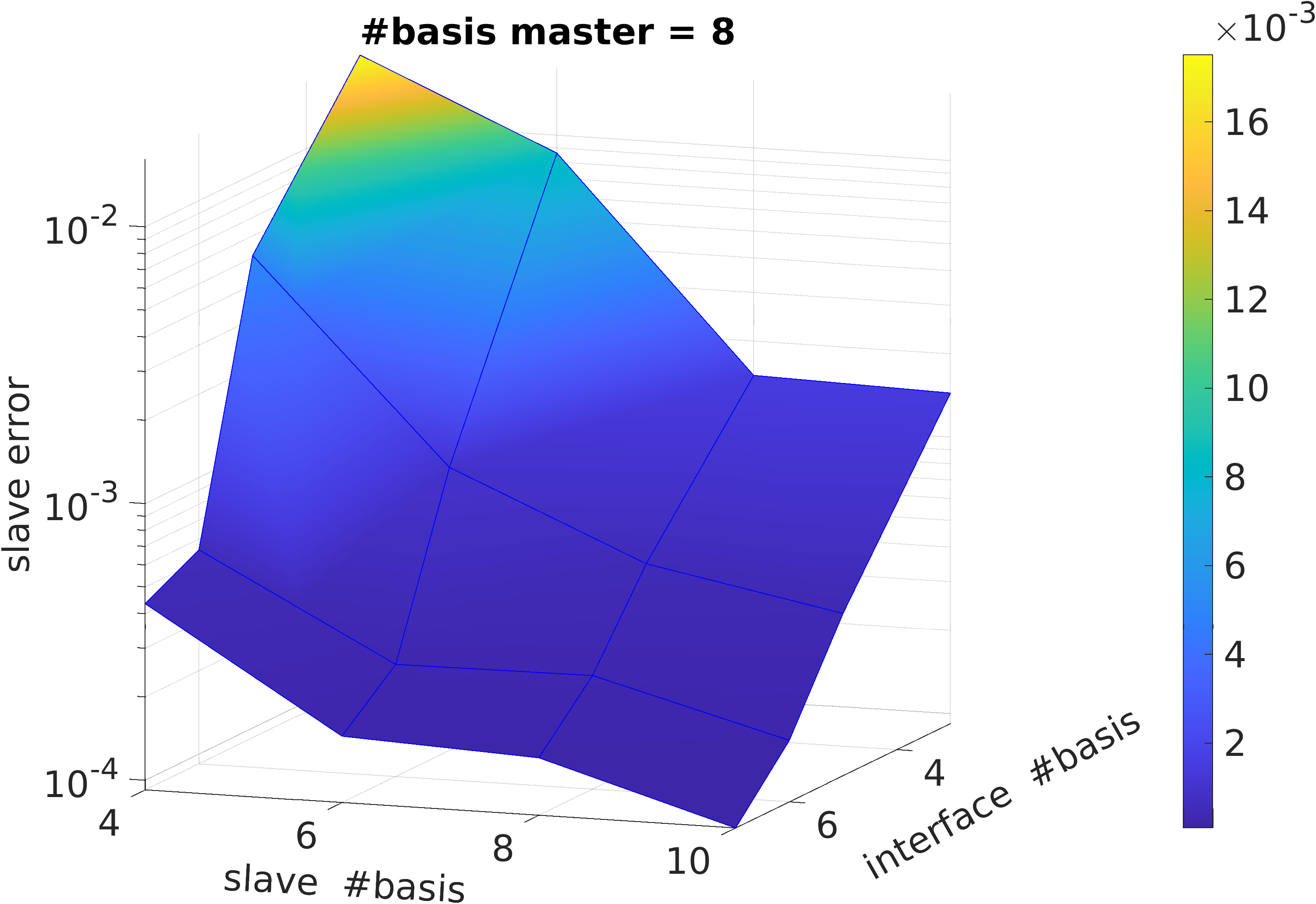} \quad
	\includegraphics[width=0.4\textwidth]{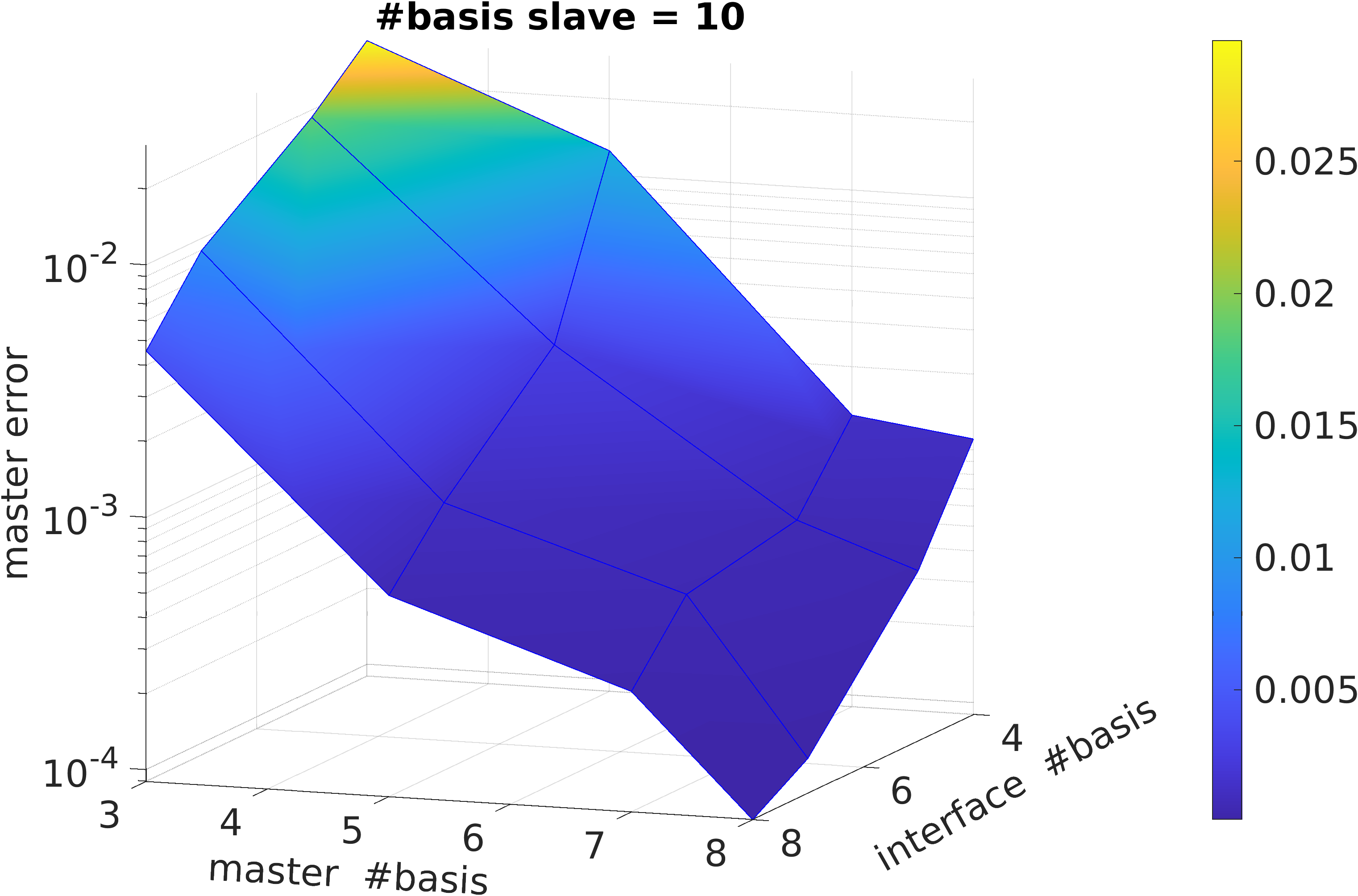}\\
	\includegraphics[width=0.4\textwidth]{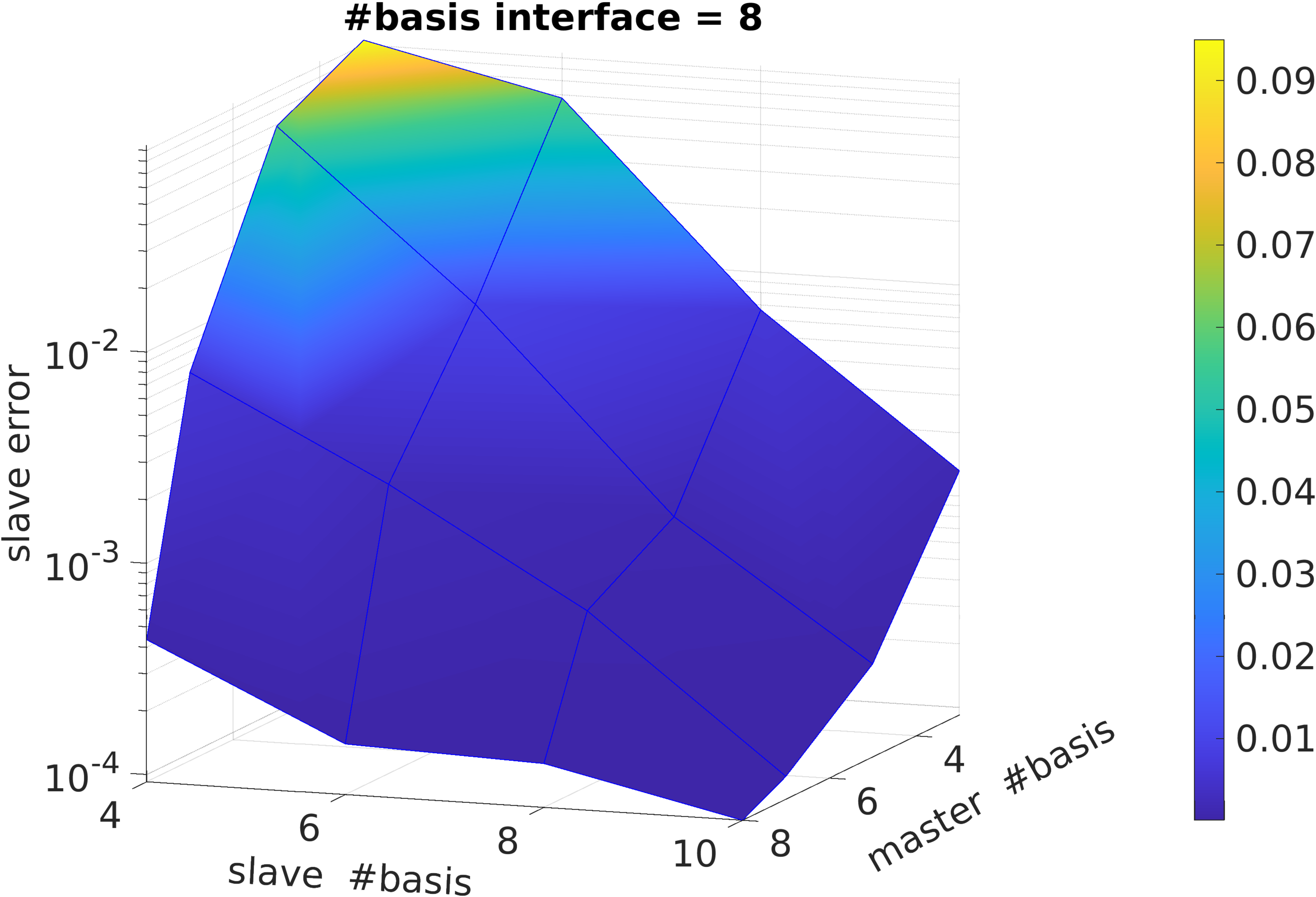} \quad
	\includegraphics[width=0.4\textwidth]{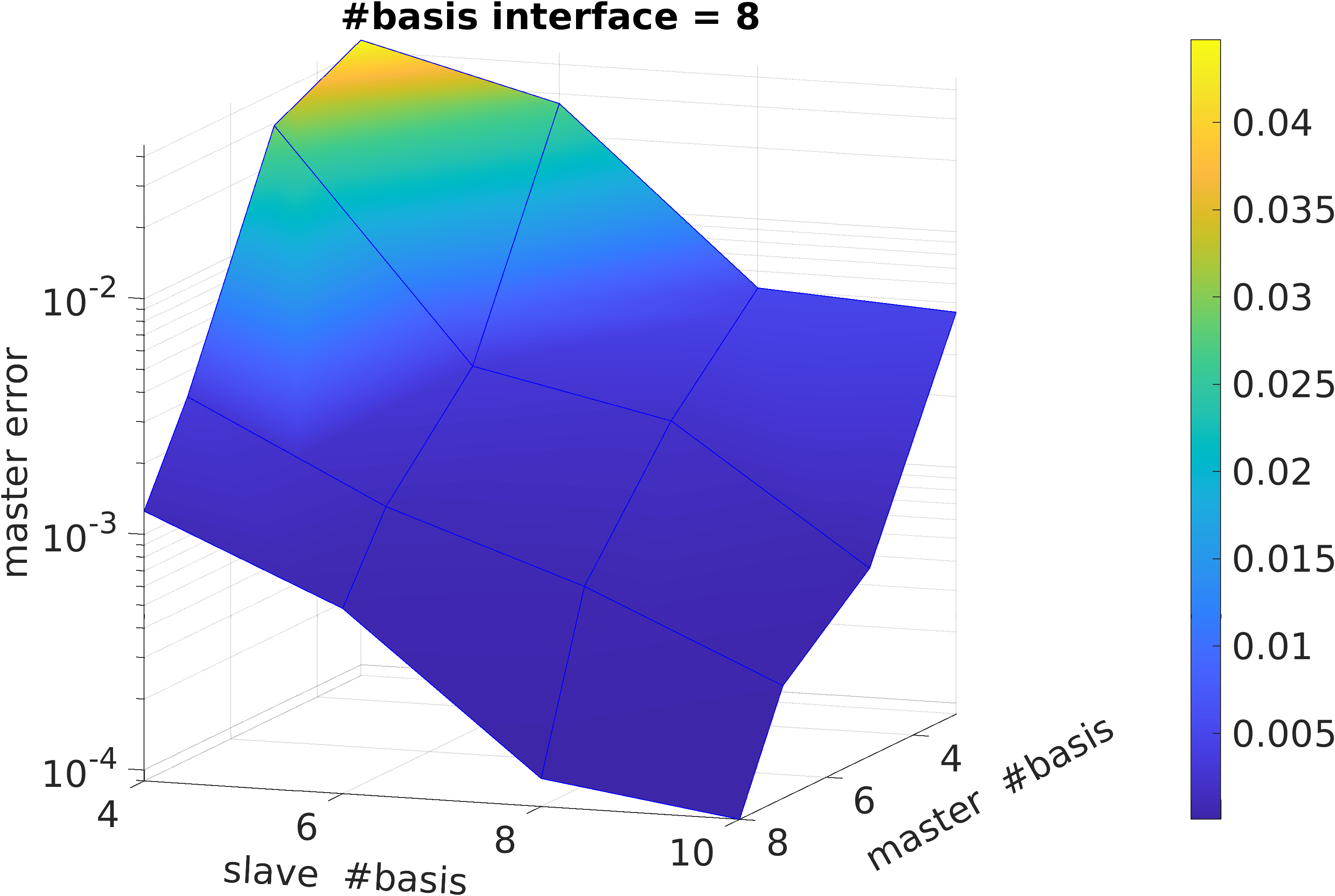}
	\caption{\emph{Test\#2.} $H_1(\Omega_i)$ mean relative error ($z$--axis) over the solution for $N_{test}=20$ different instances of the parameters between the FOM and ROM solutions varying the number of basis functions used to represent the slave and the master solution $n_1$ and $n_2$, and the interface data $M_1$ and $M_2$($x$-- and $y$--axis). On the top row, we fix the number of basis functions of the master problem to 8 (on the left) and to 10 for the slave problem (on the right), while on the bottom we fix the number of basis functions equal to 8 for the interface data representation.}
	\label{fig:error_rbf_steady}
\end{figure}

In Fig. \ref{fig:error_rbf_steady} we therefore report the approximation errors obtained by varying the basis functions for two quantities between the  master solution, slave solution, or interface data reduction, while keeping the basis functions for the third quantity constant. Specifically, the number of basis functions held constant in each graph are those required by the POD/DEIM algorithm to attain an approximation accuracy for the respective quantity of $10^{-5}$. These graphs aim at providing a comprehensive overview of the algorithm's behavior.

We observe that although the overall reduction in computational error appears to be more stable and consistent compared to the test presented in Subsection \ref{Subsect:steady_case_source} (which utilized nearest neighbor interpolation), the overall values of the approximation errors are of a similar magnitude to those discussed in the paper. Furthermore, in terms of the number of Dirichlet--Neumann iterations required for the algorithms to reach convergence, we have not observed any differences compared to the case of nearest neighbor interpolation. We believe this can be largely due to the number of FOM Dirichlet--Neumann iterations, and not to the interpolation method employed.

\subsection{Test\#3 - unsteady case: time--dependent heat equation}
In this subsection we instead investigate the effect of the RBF interpolation on the ROM algorithm on test case \#3 presented in Subsection \ref{sub:heat_problem}. 

\begin{figure}[!h]
	\centering
	\includegraphics[width=0.4\textwidth]{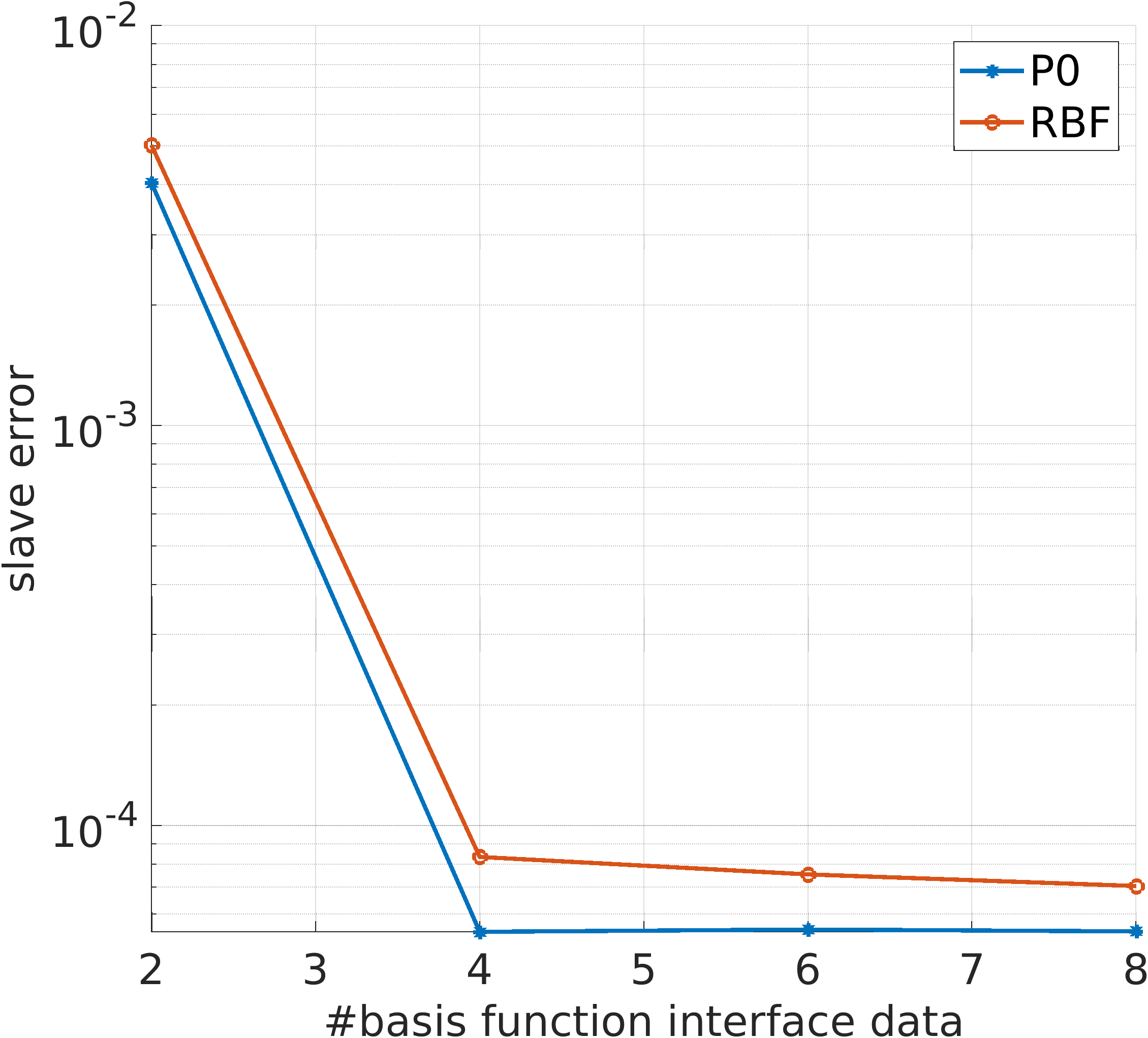}\quad
	\includegraphics[width=0.4\textwidth]{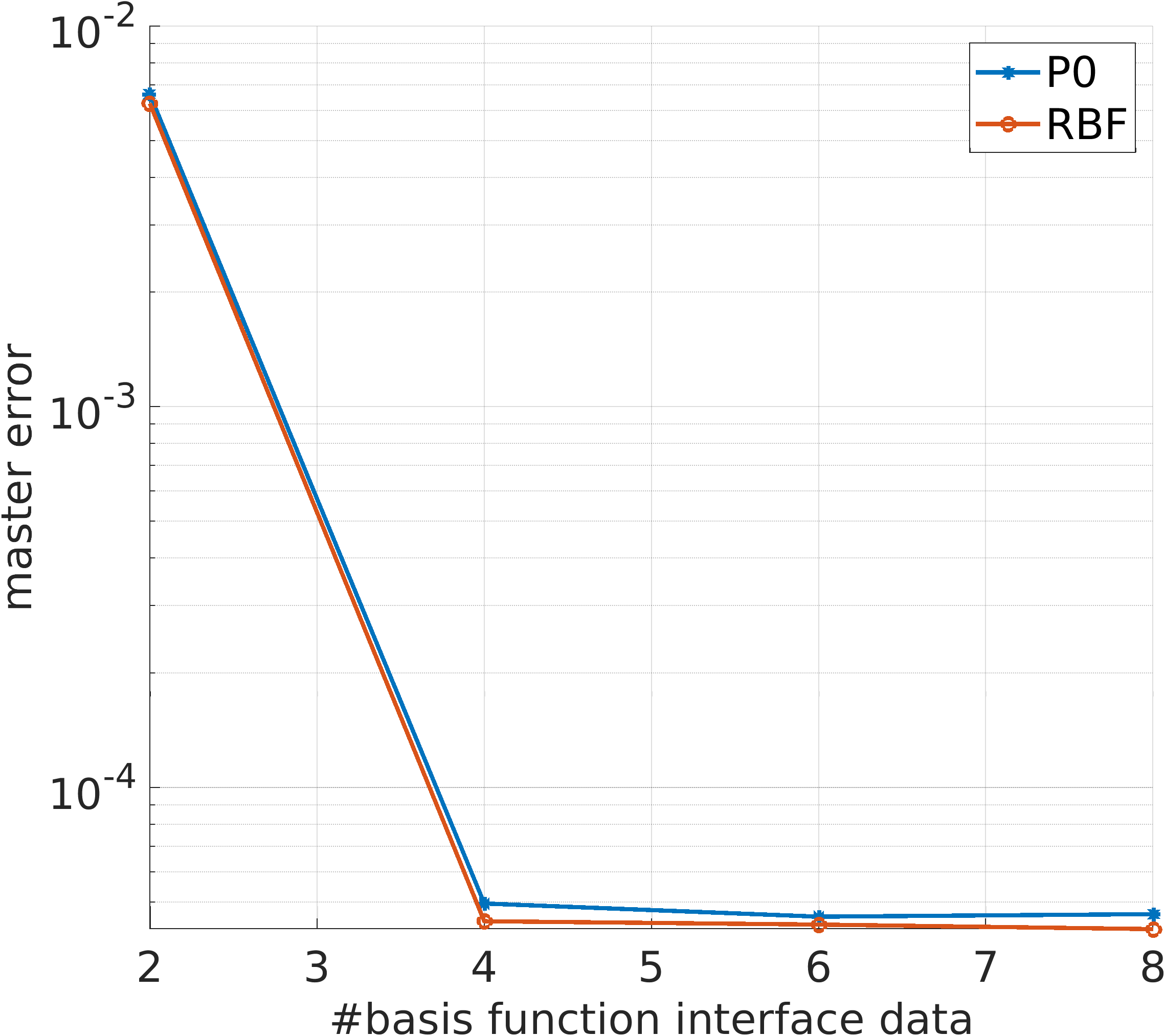}
	\caption{\emph{Test\#3.} $H_1(\Omega_i)$ mean relative error over the solution for $N_{test}=20$ different instances of the parameters between the FOM and ROM solutions varying the number of basis functions used to represent the interface data $M_1$ and $M_2$ when either nearest neighbor (P0) or Radial Basis Functions (RBF) interpolation is used to exchange the data across the interfaces. The number of basis functions for master and slave solutions is fixed to 10 and 8, respectively, to prescribe an approximation accuracy of $10^{-5}$.}
	\label{fig:comparison_rbf_heat}
\end{figure}

Considering the results obtained in Subsection \ref{App_sub:test_case_2} and recognizing that the interpolation specifically targets the interface data, thereby impacting mostly the effectiveness of the DEIM, we proceed to investigate the influence of the RBF methods, by varying the prescribed approximation accuracy only for the interface data. Therefore, differently from previous tests, here we present error graphs obtained by varying the number of basis functions employed to approximate the interface quantities, while always prescribing an accuracy of $10^{-5}$ for both the slave and master solutions. The results are illustrated in Fig. \ref{fig:comparison_rbf_heat}, where we depict the approximation error on the solution in both the slave and master domains.
We observed a remarkably similar behavior in the reduction of the approximation error when employing both the nearest neighbor approach and the RBF interpolation. It is noteworthy that the RBF proves to be more effective in approximating the master solution, while the nearest neighbor approach yields better results for the slave model. Nevertheless, these distinctions are quite minor, whereas the computational error magnitude remains consistent for each set of selected basis functions, regardless of the interpolation method employed.

Taking into account the results presented in this appendix, we observe a very similar overall outcome when using either the RBF or the nearest neighbor approach. While one might anticipate the RBF to exhibit a superior convergence rate, several factors could contribute to these findings: (i) the interface non-conformity may not be pronounced enough to favor the RBF over the simpler P0 method, (ii) the test cases considered are straightforward and do not experience significant issues with the basic P0 interpolation, (iii) the snapshots used for constructing the RBF do not incorporate interface interpolation. Each of these points requires separate investigation through a comprehensive convergence analysis, which will be the focus of future work.
\end{appendices}

\end{document}